%&amstex
\input amstex\documentstyle{amsppt}  
\pagewidth{12.5 cm}\pageheight{19 cm}\magnification\magstep1
\topmatter
\title Hecke algebras with unequal parameters\endtitle
\author G. Lusztig\endauthor
\address Department of Mathematics, M.I.T., Cambridge, MA 02139\endaddress
\thanks Supported in part by the National Science Foundation\endthanks
\endtopmatter   
\document  
\define\unb{\un{\bul}}
\define\che{\check}
\define\ZZ{\bold Z}
\define\NN{\bold N}
\define\RR{\bold R}
\define\CC{\bold C}
\define\EE{\bold E}
\define\bEE{\bar{\bold E}}
\define\mpb{\medpagebreak}
\define\frl{\forall}
\define\bul{\bullet}
\define\pe{\perp}
\define\si{\sim}
\define\wt{\widetilde}
\define\su{\sum}
\define\sqc{\sqcup}
\define\C{\text{\rm Con}}

\define\qua{\quad}
\define\hn{\hat n}

\define\hM{\hat M}
\define\uss{\un\ss}
\define\uZ{\un Z}

\define\bap{\bar p}
\define\bP{\bar P}
\define\lb{\linebreak}
\define\taa{\ti\aa}
\define\tA{\ti A}
\define\tB{\ti B}
\define\tC{\ti C}
\define\tD{\ti D}
\define\tE{\ti E}
\define\tJ{\ti J}
\define\tP{\ti P}  
\define\tS{\ti S}
\define\tZ{\ti Z}
\define\tW{\ti W}
\define\ta{\ti a}
\define\tb{\ti b}
\define\tc{\ti c}
\define\tf{\ti f}
\define\tk{\ti k}
\define\tp{\ti p}
\define\ts{\ti s}
\define\tz{\ti z}
\define\tl{\ti l}
\define\ty{\ti y}
\define\tch{\ti\ch}
\define\tca{\ti\ca}
\define\tw{\ti w}
\define\eSb{\endSb}
\define\gt{\gets}
\define\bin{\binom}
\define\tcd{\ti\cd}
\define\tcg{\ti\cg}
\define\lr{\cl\car}
\define\btZ{\bar{\ti Z}}
\define\bZ{\bar Z}

\define\bH{\bar H}
\define\op{\oplus}

\define\h{\frac{\hphantom{aa}}{\hphantom{aa}}}   
\define\da{\dag}
\redefine\sp{\spadesuit}
\define\em{\emptyset}
\define\imp{\implies}
\define\ra{\rangle}
\define\n{\notin}
\define\iy{\infty}
\define\m{\mapsto}
\define\do{\dots}
\define\la{\langle}
\define\bsl{\backslash}
\define\lras{\leftrightarrows}
\define\lra{\leftrightarrow}
\define\Lra{\Leftrightarrow}

\define\sm{\smallmatrix}
\define\esm{\endsmallmatrix}
\define\sub{\subset}
\define\bxt{\boxtimes}
\define\T{\times}
\define\ti{\tilde}
\define\nl{\newline}
\redefine\i{^{-1}}
\define\fra{\frac}
\define\un{\underline}
\define\ov{\overline}
\define\ot{\otimes}
\define\bbq{\bar{\bq}_l}

\define\Ad{\text{\rm Ad}}
\define\Hom{\text{\rm Hom}}
\define\End{\text{\rm End}}
\define\Aut{\text{\rm Aut}}

\define\sg{\text{\rm sgn}}
\define\tr{\text{\rm tr}}

\define\supp{\text{\rm supp}}

\define\bst{{}^\bigstar\!}

\define\a{\alpha}
\redefine\b{\beta}
\redefine\c{\chi}
\define\g{\gamma}
\redefine\d{\delta}
\define\e{\epsilon}
\define\et{\eta}
\define\io{\iota}
\redefine\o{\omega}
\define\p{\pi}
\define\ph{\phi}
\define\ps{\psi}
\define\r{\rho}
\define\s{\sigma}
\redefine\t{\tau}
\define\th{\theta}
\define\k{\kappa}
\redefine\l{\lambda}
\define\z{\zeta}
\define\x{\xi}

\redefine\G{\Gamma}
\redefine\D{\Delta}
\define\Om{\Omega}

\redefine\L{\Lambda}
\define\Ph{\Phi}

\redefine\aa{\bold a}

\define\boc{\bold c}

\redefine\gg{\bold g}

\define\jj{\bold j}

\define\nn{\bold n}

\redefine\ss{\bold s}

\define\bc{\bold C}

\define\be{\bold E}
\define\bF{\bold F}
\define\bg{\bold G}

\define\bk{\bold K}

\define\bn{\bold N}

\define\bq{\bold Q}
\define\br{\bold R}

\define\bz{\bold Z}

\define\ca{\Cal A}
\define\cb{\Cal B}
\define\cc{\Cal C}
\define\cd{\Cal D}
\define\ce{\Cal E}
\define\cf{\Cal F}
\define\cg{\Cal G}
\define\ch{\Cal H}
\define\ci{\Cal I}
\define\cj{\Cal J}
\define\ck{\Cal K}
\define\cl{\Cal L}
\define\cm{\Cal M}
\define\cn{\Cal N}
\define\co{\Cal O}
\define\cp{\Cal P}

\define\car{\Cal R}
\define\cs{\Cal S}
\define\ct{\Cal T}
\define\cu{\Cal U}
\define\cv{\Cal V}
\define\cw{\Cal W}
\define\cz{\Cal Z}

\define\ff{\frak f}

\define\fh{\frak h}

\define\fA{\frak A}
\define\fB{\frak B}
\define\fC{\frak C}

\define\fS{\frak S}

\define\sh{\sharp}
\define\Mod{\text{\rm Mod}}
\define\Ir{\text{\rm Irr}}
\define\In{\text{\rm Inv}}  
\define\Sy{\text{\rm Sy}}
\define\BE{Be}
\define\BO{Bo}
\define\BR{B}
\define\BM{BM}
\define\DE{De}
\define\DL{DL}
\define\EW{EW}
\define\GE{G}
\define\GP{GP}
\define\HO{H}
\define\IW{I}
\define\IM{IM}
\define\KI{Ki}
\define\KL{KL1}
\define\KLL{KL2}
\define\COX{L1} 
\define\CL{L2}  
\define\LU{L3}  
\define\SQ{L4}
\define\SING{L5}  
\define\RED{L6} 
\define\LC{L7}  
\define\CG{L8}  
\define\LCC{L9} 
\define\ICS{L10} 
\define\KYO{L11}
\define\QG{L12}
\define\LPI{L13}
\define\TENS{L14}
\define\LH{L15}
\define\LP{L16}
\define\BAR{L17}
\define\LV{LV}
\define\LX{LX}
\define\SO{So}
\define\XI{X}
\head Introduction\endhead
This is revised version of my book ``Hecke algebra with unequal parameters'',
CRM monograph series vol.18, Amer.Math.Soc. 2003. (That book was based 
on the Aisenstadt lectures given at the CRM, Universit\'e de Montr\'eal, in 
May/June 2002; and included material from lectures given at MIT during the 
Fall of 1999 \cite{\LH}.) 

This version updates the book by taking into account the recent results
of Elias and Williamson \cite{\EW} which allow us to present the results in a
more general setup. In particular the chapter on the quasisplit case (\S16) 
is more general than that in the book; it makes use of results in
the Appendix. In \S9 a discussion of 
double cosets in a Coxeter group with respect to two parabolic subgroups has
been added.
The definition of the ring J (\S18) is now done in more generality
than in the book. Since in that generality J may not have a unit element, an
imbedding of J into a ring with unit is described. A discussion of a (tensor) category
version of $J$ is given in 18.15-18.20 where it is shown how this leads to a new
construction of a simple algebraic group from an affine Weyl group which, unlike earlier
constructions, does not use perverse sheaves on an affine Grassmannian.

I am very grateful to Meinolf Geck and Darij Grinberg for supplying a 
list of misprints and minor errors. 
 
\mpb

Hecke algebras arise as endomorphism algebras of representations of groups induced 
by representations of subgroups. In these notes we are mainly interested in a 
particular kind of Hecke algebras, which arise in the representation theory of
reductive algebraic groups over finite or $p$-adic fields (see 0.3, 0.6). These 
Hecke algebras are specializations of certain algebras (known as Iwahori-Hecke 
algebras) which can be defined without reference to algebraic groups, namely by 
explicit generators and relations (see 3.2) in terms of a Coxeter group $W$ (see 
3.1) and a weight function $L:W@>>>\bz$ (see 3.1), that is, a {\it weighted
Coxeter group}. An Iwahori-Hecke algebra is completely specified by a {\it weighted
Coxeter graph}, that is, the Coxeter graph of $W$ (see 1.1) where for each vertex we
specify the value of $L$ at the corresponding simple reflection. 

A particularly simple kind of Iwahori-Hecke algebras corresponds to the case where 
the weight
function is constant on the set of simple reflections (equal parameter case). In 
this case one has the theory of the "new basis" \cite{\KL} and cells 
\cite{\KL}, \cite{\LC}, 
\cite{\LCC}. The main goal of these notes is to try to extend as much as possible 
the theory of the new basis to the general case (of not necessarily equal 
parameters). We give a number of conjectures for what should happen in the general 
case and we present some evidence for these conjectures.

We now review the contents of these notes.

\S1 introduces Coxeter groups following \cite{\BO}.
We also give a realization of the classical affine Weyl groups as periodic
permutations of
$\bz$ following an idea of \cite{\SQ}. \S2 contains some standard results on the
partial order on a Coxeter group. In \S3 we introduce the Iwahori-Hecke algebra
attached to a weighted Coxeter group. Useful references 
for this are \cite{\BO},\cite{\GP}. In \S4 we define the bar operator following
\cite{\KL}. This is used in \S5 to define the "new basis" $(c_w)$ of an 
Iwahori-Hecke algebra following \cite{\KL} for equal parameters and \cite{\LU} in
general. In \S6 we study some multiplicative properties of the new basis, following
\cite{\KL} and \cite{\LU}. In \S7 we compute explicitly the "new basis" in the case
of dihedral groups. In \S8 we define left, right and two-sided cells. In \S9 we 
study the behaviour of the new basis in relation to a given parabolic subgroup. In
\S10,\S12 we study a "basis" dual to the new basis. In \S11 we consider the case of
finite Coxeter groups. In \S13 we study the function $\aa$ on certain weighted 
Coxeter groups following an idea from \cite{\LC}. In \S14 we present a list of 
conjectures concerning cells and the function $\aa$ and we show that they can be 
deduced from a much shorter list of conjectures. These conjectures are established 
in a "split case" in \S15 (following \cite{\LCC}), in a "quasisplit case" in \S16 
and for an infinite dihedral group in \S17. Note that in the first two cases the 
proof requires arguments from the theory of Soergel modules while in the third case the 
argument is computational. In \S18, assuming the truth of these conjectures we 
develop the theory of $J$-rings in the weighted case, 
following an idea from \cite{\LCC}; we also discuss a tensor category version of a $J$-ring.
\S19,\S20,\S21 (where $W$ is assumed to be a
Weyl group) are in preparation for \S22 where the class of constructible 
representations of $W$ is introduced and studied in the weighted case 
(conjecturally these are the representations of $W$ carried by left cells), for \S23
where two-sided cells are discussed and for \S24 where certain virtual 
representations of $W$ ("virtual cells) are discussed. In \S25 we discuss the 
weighted Coxeter groups which arise in the examples 0.3 and 0.6. We formulate a 
conjecture (25.3) which relates the two-sided cells of such a weighted Coxeter group
to the two-sided cells of a larger Coxeter group with the weight function given by 
the length. In \S26 we state (following \cite{\LP}) the classification of 
irreducible representations of Hecke algebras of the type discussed in 0.6 in terms
of the geometry of the dual group. In \S27 we give a new realization of a Hecke 
algebra as in 0.3 or 0.6 as a space of functions on the rational points of an 
algebraic variety defined over $\bF_q$. This leads us to a (partly conjectural) 
geometrical interpretation of the coefficients $p_{y,w}$ of the new basis of the 
Hecke algebra in terms of intersection cohomology, generalizing the results of 
\cite{\KLL}. We expect that this geometrical interpretation should play a role in 
the proof of the conjectures in \S14 in the cases arising from algebraic groups as 
in 0.3, 0.6. In the Appendix we discuss Coxeter groups with a given automorphism
which preserve the set of simple reflections.

\subhead 0.1\endsubhead
In 0.1-0.8 we give a survey of the theory of Hecke algebras arising from reductive 
groups.

Let $\G$ be a group acting transitively on a set $X$. If $\be$ is a $\G$-equivariant
$\bc$-vector bundle over $X$ (with discrete topology) then the fibre $\be_x$ of 
$\be$ at $x\in X$ is naturally a representation of $\G_x=\{g\in\G;gx=x\}$. Moreover,
for $x\in X$, $\be\m\be_x$ is an equivalence from the category of $\G$-equivariant 
vector bundles on $X$ of finite dimension and that of finite dimensional 
$\G_x$-modules over $\bc$. 

Let $\be$ be a $\G$-equivariant $\bc$-vector bundle of finite dimension over $X$. 
Then $\G$ acts naturally on the vector space $\op_{x\in X}\be_x$. (This is the 
representation of $\G$ induced by the representation of $\G_x$ on $\be_x$, for any 
$x\in X$.) The $\bc$-algebra
$$H=H(\G,X,\be)=\End_\G(\op_{x\in X}\be_x)$$
is called the {\it Hecke algebra}. The image of the obvious imbedding 
$$H\sub\prod_{(x,x')\in X\T X}\Hom(\be_x,\be_{x'}),\qua 
\ph\m(\ph^x_{x'})_{(x,x')\in X\T X}$$
consists of all $(f^x_{x'})\in\prod_{(x,x')\in X\T X}\Hom(\be_x,\be_{x'})$ such that

for any $x\in X$ we have $f^x_{x'}=0$ for all but finitely many $x'\in X$;

for any $g\in\G$ and any $(x,x')\in X\T X$, the compositions
$\be_x@>f^x_{x'}>>\be_{x'}@>g>>\be_{gx'}$,
$\be_x@>g>>\be_{gx}@>f^{gx}_{gx'}>>\be_{gx'}$ coincide.
\nl
For any $\G$-orbit $\cc$ in $X\T X$ we set
$$H_\cc=H(\G,X,\be)_\cc=\{\ph\in H;\ph^x_{x'}\ne 0\imp(x',x)\in\cc\}.$$
Then $H_\cc=0$ unless $\cc$ is {\it finitary} in the following sense:

for some (or any) $x\in X$, the set $\{x'\in X;(x',x)\in\cc\}$ is finite,
\nl
in which case 
$$H_\cc@>\si>>\Hom_{\G_x\cap\G_{x'}}(\be_x,\be_{x'}),\qua\ph\m\ph^x_{x'}$$ 
for $(x',x)\in\cc$. Moreover, 

(a) $H=\op_{\cc \text{ finitary}}H_\cc$.

\subhead 0.2\endsubhead
To explain how Hecke algebras arise from reductive algebraic groups we need the
notion of "unipotent cuspidal representation".

Let $p$ be a prime number and let $\bF$ be an algebraic closure of the finite field
with $p$ elements. Let $q$ be a power of $p$ and let $\bF_q$ be the subfield of 
$\bF$ with $q$ elements. Let $G$ be a connected reductive algebraic group over $\bF$
with a fixed $\bF_q$ structure and let $F:G@>>>G$ be the corresponding Frobenius 
map. 

We refer to \cite{\DL} for the notion of unipotent cuspidal representation of the 
finite group $G^F=G(\bF_q)$. We will only give here the definition assuming that $q$
is sufficiently large. Let $E$ be an irreducible representation over $\bc$ of $G^F$
and let $\c_E:G^F@>>>\bc$ be its character. We say that $E$ is unipotent if, for any
$F$-stable maximal torus $T$ of $G$, the restriction of $\c_E$ to the set of regular
elements in $T^F$ is a constant, say $c_T\in\bz$. We say that $E$ is unipotent 
cuspidal if, in addition, for any $T$ as above that is contained in some proper 
$F$-stable parabolic subgroup of $G$, we have $c_T=0$.

The unipotent cuspidal representations of $G^F$ are classified  in \cite{\RED}. For
example, if $G$ is a torus times a symplectic group of rank $n\ge 0$ then $G^F$ has
(up to isomorphism) a unique unipotent cuspidal representation if $n=k^2+k$ for some
integer $k\ge 0$, and none, otherwise.

\subhead 0.3\endsubhead
Let $G,F$ be as in 0.2. For any parabolic subgroup $P$ of $G$ let $U_P$ be the 
unipotent radical of $P$ and let $\bP=P/U_P$. Let $\cp$ be an $F$-stable 
$G$-conjugacy class of parabolic subgroups of $G$ and let $\be$ be a 
$G^F$-equivariant vector bundle over $\cp^F$ (a $G^F$-homogeneous space) such that 
for some (or any) $P\in\cp^F$, the $P^F$-action on the fibre $\be_P$ of $\be$ at $P$
factors through a unipotent, cuspidal $\bP^F$-module. (To give such $\be$ is the 
same as to give, for some $P\in\cp^F$, a unipotent cuspidal representation of 
$\bP^F$.) The Hecke algebra $H(G^F,\cp^F,\be)$ is defined.

Let $W$ be the set of $G$-orbits on the set of ordered pairs of Borel subgroups in 
$G$; it is known that $W$ may be naturally regarded as a finite Coxeter group (see 
1.1) with a set $S$ of simple reflections. Now any Borel subgroup of $G$ is 
contained in a unique subgroup in $\cp$; this defines a (surjective) map from $W$ to
$G\bsl(\cp\T\cp)$, the set of $G$-orbits in $\cp\T\cp$. The inverse image of the 
diagonal orbit under this map is the subgroup $W_J$ of $W$ generated by a subset $J$
of $S$ and $W@>>>G\bsl(\cp\T\cp)$ factors through a bijection 
$$W_J\bsl W/W_J@>\si>>G\bsl(\cp\T\cp).\tag a$$
Let $\cw$ be the set of all $w\in W$ such that $wW_J=W_Jw$ and $w$ has minimal 
length in $wW_J=W_Jw$. Then $\cw$ is a subgroup of $W$. The Frobenius map $u:W@>>>W$
restricts to an isomorphism $u:\cw@>>>\cw$ whose fixed point set $\cw^u$ is 
naturally a Coxeter group with simple reflections indexed by $u\bsl(S-J)$ (set of 
orbits of $u:S-J@>>>S-J$). See 25.1(a). A $G$-orbit $\co$ on $\cp\T\cp$ is said to 
be {\it good} if for $(P,P')\in\co$ we have $(P\cap P')U_P=P$ or equivalently 
$(P\cap P')U_{P'}=P'$. Otherwise, $\co$ is said to be {\it bad}. If $\co$ is a good,
$F$-stable $G$-orbit on $\cp\T\cp$ then $\co^F$ is a $G^F$-orbit on $\cp^F\T\cp^F$ 
and $\dim H(G^F,\cp^F,\be)_{\co^F}=1$. If $\co$ is an $F$-stable bad $G$-orbit on 
$\cp\T\cp$ then $\co^F$ is a $G^F$-orbit on $\cp^F\T\cp^F$ and 
$\dim H(G^F,\cp^F,\be)_{\co^F}=0$. Now the bijection (a) restricts (via the 
imbedding $\cw^u\sub W_J\bsl W/W_J$, $w\m W_JwW_J$) to a bijection $w\m\co_w$ of 
$\cw^u$ onto the set of good, $F$-stable $G$-orbits on $\cp\T\cp$. It follows that 
0.1(a) becomes in our case
$$H(G^F,\cp^F,\be)=\op_{w\in\cw^u}H(G^F,\cp^F,\be)_{\co_w}\tag b$$
with
$$\dim H(G^F,\cp^F,\be)_{\co_w}=1 \text{ for all } w\in\cw^u.\tag c$$
Let $\t_k$ be the generator of $\cw^u$ corresponding to $k\in u\bsl(S-J)$. There is
a unique basis element $T_{\t_k}$ of $H(G^F,\cp^F,\be)_{\co_{\t_k}}$ such that
$$(T_{\t_k}+q^{-N_k/2})(T_{\t_k}-q^{N_k/2})=0\tag d$$
for some $N_k\in\bz_{>0}$. ($N_k$ is uniquely determined.) The elements 
$T_{\t_k}(k\in u\bsl(S-J))$ generate the $\bc$-algebra $H(G^F,\cp^F,\be)$. They 
satisfy identities of the form
$$T_{\t_k}T_{\t_{k'}}T_{\t_k}\do=T_{\t_{k'}}T_{\t_k}T_{\t_{k'}}\do\tag e$$
for $k\ne k'$ in $u\bsl(S-J)$; both products have a number of factors equal to the
order of $\t_k\t_{k'}$ in $\cw^u$. Now $T_{\t_k}\m T_{\t_k}$ gives an isomorphism 
from an Iwahori-Hecke algebra (see 3.2) specialized at $v=\sqrt q$ to the algebra 
$H(G^F,\cp^F,\be)$. 

The function $k\m N_k$ coincides with the function $k\m L(\t_k)$ in 25.2.

(The results in this subsection appeared in \cite{\CL,\COX}. In the special case 
where $\cp$ is the set of Borel subgroups of $G$ and $\be$ is the trivial vector 
bundle $\bc$, they were first proved by Iwahori \cite{\IW}; if, in addition, $u=1$
on $W$ then $N_k=1$ for all $k$.)

\subhead 0.4\endsubhead
Let $V$ be an $\bF_q$-vector space of dimension $n\ge 2$. Then $G=SL(\bF\ot V)$ has
a natural $\bF_q$-structure. Let $\cp$ be the set of all Borel subgroups of $G$. 
Then $\cp^F$ may be identified as a set with $G^F$-action with the set $\cf$ of all
flags \lb $V_*=(V_0\sub V_1\sub V_2\sub\do\sub V_n)$ of subspaces of $V$ 
($\dim V_i=i$ for all $i$). 

Let $V_*=(V_0\sub V_1\sub V_2\sub\do\sub V_n)$,
$V'_*=(V'_0\sub V'_1\sub V'_2\sub\do\sub V'_n)$ be flags in $\cf$. For 
$i\in[0,n],j\in[1,n]$ we set 
$d_{ij}=\dim\fra{V'_i\cap V_j}{V'_i\cap V_{j-1}}\in\{0,1\}$. For $i\in[0,n]$ we set
$X_i=\{j\in[1,n];d_{ij}=1\}$. Then 
$\em=X_0\sub X_1\sub X_2\sub\do\sub X_n=[1,n]$ and
for $i\in[1,n]$ there is a unique $a_i\in[1,n]$ such that $X_i=X_{i-1}\sqc\{a_i\}$.
Also, $i\m a_i$ is a permutation of $[1,n]$. Now $(V_*,V'_*)\m(a_i)$ defines a 
bijection of $G^F\bsl(\cp^F\T\cp^F)=G^F\bsl(\cf\T\cf)$ with the symmetric group 
$\fS_n$. Let $\be$ be the trivial $G^F$-equivariant vector bundle $\bc$ on 
$\cp^F=\cf$. Then $H(G^F,\cp^F,\be)$ is defined. In our case we have 
$W=\cw=\cw^u=\fS_n$. 

\subhead 0.5\endsubhead
Let $V,n$ be as in 0.4. Assume that $n=2m$ and that $V$ has a fixed non-degenerate
symplectic form $\la,\ra:V\T V@>>>\bF_q$. 
Then $G=Sp(\bF\ot V)$ has a natural $\bF_q$-structure. 
Assume that $m=r+k^2+k$ where $k\in\bn,r\in\bz_{>0}$. Let $\cf$ be the set of all 
flags $V_*=(V_0\sub V_1\sub V_2\sub\do\sub V_r)$ of isotropic subspaces of $V$ 
($\dim V_i=i$ for all $i$). There is a unique $G$-conjugacy class $\cp$ of parabolic
subgroups of $G$ such that, if $V_*\in\cf$, then 
$$\{g\in G;g(\bF\ot V_j)=\bF\ot V_j\qua\frl j\}\in\cp.$$
We may identify $\cp^F=\cf$ as spaces with $G^F$-action.

Let $U\m\cd_k(U)$ be a functor from the category of symplectic vector spaces of
dimension $2k^2+2k$ over $\bF_q$ (and isomorphisms between them) to the category of
$\bc$-vector spaces (and isomorphisms between them) such that for any $U$, the
$Sp(U)$-module $\cd_k(U)$ is unipotent, cuspidal. (Such a functor exists and is 
unique up to isomorphism.) Let $\be$ be the vector bundle over $\cp^F$ (or 
equivalently $\cf$) whose fibre at 
$V_*=(V_0\sub V_1\sub V_2\sub\do\sub V_r)\in\cf^F$ is $\cd_k(V_r^\pe/V_r)$. (Here
$V_s^\pe=\{x\in V;\la x,V_s\ra=0\}$.)

This vector bundle is naturally $G^F$-equivariant (since $\cd_k$ is a functor). 
Hence $H(G^F,\cp^F,\be)$ is defined.

Let $V_*=(V_0\sub V_1\sub V_2\sub\do\sub V_r)$,
$V'_*=(V'_0\sub V'_1\sub V'_2\sub\do\sub V'_r)$ be flags in $\cf$. The $G$-orbit of
the point of $\cp\T\cp$ corresponding to $(V_*,V'_*)$ is good if the following three
equivalent conditions hold:

$V_r\cap V'_r=V_r\cap V'_r{}^\pe$,

$V_r\cap V'_r=V_r^\pe\cap V'_r$,

$V_r\cap V'_r=(V_r^\pe\cap V'_r{}^\pe)\cap(V_r+V'_r)$.
\nl
If these conditions hold, we can define an isomorphism
$\psi^{V_r}_{V'_r}:V_r^\pe/V_r@>>>V'_r{}^\pe/V'_r$ by requiring that the diagram
$$\CD (V_r^\pe\cap V'_r{}^\pe)/(V_r\cap V'_r)@>=>>
(V_r^\pe\cap V'_r{}^\pe)/(V_r\cap V'_r)\\  @VVV @VVV\\
V_r^\pe/V_r@>\psi^{V_r}_{V'_r}>>V'_r{}^\pe/V'_r\endCD$$
(where the vertical maps are the isomorphisms induced by the inclusion) is
commutative. In this case, to $(V_*,V'_*)$ we associate an element $\s$
of the group $W_r$ of permutations of $\cs=\{1,2,\do,r,r',\do,2',1'\}$ that commute
with the involution $f:\cs@>>>\cs,j\m j',j'\m j$. For $j\in[1,r]$ we set
$$A_j=\{h\in[1,r]; V_{h-1}\cap V'_j\ne V_h\cap V'_j\},\qua
B_j=f\{h\in[1,r];V_{h-1}^\pe\cap V'_j\ne V_h^\pe\cap V'_j\}.$$
Then $\sh(A_j\cap B_j)=j,A_1\cap B_1\sub A_2\cap B_2\sub\do\sub A_s\cap B_s$ and
$h\in A_j\imp h'\n B_j$. For $j\in[1,r]$ define $a_j\in\cs$ by
$A_j\cup B_j=\{a_1,a_2,\do,a_j\}$. Then $\s$ is defined by the condition that
$\s(j)=a_j$ for $j\in[1,r]$. We see that in our case, $\cw=\cw^u$ may be identified
with $W_r$. In our case, the Iwahori-Hecke algebra corresponds to the weighted 
Coxeter graph
$$\bul_{2k+1} =\bul_1-\bul_1-\do-\bul_1$$
($r$ vertices); in the case where $r=1$ this should be interpreted as a graph with
one vertex marked by $2k+1$.

\subhead 0.6\endsubhead
Let $\e$ be an indeterminate. Let $\bk$ be the subfield of $\bF((\e))$ generated by
$\bF_q((\e))$ and $\bF$. Let $\bg$ be a split connected simply connected almost 
simple algebraic group over $\bk$ with a fixed $\bF_q((\e))$-rational structure. We
identify $\bg$ with its group of $\bk$-points. There is a "Frobenius map" 
$F:\bg@>>>\bg$ whose fixed point set is $\bg(\bF_q((\e)))$. Let $\cb$ be the set of
all Iwahori subgroups of $\bg$. (This concept will be illustrated in 0.7.) A 
subgroup of $\bg$ is said to be a parahoric subgroup if it is $\ne\bg$ and it 
contains some Iwahori subgroup. If $P$ is a parahoric subgroup then $P$ has a 
"pro-unipotent radical" $U_P$ and $\bP=P/U_P$ is a connected, reductive group over 
$\bF$. Let $\cp$ be an $F$-stable $\bg$-conjugacy class of parahoric subgroups of 
$\bg$ and let $\be$ be a $\bg^F$-equivariant vector bundle over $\cp^F$ (a
$\bg^F$-homogeneous space) such that for some (or any) $P\in\cp^F$, the $P^F$-action
on the fibre $\be_P$ of $\be$ at $P$ factors through a unipotent, cuspidal 
$\bP^F$-module. (To give such $\be$ is the same as to give, for some $P\in\cp^F$, a
unipotent cuspidal representation of $\bP^F$.) The Hecke algebra 
$H(\bg^F,\cp^F,\be)$ is defined.

Let $W$ be the set of $\bg$-orbits on $\cb\T\cb$; it is known that $W$ may be 
naturally regarded as a Coxeter group (more precisely, an affine Weyl group, see 
1.15) with a set $S$ of simple reflections. Now any Iwahori subgroup of $\bg$ is 
contained in a unique subgroup in $\cp$; this defines a (surjective) map from $W$ to
$\bg\bsl(\cp\T\cp)$, the set of $\bg$-orbits in $\cp\T\cp$. The inverse image of the
diagonal orbit under this map is the subgroup $W_J$ of $W$ generated by a subset $J$
of $S$ and $W@>>>\bg\bsl(\cp\T\cp)$ factors through a bijection 
$$W_J\bsl W/W_J@>\si>>\bg\bsl(\cp\T\cp).\tag a$$
Let $\cw$ be the set of all $w\in W$ such that $wW_J=W_Jw$ and $w$ has minimal 
length in $wW_J=W_Jw$. Then $\cw$ is a subgroup of $W$. The Frobenius map $u:W@>>>W$
restricts to an isomorphism $u:\cw@>>>\cw$ whose fixed point set $\cw^u$ is 
naturally an infinite Coxeter group with simple reflections indexed by $u\bsl(S-J)$
(set of orbits of $u:S-J@>>>S-J$). (We make the additional assumption that 
$\sh(u\bsl(S-J))\ge 2$.) See 25.1(a). A $\bg$-orbit $\co$ on $\cp\T\cp$ is said to 
be {\it good} if for $(P,P')\in\co$ we have $(P\cap P')U_P=P$ or equivalently 
$(P\cap P')U_{P'}=P'$. Otherwise, $\co$ is said to be {\it bad}. If $\co$ is a good,
$F$-stable $\bg$-orbit on $\cp\T\cp$ then $\co^F$ is a $\bg^F$-orbit on 
$\cp^F\T\cp^F$ and $\dim H(\bg^F,\cp^F,\be)_{\co^F}=1$. If $\co$ is an $F$-stable 
bad $\bg$-orbit on $\cp\T\cp$ then $\co^F$ is a $\bg^F$-orbit on $\cp^F\T\cp^F$ and
$\dim H(\bg^F,\cp^F,\be)_{\co^F}=0$. Now the bijection (a) restricts (via the 
imbedding $\cw^u\sub W_J\bsl W/W_J$, $w\m W_JwW_J$) to a bijection $w\m\co_w$ of 
$\cw^u$ onto the set of good, $F$-stable $\bg$-orbits on $\cp\T\cp$. It follows that
0.1(a) becomes in our case
$$H(\bg^F,\cp^F,\be)=\op_{w\in\cw^u}H(\bg^F,\cp^F,\be)_{\co_w}\tag b$$
with
$$\dim H(\bg^F,\cp^F,\be)_{\co_w}=1 \text{ for all } w\in\cw^u.\tag c$$
Let $\t_k$ be the generator of $\cw^u$ corresponding to $k\in u\bsl(S-J)$. There is
a unique basis element $T_{\t_k}$ of $H(\bg^F,\cp^F,\be)_{\co_{\t_k}}$ such that
$$(T_{\t_k}+q^{-N_k/2})(T_{\t_k}-q^{N_k/2})=0\tag d$$
for some $N_k\in\bz_{>0}$. ($N_k$ is uniquely determined.) The elements 
$T_{\t_k} (k\in u\bsl(S-J))$ generate the $\bc$-algebra $H(\bg^F,\cp^F,\be)$. They 
satisfy identities of the form
$$T_{\t_k}T_{\t_{k'}}T_{\t_k}\do=T_{\t_{k'}}T_{\t_k}T_{\t_{k'}}\do$$
for $k\ne k'$ in $u\bsl(S-J)$ with $\t_k\t_{k'}$ of order $m_{k,k'}<\iy$ in 
$\cw^u$ (both products have $m_{k,k'}$ factors). Now $T_{\t_k}\m T_{\t_k}$ 
gives an isomorphism from an Iwahori-Hecke algebra (see 3.2) specialized at 
$v=\sqrt q$ to the algebra $H(\bg^F,\cp^F,\be)$. The function $k\m N_k$ coincides 
with the function $k\m L(\t_k)$ in 25.2.

(The results in this subsection appeared in  \cite{\KYO,\LPI}. In the special case 
where $\cp=\cb$, $u=1$ and $\be$ is the trivial vector bundle $\bc$, they were first
proved by Iwahori-Matsumoto \cite{\IM}; in this case, $N_k=1$ for all $k$.)

\subhead 0.7\endsubhead
Let $V$ be an $\bF_q((\e))$-vector space of dimension $n\in[2,\iy)$ with a fixed
volume form $\o$. Then $G=SL(\bF\ot_{\bF_q}V)$ is as in 0.6. Let $\fA=\bF_q[[\e]]$.
A lattice in $V$ is an $\fA$-submodule of $V$ of rank $n$ which generates $V$. For a
lattice $\cl$ in $V$ we set $vol(\cl)=r$ where $m\in\bz$ is defined by the condition
that the $n$-th exterior power of $\cl$ (an $\fA$-module) has $\e^{-r}\o$ as basis 
element. Let $\cf$ be the set of all sequences of lattices $(\cl_j)_{j\in\bz}$ such
that 
$$\cl_{j-1}\sub\cl_j,vol(\cl_j)=j,\e\cl_j=\cl_{j-n}$$
for all $j$. We may identify $\cb^F$ with $\cf$ as sets with (transitive) 
$G^F$-action. 

Let $\cl_*=(\cl_j)_{j\in\bz},\cl'_*=(\cl'_j)_{j\in\bz}$ be elements of $\cf$. For 
$i,j\in\bz$ we set 
$$d_{ij}=\dim\fra{\cl'_i\cap\cl_j}{\cl'_i\cap\cl_{j-1}}\in\{0,1\}.$$
For $i\in\bz$ let $X_i=\{j\in\bz;d_{ij}=1\}$. Then $X_{i-1}\sub X_i$ for all $i$ and
$\sh(X_i-X_{i-1})=1$ for all $i$. Define $a_i\in\bz$ by $X_i=X_{i-1}\sqc\{a_i\}$.
Now $X_{i-n}=X_i-n$. Hence 
$$a_{i+n}=a_i+n \text{ for all } i\in\bz.\tag a$$
One can check that $i\m a_i$ is a bijection $\bz@>>>\bz$. Using the fact that 
$vol(\cl_j)=vol(\cl'_j)=j$ we see that
$$\su_{i=1}^n(a_i-i)=0.\tag b$$ 
This gives a bijection of $W$ onto the group of all bijections $\bz@>\si>>\bz$ that
satisfy (a),(b) (see 1.12).

\subhead 0.8\endsubhead
Let $V,n$ be as in 0.7. Assume that $n=2m$ and that $V$ has a fixed non-degenerate
symplectic form $\la,\ra:V\T V@>>>\bF_q((\e))$. 
Then $G=Sp(\bF\ot_{\bF_q}V)$ is as in 0.6. If $\cl$ is a
lattice in $V$ then $\cl^\sh=\{x\in V;\la x,\cl\ra\in\fA\}$ is again a lattice; 
moreover, $(\cl^\sh)^\sh=\cl$. Assume that $m=r+k^2+k+l^2+l$ where 
$k,l,r\in\bn,r\ge 1$. Let $\cn$ be the set of all integers of the form $a+2mt$ where
$t\in\bz$ and 
$$k^2+k\le a\le k^2+k+r \text{ or } -(k^2+k+r)\le a\le-(k^2+k).$$
Let $\cf$ be the set of all sequences of lattices $(\cl_j)_{j\in\cn}$ such that 

$\cl_j\sub\cl_{j'}$ if $j\le j'$ in $\cn$,

$\cl_j^\sh=\cl_{-j}$ for all $j\in\cn$, 

$\e\cl_j=\cl_{j-2m}$ for all $j\in\cn$, 

$vol(\cl_j)=j$ for all $j\in\cn$.
\nl 
Here the volume of a lattice is defined in terms of the volume form on $V$ attached
to the symplectic form. There is a unique $G$-conjugacy class $\cp$ of parahoric 
subgroups of $G$ such that, if $(\cl_j)_{j\in\cn}\in\cf$, then 
$$\{g\in G; g(\bF\ot\cl_j)=\bF\ot\cl_j\qua\frl j\in\cn\}\in\cp.$$
We may identify $\cp^F$ and $\cf$ as sets with $G^F$-action. If 
$(\cl_j)_{j\in\cn}\in\cf$, then the $\bF_q$-vector space $\cl_{k^2+k}/\cl_{-k^2-k}$
(of dimension $2k^2+2k$) has a natural non-degenerate symplectic form induced by 
$x,y\m Res\la x,y\ra$ and the $\bF_q$-vector space $\cl_{2m-k^2-k-r}/\cl_{k^2+k+r}$
(of dimension $2l^2+2l$) has a natural non-degenerate symplectic form induced by 
$x,y\m Res\la x,\e y\ra$. Here $Res:\bF_q((\e))@>>>\bF_q$ denotes residue at $0$.

Let $\cd_k$ be a functor as in 0.5 and let $\cd_l$ be an analogous functor obtained
by replacing $k$ by $l$. Let $\be$ be the vector bundle over $\cp^F$ (or 
equivalently $\cf$) whose fibre at $(\cl_j)_{j\in\cn}\in\cf$ is
$$\cd_k(\cl_{k^2+k}/\cl_{-k^2-k})\ot\cd_l(\cl_{2m-k^2-k-r}/\cl_{k^2+k+r}).$$
This vector bundle is naturally $G^F$-equivariant (since $\cd_k,\cd_l$ are 
functors). Hence $H(G^F,\cp^F,\be)$ is defined. In our case, the Iwahori-Hecke 
algebra corresponds to the weighted Coxeter graph
$$\bul_{2k+1} =\bul_1-\bul_1-\do-\bul_1=\bul_{2l+1}$$
($r+1$ vertices); in the case where $r=1$ this should be interpreted as a Coxeter
graph with 2 vertices marked by $2k+1,2l+1$, joined by a quadruple edge.

\subhead 0.9. Notation\endsubhead
We set $[a,b]=\{z\in\bz;a\le z\le b\}$, 
$[a,b)=\{z\in\bz;a\le z<b\}$. If $X$ is a subset of a group $G$, we denote by 
$\la X\ra$ the subgroup of $G$ generated by $X$. 

\head Contents\endhead
 1. Coxeter groups

 2. Partial order on $W$

 3. The algebra $\ch$

 4. The bar operator 

 5. The elements $c_w$

 6. Left or right multiplication by $c_s$

 7. Dihedral groups 

 8. Cells

 9. Cosets of parabolic subgroups

 10. Inversion

 11. The longest element for a finite $W$

 12. Examples of elements $D_w$

 13. The function $\aa$

 14. Conjectures

 15. Example: the split case

 16. Example: the quasisplit case

 17. Example: the infinite dihedral case

 18. The ring $J$

 19. Algebras with trace form  

 20. The function $\aa_E$

 21. Study of a left cell

 22. Constructible representations

 23. Two-sided cells

 24. Virtual cells
 
 25. Relative Coxeter groups

 26. Representations

 27. A new realization of Hecke algebras

 Appendix

\head 1. Coxeter groups\endhead
\subhead 1.1\endsubhead
Let $S$ be a finite set and let $(m_{s,s'})_{(s,s')\in S\T S}$ be a matrix with 
entries in $\bn\cup\{\iy\}$ such that $m_{s,s}=1$ for all $s$ and 
$m_{s,s'}=m_{s',s}\ge 2$ for all $s\ne s'$. (A {\it Coxeter matrix}.) In the case 
where $m_{s,s'}\in\{2,3,4,6,\iy\}$ for all $s\ne s'$, the matrix
$(m_{s,s'})_{(s,s')\in S\T S}$ is completely described by a graph (the {\it Coxeter
graph}) with set of vertices in bijection with $S$ where the vertices corresponding
to $s\ne s'$ are joined by an edge if $m_{s,s'}=3$, by a double edge if 
$m_{s,s'}=4$, by a triple edge if $m_{s,s'}=6$, by a quadruple edge if 
$m_{s,s'}=\iy$.

Let $W$ be the group defined by the generators $s (s\in S)$ and relations
$$(ss')^{m_{s,s'}}=1$$
for any $s,s'$ in $S$ such that $m_{s,s'}<\iy$. We say that $W,S$ is a {\it Coxeter
group}. Note that the Coxeter matrix is uniquely determined by $W,S$ (see 1.3(b) 
below). We sometimes refer to $W$ itself as a Coxeter group. In $W$ we have $s^2=1$
for all $s$. Clearly, there is a unique homomorphism 
$$\sg:W@>>>\{1,-1\}$$
such that $\sg(s)=-1$ for all $s$. ("{\it Sign} representation".)

For $w\in W$ let $l(w)$ be the smallest integer $q\ge 0$ such that $w=s_1s_2\do s_q$
with $s_1,s_2,\do,s_q$ in $S$. (We then say that $w=s_1s_2\do s_q$ is a {\it reduced
expression} and $l(w)$ is the {\it length} of $w$.) Now $l(1)=0,l(s)=1$ for 
$s\in S$. (Indeed, $s\ne 1$ in $W$ since $\sg(s)=-1,\sg(1)=1$.)

\proclaim{Lemma 1.2} Let $w\in W,s\in S$.

(a) We have either $l(sw)=l(w)+1$ or $l(sw)=l(w)-1$.

(b) We have either $l(ws)=l(w)+1$ or $l(ws)=l(w)-1$.
\endproclaim
Clearly, $\sg(w)=(-1)^{l(w)}$. Since $\sg(sw)=-\sg(w)$, we have 
$(-1)^{l(sw)}=-(-1)^{l(w)}$. Hence $l(sw)\ne l(w)$. This, together with the obvious
inequalities $l(w)-1\le l(sw)\le l(w)+1$ gives (a). The proof of (b) is similar.

\proclaim{Proposition 1.3} Let $E$ be an $\br$-vector space with basis 
$(e_s)_{s\in S}$. For $s\in S$ define a linear map $\s_s:E@>>>E$ by 
$\s_s(e_{s'})=e_{s'}+2\cos\fra{\p}{m_{s,s'}}e_s$ for all $s'\in S$.

(a) There is a unique homomorphism $\s:W@>>>GL(E)$ such that $\s(s)=\s_s$ for all
$s\in S$.

(b) If $s\ne s'$ in $S$, then $ss'$ has order $m_{s,s'}$ in $W$. In particular,
$s\ne s'$ in $W$.
\endproclaim
We have $\s_s(e_s)=-e_s$ and $\s_s$ induces the identity map on $E/\br e_s$. Hence
$\s_s^2=1$. Now let $s\ne s'$ in $S$, $m=m_{s,s'}$, $\Ph=\s_s\s_{s'}$. We have

$\Ph(e_s)=(4\cos^2\fra{\p}{m}-1)e_s+2\cos\fra{\p}{m}e_{s'}$,

$\Ph(e_{s'})=-2\cos\fra{\p}{m}e_s-e_{s'}$.
\nl
Hence $\Ph$ restricts to an endomorphism $\ph$ of $\br e_s\op\br e_{s'}$ whose 
characteristic polynomial is
$$X^2-2\cos\fra{2\p}{m}X+1=(X-e^{2\p\sqrt{-1}/m})(X-e^{-2\p\sqrt{-1}/m}).$$
It follows that, if $2<m<\iy$, then $1+\ph+\ph^2+\do+\ph^{m-1}=0$. The same holds if
$m=2$ (in this case we see directly that $\ph=-1$). Since $\Ph$ induces the identity
map on $E/(\br e_s\op\br e_{s'})$, it follows that $\Phi:E@>>>E$ has order $m$ (if 
$m<\iy$). If $m=\iy$, we have $\ph\ne 1$ and $(\ph-1)^2=0$, hence $\ph$ has infinite
order and $\Ph$ has also infinite order. Now (a), (b) follow.

\proclaim{Corollary 1.4} Let $s_1\ne s_2$ in $S$. For $k\ge 0$ let 
$1_k=s_1s_2s_1\do$ ($k$ factors), $2_k=s_2s_1s_2\do$ ($k$ factors).

(a) Assume that $m=m_{s_1,s_2}<\iy$. Then $\la s_1,s_2\ra$ consists of the elements
$1_k,2_k$ ($k=0,1,\do,m$); these elements are distinct except for the equalities 
$1_0=2_0,1_m=2_m$. For $k\in[0,m]$ we have $l(1_k)=l(2_k)=k$.

(b) Assume that $m_{s_1,s_2}=\iy$. Then $\la s_1,s_2\ra$ consists of the elements 
$1_k,2_k$ ($k=0,1,\do$); these elements are distinct except for the equality 
$1_0=2_0$. For all $k\ge 0$ we have $l(1_k)=l(2_k)=k$.
\endproclaim
This follows immediately from 1.3(b).

We identify $S$ with a subset of $W$ (see 1.3(b)) said to be the set of
{\it simple reflections}. Let 
$$T=\cup_{w\in W}wSw\i\sub W.$$

\proclaim{Proposition 1.5} Let $R=\{1,-1\}\T T$. For $s\in S$ define $U_s:R@>>>R$ by
$U_s(\e,t)=(\e(-1)^{\d_{s,t}},sts)$ where $\d$ is the Kronecker symbol. There is a 
unique homomorphism $U$ of $W$ into the group of permutations of $R$ such that
$U(s)=U_s$ for all $s\in S$.
\endproclaim
We have $U_s^2(\e,t)=(\e(-1)^{\d_{s,t}+\d_{s,sts}},t)=(\e,t)$ since the conditions 
$s=t$, $s=sts$ are equivalent. Thus, $U_s^2=1$. For $s\ne s'$ in $S$ with
$m=m_{s,s'}<\iy$ we have
$$U_sU_{s'}(\e,t)=(\e(-1)^{\d_{s',t}+\d_{s,s'ts'}},ss'ts's)$$
hence
$$(U_sU_{s'})^m(\e,t)=(\e(-1)^\k,(ss')^mt(s's)^m)=(\e(-1)^\k,t)$$
where 
$$\k=\d_{s',t}+\d_{s,s'ts'}+\d_{s',ss'ts's}+\do=
\d_{s',t}+\d_{s'ss',t}+\d_{s'ss'ss',t}+\do$$
(both sums have exactly $2m$ terms). It is enough to show that $\k$ is even, or that
$t$ appears an even number of times in the $2m$-term sequence 
$s',s'ss',s'ss'ss',\do$. This follows from the fact that in this sequence the $k$-th
term is equal to the $(k+m)$-th term for $k=1,2,\do,m$.

\proclaim{Proposition 1.6} Let $w\in W$. Let $w=s_1s_2\do s_q$ be a reduced 
expression.

(a) The elements $s_1,s_1s_2s_1,s_1s_2s_3s_2s_1,\do,s_1s_2\do s_q\do s_2s_1$ are 
distinct. 

(b) These elements form a subset of $T$ that depends only on $w$, not on the choice
of reduced expression for it.
\endproclaim
Assume that $s_1s_2\do s_i\do s_2s_1=s_1s_2\do s_j\do s_2s_1$ for some
$1\le i<j\le q$. Then $s_i=s_{i+1}s_{i+2}\do s_j\do s_{i+2}s_{i+1}$ hence
$$\align s_1s_2\do s_q&=s_1s_2\do s_{i-1}(s_{i+1}s_{i+2}\do s_j\do 
s_{i+2}s_{i+1})s_{i+1}\do s_js_{j+1}\do s_q\\&=
s_1s_2\do s_{i-1}s_{i+1}s_{i+2}\do s_{j-1}s_{j+1}\do s_q,\endalign$$
which shows that $l(w)\le q-2$, contradiction. This proves (a). 

For $(\e,t)\in R$ we have (see 1.5) $U(w\i)(\e,t)=(\e\et(w,t),w\i tw)$ where 
$\et(w,t)=\pm 1$ depends only on $w,t$. On the other hand, 
$$\align U(w\i)(\e,t)&=U_{s_q}\do U_{s_1}(\e,t)\\&=(\e(-1)^{\d_{s_1,t}+
\d_{s_2,s_1ts_1}+\do+\d_{s_q,s_{q-1}\do s_1ts_1\do s_{q-1}}},w\i tw)\\&
=(\e(-1)^{\d_{s_1,t}+\d_{s_1s_2s_1,t}+\do+\d_{s_1\do s_q\do s_1,t}},w\i tw).
\endalign$$
Thus, $\et(w,t)=(-1)^{\d_{s_1,t}+\d_{s_1s_2s_1,t}+\do+\d_{s_1\do s_q\do s_1,t}}$. 
Using (a), we see that for $t\in T$, the sum 
$\d_{s_1,t}+\d_{s_1s_2s_1,t}+\do+\d_{s_1\do s_q\do s_1,t}$ is $1$ if $t$ belongs to
the subset in (b) and is $0$, otherwise. Hence the subset in (b) is just 
$\{t\in T;\et(w,t)=-1\}$. This completes the proof.

\proclaim{Proposition 1.7} Let $w\in W,s\in S$ be such that $l(sw)=l(w)-1$. Let
$w=s_1s_2\do s_q$ be a reduced expression. Then there exists $j\in[1,q]$ such that 

$ss_1s_2\do s_{j-1}=s_1s_2\do s_j$.
\endproclaim
Let $w'=sw$. Let $w'=s'_1s'_2\do s'_{q-1}$ be a reduced expression. Then
$w=ss'_1s'_2\do s'_{q-1}$ is another reduced expression. By 1.6(b), the $q$-term 
sequences 
$$s_1,s_1s_2s_1,s_1s_2s_3s_2s_1,\do \text{ and } s,ss'_1s,ss'_1s'_2s'_1s,\do$$
coincide up to rearranging terms. In particular, $s=s_1s_2\do s_j\do s_2s_1$ for
some $j\in[1,q]$. The proposition follows.

\subhead 1.8 \endsubhead
Let $X$ be the set of all sequences $(s_1,s_2,\do,s_q)$ in $S$ such that
$s_1s_2\do s_q$ is a reduced expression in $W$. We regard $X$ as the vertices of a 
graph in which $(s_1,s_2,\do,s_q),(s'_1,s'_2,\do,s'_{q'})$ are joined if one is 
obtained from the other by replacing $m$ consecutive entries of form $s,s',s,s',\do$
by the $m$ entries $s',s,s',s,\do$; here $s\ne s'$ in $S$ are such that 
$m=m_{s,s'}<\iy$. We use the notation 
$$(s_1,s_2,\do,s_q)\si(s'_1,s'_2,\do,s'_{q'})$$
for "$(s_1,s_2,\do,s_q),(s'_1,s'_2,\do,s'_{q'})$ are in the same connected component
of $X$". (When this holds we have necessarily $q=q'$ and 
$s_1s_2\do s_q=s'_1s'_2\do s'_q$ in $W$.)

\medpagebreak

The following result is due to Matsumoto and Tits.

\proclaim{Theorem 1.9} Let $\ss=(s_1,s_2,\do,s_q),\ss'=(s'_1,s'_2,\do,s'_q)$ in $X$
be such that $s_1s_2\do s_q=s'_1s'_2\do s'_q=w\in W$. Then $\ss\si\ss'$.
\endproclaim
Let $C$ (resp. $C'$) be the connected component of $X$ that contains $\ss$ (resp. 
$\ss'$). For $i\in[1,q]$ we set

$\ss(i)=(\do,s'_1,s_1,s'_1,s_1,s_2,s_3,\do,s_i)$ (a $q$-element sequence in 
$S$),
 
$\uss(i)=\do s'_1s_1s'_1s_1s_2s_3\do s_i\in W$ (the product of this sequence).
\nl
Let $C_i$ be the connected component of $X$ that contains $\ss(i)$. Then 
$\ss=\ss(q)$. Hence $C=C_q$.

We argue by induction on $q$. The theorem is obvious for $q\in\{0,1\}$. We now 
assume that $q\ge 2$ and that the theorem is known for $q-1$. We first prove the 
following weaker statement.

$(A)$ {\it In the setup of the theorem we have either $\ss\si\ss'$ or}
$$s_1s_2\do s_q=s'_1s_1s_2\do s_{q-1} \text{ and }
(s'_1,s_1,s_2,\do,s_{q-1})\si(s'_1,s'_2,\do,s'_q).\tag a$$
We have $l(s'_1w)=l(w)-1$. By 1.7 we have $s'_1s_1s_2\do s_{i-1}=s_1s_2\do s_i$ for
some $i\in[1,q]$, so that $w=s'_1s_1s_2\do s_{i-1}s_{i+1}\do s_q$. In particular, 
$$(s'_1,s_1,s_2,\do,s_{i-1},s_{i+1},\do,s_q)\in X.$$
By the induction hypothesis, 
$(s_1,s_2,\do,s_{i-1},s_{i+1},\do,s_q)\si(s'_2,\do,s'_q)$. Hence 

$$(s'_1,s_1,s_2,\do,s_{i-1},s_{i+1},\do,s_q)\si(s'_1,s'_2,\do,s'_q).\tag b$$
Assume first that $i<q$. Then from 
$s'_1s_1s_2\do s_{i-1}s_{i+1}\do s_{q-1}=s_1s_2\do s_{q-1}$ and the induction 
hypothesis we deduce that 

$(s'_1,s_1,s_2,\do,s_{i-1},s_{i+1},\do,s_{q-1})\si(s_1,s_2,\do,s_{q-1})$, hence

$(s'_1,s_1,s_2,\do,s_{i-1},s_{i+1},\do,s_{q-1},s_q)\in C$.
\nl
Combining this with (b) we deduce that $C=C'$.

Assume next that $i=q$ so that $s_1s_2\do s_q=s'_1s_1s_2\do s_{q-1}$. Then (b) shows
that (a) holds. Thus, $(A)$ is proved.

Next we prove for $p\in[0,q-2]$ the following generalization of $(A)$.

$(A'_p)$ {\it In the setup of the theorem we have either $C=C'$ or:

for $i\in[q-p-1,q]$ we have $\ss(i)\in X,\uss(i)=w$, $C_i=C$ if $i-q\in 2\bz$,
$C_i=C'$ if $i-q\n 2\bz$.}
\nl
For $p=0$ this reduces to $(A)$. Assume now that $p>0$ and that $(A'_{p-1})$ is
already known. We prove that $(A'_p)$ holds.

If $C=C'$, then we are done. Hence by $(A'_{p-1})$ we may assume that: for 
$i\in[q-p,q]$ we have $\ss(i)\in X,\uss(i)=w$, $C_i=C$ if $i-q\in 2\bz$, $C_i=C'$ if
$i-q\n 2\bz$.

Applying $(A)$ to $\ss(q-p),\ss(q-p+1)$ (instead of $\ss,\ss'$), we see that either
$C_{q-p}=C_{q-p+1}$ or:

$\ss(q-p),\ss(q-p-1)$ are in $X$, $\uss(q-p)=\uss(q-p-1)$ and $C_{q-p-1}=C_{q-p+1}$.
\nl
In both cases, we see that $(A'_p)$ holds.

This completes the inductive proof of $(A'_p)$. In particular, $(A'_{q-2})$ holds. 
In other words, in the setup of the theorem, either $C=C'$ holds or:

(c) {\it for $i\in[1,q]$ we have $\ss(i)\in X,\uss(i)=w$, $C_i=C$ if $i-q\in 2\bz$,
$C_i=C'$ if $i-q\n 2\bz$.}
\nl
If $C=C'$, then we are done. Hence we may assume that (c) holds. In particular,
$$\ss(2)\in X,\ss(1)\in X,\uss(2)=\uss(1).\tag d$$
From $\ss(1)\in X$ and $q\ge 2$ we see that $s'_1\ne s_1$ and that $q\le m=m_{s_1,s'_1}$. 
From $\uss(2)=\uss(1)$ we see that $s_2\in\la s_1,s'_1\ra$, hence $s_2$ is either
$s_1$ or $s'_1$. In fact we cannot have $s_2=s_1$ since this would contradict
$\ss(2)\in X$. Hence $s_2=s'_1$. We see that $\ss(2)=(\do,s'_1,s_1,s'_1,s_1,s'_1)$ 
(the number of terms is $q,q\le m$). Since $\uss(2)=\uss(1)$, it follows that $q=m$,
so that $\ss(2),\ss(1)$ are joined in $X$. It follows that $C_2=C_1$. By (c), for 
some permutation $a,b$ of $1,2$ we have $C_a=C,C_b=C'$. Since $C_a=C_b$ it follows 
that $C=C'$. The theorem is proved.

\proclaim{Proposition 1.10} Let $w\in W$ and let $s,t\in S$ be such that 
$l(swt)=l(w),l(sw)=l(wt)$. Then $sw=wt$.
\endproclaim
Let $w=s_1s_2\do s_q$ be a reduced expression.

Assume first that $l(wt)=q+1$. Then $s_1s_2\do s_qt$ is a reduced expression for 
$wt$. Now $l(swt)=l(wt)-1$ hence by 1.7 there exists $i\in[1,q]$ such that
$ss_1s_2\do s_{i-1}=s_1s_2\do s_i$ or else $ss_1s_2\do s_q=s_1s_2\do s_qt$. If the 
second alternative occurs, we are done. If the first alternative occurs, we have
$sw=s_1s_2\do s_{i-1}s_{i+1}\do s_q$ hence $l(sw)\le q-1$. This contradicts
$l(sw)=l(wt)$.

Assume next that $l(wt)=q-1$. Let $w'=wt$. Then $l(sw't)=l(w'),l(sw')=l(w't)$. We 
have $l(w't)=l(w')+1$ hence the first part of the proof applies and gives $sw'=w't$.
Hence $sw=wt$. The proposition is proved.

\subhead 1.11\endsubhead
We can regard $S$ as the set of vertices of a graph in which $s,s'$ are joined if 
$m_{s,s'}>2$. We say that $W$ is {\it irreducible} if this graph is connected. It is
easy to see that in general, $W$ is naturally a product of irreducible Coxeter 
groups, corresponding to the connected components of $S$.

In the setup of 1.3, let $(,):E\T E\to\br$ be the symmetric $\br$-bilinear form 
given by $(e_s,e_{s'})=-\cos\fra{\p}{m_{s,s'}}$. Then $\s(w):E@>>>E$ preserves $(,)$
for any $w\in W$. 

We say that $W$ is {\it tame} if $(e,e)\ge 0$ for any $e\in E$. It is easy to see 
that, if $W$ is finite then $W$ is tame. 

We say that $W$ is {\it integral} if, for any $s\ne s'$ in $S$, we have
$4\cos^2\fra{\p}{m_{s,s'}}\in\bn$ (or equivalently $m_{s,s'}\in\{2,3,4,6,\iy\}$).

We will be mainly interested in the case where $W$ is tame. The tame, irreducible 
$W$ are of three kinds:

(a) finite, integral;

(b) finite, non-integral;

(c) tame, infinite (and automatically integral).

\subhead 1.12\endsubhead
For $k\in\bz$ define $\r_k:\bz@>>>\bz$ by $\r_k(z)=z+k$. Let $n\ge 2$. Let 
$\tW$ be 
the group of all permutations $\s:\bz@>>>\bz$ such that $\s\r_n=\r_n\s$. Define 
$\c:\tW@>>>\bz$ by $\c(\s)=\su_{k\in X}(\s(k)-k)$ where $X$ is a set of 
representatives for the residue classes $\mod n$ in $\bz$. One checks that $\c$ does
not depend on the choice of $X$ and $\c$ is a group homomorphism with image $n\bz$.
Now $\tW'=\ker(\c)$ is generated by $\{s_m;m\in\bz/n\bz\}$ where 
$s_m:\bz@>>>\bz$ is defined by

$s_m(z)=z+1$ if $z=m\mod n$, 

$s_m(z)=z-1$ if $z=m+1\mod n$, 

$s_m(z)=z$ for all other $z\in\bz$.
\nl
It is a Coxeter group on these generators, said to be of type $\tA_{n-1}$. (This 
description of $\tW'$ appears in \cite{\SQ}.) For $n\ge 3$, $m,m'\in\bz/n\bz$ are
joined by a single edge in the Coxeter graph if $m-m'=1\mod n$ and are not joined
otherwise. For $n=2$, $0,1\in\bz/2\bz$ are joined by a quadruple edge in
the Coxeter graph. The length function on $\tW'$ is given by
$$l(\s)=\sh(Y_\s/\t_n)$$
where, for $\s\in\tW'$, 
$$Y_\s=\{(i,j)\in\bz\T\bz;i<j,\s(i)>\s(j)\}$$
and $Y_\s/\t_n$ is the (finite) set of orbits of $\t_n:Y_\s@>>>Y_\s$, 
$(i,j)\m(i+n,j+n)$.

\subhead 1.13\endsubhead
Assume now that $n=2p\ge 4$, where $p\in\bn$. Let $W$ be the subgroup of $\tW$
consisting of all $\s\in\tW$ that commute with the involution $\bz@>>>\bz,z\m 1-z$.
We compute $\c(\s)$ for $\s\in W$, taking $X=\{-(p-1),\do,-1,0,1,2,\do,p\}$:
$$\align&\c(\s)=\su_{k\in[1,p]}(\s(k)-k)+\su_{k\in[1,p]}(\s(1-k)-(1-k))
\\&=\su_{k\in[1,p]}(\s(k)-k)+\su_{k\in[1,p]}(1-\s(k)-(1-k))=0.\endalign$$
Thus, $W$ is a subgroup of $\tW'$. Now $W$ is generated by $s'_0,s'_1,\do,s'_p$ 
where 
$$s'_0=s_0,s'_1=s_1s_{-1},s'_2=s_2s_{-2},\do,s'_{p-1}=s_{p-1}s_{1-p},s'_p=s_p.$$
It is a Coxeter group on these generators, said to be of type $\tC_p$. The Coxeter 
graph is 
$$\bul=\bul\h\bul\h\do\h\bul=\bul$$
with vertices corresponding to $0,1,2,\do,p-1,p$.

Let $\s\in W$. We have a partition $Y_\s=Y^0_\s\sqc Y^1_\s$ where
$$Y^0_\s=\{(i,j)\in\bz\T\bz;i<j,\s(i)>\s(j),i+j\ne 1\mod 2p\},$$
$$Y^1_\s=\{(i,j)\in\bz\T\bz;i<j,\s(i)>\s(j),i+j=1\mod 2p\}.$$
Now $Y^1_\s/\t_n$ is the fixed point set of the involution of $Y_\s/\t_n$ induced by
the involution $(i,j)\m(1-j,1-i)$ of $Y_\s$. Hence we have 
$\sh(Y_\s/\t_n)=2l^0(\s)+l^1(\s)$ where $l^0(\s)=\sh(Y^0_\s/\t_n)/2$,
$l^1(\s)=\sh(Y^1_\s/\t_n)$ are integers. Now $(i,j)\m(i,\fra{i+j-1}{2p})$ is a 
bijection of $Y^1_\s$ with 
$$\align%\in\bz\T\bz;i<1-i+2ph,\s(i)>\s(1-i+2ph)\}\\&
\{(i,h)\in\bz\T\bz;2i<1+2ph,2\s(i)>1+2ph\}=\{(i,h)\in\bz\T\bz;i\le ph<\s(i)\}.
\endalign$$
It follows that 

$l^1(\s)=\su_{i\in[1-p,p];i<\s(i)}f(i)$ 
\nl
where 

$f(i)=\sh(x\in p\bz; i\le x<\s(i))$. 
\nl
Let $\bz'=[1,p]+2p\bz,\bz''=[1-p,0]+2p\bz$; then 
$\bz=\bz'\sqc\bz''$. We have $l^1(\s)=l'(\s)+l''(\s)$ where
$$\align&l'(\s)=\su_{i\in[1-p,0];\s(i)\in\bz'',i<\s(i)}\fra{f(i)}{2}+ 
\su_{i\in[1,p];\s(i)\in\bz',i<\s(i)}\fra{f(i)}{2}\\&+ 
\su_{i\in[1-p,0];\s(i)\in\bz',i<\s(i)}\fra{f(i)+1}{2}+ 
\su_{i\in[1,p];\s(i)\in\bz'',i<\s(i)}\fra{f(i)-1}{2},\endalign$$
$$\align&l''(\s)=\su_{i\in[1-p,0];\s(i)\in\bz'',i<\s(i)}\fra{f(i)}{2}+ 
\su_{i\in[1,p];\s(i)\in\bz',i<\s(i)}\fra{f(i)}{2}\\&+ 
\su_{i\in[1-p,0];\s(i)\in\bz',i<\s(i)}\fra{f(i)-1}{2}+ 
\su_{i\in[1,p];\s(i)\in\bz'',i<\s(i)}\fra{f(i)+1}{2}.\endalign$$
These are integers since $f(i)$ is even if $i,\s(i)$ are in the same set $\bz'$ or
$\bz''$ and is odd 
otherwise. We see that the length of $\s$ in $\ti W'$ is $2l^0(\s)+l'(\s)+l''(\s)$.
On the other hand, the length of $\s$ in $W$ is 
$$l^0(\s)+l'(\s)+l''(\s).$$
Now $l'(\s)$ (resp. $l''(\s)$) is the number of times that $s'_0$ 
(resp. $s'_p$) appears in a reduced expression of $\s$ in $W$.

One can show that $l^0,l',l''$ are weight functions on $W$ in the sense of 3.1.

Define $\c':W@>>>\{\pm 1\}$ and $\c'':W@>>>\{\pm 1\}$ 
by $\c'(\s)=(-1)^{l'(\s)}$,
$\c''(\s)=(-1)^{l''(\s)}$. Then $\c',\c''$ are group homomorphisms.

Assuming that $p\ge 3$, let $W'=\ker(\c')$. This is the subgroup of $W$ generated by
$$s'_0s'_1s'_0,s'_1,s'_2,\do,s'_{p-1},s'_p.$$
It is a Coxeter group on these generators, said to be of type $\tB_p$. The Coxeter 
graph has vertices $\ti 1,1,2,\do,(p-1),p$ corresponding to 
$s'_0s'_1s'_0,s'_1,s'_2,\do,s'_{p-1},s'_p$ and edges 
$$\matrix \bul&\h&\bul&\h&\bul&\h&\do&\h&\bul&=&\bul&\\
            {}  & {}&\vert&{}&  {}  &{}&{} &{}&       {} &{}&{}   &\\
            {}  & {}&\bul&{}&{}&{} &{}&       {} &{}&{}   &\endmatrix$$
Assuming that $p\ge 4$, let $W''=\ker(\c')\cap\ker(\c'')$. This is the subgroup of 
$W$ (or $W'$) generated by 
$$s'_0s'_1s'_0,s'_1,s'_2,\do,s'_{p-1},s'_ps'_{p-1}s'_p.$$
It is a Coxeter group on these generators, said to be of type $\tD_p$. The Coxeter 
graph has vertices $\ti 1,1,2,\do,(p-1),\wt{(p-1)}$ and edges 
$$\matrix \bul&\h&\bul&\h&\bul&\h&\do&\h&\bul&\h&\bul&\\
            {}  &{}&\vert &{}&   {} & {}&{}& {}&  \vert  &{}&{}       &\\
            {}  &{}&\bul&{}&{}&{}&{}&{}&\bul&{}&{}  &\endmatrix
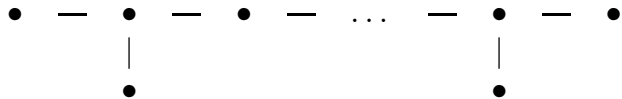
$$

\subhead 1.14\endsubhead
Let $p\ge q\ge r\ge 1$ be integers such that $p\i+q\i+r\i=1$. Then $p,q,r$ is 
$3,3,3$ or $4,4,2$ or $6,3,2$. Thus, $q$ and $r$ divide $p$. Consider the graph with
vertices
$$\{(\fra{ip}{p},0,0);i\in[1,p]\}\cup\{(0,\fra{ip}{q},0);i\in[1,q]\}
\cup\{(0,0,\fra{ip}{r});i\in[1,r]\}$$
where $(p,0,0),(0,p,0),(0,0,p)$ are identified; the edges are

$(\fra{ip}{p},0,0)\h(\fra{(i+1)p}{p},0,0),1\le i<p$, 

$(0,\fra{ip}{q},0)\h(0,\fra{(i+1)p}{q},0),1\le i<q$,

$(0,0,\fra{ip}{r})\h(0,0,\fra{(i+1)p}{r}),1\le i<r$.
\nl
The Coxeter group $W$ corresponding to this graph is said to be of type $\tE_n$ 
where $n=p+q+r-3$. Thus, $n\in\{6,7,8\}$.

Let $W$ be of type $\tE_6$. Let $W'$ be the subgroup of $W$ generated by 
$$s_{1,0,0},s_{2,0,0},s_{3,0,0},s_{0,2,0}s_{0,0,2},s_{0,1,0}s_{0,0,1}.$$
(The index of a generator of $W$ is the corresponding vertex of the graph.) Then 
$W'$ is a Coxeter group on these generators, said to be of type $\ti F_4$. The 
Coxeter graph is 
$$\bul\h\bul\h\bul=\bul\h\bul.$$
Let $W$ be of type $\tD_4$. The standard generators may be denoted by 
$s_0,s_1,s_2,s_3,s_4$ where the Coxeter graph has vertices $0,1,2,3,4$ with four
edges joining $0$ with $1,2,3,4$. Let $W'$ be the subgroup of $W$ generated by 
$s_1,s_0,s_2s_3s_4$. Then $W'$ is a Coxeter group on these generators, said to be of
type $\ti G_2$. The Coxeter graph is 
$$\bul\h\bul\equiv\bul.$$

\subhead 1.15\endsubhead
The collection of Coxeter groups of type $\tA_{n-1}(n\ge 2)$, $\tD_n(n\ge 4)$,
$\tC_n(n\ge 2)$, $\tB_n(n\ge 3)$, $\tE_n(n=6,7,8)$, $\ti F_4,\ti G_2$ (see 
1.12-1.14) coincides with the collection of infinite, tame, irreducible Coxeter 
groups (or {\it affine Weyl groups}).

\subhead 1.16\endsubhead
Let $W,S$ be an affine Weyl group. Let $\ct$ be the union of all finite conjugacy 
classes in $W$. Then $\ct$ is a normal, finitely generated free abelian subgroup of
$W$ of finite index.  Let $S_{min}$ be the set of all $s\in S$ such that the obvious
composition $\la S-\{s\}\ra@>>>W@>>>W/\ct$ is an isomorphism. (This composition is 
injective for any $s\in S$.) Now $S_{min}\ne\em$. We describe $S_{min}$ in each 
case.

If $W$ is of type $\tA_{n-1}$ we have $S_{min}=S$. 
In the setup of 1.13, $S_{min}$ corresponds to the following vertices of
the Coxeter graph: $1,p$, if $W$ is of type $\tC_p$; $1,\ti 1$, if $W$ is of 
type $\tB_p$; $1,\ti 1,(p-1),\wt{(p-1)}$, if $W$ is of type $\tD_p$.
In the setup of 1.14, if $W$ is of type $\tE_n$ then $S_{min}$ corresponds
to the following vertices of the Coxeter graph:
$(1,0,0),(0,1,0),(0,0,1)$ if $W$ is of type $\tE_6$;
$(1,0,0),(0,1,0)$ if $W$ is of type $\tE_7$;
$(1,0,0)$ if $W$ is of type $\tE_8$.
If $W$ is of type $\ti F_4$ then $S_{min}=\{s_{1,0,0}\}$; if 
$W$ is of type $\ti G_2$ then $S_{min}=\{s_1\}$.

\subhead 1.17\endsubhead
For a Coxeter group $W,S$ we denote by $A_W$ the group of all automorphisms of $W$ 
that map $S$ into itself. (This is also the group of automorphisms of the
corresponding Coxeter graph.) 

\subhead 1.18\endsubhead
Let $W,S$ be an affine Weyl group. Let $\ct\sub W$ be as in 1.16. Let $\Om$ be the
set of all $a\in A_W$ such that there exists $w\in W$ with $a(t)=wtw\i$ for all 
$t\in\ct$. Now $\Om$ is a commutative normal subgroup of $A_W$. The action of $A_W$
on $W$ restricts to an action of $\Om$ on $S_{min}$ which is simply transitive. 

\subhead 1.19\endsubhead
For any $I\sub S$, let $W_I=\la I\ra$. Then $(W_I,I)$ is a Coxeter group whose 
Coxeter matrix is a submatrix of that of $W,S$. See \S9.

\subhead 1.20\endsubhead
Let $W,S$ be an affine Weyl group. Let $s\in S_{min}$. Then $W_{S-\{s\}},S-\{s\}$ is
a finite Coxeter group.

A finite Coxeter group is said to be of type $A_{n-1}(n\ge 2)$ (resp. $C_n(n\ge 2)$,
$B_n(n\ge 3)$, $D_n(n\ge 4)$, $E_n(n=6,7,8)$, $F_4,G_2$) if it is isomorphic to
$W_{S-\{s\}}$ for some $W,S,s$ as above, where $W$ has type $\tA_{n-1}(n\ge 2)$
(resp. $\tC_n(n\ge 2)$, $\tB_n(n\ge 3)$, $\tD_n(n\ge 4)$, $\tE_n(n=6,7,8)$, 
$\ti F_4,\ti G_2$). 

The collection of Coxeter groups of type $A_{n-1}(n\ge 2)$, $C_n(n\ge 2)$, 
$B_n(n\ge 3)$, $D_n(n\ge 4)$, $E_n(n=6,7,8)$, $F_4,G_2$ coincides with the 
collection of finite, integral, irreducible Coxeter groups $\ne\{1\}$ (or {\it Weyl 
groups}). The group $W=\{1\}$ with $S=\em$ is also considered to be a Weyl group.
Note that the types $C_n$ and $B_n$ coincide for $n\ge 3$.

For a Weyl group $W,S$ we set $n(W)=2\sh(T)/\sh(S)^2$ where $T$ is as in 1.4. We list below the numbers 
$n(W)$ for $W$ of various types:

$A_{n-1}:$  $n(W)=1+\frac{1}{n-1}$

$D_n:$      $n(W)=2-\fra{2}{n}$

$B_n:$      $n(W)=2$

$E_6:$      $n(W)=2$

$E_7:$      $n(W)=2.57\do$

$F_4:$      $n(W)=3$

$G_2:$      $n(W)=3$

$E_8:$      $n(W)=3.75$.
\nl
We see that the maximum value of $n(W)$ is achieved in type $E_8$.

\subhead 1.21\endsubhead
Let $W,S$ be a Weyl group of type $E_8$. Let $W'$ be the subgroup of $W$ generated 
by 
$$s_{2,0,0}s_{0,2,0},s_{3,0,0}s_{0,4,0},s_{4,0,0}s_{6,0,0},s_{5,0,0}s_{0,0,3}.$$
(The index of a generator of $W$ is the corresponding vertex of the graph, see 
1.13.) Then $W'$ is a (finite, non-integral) Coxeter group on these generators, said
to be of type $H_4$. This description of $H_4$ appeared in \cite{\SQ}.

\head 2. Partial order on $W$\endhead
\subhead 2.1\endsubhead
Let $W,S$ be a Coxeter group. 
Let $y,w$ be two elements of $W$. Following Chevalley, we say that $y\le w$ if there
exists a sequence $y=y_0,y_1,y_2,\do,y_n=w$ in $W$ such that $l(y_k)-l(y_{k-1})=1$ 
for $k\in[1,n]$ and $y_ky_{k-1}\i\in T$ (or equivalently $y_{k-1}y_k\i\in T$, or 
$y_k\i y_{k-1}\in T$, or $y_{k-1}\i y_k\in T$) for $k\in[1,n]$. This is a partial 
order on $W$. Note that $y\le w$ implies $y\i\le w\i$ and $l(y)\le l(w)$. 
We write $y<w$ or $w>y$ instead of $y\le w, y\ne w$. If $w\in W,s\in S$ then, 

$sw<w$ if and only if $l(sw)=l(w)-1$;

$sw>w$ if and only if $l(sw)=l(w)+1$.

\proclaim{Lemma 2.2} Let $w=s_1s_2\do s_q$ be a reduced expression in $W$ and let
$t\in T$. The following are equivalent:

(i) $U(w\i)(\e,t)=(-\e,w\i tw)$ for $\e=\pm 1$;

(ii) $t=s_1s_2\do s_i\do s_2s_1$ for some $i\in[1,q]$;

(iii) $l(tw)<l(w)$.
\endproclaim
The equivalence of (i),(ii) has been proved in 1.6.

Proof of (ii)$\imp$(iii). Assume that (ii) holds. Then 
$tw=s_1\do s_{i-1}s_{i+1}\do s_q$ hence $l(tw)<q$ and (iii) holds.

Proof of (iii)$\imp$(i). First we check that 
$$U(t)(\e,t)=(-\e,t).\tag a$$
If $t\in S$, (a) is clear. If (a) is true for $t$ then it is also true for $sts$ 
where $s\in S$. Indeed, 
$$\align U(sts)(\e,sts)&=U_sU(t)U_s(\e,sts)=U_sU(t)(\e(-1)^{\d_{s,sts}},t)
=U_s(-\e(-1)^{\d_{s,sts}},t)\\&=(-\e(-1)^{\d_{s,sts}+\d_{s,t}},sts)=(-\e,sts);\endalign
$$
(a) follows. Assume now that (i) does not hold; thus, $U(w\i)(\e,t)=(\e,w\i tw)$.
Then
$$\align U((tw)\i)(\e,t)&=U(w\i)U(t)(\e,t)=U(w\i)(-\e,t)=(-\e,w\i tw)\\&=
(-\e,(tw)\i t(tw)).\endalign$$
Since (i)$\imp$(iii) we deduce that $l(w)<l(tw)$; thus, (iii) does not hold. The
lemma is proved.

\proclaim{Lemma 2.3} Let $y,z\in W$ and let $s\in S$. If $sy\le z<sz$, then
$y\le sz$.
\endproclaim
We argue by induction on $l(z)-l(sy)$. If 
$l(z)-l(sy)=0$ then $z=sy$ and the result is clear. Now assume that $l(z)>l(sy)$. 
Then $sy<z$. We can assume that $sy<y$ (otherwise the result is trivial). We can 
find $t\in T$ such that $sy<tsy\le z$ and $l(tsy)=l(sy)+1$. If $t=s$, then $y\le z$
and we are done. Hence we may assume that $t\ne s$. We show that 
$$y<stsy.\tag a$$
Assume that (a) does not hold. Then $y,tsy,sy,stsy$ have lengths $q+1,q+1,q,q$. We
can find a reduced expression $y=ss_1s_2\do s_q$. Since $l(stsy)<l(y)$, we see from
2.2 that either $sts=ss_1\do s_i\do s_1s$ for some $i\in[1,q]$ or $sts=s$. (This 
last case has been excluded.) It follows that
$$tsy=s_1\do s_i\do s_1sss_1s_2\do s_q=s_1\do s_{i-1}s_{i+1}\do s_q.$$
Thus, $l(tsy)\le q-1$, a contradiction. Thus, (a) holds. Let $y'=stsy$. We have
$sy'\le z<sz$ and $l(z)-l(sy')<l(z)-l(sy)$. By the induction hypothesis, 
$y'\le sz$. We have $y<y'$ by (a), hence $y\le sz$. The lemma is proved.

\proclaim{Proposition 2.4} The following three conditions on $y,w\in W$ are
equivalent:

(i) $y\le w$;

(ii) for any reduced expression $w=s_1s_2\do s_q$ there exists a subsequence
$i_1<i_2<\do<i_r$ of $1,2,\do,q$ such that $y=s_{i_1}s_{i_2}\do s_{i_r}$, $r=l(y)$;

(iii) there exists a reduced expression $w=s_1s_2\do s_q$ and a subsequence
$i_1<i_2<\do<i_r$ of $1,2,\do,q$ such that $y=s_{i_1}s_{i_2}\do s_{i_r}$.
\endproclaim
Proof of (i)$\imp$(ii). We may assume that $y<w$. Let $y=y_0,y_1,y_2,\do,y_n=w$ be
as in 2.1. Let $w=s_1s_2\do s_q$ be a reduced expression. Since $y_{n-1}y_n\i\in T$,
$l(y_{n-1})=l(y_n)-1$, we see from 2.2 that there exists $i\in[1,q]$ such that 
$y_{n-1}y_n\i=s_1s_2\do s_i\do s_2s_1$ hence
$y_{n-1}=s_1s_2\do s_{i-1}s_{i+1}\do s_q$ (a reduced expression). Similarly, since 
$y_{n-2}y_{n-1}\i\in T$, $l(y_{n-2})=l(y_{n-1})-1$, we see from 2.2 (applied to 
$y_{n-1}$) that there exists $j\in[1,q]-\{i\}$ such that $y_{n-2}$ equals
$$s_1s_2\do s_{i-1}s_{i+1}\do s_{j-1}s_{j+1}\do s_q \text{ or }
s_1s_2\do s_{j-1}s_{j+1}\do s_{i-1}s_{i+1}\do s_q$$
(depending on whether $i<j$ or $i>j$). Continuing in this way we see that $y$ is of
the required form.

The proof of (ii)$\imp$(iii) is trivial.

Proof of (iii)$\imp$(i). Assume that $w=s_1s_2\do s_q$ (reduced expression) and
$y=s_{i_1}s_{i_2}\do s_{i_r}$ where $i_1<i_2<\do<i_r$ is a subsequence of 
$1,2,\do,q$. We argue by induction on $q$. If $q=0$ there is nothing to prove. Now 
assume $q>0$. 

If $i_1>1$, then the induction hypothesis is applicable to $y,w'=s_2\do s_q$ and 
yields $y\le w'$. But $w'\le w$ hence $y\le w$. If $i_1=1$ then the induction 
hypothesis is applicable to $y'=s_{i_2}\do s_{i_r},w'=s_2\do s_q$ and yields 
$y'\le w'$. Thus, $s_1y\le s_1w<w$. By 2.3 we then have $y\le w$. The proposition is
proved.

\proclaim{Corollary 2.5} Let $y,z\in W$ and let $s\in S$. 

(a) Assume that $sz<z$. Then $y\le z\Lra sy\le z$.

(b) Assume that $y<sy$. Then $y\le z\Lra y\le sz$.
\endproclaim
We prove (a). We can find a reduced expression of $z$ of form $z=ss_1s_2\do s_q$. 
Assume that $y\le z$. By 2.4 we can find a subsequence $i_1<i_2<\do<i_r$ of 
$1,2,\do,q$ such that either $y=s_{i_1}s_{i_2}\do s_{i_r}$ or 
$y=ss_{i_1}s_{i_2}\do s_{i_r}$. In the first case we have 
$sy=ss_{i_1}s_{i_2}\do s_{i_r}$ and in the second case we have 
$sy=s_{i_1}s_{i_2}\do s_{i_r}$. In both cases, $sy\le z$ by 2.4. The same argument 
shows that, if $sy\le z$ then $y\le z$. This proves (a).

We prove (b). Assume that $y\le z$. We must prove that $y\le sz$. If $z<sz$, this is
clear. Thus we may assume that $sz<z$. We can find a reduced expression of $z$ of 
form $z=ss_1s_2\do s_q$. By 2.4 we can find a subsequence $i_1<i_2<\do<i_r$ of 
$1,2,\do,q$ such that either
$$y=s_{i_1}s_{i_2}\do s_{i_r},l(y)=r \text{ or }
y=ss_{i_1}s_{i_2}\do s_{i_r},l(y)=r+1.$$
In the second case we have $l(sy)=r<l(y)$, contradicting $y<sy$. Thus we are in the
first case. Hence $y$ is the product of a subsequence of $s_1,s_2,\do,s_q$ and using
again 2.4, we deduce that $y\le sz$ (note that $sz=s_1s_2\do s_q$ is a reduced 
expression). The lemma is proved.

\head 3. The algebra $\ch$\endhead
\subhead 3.1\endsubhead
Let $W,S$ be a Coxeter group. 
A map $L:W@>>>\bz$ is said to be a {\it weight function} for $W$ if 
$L(ww')=L(w)+L(w')$ for any $w,w'\in W$ such that $l(ww')=l(w)+l(w')$. We will 
assume that a weight function $L:W@>>>\bz$ is fixed; we then say that $W,L$ is a
{\it weighted Coxeter group}. (For example we could take $L=l$; in that case we say
that we are in the {\it split case}.) Note that $L$ is determined by its values on 
$S$ which are subject only to the condition

$L(s)=L(s')$ for any $s\ne s'$ in $S$ such that $m_{s,s'}$ is finite and odd.
\nl
We have $L(w)=L(w\i)$ for all $w\in W$.

Let $\ca=\bz[v,v\i]$ where $v$ is an indeterminate. For $s\in S$ we set 
$v_s=v^{L(s)}\in\ca$.

\subhead 3.2\endsubhead
Let $\ch$ be the $\ca$-algebra with $1$ defined by the generators $T_s(s\in S)$ and the 
relations
$$(T_s-v_s)(T_s+v_s\i)=0 \text{ for } s\in S \tag a$$
$$T_sT_{s'}T_s\do=T_{s'}T_sT_{s'}\do\tag b$$ 
(both products have $m_{s,s'}$ factors) for any $s\ne s'$ in $S$ such that 
$m_{s,s'}<\iy$; $\ch$ is called the {\it Iwahori-Hecke algebra}. 

For $w\in W$ we define $T_w\in\ch$ by $T_w=T_{s_1}T_{s_2}\do T_{s_q}$, where
$w=s_1s_2\do s_q$ is a reduced expression in $W$. By (b) and 1.9, $T_w$ is 
independent of the choice of reduced expression. We have $T_1=1$. From the definitions it is clear 
that for $s\in S,w\in W$ we have
$$T_sT_w=T_{sw} \text{ if } l(sw)=l(w)+1,$$
$$T_sT_w=T_{sw}+(v_s-v_s\i)T_w \text{ if } l(sw)=l(w)-1.$$
In particular, the $\ca$-submodule of $\ch$ generated by $\{T_w;w\in W\}$ is a left
ideal of $\ch$. It contains $1=T_1$ hence it is the whole of $\ch$. Thus 
$\{T_w;w\in W\}$ generates the $\ca$-module $\ch$.

\proclaim{Proposition 3.3} $\{T_w;w\in W\}$ is an $\ca$-basis of $\ch$.
\endproclaim
We follow the lines of the proof in \cite{\BO, Ex.23, p.55}.

We consider the free $\ca$-module $\ce$ with basis $(e_w)_{w\in W}$. For any 
$s\in S$ we define $\ca$-linear maps $P_s,Q_s:\ce@>>>\ce$ by

$P_s(e_w)=e_{sw}$ if $l(sw)=l(w)+1$,

$P_s(e_w)=e_{sw}+(v_s-v_s\i)e_w$ if $l(sw)=l(w)-1$;

$Q_s(e_w)=e_{ws}$ if $l(ws)=l(w)+1$,

$Q_s(e_w)=e_{ws}+(v_s-v_s\i)e_w$ if $l(ws)=l(w)-1$.
\nl
We shall continue the proof assuming that

(a) $P_sQ_t=Q_tP_s$ {\it for any $s,t$ in} $S$.
\nl
Let $\fA$ be the $\ca$-subalgebra with $1$ of $\End(\ce)$ generated by 
$\{P_s;s\in S\}$. The map $\fA@>>>\ce$ given by $\p\m\p(e_1)$ is surjective. Indeed,
if $w=s_1s_2\do s_q$ is a reduced expression in $W$, then 
$e_w=P_{s_1}\do P_{s_q}e_1$. Assume now that $\p\in\fA$ satisfies $\p(e_1)=0$. Let 
$\p'=Q_{s_q}\do Q_{s_1}$. By (a) we have $\p\p'=\p'\p$ hence
$$0=\p'\p(e_1)=\p\p'(e_1)=\p(Q_{s_q}\do Q_{s_1}(e_1))=\p(e_w).$$
Since $w$ is arbitrary, it follows that $\p=0$. We see that the map $\fA@>>>\ce$ is
injective, hence an isomorphism of $\ca$-modules. Using this isomorphism we 
transport the algebra structure of $\fA$ to an algebra structure on $\ce$ with unit
element $e_1$. For this algebra structure we have $P_s(e_1)\p(e_1)=P_s(\p(e_1))$ for
$s\in S,\p\in\fA$. Hence $e_se_w=P_s(e_w)$ for any $w\in W,s\in S$. It follows that

(b) $e_se_w=e_{sw}$ if $l(sw)=l(w)+1$, 

(c) $e_se_w=e_{sw}+(v_s-v_s\i)e_w$ if $l(sw)=l(w)-1$.
\nl
From (b) it follows that, if $w=s_1s_2\do s_q$ is a reduced expression, then
$e_w=e_{s_1}e_{s_2}\do e_{s_q}$. In particular, if $s\ne s'$ in $S$ are such that 
$m=m_{s,s'}<\iy$ then $e_se_{s'}e_s\do=e_{s'}e_se_{s'}\do$ (both products have $m$ 
factors); indeed, this follows from the equality $e_{ss's\do}=e_{s'ss'\do}$ (see 
1.4). From (c) we deduce that $e_s^2=1+(v_s-v_s\i)e_s$ for $s\in S$, or that 
$(e_s-v_s)(e_s+v_s\i)=0$. We see that there is a unique algebra homomorphism 
$\ch@>>>\ce$ preserving $1$ such that $T_s\m e_s$ for all $s\in S$. It takes $T_w$ 
to $e_w$ for any $w\in W$. Assume now that $a_w\in\ca$ ($w\in W$) are zero for all 
but finitely many $w$ and that $\su_wa_wT_w=0$ in $\ch$. Applying $\ch@>>>\ce$ we 
obtain $\su_wa_we_w=0$. Since $(e_w)$ is a basis of $\ce$, it follows that $a_w=0$
for all $w$. Thus, $\{T_w;w\in W\}$ is an $\ca$-basis of $\ch$. This completes the 
proof, modulo the verification of (a).

We prove (a). Let $w\in W$. We distinguish six cases. 
\nl
{\it Case} 1. $swt,sw,wt,w$ have lengths $q+2,q+1,q+1,q$. Then 
$$P_sQ_t(e_w)=Q_tP_s(e_w)=e_{swt}.$$
{\it Case} 2. $w,sw,wt,swt$ have lengths $q+2,q+1,q+1,q$. Then 
$$\align&P_sQ_t(e_w)=Q_tP_s(e_w)\\&=e_{swt}+(v_t-v_t\i)e_{sw}+(v_s-v_s\i)e_{wt}
+(v_t-v_t\i)(v_s-v_s\i)e_w.\endalign$$
{\it Case} 3. $wt,swt,w,sw$ have lengths $q+2,q+1,q+1,q$. Then 
$$P_sQ_t(e_w)=Q_tP_s(e_w)=e_{swt}+(v_s-v_s\i)e_{wt}.$$
{\it Case} 4. $sw,swt,w,wt$ have lengths $q+2,q+1,q+1,q$. Then 
$$P_sQ_t(e_w)=Q_tP_s(e_w)=e_{swt}+(v_t-v_t\i)e_{sw}.$$
{\it Case} 5. $swt,w,wt,sw$ have lengths $q+1,q+1,q,q$. Then 
$$P_sQ_t(e_w)=e_{swt}+(v_t-v_t\i)e_{sw}+(v_t-v_t\i)(v_s-v_s\i)e_w,$$
$$Q_tP_s(e_w)=e_{swt}+(v_s-v_s\i)e_{wt}+(v_t-v_t\i)(v_s-v_s\i)e_w.$$
{\it Case} 6. $sw,wt,w,swt$ have lengths $q+1,q+1,q,q$. Then 
$$P_sQ_t(e_w)=e_{swt}+(v_s-v_s\i)e_{wt},$$
$$Q_tP_s(e_w)=e_{swt}+(v_t-v_t\i)e_{sw}.$$
In cases 5, 6 we have $sw=wt$ by 1.10. 
In case 5 we have $L(t)+L(wt)=L(w)=L(swt)=L(s)+L(wt)$ hence $L(t)=L(s)$ and 
$v_s=v_t$. In case 6 we have $L(t)+L(swt)=L(sw)=L(wt)=L(s)+L(swt)$, hence 
$L(t)=L(s)$ and $v_s=v_t$. 
Hence 
$P_sQ_t(e_w)=Q_tP_s(e_w)$ in each case. The proposition is proved.

\subhead 3.4\endsubhead
There is a unique involutive antiautomorphism $h\m h^\flat$ of the algebra $\ch$ 
which carries $T_s$ to $T_s$ for any $s\in S$. It carries $T_w$ to $T_{w\i}$ for 
any $w\in W$. 

\subhead 3.5\endsubhead
For $s\in S$, the element $T_s\in\ch$ is invertible: $T_s\i=T_s-(v_s-v_s\i)$. It 
follows that $T_w$ is invertible for each $w\in W$; if 
$w=s_1s_2\do s_q$ is a reduced expression in $W$, then
$T_w\i=T_{s_q}\i\do T_{s_2}\i T_{s_1}\i$.

There is a unique algebra involution of $\ch$ denoted $h\m h^\da$ such that 
$T_s^\da=-T_s\i$ for any $s\in S$. We have $T_w^\da=\sg(w)T_{w\i}\i$ for any 
$w\in W$. 

\head 4. The bar operator \endhead
\subhead 4.1\endsubhead
We preserve the setup of 3.1.
Let $\bar{}:\ca@>>>\ca$ be the ring involution which takes $v^n$ to $v^{-n}$ for any
$n\in\bz$.

\proclaim{Lemma 4.2}(a) There is a unique ring homomorphism $\bar{}:\ch@>>>\ch$
which is $\ca$-semilinear with respect to $\bar{}:\ca@>>>\ca$ and satisfies
$\ov T_s=T_s\i$ for any $s\in S$.

(b) This homomorphism is involutive. It takes $T_w$ to $T_{w\i}\i$ for any $w\in W$.
\endproclaim
The following identities can be deduced easily from 3.2(a),(b),(d):

$(T_s\i-v_s\i)(T_s\i+v_s)=0$ for $s\in S$,

$T_s\i T_{s'}\i T_s\i\do=T_{s'}\i T_s\i T_{s'}\i\do$ 
\nl
(both products have $m_{s,s'}$ factors) for any $s\ne s'$ in $S$ such that 
$m_{s,s'}<\iy$; (a) follows.

We prove (b). Let $s\in S$. Applying $\bar{}$ to $T_s\bar{T_s}=1$ gives
$\bar{T_s}\bar{\bar{T_s}}=1$. We have also $\bar{T_s}T_s=1$ hence 
$\bar{\bar{T_s}}=T_s$. It follows that the square of $\bar{}$ is $1$. The second 
assertion of (b) is immediate. The lemma is proved.

\subhead 4.3\endsubhead
For any $w\in W$ we can write uniquely $\ov T_w=\su_{y\in W}\ov{r}_{y,w}T_y$ where
$r_{y,w}\in\ca$ are zero for all but finitely many $y$. Note that $r_{w,w}=1$.

\proclaim{Lemma 4.4} Let $w\in W$ and $s\in S$ be such that $w>sw$. For $y\in W$ we
have

$r_{y,w}=r_{sy,sw}$ if $sy<y$,

$r_{y,w}=r_{sy,sw}+(v_s-v_s\i)r_{y,sw}$ if $sy>y$.
\endproclaim
We have 
$$\align\ov T_w&=T_s\i\ov T_{sw}=(T_s-(v_s-v_s\i))\su_y\ov r_{y,sw}T_y\\&=
\su_y\ov r_{y,sw}T_{sy}-\su_y(v_s-v_s\i)\ov r_{y,sw}T_y+\su_{y;sy<y}(v_s-v_s\i)
\ov r_{y,sw}T_y\\&=\su_y\ov r_{sy,sw}T_y-\su_{y;sy>y}(v_s-v_s\i)\ov r_{y,sw}T_y.
\endalign$$
The lemma follows.

\proclaim{Lemma 4.5} For any $y,w$ we have $\ov r_{y,w}=\sg(yw)r_{y,w}$.
\endproclaim
We argue by induction on $l(w)$. If $l(w)=0$, then $w=1$ and the result is obvious.
Assume now that $l(w)\ge 1$. We can find $s\in S$ such that $w>sw$. Assume first 
that $sy<y$. From 4.4 we see, using the induction hypothesis, that
$$\ov r_{y,w}=\ov r_{sy,sw}=\sg(sysw)r_{sy,sw}=\sg(yw)r_{y,w}.$$
Assume next that $sy>y$. From 4.4 we see, using the induction hypothesis, that
$$\align\ov r_{y,w}&=\ov r_{sy,sw}+(v_s\i-v_s)\ov r_{y,sw}
=\sg(sysw)r_{sy,sw}+(v_s\i-v_s)\sg(ysw)r_{y,sw}\\&=\sg(yw)
(r_{sy,sw}+(v_s-v_s\i)r_{y,sw})=\sg(yw)r_{y,w}.\endalign$$
The lemma is proved.

\proclaim{Lemma 4.6}For any $x,z\in W$ we have $\su_y\ov r_{x,y}r_{y,z}=\d_{x,z}$.
\endproclaim
Using the fact that $\bar{}$ is an involution, we have
$$T_z=\bar{\bar T}_z=\ov{\su_y\ov r_{y,z}T_y}=\su_yr_{y,z}\ov T_y
=\su_y\su_xr_{y,z}\ov r_{x,y}T_x.$$
We now compare the coefficients of $T_x$ on both sides. The lemma follows.

\proclaim{Proposition 4.7} Let $y,w\in W$.

(a) If $r_{y,w}\ne 0$, then $y\le w$.

(b)  Assume that $L(s)>0$ for all $s\in S$. If $y\le w$, then 

$r_{y,w}=v^{L(w)-L(y)} \mod v^{L(w)-L(y)-1}\bz[v\i]$,

$r_{y,w}=\sg(yw)v^{-L(w)+L(y)} \mod v^{-L(w)+L(y)+1}\bz[v]$.

(c) Without assumption on $L$, $r_{y,w}\in v^{L(w)-L(y)}\bz[v^2,v^{-2}]$.
\endproclaim
We prove (a) by induction on $l(w)$. If $l(w)=0$ then $w=1$ and the result is 
obvious. Assume now that $l(w)\ge 1$. We can find $s\in S$ such that $w>sw$. Assume
first that $sy<y$. From 4.4 we see that $r_{sy,sw}\ne 0$ hence, by the induction 
hypothesis, $sy\le sw$. Thus $sy\le sw<w$ and, by 2.3, we deduce $y\le w$. Assume 
next that $sy>y$. From 4.4 we see that either $r_{sy,sw}\ne 0$ or $r_{y,sw}\ne 0$ 
hence, by the induction hypothesis, $sy\le sw$ or $y\le sw$. Combining this with 
$y<sy$ and $sw<w$ we see that $y\le w$. This proves (a).

We prove the first assertion of (b) by induction on $l(w)$. If $l(w)=0$ then $w=1$ 
and the result is obvious. Assume now that $l(w)\ge 1$. We can find 
$s\in S$ such that $w>sw$. Assume first that $sy<y$. Then we have also $sy<w$ and, 
using 2.5(b), we deduce $sy\le sw$. By the induction hypothesis, 
$$\align r_{sy,sw}&=v^{L(sw)-L(sy)}+\text{strictly lower powers}\\&
=v^{L(w)-L(y)}+\text{strictly lower powers}.\endalign$$
But $r_{y,w}=r_{sy,sw}$ and the result follows. Assume next that $sy>y$. From 
$y<sy,y\le w$ we deduce using 2.5(b) that $y\le sw$. By the induction hypothesis, we
have $r_{y,sw}=v^{L(sw)-L(y)}+\text{strictly lower powers}$. Hence
$$\align&(v_s-v_s\i)r_{y,sw}=v^{L(s)}v^{L(sw)-L(y)}+\text{strictly lower powers }\\&
=v^{L(w)-L(y)}+\text{strictly lower powers}.\endalign$$
On the other hand, if $sy\le sw$, then by the induction hypothesis,
$$\align&r_{sy,sw}=v^{L(sw)-L(sy)}+\text{strictly lower powers}\\&
=v^{L(w)-L(y)-2L(s)}+\text{strictly lower powers}\endalign$$
while if $sy\not\le sw$ then $r_{sy,sw}=0$ by (a). Thus, in 
$r_{y,w}=r_{sy,sw}+(v_s-v_s\i)r_{y,sw}$, the term $r_{sy,sw}$ contributes only 
powers of $v$ which are strictly smaller than $L(w)-L(y)$ hence 
$r_{y,w}=v^{L(w)-L(y)}+\text{strictly lower powers}$. This proves the first
assertion of (b). The second assertion of (b) follows from the first using 4.5.

We prove (c) by induction on $l(w)$. If $l(w)=0$ then $w=1$ and the result is 
obvious. Assume now that $l(w)\ge 1$. We can find $s\in S$ such that $w>sw$. Assume
first that $sy<y$. By the induction hypothesis, 
$$r_{y,w}=r_{sy,sw}\in v^{L(sw)-L(sy)}\bz[v^2,v^{-2}]=v^{L(w)-L(y)}\bz[v^2,v^{-2}]$$
as required. Assume next that $sy>y$. By the induction hypothesis, 
$$\align&r_{y,w}=r_{sy,sw}+(v_s-v_s\i)r_{y,sw}\\&\in v^{L(sw)-L(sy)}
\bz[v^2,v^{-2}]+v^{L(s)}v^{L(sw)-L(y)}\bz[v^2,v^{-2}]=v^{L(w)-L(y)}\bz[v^2,v^{-2}],
\endalign$$
as required. The proposition is proved. 

\proclaim{Proposition 4.8}For $x<z$ in $W$ we have
$\su_{y;x\le y\le z}\sg(y)=0$ (D.N.Verma).
\endproclaim
Using 4.5 we can rewrite 4.6 (in our case) in the form
$$\su_y\sg(xy)r_{x,y}r_{y,z}=0.\tag a$$
Here we may restrict the summation to $y$ such that $x\le y\le z$. In the rest of 
the proof we shall take $L=l$. Then 4.7(b) holds and we see that if $x\le y\le z$, 
then 

$r_{x,y}r_{y,z}=v^{l(y)-l(x)}v^{l(z)-l(y)}+$ strictly lower powers of $v$.
\nl 
Hence (a) states that

$\su_{y;x\le y\le z}\sg(xy)v^{l(z)-l(x)}+$ strictly lower powers of $v$ is $0$.
\nl
In particular $\su_{y;x\le y\le z}\sg(xy)=0$. The proposition is proved.

\subhead 4.9\endsubhead
Now $\bar{}:\ch@>>>\ch$ commutes with $h\m h^\flat$. Hence

(a) $r_{y\i,w\i}=r_{y,w}$
\nl
for any $y,w\in W$. On the other hand, it is clear that $\bar{}:\ch@>>>\ch$ and 
${}^\da:\ch@>>>\ch$ commute.

\head 5. The elements $c_w$\endhead
\subhead 5.1\endsubhead
We preserve the setup of 3.1. For any $n\in\bz$ let 
$$\ca_{\le n}=\op_{m;m\le n}\bz v^m,\ca_{\ge n}=\op_{m;m\ge n}\bz v^m,$$
$$\ca_{<n}=\op_{m;m<n}\bz v^m, \ca_{>n}=\op_{m;m>n}\bz v^m,$$
$$\ch_{\le 0}=\op_w\ca_{\le 0}T_w, \ch_{<0}=\op_w\ca_{<0}T_w.$$
We have $\ca_{\le 0}=\bz[v\i]$, $\ch_{<0}\sub\ch_{\le 0}\sub\ch$. 

\proclaim{Theorem 5.2} (a) Let $w\in W$. There exists a unique element 
$c_w\in\ch_{\le 0}$ such that $\ov c_w=c_w$ and $c_w=T_w\mod\ch_{<0}$.

(b) $\{c_w;w\in W\}$ is an $\ca_{\le 0}$-basis of $\ch_{\le 0}$ and an $\ca$-basis 
of $\ch$.
\endproclaim
We prove the existence part of (a). We will construct, for any $x$ such that
$x\le w$, an element $u_x\in\ca_{\le 0}$ such that

(c) $u_w=1$,

(d) $u_x\in \ca_{<0},\ov u_x-u_x=\su_{y;x<y\le w}r_{x,y}u_y$ for any $x<w$.
\nl
We argue by induction on $l(w)-l(x)$. If $l(w)-l(x)=0$ then $x=w$ and we set
$u_x=1$. Assume now that $l(w)-l(x)>0$ and that $u_z$ is already defined whenever 
$z\le w,l(w)-l(z)<l(w)-l(x)$ so that (c) holds and (d) holds if $x$ is replaced by 
any such $z$. Then the right hand side of the equality in (d) is defined. We denote
it by $a_x\in\ca$. We have
$$\align a_x+\ov a_x &=\su_{y;x<y\le w}r_{x,y}u_y+\su_{y;x<y\le w}\ov r_{x,y}
\ov u_y\\&
=\su_{y;x<y\le w}r_{x,y}u_y+\su_{y;x<y\le w}\ov r_{x,y}(u_y+\su_{z;y<z\le w}r_{y,z}u_z)\\&
=\su_{z;x<z\le w}r_{x,z}u_z+\su_{z;x<z\le w}\ov r_{x,z}u_z
+\su_{z;x<z\le w}\su_{y; x<y<z}\ov r_{x,y}r_{y,z}u_z\\&
=\su_{z;x<z\le w}\su_{y; x\le y\le z}\ov r_{x,y}r_{y,z}u_z
=\su_{z;x<z\le w}\d_{x,z}u_z=0.\endalign$$
(We have used 4.6 and the equality $r_{y,y}=1$.) Since $a_x+\ov a_x=0$, we have 
$a_x=\su_{n\in\bz}\g_nv^n$ (finite sum) where $\g_n\in\bz$ satisfy $\g_n+\g_{-n}=0$
for all $n$ and in particular, $\g_0=0$. Then $u_x=-\su_{n<0}\g_nv^n\in\ca_{<0}$ 
satisfies $\ov u_x-u_x=a_x$. This completes the inductive construction of the 
elements $u_x$. We set $c_w=\su_{y;y\le w}u_yT_y\in\ch_{\le 0}$. It is clear that 
$c_w=T_w\mod\ch_{<0}$. We have
$$\align&\ov c_w=\su_{y;y\le w}\ov u_y\ov T_y=\su_{y;y\le w}\ov u_y\su_{x;x\le y}
\ov r_{x,y}T_x=\su_{x;x\le w}(\su_{y;x\le y\le w}\ov r_{x,y}\ov u_y)T_x\\&
=\su_{x;x\le w}u_xT_x=c_w.\endalign$$
(We have used the fact that $r_{x,y}\ne 0$ implies $x\le y$, see 4.7, and (d).) The
existence of the element $c_w$ is established.

To prove uniqueness, it suffices to verify the following statement:

(e) {\it If $h\in\ch_{<0}$ satisfies $\ov h=h$ then} $h=0$.
\nl
We can write uniquely $h=\su_{y\in W}f_yT_y$ where $f_y\in\ca_{<0}$ are zero for all
but finitely many $y$. Assume that not all $f_y$ are $0$. Then we can find 
$l_0\in\bn$ such that 
$$Y_0:=\{y\in W;f_y\ne 0,l(y)=l_0\}\ne\em \text{ and }
\{y\in W;f_y\ne 0,l(y)>l_0\}=\em.$$
Now $\su_yf_yT_y=\ov{\su_yf_yT_y}$ implies
$$\su_{y\in Y_0}f_yT_y=\su_{y\in Y_0}\ov{f_y}T_y\mod\su_{y;l(y)<l_0}\ca T_y$$
hence $\ov{f_y}=f_y$ for any $y\in Y_0$.  Since $f_y\in\ca_{<0}$, it follows that 
$f_y=0$ for any $y\in Y_0$, a contradiction. We have proved that $f_y=0$ for all 
$y$; (e) is verified and (a) is proved.

The elements $c_w$ constructed in (a) (for various $w$) are related to the basis 
$T_w$ by a triangular matrix (with respect to $\le$) with $1$ on the diagonal. Hence
these elements satisfy (b). The theorem is proved.

\subhead 5.3\endsubhead
For any $w\in W$ we set $c_w=\su_{y\in W}p_{y,w}T_y$ where $p_{y,w}\in\ca_{\le 0}$.
By the proof of 5.2 we have 

$p_{y,w}=0$ unless $y\le w$, 

$p_{w,w}=1$,

$p_{y,w}\in\ca_{<0}$ if $y<w$.
\nl
Moreover, for any $x\le w$ in $W$ we have

$\bap_{x,w}=\su_{y;x\le y\le w}r_{x,y}p_{y,w}$.

\proclaim{Proposition 5.4}(a) Assume that $L(s)>0$ for all $s\in S$. If $x\le w$, 
then 

$p_{x,w}=v^{-L(w)+L(x)}\mod v^{-L(w)+L(x)+1}\bz[v]$.

(b) Without assumption on $L$, for $x\le w$ we have
$p_{x,w}\in v^{L(w)-L(x)}\bz[v^2,v^{-2}]$.
\endproclaim
We prove (a) by induction on $l(w)-l(x)$. If $l(w)-l(x)=0$ then $x=w$, $p_{x,w}=1$ 
and the result is obvious. Assume now that $l(w)-l(x)>0$. Using 4.7(b) and the 
induction hypothesis, we see that $\su_{y;x<y\le w}r_{x,y}p_{y,w}$ is equal to
$$\su_{y;x<y\le w}\sg(x)\sg(y)v^{-L(y)+L(x)}v^{-L(w)+L(y)}=
\su_{y;x<y\le w}\sg(x)\sg(y)v^{-L(w)+L(x)}$$
plus strictly higher powers of $v$. Using 4.8, we see that this is $-v^{-L(w)+L(x)}$
plus strictly higher powers of $v$. Thus, 

$\bap_{x,w}-p_{x,w}=-v^{-L(w)+L(x)}+$ strictly higher powers of $v$.
\nl
Since $\bap_{x,w}\in v\bz[v]$, it is in particular a $\bz$-linear combination of 
powers of $v$ strictly higher than $-L(w)+L(x)$. 
(We use that, if $a\le b$ in $W$ then $L(a)\le L(b)$ which follows from 2.4.)
Hence 

$-p_{x,w}=-v^{-L(w)+L(x)}+$ strictly higher powers of $v$. 
\nl
This proves (a).

We prove (b) by induction on $l(w)-l(x)$. If $l(w)-l(x)=0$, then $x=w$, $p_{x,w}=1$
and the result is obvious. Assume now that $l(w)-l(x)>0$. Using 4.7(c) and the 
induction hypothesis, we see that
$$\su_{y;x<y\le w}r_{x,y}p_{y,w}\in\su_{y;x<y\le w}v^{L(y)-L(x)}
v^{L(w)-L(y)}\bz[v^2,v^{-2}]\sub v^{L(w)-L(x)}\bz[v^2,v^{-2}].$$
Thus, $\bap_{x,w}-p_{x,w}\in v^{L(w)-L(x)}\bz[v^2,v^{-2}]$. Hence
$p_{x,w}\in v^{L(w)-L(x)}\bz[v^2,v^{-2}]$. The proposition is proved.

\subhead 5.5\endsubhead
Let $s\in S$. From $T_s\i=T_s-(v_s-v_s\i)$ we see that $r_{1,s}=v_s-v_s\i$. We also
see that

$\ov{T_s+v_s\i}=T_s-(v_s-v_s\i)+v_s=T_s+v_s\i$,

$\ov{T_s-v_s}=T_s-(v_s-v_s\i)-v_s\i=T_s-v_s$.
\nl
If $L(s)=0$ we have $T_s\i=T_s$. Hence, 

$c_s=T_s+v_s\i$ if $L(s)>0$, 

$c_s=T_s-v_s$ if $L(s)<0$, 

$c_s=T_s$ if $L(s)=0$.

\subhead 5.6\endsubhead
Now $h\m h^\flat$ carries $\ch_{\le 0}$ into itself; moreover, it commutes with 
$\bar{}:\ch@>>>\ch$ (as pointed out in 4.9). Hence it carries $c_w$ to $c_{w\i}$ for
any $w\in W$. It follows that

(a) $p_{y\i,w\i}=p_{y,w}$
\nl
for any $y,w\in W$.

\head 6. Left or right multiplication by $c_s$\endhead
\subhead 6.1\endsubhead 
We preserve the setup of 3.1 and we fix $s\in S$. Assume first that $L(s)=0$. In 
this case we have $c_s=T_s$; moreover, for any $y\in W$ we have $T_sT_y=T_{sy}$. 
Hence for $w\in W$ we have
$$c_sc_w=\su_yp_{y,w}T_sT_y=\su_yp_{y,w}T_{sy}=\su_yp_{sy,w}T_y.$$
We see that $c_sc_w\in\ch_{\le 0}$ and $c_sc_w=T_{sw}\mod\ch_{<0}$. Since
$\ov{c_sc_w}=c_sc_w$, it follows that, in this case, $c_sc_w=c_{sw}$. Similarly we 
have $c_wc_s=c_{ws}$.

\subhead 6.2\endsubhead 
{\it In the remainder of this chapter (except in 6.8) we assume that} $L(s)>0$.

\proclaim{Proposition 6.3} To any $y,w\in W$ such that $sy<y<w<sw$ one can assign 
uniquely an element $\mu_{y,w}^s\in\ca$ so that

(i) $\ov{\mu}_{y,w}^s=\mu_{y,w}^s$ and

(ii) $\su_{z;y\le z<w;sz<z}p_{y,z}\mu_{z,w}^s-v_sp_{y,w}\in\ca_{<0}$
\nl
for any $y,w\in W$ such that $sy<y<w<sw$.
\endproclaim
Let $y,w$ be as above. We may assume that $\mu_{z,w}^s$ are already defined for all
$z$ such that $y<z<w;sz<z$. Then condition (ii) is of the form:

$\mu_{y,w}^s$ {\it equals a known element of $\ca$ modulo} $\ca_{<0}$.
\nl
This condition determines uniquely the coefficients of $v^n$ with $n\ge 0$ in
$\mu_{y,w}^s$. Then condition (i) determines uniquely the coefficients of $v^n$ with
$n<0$ in $\mu_{y,w}^s$. The proposition is proved.

\proclaim{Proposition 6.4}Let $y,w\in W$ be such that $sy<y<w<sw$. Then 
$\mu_{y,w}^s$ is a $\bz$-linear combination of powers $v^n$ with
$-L(s)+1\le n\le L(s)-1$ and $n=L(w)-L(y)-L(s)\mod 2$.
\endproclaim
We may assume that this is already known for all $\mu_{z,w}^s$ with $z$ such that 
$y<z<w;sz<z$. Using 6.3(ii) and 5.4(b), we see that $\mu_{y,w}^s$ is a $\bz$-linear 
combination of powers $v^n$ such that, whenever $n\ge 0$, we have $n\le L(s)-1$ and
$n=L(w)-L(y)-L(s)\mod 2$. Using now 6.3(i), we deduce the remaining assertions of 
the proposition.

\proclaim{Corollary 6.5} Assume that $L(s)=1$. Let $y,w\in W$ be such that 
$sy<y<w<sw$. Then $\mu_{y,w}^s$ is an integer, equal to the coefficient of $v\i$ in
$p_{y,w}$. In particular, it is $0$ unless $L(w)-L(y)$ is odd.
\endproclaim
In this case, the inequalities of 6.4 become $0\le n\le 0$. They imply $n=0$. Thus,
$\mu_{y,w}^s\in\bz$. Picking up the coefficient of $v^0$ in the two sides of
6.3(ii), we see that $\mu_{y,w}^s$ is equal to the coefficient of $v\i$ in 
$p_{y,w}$. The last assertion follows from 5.4.

\proclaim{Theorem 6.6} Let $w\in W$.

(a) If $w<sw$, then $c_sc_w=c_{sw}+\su_{z;sz<z<w}\mu_{z,w}^sc_z$.

(b) If $sw<w$, then $c_sc_w=(v_s+v_s\i)c_w$.
\endproclaim
Since $c_s=T_s+v_s\i$ (see 5.5), we see that (b) is equivalent to $(T_s-v_s)c_w=0$,
or to

(c) $p_{x,w}=v_s\i p_{sx,w}$
\nl
(where $sw<w$ and $x<sx$). We prove the theorem by induction on $l(w)$. If $l(w)=0$, then
$w=1$ and the result is obvious. Assume now that $l(w)\ge 1$ and that the result
holds when $w$ is replaced by $w'$ with $l(w')<l(w)$. 

{\it Case} 1. Assume that $w<sw$. Using $c_s=T_s+v_s\i$, we see that the coefficient
of $T_y$ in the left hand side minus the right hand side of (a) is
$$f_y=v_s^\s p_{y,w}+p_{sy,w}-p_{y,sw}-\su_{z;y\le z<w;sz<z}p_{y,z}\mu_{z,w}^s$$
where $\s=1$ if $sy<y$ and $\s=-1$ if $sy>y$. We must show that $f_y=0$. We first
show that 

(d) $f_y\in\ca_{<0}$. 
\nl
If $sy<y$ this follows from 6.3(ii). (The contribution of $p_{sy,w}-p_{y,sw}$ is in
$\ca_{<0}$ if $sy\ne w$ and is $1-1=0$ if $sy=w$.) 

If $sy>y$ then, by (c) (applied to $z$ in the sum, instead of $w$), we have
$$\align f_y&=v_s\i p_{y,w}+p_{sy,w}-p_{y,sw}-\su_{z;y\le z<w;sz<z}v_s\i p_{sy,z}
\mu_{z,w}^s\\&=v_s\i f_{sy}+v_s\i p_{sy,sw}-p_{y,sw}\endalign$$
(the second equality holds by 2.5(a)) and this is in $\ca_{<0}$ since 
$f_{sy}\in\ca_{<0}$ (by the previous paragraph), $v_s\i\in\ca_{<0}$ and since 
$y\ne sw$. Thus, (d) is proved.

Since both sides of (a) are fixed by $\bar{}$, the sum $\su_yf_yT_y$ is fixed by 
$\bar{}$. From (d) and 5.2(e) we see that $f_y=0$ for all $y$, as required.

{\it Case} 2. Assume that $w>sw$. Then case 1 is applicable to $sw$ (by the 
induction hypothesis). We see that
$$c_w=(T_s+v_s\i)c_{sw}-\su_{z;sz<z<sw}\mu_{z,sw}^sc_z.$$
Now $(T_s-v_s)(T_s+v_s\i)=0$ and $(T_s-v_s)c_z=0$ for each $z$ in the sum (by the 
induction hypothesis). Hence $(T_s-v_s)c_w=0$. The theorem is proved.

\proclaim{Corollary 6.7} Let $w\in W$.

(a) If $w<ws$, then $c_wc_s=c_{ws}+\su_{z;zs<z<w}\mu_{z\i,w\i}^sc_z$.

(b) If $ws<w$, then $c_wc_s=(v_s+v_s\i)c_w$.
\endproclaim
We write the equalities in 6.6(a),(b) for $w\i$ instead of $w$ and we apply to these
equalities $h\m h^\flat$ which carries $c_w$ to $c_{w\i}$; the corollary follows.

\subhead 6.8\endsubhead
Now 6.3, 6.6, 6.7 remain valid when $L(s)<0$ provided that we replace in their 
statements and proofs $v_s$ by $-v_s\i$.

\head 7. Dihedral groups \endhead
\subhead 7.1\endsubhead
We preserve the setup of 3.1; we assume that $S$ consists of two elements $s_1,s_2$.
For $i=1,2$, let $L_i=L(s_i),T_i=T_{s_i},c_i=c_{s_i}$. We assume that $L_1>0,L_2>0$.
Let $\z=v^{L_1-L_2}+v^{L_2-L_1}\in\ca$. Let $m=m_{s_1,s_2}$. Let $1_k,2_k$ be as in
1.4. For $w\in W$ we set

$\G_w=\su_{y;y\le w}v^{-L(w)+L(y)}T_y$.

\proclaim{Lemma 7.2} We have

$c_1\G_{2_k}=\G_{1_{k+1}}+v^{L_1-L_2}\G_{1_{k-1}}$ if $k\in[2,m)$,

$c_2\G_{1_k}=\G_{2_{k+1}}+v^{-L_1+L_2}\G_{2_{k-1}}$ if $k\in[2,m)$,

$c_1\G_{2_k}=\G_{1_{k+1}}$ if $k=0,1$,

$c_2\G_{1_k}=\G_{2_{k+1}}$ if $k=0,1$.
\endproclaim
Since $c_i=T_i+v^{-L_i}$, the proof is an easy exercise.

\proclaim{Proposition 7.3} Assume that $L_1=L_2$. For any $w\in W$ we have 
$c_w=\G_w$.
\endproclaim
This is clear when $l(w)\le 1$. In the present case Lemma 7.2 gives
$$\G_{1_{k+1}}=c_1\G_{2_k}-\G_{1_{k-1}},\qua\G_{2_{k+1}}=c_2\G_{1_k}-\G_{2_{k-1}}
\tag c$$
for $k\in[2,m)$. This and 7.2 shows by induction on $k$ that $\bar\G_w=\G_w$ for all
$w\in W$. Clearly, $\G_w=T_w\mod\ch_{<0}$. The lemma follows.

\subhead 7.4\endsubhead
{\it In 7.4-7.6 we assume that $L_2>L_1$.} In this case, if $m<\iy$, then $m$ is 
even. (See 3.1.) For $2k+1\in[1,m)$ we set 
$$\align&\G'_{2_{2k+1}}=\su_{s\in[0,k-1]}(1-v^{2L_1}+v^{4L_1}-\do+(-1)^s
v^{2sL_1})v^{-sL_1-sL_2}\\&\T(T_{2_{2k-2s+1}}+v^{-L_2}T_{2_{2k-2s}}+v^{-L_2}
T_{1_{2k-2s}}+v^{-2L_2}T_{1_{2k-2s-1}})\\&+(1-v^{2L_1}+v^{4L_1}-\do+(-1)^k
v^{2kL_1})v^{-kL_1-kL_2}(T_{2_1}+v^{-L_2}T_{2_0}).\endalign$$
For $2k+1\in[3,m)$ we set
$$\align&\G'_{1_{2k+1}}=T_{1_{2k+1}}+v^{-L_1}T_{1_{2k}}+v^{-L_1}T_{2_{2k}}+
v^{-2L_1}T_{2_{2k-1}}\\&+\su\Sb y\\y\le 1_{2k-1}\eSb v^{-L(w)+L(y)}(1+v^{2L_1})T_y
\endalign$$
where $w=1_{2k+1}$. For $w$ such that $l(w)$ is even and for $w=1_1$ we set 
$\G'_w=\G_w$.

\proclaim{Lemma 7.5} We have

(a) $c_1\G'_{2_{k'}}=\G'_{1_{k'+1}}$, if $k'\in[0,m)$;

(b) $c_2\G'_{1_{k'}}=\G'_{2_{k'+1}}+\z\G'_{2_{k'-1}}+\G'_{2_{k'-3}}$, if
$k'\in[4,m)$;

(c) $c_2\G'_{1_{k'}}=\G'_{2_{k'+1}}+\z\G'_{2_{k'-1}}$, if $k'=2,3, k'<m$;

(d) $c_2\G'_{1_{k'}}=\G'_{2_{k'+1}}$ if $k'=0,1$.
\endproclaim
From the definitions we have

(e) $\G'_{2_{2k+1}}=\su_{s\in[0,k]}(-1)^sv^{s(L_1-L_2)}\G_{2_{2k-2s+1}}$ if
$2k+1\in[1,m)$, 

(f) $\G'_{1_{2k+1}}=\G_{1_{2k+1}}+v^{L_1-L_2}\G_{1_{2k-1}}$ if $2k+1\in[3,m)$.
\nl
We prove (a) for $k'=2k+1$. The left hand side can be computed using (e) and 7.2:
$$\align&c_1\G'_{2_{2k+1}}=c_1(\G_{2_{2k+1}}-v^{L_1-L_2}\G_{2_{2k-1}}+
v^{2L_1-2L_2}\G_{2_{2k-3}}+\do)\\&=\G_{1_{2k+2}}+v^{L_1-L_2}\G_{1_{2k}}
-v^{L_1-L_2}\G_{1_{2k}}-v^{2L_1-2L_2}\G_{1_{2k-2}}\\&+v^{2L_1-2L_2}\G_{1_{2k-2}}
-v^{3L_1-3L_2}\G_{1_{2k-4}}+\do=\G_{1_{2k+2}}=\G'_{1_{2k+2}}.\endalign$$
This proves (a) for $k'=2k+1$. Now (a) for $k'=0$ is trivial. We prove (a) for
$k'=2k\ge 2$. The left hand side can be  computed using 7.2 and (f):
$$c_1\G'_{2_{2k}}=c_1\G_{2_{2k}}=\G_{1_{2k+1}}+v^{L_1-L_2}\G_{1_{2k-1}}
=\G'_{1_{2k+1}}.$$
This proves (a) for $k'=2k$. We prove (b) for $k'=2k$. The left hand side can be 
computed using 7.2:
$$c_2\G'_{1_{2k}}=c_2\G_{1_{2k}}=\G_{2_{2k+1}}+v^{-L_1+L_2}\G_{2_{2k-1}}.$$
The right hand side of (b) is (using (e)):
$$\align&\G_{2_{2k+1}}-v^{L_1-L_2}\G_{2_{2k-1}}+v^{2L_1-2L_2}\G_{2_{2k-3}}+\do\\&+
\z\G_{2_{2k-1}}-v^{L_1-L_2}\z\G_{2_{2k-3}}+v^{2L_1-2L_2}\z\G_{2_{2k-5}}+\do\\&+
\G_{2_{2k-3}}-v^{L_1-L_2}\G_{2_{2k-5}}+v^{2L_1-2L_2}\G_{2_{2k-7}}+\do=\G_{2_{2k+1}}+
v^{-L_1+L_2}\G_{2_{2k-1}}.\endalign$$
This proves (b) for $k'=2k$. We prove (b) for $k'=2k+1$. The left hand side can be 
computed using (f) and 7.2:
$$\align& c_2\G'_{1_{2k+1}}=c_2(\G_{1_{2k+1}}+v^{L_1-L_2}\G_{1_{2k-1}})\\&
=\G_{2_{2k+2}}+v^{-L_1+L_2}\G_{2_{2k}}+v^{L_1-L_2}\G_{2_{2k}}+\G_{2_{2k-2}}
=\G'_{2_{2k+2}}+\z\G'_{2_{2k}}+\G'_{2_{2k-2}}.\endalign$$
This proves (b) for $k'=2k+1$. The proof of (c),(d) is similar to that of (b). This
completes the proof.

\proclaim{Proposition 7.6} For any $w\in W$ we have $c_w=\G'_w$.
\endproclaim
Clearly, $\G'_w=T_w\mod\ch_{<0}$. From the formulas in 7.5 we see by induction on 
$l(w)$ that $\bar\G'_w=\G'_w$ for all $w$. The proposition is proved. (This was 
proved for $m=4$ in \cite{\LC}, for $m=6$ in \cite{\XI}, for general $m$ 
independently in \cite{\LH} and \cite{\GP, p.396}.)

\proclaim{Proposition 7.7} Assume that $m=\iy$. For $a\in\{1,2\}$, let
$f_a=v^{L(a)}+v^{-L(a)}$.

(a) If $L_1=L_2$ and $k,k'\ge 0$ then
$c_{a_{2k+1}}c_{a_{2k'+1}}=f_a\su_{u\in[0,\min(2k,2k')]}c_{a_{2k+2k'+1-2u}}$.

(b) If $L_2>L_1$ and $k,k'\ge 0$ then
$c_{2_{2k+1}}c_{2_{2k'+1}}=f_2\su_{u\in[0,\min(k,k')]}c_{2_{2k+2k'+1-4u}}$.

(c) If $L_2>L_1$ and $k,k'\ge 1$ then
$$c_{1_{2k+1}}c_{1_{2k'+1}}=f_1\su_{u\in[0,\min(k-1,k'-1)]}p_uc_{1_{2k+2k'+1-2u}}$$
where $p_u=\z$ for $u$ odd, $p_u\in\bz$ for $u$ even.
\endproclaim
We prove (a). For $k=k'=0$ the equality in (a) is clear. Assume now that 
$k=0,k'\ge 1$. Using 7.2, 7.3, we have
$$\align&c_2c_{2_{2k'+1}}=c_2(c_2c_{1_{2k'}}-c_{2_{2k'-1}})=
f_2c_2c_{1_{2k'}}-f_2c_{2_{2k'-1}}\\&=
f_2c_{2_{2k'+1}}+f_2c_{2_{2k'-1}}-f_2c_{2_{2k'-1}}=f_2c_{2_{2k'+1}},\endalign$$
as required. We now prove the equality in (a) for fixed $k'$, by induction on $k$. 
The case $k=0$ is already known. Assume now that $k=1$. From 7.2, 7.3 we have 
$c_{2_3}=c_2c_1c_2-c_2$. Using this and 7.2, 7.3, we have
$$\align&c_{2_3}c_{2_{2k'+1}}=c_2c_1c_2c_{2_{2k'+1}}-c_2c_{2_{2k'+1}}=f_2c_2
c_{1_{2k'+2}}+f_2c_2c_{1_{2k'}}-f_2c_{2_{2k'+1}}\\&=f_2c_{2_{2k'+3}}+f_2
c_{2_{2k'+1}}+f_2c_{2_{2k'+1}}+(1-\d_{k',0})f_2c_{2_{2k'-1}}-f_2c_{2_{2k'+1}}
\\&=f_2c_{2_{2k'+3}}+f_2c_{2_{2k'+1}}+f_2(1-\d_{k',0})c_{2_{2k'-1}},\endalign$$
as required. Assume now that $k\ge 2$. From 7.2,7.3 we have
$$c_{2_{2k+1}}=c_2c_1c_{2_{2k-1}}-2c_{2_{2k-1}}-c_{2_{2k-3}}.$$
Using this and the induction hypothesis we have
$$\align&c_{2_{2k+1}}c_{2_{2k'+1}}=c_2c_1c_{2_{2k-1}}c_{2_{2k'+1}}
-2c_{2_{2k-1}}c_{2_{2k'+1}}-c_{2_{2k-3}}c_{2_{2k'+1}}\\&=f_2c_1c_2
\su_{u\in[0,\min(2k-2,k')]}c_{2_{2k+2k'-1-2u}}-f_2\su_{u\in[0,\min(2k-2,k')]}
c_{2_{2k+2k'-1-2u}}\\&-f_2\su_{u\in[0,\min(2k-4,k')]}c_{2_{2k+2k'-3-2u}}.\endalign$$
We now use 7.2,7.3 and (a) follows (for $a=2$). The case $a=1$ is similar.

We prove (b). For $k=k'=0$ the equality in (b) is clear. Assume now that $k=0,k'=1$.
Using 7.5, 7.6, we have
$$c_2c_{2_3}=c_2(c_2c_{1_2}-\z c_{2_1})=f_2c_2c_{1_2}-f_2\z c_{2_1}
=f_2c_{2_3}+f_2\z c_{2_1}-f_2\z c_{2_1}=f_2c_{2_3},$$
as required. Assume next that $k=0,k'\ge 2$. Using 7.5, 7.6, we have
$$\align&c_2c_{2_{2k'+1}}=c_2(c_2c_{1_{2k'}}-\z c_{2_{2k'-1}}-c_{2_{2k'-3}})
=f_2c_2c_{1_{2k'}}-f_2\z c_{2_{2k'-1}}-f_2c_{2_{2k'-3}}\\&
=f_2c_{2_{2k'+1}}+f_2\z c_{2_{2k'-1}}+f_2c_{2_{2k'-3}}
-f_2\z c_{2_{2k'-1}}-f_2\z c_{2_{2k'-3}}=f_2c_{2_{2k'+1}},\endalign$$
as required. We now prove the equality in (a) for fixed $k'$, by induction on $k$. 
The case $k=0$ is already known. Assume now that $k=1$. From 7.5, 7.6 we have 
$c_{2_3}=c_2c_1c_2-\z c_2$. Using this and 7.5, 7.6, we have
$$\align&c_{2_3}c_{2_{2k'+1}}=c_2c_1c_2c_{2_{2k'+1}}-\z c_2c_{2_{2k'+1}}=
f_2c_2c_{1_{2k'+2}}-f_2\z c_{2_{2k'+1}}\\&=f_2c_{2_{2k'+3}}+f_2\z 
c_{2_{2k'+1}}+(1-\d_{k',0})f_2c_{2_{2k'-1}}-f_2\z c_{2_{2k'+1}}\\&=f_2
c_{2_{2k'+3}}+(1-\d_{k',0})f_2c_{2_{2k'-1}}\endalign$$
as required. Assume now that $k\ge 2$. From 7.5, 7.6 we have
$$c_{2_{2k+1}}=c_2c_1c_{2_{2k-1}}-\z c_{2_{2k-1}}-c_{2_{2k-3}}.$$
Using this and the induction hypothesis we have
$$\align&c_{2_{2k+1}}c_{2_{2k'+1}}=c_2c_1c_{2_{2k-1}}c_{2_{2k'+1}}
-\z c_{2_{2k-1}}c_{2_{2k'+1}}-c_{2_{2k-3}}c_{2_{2k'+1}}\\&
=f_2c_2c_1\su_{u\in[0,\min(k-1,k')]}c_{2_{2k+2k'-1-4u}}
-f_2\z\su_{u\in[0,\min(k-1,k')]}c_{2_{2k+2k'-1-4u}}\\&
-f_2\su_{u\in[0,\min(k-2,k')]}c_{2_{2k+2k'-3-4u}}.\endalign$$
We now use 7.5, 7.6 and (b) follows.

The proof of (c) is similar to that of (b). This completes the proof.

\proclaim{Proposition 7.8} Assume that $4\le m<\iy$ and $L_2>L_1$, so that $m=2k+2$ 
with $k\ge 1$. Let

$p_0
=(-1)^k(v^{L_2}+v^{-L_2})(v^{k(L_2-L_1)}+v^{(k-2)(L_2-L_1)}+\do+v^{-k(L_2-L_1)}).$
\nl
Then
$$c_{2_{m-1}}c_{2_{m-1}}=pc_{2_{m-1}}+qc_{2_m},\tag a$$
for some $p,q\in\ca$. Moreover, $p=p_0$.
\endproclaim
From 7.5, 7.6, we see that $\ca c_{2_{m-1}}+\ca c_{2_m}$ is a two-sided ideal of
$\ch$. Hence (a) holds for some (unknown) $p,q\in\ca$. It remains to compute $p$. 
Define an algebra homomorphism $\c:\ch@>>>\ca$ by 
$\c(T_1)=-v^{-L_1},\c(T_2)=v^{L_2}$. Since $c_{2_m}=(T_1+v^{-L_1})h$ for some 
$h\in\ch$ (see 7.5,7.6) we see that $\c(c_{2_m})=0$. Hence applying $\c$ to (a)
gives $\c(c_{2_{m-1}})^2=p\c(c_{2_{m-1}})$. It is thus enough to show that 
$\c(c_{2_{m-1}})=p_0$. We verify this for $m=4$:
$$\align&\c(T_2T_1T_2+v^{-L_2}T_2T_1+v^{-L_2}T_1T_2+v^{-2L_2}T_1\\&+
(v^{-L_1-L_2}-v^{L_1-L_2})T_2+(v^{-L_1-2L_2}-v^{L_1-2L_2}))\\&=-v^{-L_1+2L_2}
-2v^{-L_1}-v^{-L_1-2L_2}+(v^{-L_1-L_2}-v^{L_1-L_2})v^{L_2}\\&+v^{-L_1-2L_2}-
v^{L_1-2L_2}=-v^{-L_1+2L_2}-v^{-L_1}-v^{L_1}-v^{L_1-2L_2}\endalign$$
and for $m=6$:
$$\align&\c(T_2T_1T_2T_1T_2+v^{-L_2}T_2T_1T_2T_1+v^{-L_2}T_1T_2T_1T_2+
v^{-2L_2}T_1T_2T_1\\&+(v^{-L_1-L_2}-v^{L_1-L_2})T_2T_1T_2+(v^{-L_1-2L_2}-
v^{L_1-2L_2})T_1T_2\\&+(v^{-L_1-2L_2}-v^{L_1-2L_2})T_2T_1+(v^{-L_1-3L_2}-
v^{L_1-3L_2})T_1\\&+(v^{-2L_1-2L_2}-v^{-2L_2}+v^{2L_1-2L_2})T_2
+(v^{-2L_1-3L_2}-v^{-3L_2}+v^{2L_1-3L_2}))\\&=v^{-2L_1+3L_2}+2v^{-2L_1+L_2}+
v^{-2L_1-L_2}-v^{-2L_1+L_2}-2v^{-2L_1-L_2}-v^{-2L_1-3L_2}+v^{L_2}\\&+2v^{-L_2}
+v^{-3L_2}+v^{-2L_1-L_2}-v^{-L_2}+v^{2L_1-L_2}+v^{-2L_1-3L_2}-v^{-3L_2}+
v^{2L_1-3L_2}\\&=v^{-2L_1+3L_2}+v^{-2L_1+L_2}+v^{L_2}+v^{-L_2}+v^{2L_1-L_2}+
v^{2L_1-3L_2}.\endalign$$
Analogous computations can be carried out for any even $m$. The proposition is
proved.

\head 8. Cells\endhead
\subhead 8.1\endsubhead
We preserve the setup of 3.1. For $z\in W$ define $D_z\in\Hom_\ca(\ch,\ca)$ by 
$D_z(c_w)=\d_{z,w}$ for all $w\in W$. For $w,w'\in W$ we write 

$w\gt_\cl w'$ (or $w'\to_\cl w$) if $D_w(c_sc_{w'})\ne 0$ for some $s\in S$;

$w\gt_\car w'$ (or $w'\to_\car w$) if $D_w(c_{w'}c_s)\ne 0$ for some $s\in S$.
\nl
If $w,w'\in W$, we say that $w\le_\cl w'$ (resp. $w\le_\car w'$) if there exist
$w=w_0,w_1,\do,w_n=w'$ in $W$ such that for any $i\in[0,n-1]$ we have 
$w_i\gt_\cl w_{i+1}$ (resp. $w_i\gt_\car w_{i+1}$).

If $w,w'\in W$, we say that $w\le_{\lr}w'$ if there exist $w=w_0,w_1,\do,w_n=w'$ in
$W$ such that for any $i\in[0,n-1]$ we have either $w_i\gt_\cl w_{i+1}$ or 
$w_i\gt_\car w_{i+1}$.

Clearly $\le_\cl,\le_\car,\le_{\lr}$ are preorders on $W$. Let
$\si_\cl,\si_\car,\si_{\lr}$ be the associated equivalence relations. (For example,
we have $w\si_\cl w'$ if and only if $w\le_\cl w'$ and $w'\le_\cl w$.) The
equivalence classes on $W$ for $\si_\cl,\si_\car,\si_{\lr}$ are called respectively
{\it left cells, right cells, two-sided cells} of $W$. They depend on the weight 
function $L$.

If $w,w'\in W$, we say that $w<_\cl w'$ (resp. $w<_\car w'$; $w<_{\lr}w'$) if
$w\le_\cl w'$ and $w\not\si_\cl w'$ (resp. $w\le_\car w'$ and $w\not\si_\car w'$; 
$w\le_{\lr}w'$ and $w\not\si_{\lr}w'$). 

For $w,w'\in W$ we have $w\le_\cl w'\Lra w\i\le_\car w'{}\i$ and
$w\le_{\lr}w'\Lra w\i\le_{\lr}w'{}\i$.
\nl
Hence $w\m w\i$ carries left cells to right cells, right cells to left cells and 
two-sided cells to two-sided cells.

\proclaim{Lemma 8.2} Let $w'\in W$.

(a) $\ch_{\le_\cl w'}=\op_{w;w\le_\cl w'}\ca c_w$ is a left ideal of $\ch$.

(b) $\ch_{\le_\car w'}=\op_{w;w\le_\car w'}\ca c_w$ is a right ideal of $\ch$.

(c) $\ch_{\le_{\lr}w'}=\op_{w;w\le_{\lr}w'}\ca c_w$ is a two-sided ideal of $\ch$.
\endproclaim
This follows from the definitions since $c_s(s\in S)$ generate $\ch$ as an 
$\ca$-algebra.

\subhead 8.3\endsubhead
Let $Y$ be a left cell of $W$. From 8.3(a) we see that for $y\in Y$,
$$\op_{w;w\le_\cl y}\ca c_w/\op_{w;w<_\cl y}\ca c_w$$
is a quotient of two left ideals of $\ch$ (independent of the choice of $y$) hence 
it is naturally a left $\ch$-module; it has an $\ca$-basis consisting of the images
of $c_w (w\in Y)$.

Similarly, if $Y'$ is a right cell of $W$ then, for $y'\in Y'$,
$$\op_{w;w\le_\car y'}\ca c_w/\op_{w;w<_\car y'}\ca c_w$$
is a quotient of two right ideals of $\ch$ (independent of the choice of $y'$) hence
it is naturally a right $\ch$-module; it has an $\ca$-basis consisting of the images
of $c_w (w\in Y')$.

If $Y''$ is a two-sided cell of $W$ then, for $y''\in Y''$,
$$\op_{w;w\le_{\lr}y''}\ca c_w/\op_{w;w<_{\lr}y''}\ca c_w$$
is a quotient of two two-sided ideals of $\ch$ (independent of the choice of $y''$)
hence it is naturally a $\ch$-bimodule; it has an $\ca$-basis consisting of the 
images of $c_w (w\in Y'')$.

\proclaim{Lemma 8.4} Let $s\in S$. Assume that $L(s)>0$. Let 
$\ch^s=\op_{w;sw<w}\ca c_w$, ${}^s\ch=\op_{w;ws<w}\ca c_w$.

(a) $\{h\in\ch;(c_s-v_s-v_s\i)h=0\}=\ch^s$. Hence $\ch^s$ is a right ideal of $\ch$.

(b) $\{h\in\ch;h(c_s-v_s-v_s\i)=0\}={}^s\ch$. Hence ${}^s\ch$ is a left ideal of 
$\ch$.
\endproclaim
We prove the equality in (a). If $h\in\ch^s$ then $(c_s-v_s-v_s\i)h=0$ by 6.6(b). 
Conversely, by 6.6, we have $c_sh\in\ch^s$ for any $h\in\ch$. Hence, if $h\in\ch$ is
such that $c_sh=(v_s+v_s\i)h$, then $(v_s+v_s\i)h\in\ch^s$ so that $h\in\ch^s$ 
(since $\ch/\ch^s$ is a free $\ca$-module). This proves (a). The proof of (b) is 
entirely similar. The lemma is proved.

\subhead 8.5\endsubhead
For $w\in W$ we set $\cl(w)=\{s\in S;sw<w\},\car(w)=\{s\in S;ws<w\}$.

\proclaim{Lemma 8.6} Let $w,w'\in W$. Assume that $L(s)>0$ for all $s\in S$.

(a) If $w\le_\cl w'$, then $\car(w')\sub\car(w)$. If $w\si_\cl w'$, then 
$\car(w')=\car(w)$. 

(b) If $w\le_\car w'$, then $\cl(w')\sub\cl(w)$. If $w\si_\car w'$, then
$\cl(w')=\cl(w)$.
\endproclaim
To prove the first assertion of (a), we may assume that $D_w(c_sc_{w'})\ne 0$ 
for some $s\in S$. In this case, let $t\in\car(w')$. We must prove that 
$t\in\car(w)$. We have $c_{w'}\in{}^t\ch$. By 8.4, ${}^t\ch$ is a left ideal of 
$\ch$. Hence $c_sc_{w'}\in{}^t\ch$. By the definition of ${}^t\ch$, for 
$h\in{}^t\ch$ we have $D_w(h)=0$ unless $wt<w$. Hence from $D_w(c_sc_{w'})\ne 0$ we
deduce $wt<w$, as required. This proves the first assertion of (a). The second 
assertion of (a) follows immediately from the first. The proof of (b) is entirely 
similar to that of (a). The lemma is proved.

\subhead 8.7\endsubhead
In the remainder of this chapter we write $\gt,\to$ instead of $\gt_\cl,\to_\cl$.
We describe the left cells of $W$ in the setup of 7.3. From 7.2 and 7.3 we can
determine all pairs $y\ne w$ such that $y\gt w$ 

$1_0\to 2_1\lras 1_2\lras 2_3\lras\do$,

$2_0\to 1_1\lras 2_2\lras 1_3\lras\do$,
\nl
if $m=\iy$,

$1_0\to 2_1\lras 1_2\lras 2_3\lras\do\lras 2_{m-1}\to 1_m$,

$2_0\to 1_1\lras 2_2\lras 1_3\lras\do\lras 1_{m-1}\to 2_m$,
\nl
if $m<\iy, m$ even,

$1_0\to 2_1\lras 1_2\lras 2_3\lras\do\lras 1_{m-1}\to 2_m$,

$2_0\to 1_1\lras 2_2\lras 1_3\lras\do\lras 2_{m-1}\to 1_m$,
\nl
if $m<\iy, m$ odd. Hence the left cells are 

$\{1_0\},\{2_1,1_2,2_3,\do\},\{1_1,2_2,1_3,\do\}$,
\nl
if $m=\iy$, 

$\{1_0\},\{2_1,1_2,2_3,\do,2_{m-1}\},\{1_1,2_2,1_3,\do,1_{m-1}\},\{2_m\}$,
\nl
if $m<\iy, m$ even,

$\{1_0\},\{2_1,1_2,2_3,\do,1_{m-1}\},\{1_1,2_2,1_3,\do,2_{m-1}\},\{2_m\}$,
\nl
if $m<\iy, m$ odd.

The two-sided cells are $\{1_0\},W-\{1_0\}$ if $m=\iy$ and
$\{1_0\},\{2_m\},W-\{1_0,2_m\}$ if $m<\iy$.

\subhead 8.8\endsubhead
We describe the left cells of $W$ in the setup of 7.4. From 7.5 and 7.6 we can
determine all pairs $y\ne w$ such that $y\gt w$. If $m=\iy$, these pairs are:
$$1_0\to 2_1\lras 1_2\to 2_3\lras 1_4\to \do,\qua
2_0\to 1_1\to 2_2\lras 1_3\to 2_4\lras\do,$$
and $2_1\gt 1_4, 2_2\gt 1_5, 2_3\gt 1_6,\do$.

If $m=4$, these pairs are:
$$1_0\to 2_1\lras 1_2\to 2_3\to 1_4,\qua 2_0\to 1_1\to 2_2\lras 1_3\to 2_4.$$
If $m=6$, these pairs are:
$$1_0\to 2_1\lras 1_2\to 2_3\lras 1_4\to 2_5\to 1_6,\qua
2_0\to 1_1\to 2_2\lras 1_3\to 2_4\lras 1_5\to 2_6,$$
and $2_1\gt 1_4, 2_2\gt 1_5$. An analogous pattern holds for any even $m$.

Hence the left cells are 
$$\{1_0\},\{2_1,1_2,2_3,1_4,\do\},\{1_1\},\{2_2,1_3,2_4,1_5,\do\},$$
if $m=\iy$,
$$\{1_0\},\{2_1,1_2,2_3,\do, 1_{m-2}\},\{2_{m-1}\},\{1_1\},
\{2_2,1_3,2_4,\do,1_{m-1}\},\{2_m\},$$
if $m<\iy$.

The two-sided cells are 

$\{1_0\},\{1_1\}, W-\{1_0,1_1\}$, if $m=\iy$ and

$\{1_0\},\{1_1\},\{2_{m-1}\},\{2_m\},W-\{1_0,1_1,2_{m-1},2_m\}$, if $m<\iy$.

\subhead 8.9\endsubhead
For further examples of cells (in the case where $L$ is non-costant) see \cite{\LU},
\cite{\BR}, \cite{\GE}.

\head 9. Cosets of parabolic subgroups\endhead
We preserve the setup of 3.1.

\proclaim{Lemma 9.1} Let $w\in W$. Assume that $w=s_1s_2\do s_q$ with $s_i\in S$.
We can find a subsequence $i_1<i_2<\do<i_r$ of $1,2,\do,q$ such that 
$w=s_{i_1}s_{i_2}\do s_{i_r}$ is a reduced expression in $W$.
\endproclaim
We argue by induction on $q$. If $q=0$ the result is obvious. Assume that $q>0$. 
Using the induction hypothesis we can assume that $s_2\do s_q$ is a reduced 
expression. If $s_1s_2\do s_q$ is a reduced expression, we are done. Hence we may 
assume that $s_1s_2\do s_q$ is not a reduced expression. Then $l(w)\le q-1$. By 1.7, we
can find $j\in [2,q]$ such that $s_1s_2\do s_{j-1}=s_2s_3\do s_j$. Then
$w=s_2s_3\do s_{j-1}s_{j+1}\do s_q$ is a reduced expression. The lemma is proved.

\subhead 9.2\endsubhead
Let $w\in W$. Let $w=s_1s_2\do s_q$ be a reduced expression of $w$. Using 1.9, we
see that the set $\{s\in S;s=s_i \text{ for some } i\in[1,q]\}$ is independent of
the choice of reduced expression. We denote it by $S_w$.

\subhead 9.3\endsubhead
In the remainder of this chapter we fix $I\sub S$. Recall that $W_I=\la I\ra$. 

If $w\in W_I$ then we can find a reduced expression $w=s_1s_2\do s_q$ in $W$ with
all $s_i\in I$ (we first write $w=s_1s_2\do s_q$, a not necessarily reduced 
expression with all $s_i\in I$, and then we apply 9.1). Thus, $S_w\sub I$. 
Conversely, it is clear that if $w'\in W$ satisfies $S_{w'}\sub I$ then 
$w'\in W_I$. It follows that

$W_I=\{w\in W;S_w\sub I\}$.

\subhead 9.4\endsubhead
Replacing $S,(m_{s,s'})_{(s,s')\in S\T S}$ by $I,(m_{s,s'})_{(s,s')\in I\T I}$ in 
the definition of $W$ we obtain a Coxeter group denoted by $W^*_I$. We have an 
obvious homomorphism $f:W^*_I@>>>W_I$ which takes $s$ to $s$ for $s\in I$.

\proclaim{Proposition 9.5}$f:W^*_I@>>>W_I$ is an isomorphism.
\endproclaim
We define $f':W_I@>>>W_I^*$ as follows: for $w\in W_I$ we choose a reduced
expression $w=s_1s_2\do s_q$ in $W$; then $s_i\in I$ for all $i$ (see 9.3) and we 
set $f'(w)=s_1s_2\do s_q$ (product in $W_I^*$). This map is well defined. Indeed, if
$s'_1s'_2\do s'_q$ is another reduced expression for $w$ with all $s'_i\in I$, then 
we can pass from $(s_1,s_2,\do,s_q)$ to $(s'_1,s'_2,\do,s'_q)$ by moving along edges
of the graph $X$ (see 1.9); but each edge involved in this move will necessarily 
involve only pairs $(s,s')$ in $I$, hence the equation 
$s_1s_2\do s_q=s'_1s'_2\do s'_q$ must hold in $W_I^*$. It is clear that $ff'(w)=w$ 
for all $w\in W_I$. Hence $f'$ is injective. 

We show that $f'$ is a group homomorphism. It suffices to show that 
$f'(sw)=f'(s)f'(w)$ for any $w\in W_I,s\in I$. This is clear if $l(sw)=l(w)+1$ (in 
$W$). Assume now that $l(sw)=l(w)-1$ (in $W$). Let $w=s_1s_2\do s_q$ be a 
reduced expression in $W$. Then $s_i\in I$ for all $i$. By 1.7 we have (in $W$)
$sw=s_1s_2\do s_{i-1}s_{i+1}\do s_q$ for some $i\in[1,q]$. Since 
$ss_1s_2\do s_{i-1}s_{i+1}\do s_q$ is a reduced expression for $w$ in $W$, we have
$f'(w)=ss_1s_2\do s_{i-1}s_{i+1}\do s_q$ (product in $W_I^*$). We also have 
$f'(w)=s_1s_2\do s_q$ (product in $W_I^*$). Hence
$ss_1s_2\do s_{i-1}s_{i+1}\do s_q=s_1s_2\do s_q$ (in $W_I^*$). Hence
$s_1s_2\do s_{i-1}s_{i+1}\do s_q=ss_1s_2\do s_q$ (in $W_I^*$). Hence
$f'(sw)=f'(s)f'(w)$, as required.

Since the image of $f'$ contains the generators $s\in I$ of $W_I^*$ and $f'$ is a 
group homomorphism, it follows that $f'$ is surjective. Hence $f'$ is bijective. 
Since $ff'=1$ it follows that $f$ is bijective. The proposition is proved.

\subhead 9.6\endsubhead
We identify $W_I^*$ and $W_I$ via $f$. Thus, $W_I$ is naturally a Coxeter group. Let
$l_I:W_I@>>>\bn$ be the length function of this Coxeter group. Let $w\in W_I$. Let 
$w=s_1s_2\do s_q$ be a reduced expression of $w$ (in $W$). Then $s_i\in I$ for all
$i$ (see 9.3). Hence $l_I(w)\le l(w)$. The reverse inequality $l(w)\le l_I(w)$ is 
obvious. Hence $l_I(w)=l(w)$.

From 2.4 we see that the partial order on $W_I$ defined in the same way as $\le$ on
$W$ is just the restriction of $\le$ from $W$ to $W_I$.

\proclaim{Lemma 9.7} Let $W_Ia$ be a coset in $W$. 

(a) This coset has a unique element $w$ of minimal length. 

(b) If $y\in W_I$ then $l(yw)=l(y)+l(w)$.

(c) $w$ is characterized by the property that $l(sw)>l(w)$ for all $s\in I$.
\endproclaim
Let $w$ be an element of minimal length in the coset. Let $w=s_1s_2\do s_q$ be a 
reduced expression. Let $y\in W_I$ and let $y=s'_1s'_2\do s'_p$ be a reduced 
expression in $W_I$. Then $yw=s'_1s'_2\do s'_ps_1s_2\do s_q$. By 9.1 we can drop 
some of the factors in the last product so that we are left with a reduced 
expression for $yw$. The factors dropped cannot contain any among the last $q$ since
otherwise we would find an element in $W_Ia$ of strictly smaller length than $w$. 
Thus, we can find a subsequence $i_1<i_2<\do<i_r$ of $1,2,\do,p$ such that 
$yw=s'_{i_1}s'_{i_2}\do s'_{i_r}s_1s_2\do s_q$ is a reduced expression. It follows 
that $y=s'_1s'_2\do s'_p=s'_{i_1}s'_{i_2}\do s'_{i_r}$. Since $p=l(y)$, we must have
$r=p$ so that $s'_1s'_2\do s'_ps_1s_2\do s_q$ is a reduced expression and 
$l(yw)=p+q=l(y)+l(w)$.

If now $w'$ is another element of minimal length in $W_Ia$ then $w'=yw$ for some 
$y\in W_I$. We have $l(w)=l(w')=l(y)+l(w)$ hence $l(y)=0$ hence $y=1$ and $w'=w$. 
This proves (a). Now (b) is already proved. Note that by (b), $w$ has the property 
in (c). Conversely, let $w'\in W_Ia$ be an element such that $l(sw')>l(w')$ for all
$s\in I$. We have $w'=yw$ for some $y\in W_I$. If $y\ne 1$ then for some $s\in I$ we
have $l(y)=l(sy)+1$. By (b) we have $l(w')=l(y)+l(w)$, $l(sw')=l(sy)+l(w)$. Thus 
$l(w')-l(sw')=l(y)-l(sy)=1$, a contradiction. Thus $y=1$ and $w'=w$. The lemma is 
proved.

\mpb

We shall denote by ${}^IW$ (resp. $W^I$) the set of all $w\in W$ such that $w$ has minimal length in $W_Iw$
(resp. in $wW_I$). If $I\sub J\sub S$ we write ${}^IW_J,W^I_J$ instead of ${}^I(W_J),(W_J)^I$.

\proclaim{Lemma 9.8} Let $W_Ia$ be a coset in $W$. 

(a) If $W_I$ is finite, this coset has a unique element $w$ of maximal length. If 
$W_I$ is infinite, this coset has no element of maximal length.

(b) Assume that $W_I$ is finite. If $y\in W_I$ then $l(yw)=l(w)-l(y)$.

(c) Assume that $W_I$ is finite. Then $w$ is characterized by the property that
$l(sw)<l(w)$ for all $s\in I$.
\endproclaim
Assume that $w$ has maximal length in $W_Ia$. We show that for any $y\in W_I$ we 
have 

(d) $l(yw)=l(w)-l(y)$. 
\nl
We argue by induction on $l(y)$. If $l(y)=0$, the result is clear. Assume now that 
$l(y)=p+1\ge 1$. Let $y=s_1\do s_ps_{p+1}$ be a reduced expression. By the induction
hypothesis, $l(w)=l(s_1s_2\do s_pw)+p$. Hence we can find a reduced expression of 
$w$ of the form $s_p\do s_2s_1s'_1s'_2\do s'_q$. Since $s_{p+1}\in I$, by our 
assumption on $w$ we have $l(s_{p+1}w)=l(w)-1$. Using 1.7, we deduce that either

(1) $s_{p+1}s_p\do s_{j+1}=s_p\do s_{j+1}s_j$ for some $j\in[1,p]$ or

(2) $s_{p+1}s_p\do s_2s_1s'_1s'_2\do s'_{i-1}=s_p\do s_2s_1s'_1s'_2\do s'_{i-1}s'_i$
for some $i\in[1,q]$.
\nl
In case (1) it follows that $y=s_1\do s_ps_{p+1}=s_1s_2\do s_{j-1}s_{j+1}\do s_p$,
contradicting $l(y)=p+1$. Thus, we must be in case (2). We have
$$yw=s'_1s'_2\do s'_{i-1}s'_{i+1}\do s'_q \text{ and }
l(yw)\le q-1=l(w)-p-1=l(w)-l(y).$$
Thus, $l(w)\ge l(yw)+l(y)$. The reverse inequality is obvious. Hence 
$l(w)=l(yw)+l(y)$. This completes the induction.

From (d) we see that $l(y)\le l(w)$. Thus $l:W_I@>>>\bn$ is bounded above. Hence 
there exists $y\in W_I$ of maximal length in $W_I$. Applying (d) to $y,W_I$ instead
of $w,W_Ia$ we see that 
$$l(y)=l(y'{}\i)+l(y'{}\i y)=l(y')+l(y'{}\i y)$$
for any $y'\in W_I$. Hence a reduced expression of $y'$ followed by a reduced
expression of $y'{}\i y$ gives a reduced expression of $y$. In particular $y'\le y$.
Since the set $\{y'\in W;y'\le y\}$ is finite, we see that $W_I$ is finite. 
Conversely, if $W_I$ is finite then $W_Ia$ clearly has some element of maximal 
length.

If $w'$ is another element of maximal length in $W_Ia$ then $w'=yw$ for some 
$y\in W_I$. We have $l(w)=l(w')=l(w)-l(y)$ hence $l(y)=0$ hence $y=1$ and $w'=w$. 
This proves (a) and (b). The proof of (c) is entirely similar to that of 9.7(c). The
lemma is proved.

\subhead 9.9\endsubhead 
Replacing $W,L$ by $W_I,L|_{W_I}$ in the definition of $\ch$ we obtain an
$\ca$-algebra $\ch_I$ (naturally a subalgebra of $\ch)$; instead of
$r_{x,y},p_{x,y},c_y,\mu^s_{x,y}$ we obtain for $x,y\in W_I$ elements 
$r_{x,y}^I\in\ca,p_{x,y}^I\in\ca_{\le 0},c_y^I\in\ch_I,\mu^{s,I}_{x,y}\in\ca$.

\proclaim{Lemma 9.10} Let $z\in W$ be such that $z$ is the element of minimal length
of $W_Iz$. Let $x,y\in W_I$. We have

(a) $\{u'\in W;xz\le u'\le yz\}=\{u\in W_I;x\le u\le y\}z$;

(b) $r_{xz,yz}=r^I_{x,y}$;

(c) $p_{xz,yz}=p^I_{x,y}$;

(d) $c^I_y=c_y$.

(e) If in addition, $s\in I$, $L(s)>0$ and $sx<x<y<sy$, then $sxz<xz<yz<syz$ and
$\mu^{s,I}_{x,y}=\mu^s_{xz,yz}$.
\endproclaim
We first prove the following statement.

{\it Assume that $z_1,z_2$ have minimal length in $W_Iz_1,W_Iz_2$ respectively, that
$u_1,u_2\in W_I$ and that $u_1z_1\le u_2z_2$. Then

(f) $z_1\le z_2$; if in addition, $z_1=z_2$ then $u_1\le u_2$. }
\nl
Indeed, using 2.4 we see that there exist $u'_1,z'_1$ such that

$u_1z_1=u'_1z'_1, u'_1\le u_2,z'_1\le z_2$.
\nl
Then $u'_1\in W_I$ and $z'_1\in W_Iz_1$ hence $z'_1=wz_1$ where $w\in W_I$, 
$l(z'_1)=l(w)+l(z_1)$. Hence $z_1\le z'_1$. Since $z'_1\le z_2$, we see that
$z_1\le z_2$. If we know that $z_1=z_2$, then $z'_1=z_1$ hence $u_1=u'_1$. Since 
$u'_1\le u_2$, it follows that $u_1\le u_2$ and (f) is proved.

We prove (a). If $u\in W_I$ and $x\le u\le y$, then $xz\le uz\le yz$ by 2.4 and
9.7(b). Conversely, assume that $u'\in W$ satisfies $xz\le u'\le yz$. Then $u'=uz_1$
where $z_1$ has minimal length in $W_Iu'$ and $u\in W_I$. Applying (f) to
$xz\le uz_1$ and to $uz_1\le yz$ we deduce $z\le z_1\le z$. Hence $z=z_1$. Applying
the second part of (f) to $xz\le uz$ and to $uz\le yz$ we deduce $x\le u\le y$. This
proves (a).

We prove (b) by induction on $l(y)$. Assume first that $y=1$. Then
$r^I_{x,y}=\d_{x,1}$. Now $r_{xz,z}=0$ unless $xz\le z$ (see 4.7(a)) in which case
$x=1$ and $r_{z,z}=1$. Thus, (b) holds for $y=1$. Assume now that $l(y)\ge 1$. We 
can find $s\in I$ such that $l(sy)=l(y)-1$. We have 
$$l(syz)=l(sy)+l(z)=l(y)-1+l(z)=l(yz)-1.$$
If $sx<x$ then we have (as above) $sxz<xz$. Using 4.4 and the induction hypothesis,
we have
$$r_{xz,yz}=r_{sxz,syz}=r^I_{sx,sy}=r^I_{x,y}.$$
If $sx>x$ then we have (as above) $sxz>xz$. Using 4.4 and the induction hypothesis,
we have
$$r_{xz,yz}=r_{sxz,syz}+(v_s-v_s\i)r_{xz,syz}=r^I_{sx,sy}+(v_s-v_s\i)r^I_{x,sy}
=r^I_{x,y}.$$
This completes the proof of (b).

We prove (c). Using (a), we may assume that $x\le y$ (otherwise, both sides are
zero). We argue by induction on $l(y)-l(x)\ge 0$. If $y=x$, the result is clear
(both sides are $1$). Assume now that $l(y)-l(x)\ge 1$. Using 5.3, then (a),(b) and
the induction hypothesis, we have
$$\align &\bap_{xz,yz}=\su_{u';xz\le u'\le yz}r_{xz,u'}p_{u',yz}=\su_{u\in W_I;
x\le u\le y}r_{xz,uz}p_{uz,yz}\\&=\su_{u\in W_I;x\le u\le y}r^I_{x,u}p_{uz,yz}=
\su_{u\in W_I;x<u\le y}r^I_{x,u}p^I_{u,y}+p_{xz,yz}.\endalign$$
Using 5.3 for $W_I$ we have $\bap^I_{x,y}=\su_{y;x\le u\le y}r^I_{x,u}p^I_{u,y}$.
Comparison with the previous equality gives
$$\bap_{xz,yz}-\bap^I_{x,y}=p_{xz,yz}-p^I_{x,y}.$$
The right hand side of this equality is in $\ca_{<0}$. Since it is fixed by
$\bar{}$, it must be $0$. This proves (c). Now (d) is an immediate consequence of 
(c) (with $z=1$). 

We prove (e). By 6.3(ii),
$$\su_{u';xz\le u'<yz;su'<u'}p_{xz,u'}\mu_{u',yz}^s-v_sp_{xz,yz}\in\ca_{<0}.$$
We rewrite this using (a): 
$$\su_{u\in W_I;x\le u<y;su<u}p_{xz,uz}\mu_{uz,yz}^s-v_sp_{xz,yz}\in\ca_{<0}.$$
We may assume that for all $u$ in the sum, other than $u=x$, we have
$\mu_{uz,yz}^s=\mu_{u,y}^{s,I}$. Using this and (d), we obtain
$$\mu^s_{xz,yz}+\su_{u\in W_I;x<u<y;su<u}p^I_{x,u}\mu_{u,y}^{s,I}-v_sp^I_{x,y}\in
\ca_{<0}.$$
By 6.3(ii) for $W_I$ we have
$$\mu^{s,I}_{x,y}+\su_{u\in W_I;x<u<y;su<u}p^I_{x,u}\mu_{u,y}^{s,I}-v_sp^I_{x,y}\in
\ca_{<0}.$$
It follows that $\mu^s_{xz,yz}-\mu^{s,I}_{x,y}\in\ca_{<0}$. On the other hand,
$\mu^s_{xz,yz}-\mu^{s,I}_{x,y}$ is fixed by $\bar{}$ (see 6.3(i)) hence it is $0$.
This proves (e). The lemma is proved.

\proclaim{Proposition 9.11} Assume that $L(s)>0$ for all $s\in I$.

(a) Let $z\in W$ be such that $z$ is the element of minimal length of $W_Iz$. If 
$x,y$ in $W_I$ satisfy $x\le_\cl y$ (relative to $W_I$), then $xz\le_\cl yz$ (in 
$W$). If $x,y$ in $W_I$ satisfy $x\si_\cl y$ (relative to $W_I$), then 
$xz\si_\cl yz$ (in $W$). 

(b)  Let $z\in W$ be such that $z$ is the element of minimal length of $zW_I$. If
$x,y$ in $W_I$ satisfy $x\le_\car y$ (relative to $W_I$), then $zx\le_\car zy$ (in 
$W$). If $x,y$ in $W_I$ satisfy $x\si_\car y$ (relative to $W_I$), then 
$zx\si_\car zy$ (in $W$). 
\endproclaim
We prove the first assertion of (a). We may assume that $x\gt_\cl y$ (relative to
$W_I$) and $x\ne y$. Thus, there exists $s\in I$ such that $sy>y$, $sx<x$ and we 
have either $x=sy$ or $x<y$ and $\mu^s_{x,y}\ne 0$. If $x=sy$, then $sxz<xz=syz>yz$,
hence $xz\gt_\cl yz$ (in $W$). Thus, we may assume that $x<y$ and 
$\mu^s_{x,y}\ne 0$. By 9.10(e) we then have $\mu^s_{xz,yz}\ne 0$, hence 
$xz\gt_\cl yz$ (in $W$). The first assertion of (a) is proved. The second assertion
of (a) follows from the first. (b) follows by applying (a) to $z\i,x\i,y\i$ instead
of $z,x,y$.

\subhead 9.12\endsubhead 
Assume that $z\in W$ is such that $W_Iz=zW_I$ and $z$ is the element of minimal
length of $W_Iz=zW_I$. Then $y\m z\i yz$ is an automorphism of $W_I$. If $s\in I$ 
then, by 9.7, we have $l(sz)=l(s)+l(z)=1+l(z)$; by 9.7 applied to $W_Iz\i$ instead 
of $W_Iz$ we have $l((z\i sz)z\i)=l(z\i sz)+l(z\i)$ hence 
$l(z\i s)=l(z\i sz)+l(z\i)$; since $l(z\i s)=l(sz)$ and $l(z\i)=l(z)$, it follows 
that $l(z\i sz)+l(z\i)=1+l(z)$, hence $l(z\i sz)=1$. We see that $y\m z\i yz$ maps 
$I$ onto itself hence it is an automorphism of $W_I$ as a Coxeter group. This 
automorphism preserves the function $L|_{W_I}$. Indeed, if $y\in W_I$, then
$$l(zyz\i)+l(z)=l((zyz\i)z)=l(zy)=l(y\i z\i)=l(y\i)+l(z\i)=l(y)+l(z)$$
(by 9.7 applied to $W_Iz$ and to $W_Iz\i$) hence 
$$\align&L(zyz\i)+L(z)=L((zyz\i)z)=L(zy)=L(y\i z\i)=L(y\i)+L(z\i)\\&=L(y)+L(z),
\endalign$$
so that $L(zyz\i)=L(y)$. In particular, this automorphism respects the preorders 
$\le_\cl,\le_\car,\le_{\lr}$ of $W_I$ (defined in terms of $L|_{W_I}$) and the 
associated equivalence relations.

\proclaim{Proposition 9.13} Assume that $L(s)>0$ for all $s\in I$. Let $z$ be as in
9.12. If $x,y$ in $W_I$ satisfy $x\le_{\lr}y$ (relative to $W_I$), then
$xz\le_{\lr}yz$ (in $W$). If $x,y$ in $W_I$ satisfy $x\si_{\lr}y$ (relative to
$W_I$), then $xz\si_{\lr}yz$ (in $W$).
\endproclaim
We prove the first assertion. We may assume that either $x\le_\cl y$ (in $W_I)$ or
$x\le_\car y$ (in $W_I$). In the first case, by 9.11(a), we have $xz\le_\cl yz$ (in
$W$) hence $xz\le_{\lr}yz$ (in $W$). In the second case, by 9.12, we have 
$z\i xz\le_\car z\i yz$. Applying 9.11(b) to $z\i xz,z\i yz$ instead of $x,y$ we see
that $xz\le_\car yz$ (in $W$) hence $xz\le_{\lr}yz$ (in $W$). This proves the first
assertion. The second assertion follows from the first.

\subhead 9.14\endsubhead
In the remainder of this chapter we fix two subsets $K,K'$ of $S$ and a $(W_K,W_{K'})$-double coset $\Om$ in
$W$. We have the following result.

\proclaim{Proposition 9.15} (a) $\Om$ contains a unique element $b$ of minimal length.

(b) Setting $J=K\cap bK'b\i$, we have $W_K\cap bW_{K'}b\i=W_J$.

(c) The map $W_K^J\T W_{K'}@>>>\Om$, $(a,c)\m abc$ is a bijection. 

(d) For any $a\in W_K^J$, $c\in W_{K'}$, we have $l(abc)=l(a)+l(b)+l(c)$.   

(e) If $W_K$ and $W_{K'}$ are finite then $\Om$ contains a unique element $\tb$ of maximal length. We have 
$\tb=w_0^Kw_0^Jbw_0^{K'}$ where $w_0^K,w_0^J,w_0^{K'}$ is the unique element of maximal length of 
$W_K,W_J,W_{K'}$
respectively.
\endproclaim
Note that (a) is stated in \cite{\BO, Ch.IV,\S1, Ex.3}; (b)-(d) are due to Kilmoyer \cite{\KI} under the
assumption that $W$ is a Weyl group. The proof of the proposition is given in 9.16.

\subhead 9.16\endsubhead
Let $b$ be an element of minimal length in $\Om$ (it clearly exists). We fix a reduced expression 
$s_1s_2\do s_q$ for $b$. We show:

(a) {\it If $w\in\Om$, then $w$ has a reduced expression of the form
$$s'_1s'_2\do s'_ps_1s_2\do s_q\ts_1\ts_1\do\ts_r$$ 
where $s'_i\in K$ for $i\in[1,p]$, $\ts_i\in K'$ for 
$i\in[1,r]$.}
\nl
We can find a not necessarily reduced expression $w=s'_1s'_2\do s'_ps_1s2\do s_q\ts_1\ts_2\do\ts_r$ with 
$s'_i,\ts_i$ as in the Lemma. Using 9.1, we can drop some of the simple reflections in this expression so 
that the resulting expression is a reduced expression for $w$; none of the dropped reflections can be among 
$s_1,s_2,\do,s_p$ (otherwise the minimality of $l(b)$ would be contradicted). The result follows.

From (a) we see that $b$ above is unique (hence 9.15(a) holds) and any $w\in\Om$ is of the form (b) $w=abc$ with $a\in W_K$, $c\in W_{K'}$, $l(w)=l(a)+l(b)+l(c)$.
\nl
We show:

(c) {\it Let $y\in W_K\cap bW_{K'}b\i$ and let $y=s'_1s'_2\do s'_p$ be a reduced expression in $W_K$. Then 
for any $i\in[1,p]$ we have $b\i s'_ib\in W_{K'}$.}
\nl
We argue by induction on $p$. When $p=0$ there is nothing to prove. We now assume that $p\ge1$.
Let $\ty=b\i yb\in W_{K'}$. Since $b\in W^{K'}$ and $\ty\in W_{K'}$, we see from 
9.7 that $l(b\ty)=l(b)+l(\ty)$. Similarly, since $b\in{}^KW$ and 
$y\in W_K$, we have $l(yb)=l(y)+l(b)$. Since $yb=b\ty$ it follows that $l(y)+l(b)=l(b)+l(\ty)$ hence 
$l(y)=l(\ty)$.
Let $y=s'_1s'_2\do s'_p$ be a reduced expression in $W_K$ (it is also a reduced expression in $W$).
Let $\ty=\ts_1\ts_2\do\ts_p$ be a reduced expression in $W_{K'}$ (it is also a reduced expression in $W$).
Then $s'_1s'_2\do s'_ps_1s_2\do s_q$ and
$s_1s_2\do s_q\ts_1\ts_2\do\ts_p$ are reduced expressions in $W$ for the same element $yb=b\ty$.
Now $$s'_1(s_1s_2\do s_q\ts_1\ts_2\do\ts_p)<s_1s_2\do s_q\ts_1\ts_2\do\ts_p$$ hence by 1.7
we have either $$s'_1s_1s_2\do s_{i-1}=s_1s_2\do s_{i-1}s_i$$
 for some $i\in[1,q]$
or 
$$s'_1s_1s_2\do s_q\ts_1\ts_2\do\ts_{j-1}=s_1s_2\do s_q\ts_1\ts_2\do\ts_{j-1}s_j$$ for some $j\in[1,p]$.
In the first case we have 
$$s'_2\do s'_ps_1s_2\do s_q=s_1\do s_{i-1}s_{i+1}\do s_q\ts_1\ts_2\do\ts_p$$,
so that $\Om$ contains the element $s_1\do s_{i-1}s_{i+1}\do s_q$ of length $<q$, contradicting the 
definition of $b$. Thus we must be in the second case so that
$$s'_2\do s'_ps_1s_2\do s_q=s_1s_2\do s_q\ts_1\do \ts_{j-1}\ts_{j+1}\do\ts_p.$$
 We then have
$s'_2\do s'_pb\in bW_{K'}$. We set $y'=s'_2\do s'_p$. We have $y'\in W_K\cap bW_{K'}b\i$ and $l(y')=p-1$. 
By the induction hypothesis we have $b\i s'_ib\in W_{K'}$ for $i\in[2,p]$. We have
$$b\i s'_1b=(b\i yb)(b\i s'_pb)\do(b\i s'_22b)\in W_{K'}.$$
 This proves (c).

We show:

(d) {\it Let $s'\in K$ be such that $b\i s'b\in W_{K'}$. Then $b\i s'b\in K'$.}
\nl
We have $s'b=bz$ where $z\in W_{K'}$. The proof of the equality $l(y)=l(\ty)$ in (c) can be repeated with 
$y,\ty$ replaced by $s',z$; it yields $l(z)=l(s')=1$ hence $z\in K'$ as required.

Now 9.15(b) follows immediately from (c),(d).

We show:

(e) {\it Let $w\in\Om$. We can find $a\in W_K^J,c\in W_{K'}$ so that (b) holds.}
\nl
By (b) we can write $w=\ta b\tc$ with $\ta\in W_K,\tc\in W_{K'}$, $l(w)=l(\ta)+l(b)+l(\tc)$. Using
9.7 with $W,W_I$ replaced by $W_K,W_J$, we see that we can write $\ta=az$ where
$a\in W_K^J$, $z\in W_J$ and we have $l(\ta)=l(a)+l(z)$. We have
$w=azb\tc$ where $b\i zb\in W_{K'}$. Moreover, since $b\i Jb\sub K'$ we have $l(z)=l(b\i zb)$
and $l(w)=l(a)+l(z)+l(b)+l(\tc)=l(a)+l(b)+l(b\i zb)+l(\tc)$. From $l(w)=l(a)+l(b)+l(b\i zb)+l(\tc)$, 
$w=ab(b\i zb)\tc$, we deduce that $l(b\i zb)+l(\tc)=l(c)$ where $c=(b\i zb)\tc$. We have $c\in W_{K'}$ and 
$w=abc$, $l(w)=l(a)+l(b)+l(c)$. This proves (e).

We show:

(f) {\it Let $a,a'\in W_K$, $c,c'\in W_{K'}$ be such that $a\in W_K^J,a'\in W_K^J$ and $abc=a'bc'$. Then 
$a=a'$ and $c=c'$.}
\nl
We have $a\i a'=bcc'{}\i b\i\in W_K\cap bW_{K'}b\i=W_J$ (see 9.15(b)). Thus $a'\in aW_J$. Since $a',a$ 
belong to $W_K^J$, we have $a=a'$ (we use 9.7(a)). It follows that $bc=bc'$ hence $c=c'$. This proves (f).

From (e),(f) we see that 9.15(c) holds.

We show:

(g) {\it Let $a\in W_K$, $c\in W_{K'}$ be such that $a\in W_K^J$. Let $w=abc$. Then $l(w)=l(a)+l(b)+l(c)$.} 
\nl
By (e) we can write $w=a'bc'$ where $a'\in W_K^J,c'\in W_{K'}$ and $l(w)=l(a')+l(b)+l(c')$. We have 
$abc=a'bc'$. By (f) we have $a=a'$ and $c=c'$. This proves (g). 

We prove 9.15(e). Let $w\in\Om$. Write $w=abc$ with $a\in W_K^J$, $c\in W_{K'}$. 
We have $l(w_0^K)\ge l(aw_0^J)=l(a)+l(w_0^J)$ hence $l(a)\le l(w_0^K)-l(w_0^J)=l(w_0^Kw_0^J)$ and 
$l(c)\le l(w_0^{K'})$
hence $l(w)=l(a)+l(b)+l(c)\le l(w_0^Kw_0^J)+l(b)+l(w_0^{K'})=
l(w_0^Kw_0^Jbw_0^{K'})$. Moreover if $l(w)=l(w_0^Kw_0^Jbw_0^{K'})$ 
then $w_0^K=aw_0^J$, $c=w_0^{K'}$. This proves 9.15(e). Proposition 9.15 is proved.

\head 10. Inversion \endhead
\subhead 10.1\endsubhead
We preserve the setup of 3.1. For $y,w\in W$ we set
$$q'_{y,w}=\su(-1)^np_{z_0,z_1}p_{z_1,z_2}\do p_{z_{n-1},z_n}\in\ca$$
(sum over all sequences $y=z_0<z_1<z_2<\do<z_n=w$ in $W$) and
$$q_{y,w}=\sg(y)\sg(w)q'_{y,w}.$$
We have

$q_{w,w}=1$,

$q_{y,w}\in\ca_{<0}$ if $y\ne w$,

$q_{y,w}=0$ unless $y\le w$.

\proclaim{Proposition 10.2} For any $y,w \in W$ we have
$\ov q_{y,w}=\su_{z;y\le z\le w}q_{y,z}r_{z,w}$.
\endproclaim
The (triangular) matrices $Q'=(q'_{y,w}),P=(p_{y,w}),R=(r_{y,w})$ are related by

(a) $Q'P=PQ'=1$, $\ov P=RP$, $\ov RR=R\ov R=1$
\nl
where $\bar{}$ over a matrix is the matrix obtained by applying $\bar{}$ to each 
entry. (Although the matrices may be infinite, the products are well defined as each
entry of a product is obtained by finitely many operations.) The last three 
equations in (a) are obtained from 5.3, 4.6; the equations involving $Q'$ follow 
from the definition. From (a) we deduce $Q'P=1=\ov Q'\ov P=\ov Q'RP$. Hence 
$Q'P=\ov Q'RP$. Multiplying on the right by $Q'$ and using $PQ'=1$ we deduce 
$Q'=\ov Q'R$. Multiplying on the right by $\ov R$ gives 

(b) $\ov Q'=Q'\ov R$.
\nl
Let $\ss$ be the matrix whose $y,w$ entry is $\sg(y)\d_{y,w}$. We have $\ss^2=1$.
Let $Q$ be the triangular matrix $(q_{y,w})$. Note that $Q=\ss Q'\ss$. By 4.5 we 
have $\ov R=\ss R\ss$. Hence by multiplying the two sides of (b) on the left and
right by $\ss$ we obtain $\ov Q=QR$. The proposition is proved.

\subhead 10.3\endsubhead
Define an $\ca$-linear map $\t:\ch@>>>\ca$ by $\t(T_w)=\d_{w,1}$ for $w\in W$. 

\proclaim{Lemma 10.4} (a) For $x,y\in W$ we have $\t(T_xT_y)=\d_{xy,1}$.

(b) For $h,h'\in\ch$ we have $\t(hh')=\t(h'h)$.

(c) Assume that $L(s)\ge0$ for all $s\in S$. Let $x,y,z\in W$ and let $M=\min(L(x),L(y),L(z))$. We 
have $\t(T_xT_yT_z)\in v^M\bz[v\i]$.
\endproclaim
We prove (a) by induction on $l(y)$. If $l(y)=0$, the result is clear. Assume now 
that $l(y)\ge 1$. If $l(xy)=l(x)+l(y)$ then $T_xT_y=T_{xy}$ and the result is clear.
Hence we may assume that $l(xy)\ne l(x)+l(y)$. Then $l(xy)<l(x)+l(y)$. Let 
$y=s_1s_2\do s_q$ be a reduced expression. We can find $i\in[1,q]$ such that 
$$l(x)+i-1=l(xs_1s_2\do s_{i-1})>l(xs_1s_2\do s_{i-1}s_i).\tag d$$
We show that
$$xs_1s_2\do s_{i-1}s_{i+1}\do s_q\ne 1.\tag e$$
If (e) does not hold, then $x=s_q\do s_{i+1}s_{i-1}\do s_1$, so that
$$l(xs_1s_2\do s_{i-1})=l(s_q\do s_{i+1}s_{i-1}\do s_1s_1\do s_{i-1})=
l(s_q\do s_{i+1})=q-i,$$
$$l(xs_1s_2\do s_{i-1}s_i)=l(s_q\do s_{i+1}s_{i-1}\do s_1s_1\do s_i)=
l(s_q\do s_{i+1}s_i)=q-i+1,$$
contradicting (d). Thus (e) holds. We have 
$$\align&\t(T_xT_y)=\t(T_xT_{s_1}T_{s_2}\do T_{s_q})=
\t(T_{xs_1s_2\do s_{i-1}}T_{s_i}T_{s_{i+1}\do s_q})\\&=  
\t(T_{xs_1s_2\do s_{i-1}s_i}T_{s_{i+1}\do s_q})+  
(v_s-v_s\i)\t(T_{xs_1s_2\do s_{i-1}}T_{s_{i+1}\do s_q}).\endalign$$
By the induction hypothesis and (e), this equals
$$\d_{xs_1s_2\do s_{i-1}s_is_{i+1}\do s_q,1}+  
(v_s-v_s\i)\d_{xs_1s_2\do s_{i-1}s_{i+1}\do s_q,1}=\d_{xy,1}.$$
This completes the proof of (a). To prove (b), we may assume that $h=T_x,h'=T_y$ for
$x,y\in W$; we then use (a) and the obvious equality $\d_{xy,1}=\d_{yx,1}$. 

We prove (c). Using (b) we see that $\t(T_xT_yT_z)=\t(T_yT_zT_x)=\t(T_zT_xT_y)$. 
Hence it is enough to show that, for any $x,y,z$ we have
$$\t(T_xT_yT_z)\in v^{L(x)}\bz[v\i].$$
We argue by induction on $l(x)$. If $l(x)=0$, then $x=1$ and the result follows from
(a). Assume now that $l(x)\ge 1$. We can find $s\in S$ such that $xs<x$. If $sy>y$,
then by the induction hypothesis,
$$\t(T_xT_yT_z)=\t(T_{xs}T_{sy}T_z)\in v^{L(x)-L(s)}\bz[v\i]\sub v^{L(x)}\bz[v\i].$$
If $sy<y$, then by the induction hypothesis,
$$\align&\t(T_xT_yT_z)=\t(T_{xs}T_{sy}T_z)+(v_s-v_s\i)\t(T_{xs}T_yT_z)\\&
\in v^{L(x)-L(s)}\bz[v\i]+v_sv^{L(x)-L(s)}\bz[v\i]\sub v^{L(x)}\bz[v\i].\endalign$$
The lemma is proved.

\subhead 10.5\endsubhead
Let $\ch'=\Hom_\ca(\ch,\ca)$. We regard $\ch'$ as a left $\ch$-module where, for 
$h\in\ch,\ph\in\ch'$ we have $(h\ph)(h_1)=\ph(h_1h)$ for all $h_1\in\ch$ and as a
right $\ch$-module where, for $h\in\ch,\ph\in\ch'$, we have 
$(\ph h)(h_1)=\ph(hh_1)$ for all $h_1\in\ch$. 

\subhead 10.6\endsubhead
We sometimes identify $\ch'$ with the set of all formal sums $\su_{x\in W}a_xT_x$ 
with $a_x\in\ca$; to $\ph\in\ch'$ corresponds the formal sum
$\su_{x\in W}\ph(T_{x\i})T_x$. Since $\ch$ is contained in the set of such formal
sums (it is the set of sums such that $a_x=0$ for all but finitely many $x$), we see
that $\ch$ is naturally a subset of $\ch'$. Using 10.4(a) we see that the imbedding
$\ch\sub\ch'$ is an imbedding of $\ch$-bimodules; it is an equality if $W$ is 
finite.

\subhead 10.7\endsubhead
Let $z\in W$. Recall that in 8.1 we have defined $D_z\in\ch'$ by $D_z(c_w)=\d_{z,w}$
for all $w$. An equivalent definition is

(a) $D_z(T_y)=q'_{z,y}$ 
\nl
for all $y\in W$. Indeed, assuming that (a) holds, we have
$$D_z(c_w)=\su_yq'_{z,y}p_{y,w}=\d_{z,w}.$$

\proclaim{Proposition 10.8} Let $z\in W,s\in S$. Assume that $L(s)>0$.

(a) If $zs<z$, then $c_sD_z=(v_s+v_s\i)D_z+D_{zs}+\su_{u;z<u<us}\mu^s_{z\i,u\i}D_u$.

(b) If $zs>z$, then $c_sD_z=0$.
\endproclaim
For $a,b\in W$ we define $\d_{a<b}$ to be $1$ if $a<b$ and $0$ otherwise. Let
$w\in W$. If $ws>w$, then by 6.7(a), we have
$$\align&(c_sD_z)(c_w)=D_z(c_wc_s)=D_z(c_{ws}+\su\Sb x\\xs<x<w\eSb
\mu_{x\i,w\i}^sc_x)\\&=\d_{z,ws}+\su\Sb x\\xs<x<w\eSb\mu_{x\i,w\i}^s\d_{z,x}.\tag c
\endalign$$
If $ws<w$, then by 6.7(b), we have
$$(c_sD_z)(c_w)=D_z(c_wc_s)=(v_s+v_s\i)D_z(c_w)=(v_s+v_s\i)\d_{z,w}.\tag d$$
If $zs<z,ws>w$, then by (c):
$$(c_sD_z)(c_w)=\d_{zs,w}+\d_{z<w}\mu_{z\i,w\i}^s
=(v_s+v_s\i)D_z+D_{zs}+\su\Sb u\\z<u<us\eSb \mu^s_{z\i,u\i}D_u)(c_w).$$
If $zs<z, ws<w$, then by (d):
$$(c_sD_z)(c_w)=(v_s+v_s\i)\d_{z,w}=(v_s+v_s\i)D_z+D_{zs}+
\su\Sb u\\ z<u<us\eSb\mu^s_{z\i,u\i}D_u)(c_w).$$
If $zs>z, ws>w$, then by (c), we have $(c_sD_z)(c_w)=0$. If $zs>z,ws<w$, then by 
(d), we have $(c_sD_z)(c_w)=0$. Since $(c_w)$ is an $\ca$-basis of $\ch$, the 
proposition follows.

\subhead 10.9\endsubhead
We show that
$$\{\pm c_w; w\in W\}=\{h\in\ch;\t(hh^\flat)\in 1+\ca_{<0};\bar h=h\}.\tag a$$
For any $w,w'\in W$ we have 
$$\t(c_wc_{w'}^\flat)=\t(\su_{y,y'}p_{y,w}p_{y',w'}T_yT_{y'{}\i})=
\su_{y,y'}p_{y,w}p_{y',w'}\d_{y,y'}=\d_{w,w'}+z_{w,w'}\tag b$$
where $z_{w,w'}\in\ca_{<0}$. In particular, $\t(c_wc_w^\flat)\in 1+\ca_{<0}$. 

Conversely, assume that $h\in\ch$ satisfies $\t(hh^\flat)\in 1+\ca_{<0},\bar h=h$. 
We have $h=\su_{w\in W}x_wc_w$ where $x_w\in\ca$ are $0$ for all but finitely many 
$w$. We can find $t\in\bz$ such that $x_w=b_wv^t\mod\ca_{<t}$ where $b_w\in\bz$ for
all $w$ and $b_w\ne 0$ for some $w$. Using (b), we have 
$$\align&\t(hh^\flat)=\t(\su_{w,w'}x_wx_{w'}c_wc_{w'}^\flat)=\su_{w,w'}x_wx_{w'}
(\d_{w,w'}+z_{w,w'})\\&=\su_wx_w^2+\su_{w,w'}x_wx_{w'}z_{w,w'}.\endalign$$
This equals $\su_wb_w^2v^{2t}$ modulo $\ca_{<2t}$ and also equals $1$ modulo 
$\ca_{<0}$. It follows that $t=0$ and $\su_wb_w^2=1$. Since $b_w$ are integers,
there exists $u\in W$ such that $b_u=\pm 1$ and $b_w=0$ for $w\ne u$. Thus 
$x_w\in\ca_{\le 0}$ for all $w$. Since $\bar h=h$ we have $\bar x_w=x_w$ for all 
$w$. It follows that $x_w=b_w$ for all $w$. Thus $h=\pm c_u$. This proves (a).

\subhead 10.10\endsubhead
The interest of 10.9(a) is that it provides a definition of $c_w$ (up to sign)
without using the basis $T_w$ of $\ch$ (instead, it uses $\t:\ch@>>>\ca$). The 
equality 10.9(a) could be used to give a definition of $c_w$ (up to sign) in more 
general situations than that considered above, when the basis $(T_w)$ is not defined
but $\t:\ch@>>>\ca$ is defined (it is known that $\t$ is defined when $W$ is 
replaced by certain complex reflection groups, see \cite{\BM}). This should lead to
a definition of cells for complex reflection groups.

\head 11. The longest element for a finite $W$ \endhead
\subhead 11.1\endsubhead
We preserve the setup of 3.1. Let $I\sub S$ be such that $W_I$ is finite. By 9.8, 
there is a unique element of maximal length of $W_I$. We denote it by $w^I_0$. If 
$w_1$ has minimal length in $W_Ia$ then $w^I_0w_1$ has maximal length in $W_Ia$.

\subhead 11.2\endsubhead   
In the remainder of this chapter we assume that $W$ is finite. Then
$w_0:=w^S_0$, the unique element of maximal length of $W$, is well defined. Since 
$l(w_0\i)=l(w_0)$, we must have $w_0\i=w_0$. By the argument in the proof of 9.8 we
have $w\le w_0$ for any $w\in W$. By 9.8 we have 

(a) $l(ww_0)=l(w_0)-l(w)$
\nl
for any $w\in W$. Applying this to $w\i$ and using the equalities
$l(w\i w_0)=l(w_0\i w)=l(w_0w),l(w\i)=l(w)$, we deduce that

(b) $l(w_0w)=l(w_0)-l(w)$.
\nl
We can rewrite (a),(b) as $l(w_0)=l(w\i)+l(ww_0)$, $l(w_0)=l(w_0w)+l(w\i)$. Using 
this and the definition of $L$ we deduce that

$L(w\i)+L(ww_0)=L(w_0)=L(w_0w)+L(w\i)$, 
\nl
hence $L(ww_0)=L(w_0w)$. This implies $L(w_0ww_0)=L(w)$ for all $w$. Replacing $L$ 
by $l$ gives $l(w_0ww_0)=l(w)$. Thus, the involution $w\m w_0ww_0$ of $W$ maps $S$ 
into itself hence is a Coxeter group automorphism preserving the function $L$. 

\proclaim{Lemma 11.3} Let $y,w\in W$. We have

(a) $y\le w\Lra w_0w\le w_0y\Lra ww_0\le yw_0$;

(b) $r_{y,w}=r_{ww_0,yw_0}=r_{w_0w,w_0y}$;

(c) $\bap_{ww_0,yw_0}=\su_{z; y\le z\le w}p_{zw_0,yw_0}r_{z,w}$.
\endproclaim
We prove (a). To prove that $y\le w\imp ww_0\le yw_0$, we may assume that
$l(w)-l(y)=1,yw\i\in T$. Then 
$$l(yw_0)-l(ww_0)=l(w_0)-l(y)-(l(w_0)-l(w))=l(w)-l(y)=1$$
and $(ww_0)(yw_0)\i=wy\i\in T$. Hence $ww_0\le yw_0$. The opposite implication is
proved in the same way. The second equivalence in (a) follows from the last sentence
in 11.2.

We prove the first equality in (b) by induction on $l(w)$. If $l(w)=0$ then $w=1$. 
We have $r_{y,1}=\d_{y,1}$. Now $r_{w_0,yw_0}$ is zero unless $w_0\le yw_0$ (see 
4.7). On the other hand we have $yw_0\le w_0$ (see 11.2). Hence $r_{w_0,yw_0}$ is 
zero unless $yw_0=w_0$, that is unless $y=1$ in which case it is $1$. Thus the 
desired equality holds when $l(w)=0$. Assume now that $l(w)\ge 1$. We can find 
$s\in S$ such that $sw<w$. Then $sww_0>ww_0$ by (a).

Assume first that $sy<y$ (hence $syw_0>yw_0$). By 4.4 and the induction hypothesis 
we have
$$r_{y,w}=r_{sy,sw}=r_{sww_0,syw_0}=r_{ww_0,yw_0}.$$
Assume next that $sy>y$ (hence $syw_0<yw_0$.) By 4.4 and the induction hypothesis we
have
$$\align&r_{y,w}=r_{sy,sw}+(v_s-v_s\i)r_{y,sw}
=r_{sww_0,syw_0}+(v_s-v_s\i)r_{sww_0,yw_0}\\&
=r_{sww_0,syw_0}+(v_s-v_s\i)r_{ww_0,syw_0}=r_{ww_0,yw_0}.\endalign$$
This proves the first equality in (b). The second equality in (b) follows from the 
last sentence in 11.2. 

We prove (c). We may assume that $y\le w$. By 5.3 (for $ww_0,yw_0$ instead of $y,w$)
we have $\bap_{ww_0,yw_0}=\su_{z;y\le z\le w}r_{ww_0,zw_0}p_{zw_0,yw_0}$ (we have
used (a)). Here we substitute $r_{ww_0,zw_0}=r_{z,w}$ (see (b)) and the result 
follows. The lemma is proved.

\proclaim{Proposition 11.4} For any $y,w\in W$ we have $q_{y,w}=p_{ww_0,yw_0}=
p_{w_0w,w_0y}$.
\endproclaim
The second equality follows from the last sentence in 11.2. We prove the first
equality. We may assume that $y\le w$. We argue by induction on $l(w)-l(y)\ge 0$. If
$l(w)-l(y)=0$ we have $y=w$ and both sides are $1$. Assume now that 
$l(w)-l(y)\ge 1$. Subtracting the identity in 11.3(c) from that in 10.2 and using 
the induction hypothesis, we obtain
$$\ov q_{y,w}-\bap_{ww_0,yw_0}=q_{y,w}-p_{ww_0,yw_0}.$$
The right hand side is in $\ca_{<0}$; since it is fixed by $\bar{}$, it is $0$. The
proposition is proved.

\proclaim{Proposition 11.5} We identify $\ch=\ch'$ as in 10.6. If $z\in W$, then 
$D_{z\i}\in\ch'$ (see 10.7) becomes an element of $\ch$. We have
$D_{z\i}T_{w_0}\i=\sg(zw_0)c_{zw_0}^\da$, ${}^\da$ as in 3.5. 
\endproclaim
By definition, $D_{z\i}\in\ch$ is characterized by
$$\t(D_{z\i}T_{y\i})=q'_{z\i,y\i}$$
for all $y\in W$. Here $\t$ is as in 10.3. Hence, by 10.4(a), we have 
$D_{z\i}=\su_yq'_{z\i,y\i}T_y$. Using 11.4, we deduce
$$D_{z\i}=\su_y\sg(yz)p_{w_0y\i,w_0z\i}T_y.$$
Multiplying on the right by $T_{w_0}\i$ gives
$$D_{z\i}T_{w_0}\i=\su_y\sg(yz)p_{w_0y\i,w_0z\i}T_{w_0y\i}\i$$
since $T_{w_0y\i}T_y=T_{w_0}$. On the other hand,
$$\sg(zw_0)c_{zw_0}^\da=\su_x\sg(zw_0x)p_{x,zw_0}T_{x\i}\i
=\su_y\sg(zw_0yw_0)p_{yw_0,zw_0}T_{w_0y\i}\i.$$
We now use the identity $p_{yw_0,zw_0}=p_{w_0y\i,w_0z\i}$. The proposition follows.

\proclaim{Proposition 11.6} Let $u,z\in W, s\in S$ be such that $sz<z<u<su$. Assume
that $L(s)>0$. Then $suw_0<uw_0<zw_0<szw_0$ and
$\mu_{uw_0,zw_0}^s=-\sg(uz)\mu^s_{z,u}$.
\endproclaim
Using 10.8(a), we see that
$$(c_s-(v_s+v_s\i))D_{z\i}T_{w_0}\i=
D_{z\i s}T_{w_0}\i+\su_{u;z\i<u\i<u\i s}\mu^s_{z,u}D_{u\i}T_{w_0}\i,$$
hence, using 11.5, we have
$$(c_s-(v_s+v_s\i))\sg(zw_0)c_{zw_0}^\da=\sg(zsw_0)c_{szw_0}^\da+
\su_{u;z<u<su}\mu^s_{z,u}\sg(uw_0)c_{uw_0}^\da.$$
Applying ${}^\da$ to both sides and using $(c_s-(v_s+v_s\i))^\da=-c_s$ gives
$$-c_sc_{zw_0}=-c_{szw_0}+\su_{u;z<u<su}\mu^s_{z,u}\sg(uz)c_{uw_0}.\tag a$$
Since $szw_0>zw_0$, we can apply 6.6(a) and we get
$$c_sc_{zw_0}=c_{szw_0}+\su_{u';su'<u'<zw_0}\mu_{u',zw_0}^sc_{u'}$$
or equivalently
$$c_sc_{zw_0}=c_{szw_0}+\su_{u;z<u<su}\mu_{uw_0,zw_0}^sc_{uw_0}.$$
Comparison with (a) gives
$$-\su_{u;z<u<su}\mu^s_{z,u}\sg(uz)c_{uw_0}=\su_{u;z<u<su}\mu_{uw_0,zw_0}^sc_{uw_0};
$$
the proposition follows.

\proclaim{Corollary 11.7}Assume that $L(s)>0$ for all $s\in S$. Let $y,w\in W$.

(a) $y\le_\cl w\Lra ww_0\le_\cl yw_0\Lra w_0w\le_\cl w_0y$;

(b) $y\le_\car w\Lra ww_0\le_\car yw_0\Lra w_0w\le_\car w_0y$;

(c) $y\le_{\lr}w\Lra ww_0\le_{\lr}yw_0\Lra w_0w\le_{\lr}w_0y$.

(d) Left multiplication by $w_0$ carries left cells to left cells, right cells to 
right cells, two-sided cells to two-sided cells. The same holds for right
multiplication by $w_0$.
\endproclaim
We prove the first equivalence in (a). It is enough to show that
$y\le_\cl w\imp ww_0\le_\cl yw_0$. We may assume that $y\gt_\cl w$ and $y\ne w$.
Then there exists $s\in S$ such that $sw>w,sy<y$ and $D_y(c_sc_w)\ne 0$. We have 
$syw_0>yw_0, sww_0<ww_0$. From 6.6 we see that either $y=sw$ or $y<w$ and 
$\mu^s_{y,w}\ne 0$. In the first case we have $ww_0=syw_0$; in the second case we 
have $ww_0<yw_0$ and $\mu^s_{ww_0,yw_0}\ne 0$ (see 11.6). In both cases, 6.6 shows 
that $D_{ww_0}(c_sc_{yw_0})\ne 0$. Hence $ww_0\le_\cl yw_0$. Thus, the first 
equivalence in (a) is established. The second equivalence in (a) follows from the 
last sentence in 11.2.

Now (b) follows by applying (a) to $y\i,w\i$ instead of y,w; (c) follows from (a) 
and (b); (d) follows from (a),(b),(c). The corollary is proved.

\head 12. Examples of elements $D_w$\endhead
We preserve the setup of 3.1.

\proclaim{Proposition 12.1} Assume that $L(s)>0$ for all $s\in S$. For any $y\in W$
we have $D_1(T_y)=\sg(y)v^{-L(y)}$. Equivalently, with the identification in 10.6,
we have

$D_1=\su_{y\in W}\sg(y)v^{-L(y)}T_y$.
\endproclaim
An equivalent statement is that $q'_{1,y}=\sg(y)v^{-L(y)}$. Since $q'_{1,y}$ are 
determined by the equations $\su_yq'_{1,y}p_{y,w}=\d_{1,w}$ (see 10.2(a)) it is 
enough to show that
$$\su_y\sg(y)v^{-L(y)}p_{y,w}=\d_{1,w}$$
for all $w\in W$. If $w=1$ this is clear. Assume now that $w\ne 1$. We can find
$s\in S$ such that $sw<w$. We must prove that
$$\su_{y;y<sy}\sg(y)v^{-L(y)}(p_{y,w}-v_s\i p_{sy,w})=0.$$
Each term of the last sum is $0$, by 6.6(c). The proposition is proved.

\proclaim{Corollary 12.2} Assume that $W$ is finite and that $L(s)>0$ for all 
$s\in S$. Then 

$c_{w_0}=\su_{y\in W}v^{-L(yw_0)}T_y$.
\endproclaim
This follows immediately from 12.1 and 11.5. Alternatively, we can argue as follows.
We prove that $p_{y,w_0}=v^{-L(yw_0)}$ for all $y$, by descending induction on 
$l(y)$. If $l(y)$ is maximal, that is $y=w_0$, then $p_{y,w_0}=1$. Assume now that 
$l(y)<l(w_0)$. We can find $s\in S$ such that $l(sy)=l(y)+1$. By the induction
hypothesis we have $p_{sy,w_0}=v^{-L(syw_0)}$. By 6.6(c), we have
$$\align&p_{y,w_0}\\&=v_s\i p_{sy,w_0}=v^{-L(s)-L(syw_0)}
=v^{-L(s)-L(w_0)+L(sy)}=v^{-L(w_0)+L(y)}=v^{-L(yw_0)}.\endalign$$
The corollary is proved.

\subhead 12.3\endsubhead
From 11.5 we see that $D_{z\i}$ can be explicitly computed when $W$ is finite and 
$c_{zw_0}$ is known. In particular, in the setup of 7.4 with $m=2k+2<\iy$, we can 
compute explicitly all $D_{z\i}$ using 7.6(a). For example:
$$\align&D_{s_1}=\su_{s\in[0,k-1]}(1-v^{2L_1}+v^{4L_1}-\do+(-1)^sv^{2sL_1})
v^{-sL_1-sL_2}\\&\T (T_{1_{2s+1}}-v^{-L_2}T_{2_{2s+2}}-v^{-L_2}T_{1_{2s+2}}
+v^{-2L_2}T_{2_{2s+3}})\\&+(1-v^{2L_1}+v^{4L_1}-\do+(-1)^kv^{2kL_1})
v^{-kL_1-kL_2}(T_{1_{2k+1}}-v^{-L_2}T_{2_{2k+2}}).\endalign$$
Using this (for larger and larger $m$) one can deduce that an analogous formula
holds in the setup of 7.4 with $m=\iy$:
$$\align&D_{s_1}=\su_{s\ge 0}(1-v^{2L_1}+v^{4L_1}-\do+(-1)^sv^{2sL_1})
v^{-sL_1-sL_2}\\&\T(T_{1_{2s+1}}-v^{-L_2}T_{2_{2s+2}}-v^{-L_2}T_{1_{2s+2}}
+v^{-2L_2}T_{2_{2s+3}})\in\ch'.\tag a\endalign$$
(We use the identification in 10.6.)

\head 13. The function $\aa$\endhead
\subhead 13.1\endsubhead
We preserve the setup of 3.1.

{\it In the remainder of these notes we assume that $L(s)>0$ for all $s\in S$.} 

For $x,y,z$ in $W$ we define $f_{x,y,z}\in\ca,f'_{x,y,z}\in\ca,h_{x,y,z}\in\ca$ by
$$T_xT_y=\su_{z\in W}f_{x,y,z}T_z=\su_{z\in W}f'_{x,y,z}c_z,$$
$$c_xc_y=\su_{z\in W}h_{x,y,z}c_z.$$
We have

(a) $f_{x,y,z}=\su_up_{z,u}f'_{x,y,u}$

(b) $f'_{x,y,z}=\su_uq'_{z,u}f_{x,y,u}$,

(c) $h_{x,y,z}=\su_{x',y'}p_{x',x}p_{y',y}f'_{x',y',z}$.
\nl
All sums in (a)-(c) are finite. (a),(c) follow from the definitions; (b) follows 
from (a) using 10.2(a). 

From 8.2, 5.6, we see that

(d) $h_{x,y,z}\ne 0\imp z\le_\car x,z\le_\cl y$,

(e) $h_{x,y,z}=h_{y\i,x\i,z\i}$. 
\nl
We show:

(f) $f_{x,y,z}\ne 0\imp l(z)\le l(x)+l(y)$;

(g) $f'_{x,y,z}\ne 0\imp l(z)\le l(x)+l(y)$;

(h) $h_{x,y,z}\ne 0\imp l(z)\le l(x)+l(y)$.
\nl
Now (f) follows from the definition, by induction on $l(x)$; (g) follows from (b) and (f) using that
$q'_{z,z'}=0$ unless $z\le z'$; (h) follows from (c) and (g) using that $p_{x',x}=0$ unless $x'\le x$.

\subhead 13.2\endsubhead
We say that $N\in\bn$ is a {\it bound} for $W,L$ if $v^{-N}f_{x,y,z}\in\ca_{\le 0}$
for all $x,y,z$ in $W$. We say that $W,L$ is {\it bounded} if there exists $N\in\bn$
such that $N$ is a bound for $W,L$. 

\proclaim{Lemma 13.3} If $W$ is finite, then $N=L(w_0^S)$ is a bound for $W,L$.
\endproclaim
By 10.4(a) we have $f_{x,y,z}=\t(T_xT_yT_{z\i})$. By 10.4(c) we have
$\t(T_xT_yT_{z\i})\in v^{L(w_0^S)}\bz[v\i]$. The lemma is proved.

\subhead 13.4\endsubhead
{\it Conjecture. In the general case $W,L$  admits a 
bound $N=\max_IL(w_0^I)$ where $I$ runs over the subsets of $S$ such that $W_I$ is finite.}

For $W$ is tame this is proved in \cite{\LC, 7.2} assuming that $L=l$, but the same proof remains valid 
without the assumption $L=l$.

We illustrate this in the setup of 7.1 with $m=\iy$. For $a,b\in\{1,2\}$ and $k>0,k'>0$, we have

$T_{a_k}T_{b_{k'}}=T_{a_{k+k'}}$ if $b=a+k\mod 2$,

$T_{a_k}T_{b_{k'}}=T_{a_kb_{k'}}+\su_{u\in[1,\min(k,k')]}
\xi_{b+u-1}T_{a_{k+k'-2u+1}}$ if $b=a+k+1\mod 2$;
\nl
here, for $n\in\bz$ we set $\xi_n=v^{L_1}-v^{-L_1}$ if $n$ is odd and 
$\xi_n=v^{L_2}-v^{-L_2}$ if $n$ is even. We see that, in this case, $\max(L_1,L_2)$
is a bound for $W,L$.

\proclaim{Lemma 13.5} Assume that $W,L$ is bounded; let $N$ be a bound for $W,L$. Then, for any $x,y,z$ in
$W$ we have

(a) $v^{-N}f'_{x,y,z}\in\ca_{\le 0}$,

(b) $v^{-N}h_{x,y,z}\in\ca_{\le 0}$.
\endproclaim
(a) follows from 13.1(b) since $q'_{z,z'}\in\ca_{\le 0}$. (b) follows from (a) and 
13.1(c) since $p_{x',x}\in\ca_{\le 0}$, $p_{y',y}\in\ca_{\le 0}$.

\subhead 13.6\endsubhead
{\it In the remainder of this chapter we assume that $W,L$ is bounded.} By 13.5(b),
for any $z\in W$ there exists a unique integer $\aa(z)\ge0$ such that

(a) $h_{x,y,z}\in v^{\aa(z)}\bz[v\i]$ for all $x,y\in W$,

(b) $h_{x,y,z}\n v^{\aa(z)-1}\bz[v\i]$ for some $x,y\in W$.
\nl
(We use that $h_{1,z,z}=1$.) We then have for any $x,y,z$:

(c) $h_{x,y,z}=\g_{x,y,z\i}v^{\aa(z)}\mod v^{\aa(z)-1}\bz[v\i]$
\nl
where $\g_{x,y,z\i}\in\bz$ is well defined; moreover, for any $z\in W$ there exists
$x,y$ such that $\g_{x,y,z\i}\ne 0$. 

For any $x,y,z$ we have 

(d) $f'_{x,y,z}=\g_{x,y,z\i}v^{\aa(z)}\mod v^{\aa(z)-1}\bz[v\i]$.
\nl
This is proved (for fixed $z$) by induction on $l(x)+l(y)$ using (c) and 13.1(c). 
(Note that $p_{x',x}p_{y',y}$ is $1$ if $x'=x,y'=y$ and is in $\ca_{<0}$ otherwise.)

\proclaim{Proposition 13.7} (a) $\aa(1)=0$.

(b) If $z\in W-\{1\}$, then $\aa(z)\ge\min_{s\in S}L(s)>0$.
\endproclaim
We prove (a). Let $x,y\in W$. Assume first that $y\ne 1$. We can find $s\in S$ such
that $ys<y$. Then $c_y\in{}^s\ch$. Since ${}^s\ch$ is a left ideal (see 8.4) we have
$c_xc_y\in{}^s\ch$. Since $s1>1$, from the definition of ${}^s\ch$ it then follows
that $h_{x,y,1}=0$. 

Similarly, if $x\ne 1$, then $h_{x,y,1}=0$. Since $h_{1,1,1}=1$, (a) follows.

In the setup of (b) we can find $s\in S$ such that $sz<z$. By 6.6(b) we have
$h_{s,z,z}=v_s+v_s\i$. This shows that $\aa(z)\ge L(s)>0$. The proposition is 
proved.

\proclaim{Proposition 13.8}Assume that $W$ is finite. 

(a) We have $\aa(w_0)=L(w_0)$.

(b) For any $w\in W-\{w_0\}$ we have $\aa(w)<L(w_0)$. 
\endproclaim
From 13.5, 13.3, for any $w\in W$ we have $\aa(w)\le L(w_0)$.

We prove (a). From 6.6(b) we see that $T_sc_{w_0}=v_sc_{w_0}$ for any $s\in S$.
Using this and 12.2, we see that 

$c_{w_0}c_{w_0}=\su_{y\in W}v^{-L(yw_0)}v^{L(y)}c_{w_0}$,
\nl
hence 
$$h_{w_0,w_0,w_0}=\su_{y\in W}v^{-L(w_0)}v^{2L(y)}\in v^{L(w_0)}\mod
v^{L(w_0)-1}\bz[v\i].$$
It follows that $\aa(w_0)\ge L(w_0)$. Hence $\aa(w_0)=L(w_0)$. This proves (a).

We prove (b). Let $z\in W$ be such that $\aa(z)=L(w_0)$. We must prove that $z=w_0$.
By 13.6(d), we can find $x,y$ such that 
$$f'_{x,y,z}=bv^{L(w_0)}+\text{strictly smaller powers of } v$$
where $b\in\bz-\{0\}$. For any $z'\ne z$ we have $f'_{x,y,z'}\in v^{L(w_0)}\bz[v\i]$
(by 13.6 and the first sentence in the proof). Since $p_{z,z'}=1$ for $z=z'$ and 
$p_{z,z'}\in\ca_{<0}$ for $z'<z$, we see that the equality 
$f_{x,y,z}=\su_{z'}p_{z,z'}f'_{x,y,z'}$ (see 13.1(a)) implies that 
$$f_{x,y,z}=bv^{L(w_0)}+\text{strictly smaller powers of } v$$
with $b\ne 0$. Now $f_{x,y,z}=\t(T_xT_yT_{z\i})$. Using now 10.4(c) we see that 
$$\min(L(x),L(y),L(z\i))=L(w_0).$$
It follows that $x=y=z\i=w_0$. The proposition is proved.

\proclaim{Proposition 13.9} (a) For any $z\in W$ we have $\aa(z)=\aa(z\i)$.

(b) For any $x,y,z\in W$ we have $\g_{x,y,z}=\g_{y\i,x\i,z\i}$.
\endproclaim
(a),(b) follow from 13.1(e). 

\subhead 13.10\endsubhead
We show that, in the setup of 7.1 with $m=\iy$ and $L_2\ge L_1$, the function 
$\aa:W@>>>\bn$ is given as follows:

(a) $\aa(1)=0$,

(b) $\aa(1_1)=L_1, \aa(2_1)=L_2$,
 
(c) $\aa(1_k)=\aa(2_k)=L_2$ if $k\ge 2$.
\nl
Now (a) is contained in 13.7(a). If $s_2z<z$ then, by the proof of 13.7(b) we have 
$\aa(z)\ge L_2$. By 13.4, $L_2$ is a bound for $W,L$ hence $\aa(z)\le L_2$ so that 
$\aa(z)=L_2$. If $zs_2<z$ then the previous argument is applicable to $z\i$. Using 
13.9, we see that $\aa(z)=\aa(z\i)=L_2$. 

Assume next that $z=1_{2k+1}$ where $k\ge 1$. By 7.5, 7.6, we have 
$$c_{1_2}c_{2_{2k}}=c_1c_2c_{2_{2k}}=(v^{L_2}+v^{-L_2})c_1c_{2_{2k}}
=(v^{L_2}+v^{-L_2})c_{1_{2k+1}},$$
hence $h_{1_2,2_{2k},z}=v^{L_2}+v^{-L_2}$. Thus, $\aa(z)\ge L_2$. By 13.4 we have
$\aa(z)\le L_2$ hence $\aa(z)=L_2$.

It remains to consider the case where $z=s_1$. Assume first that $L_1=L_2$. Then 
$\aa(s_1)\le L_1$ by 13.4 and $\aa(s_1)\ge L_1$ by 13.7(b). Hence $\aa(s_1)=L_1$.

Assume next that $L_1<L_2$. Then $\ci=\su_{w\in W-\{1,s_1\}}\ca c_w$ is a two-sided
ideal $\ci$ of $\ch$ (see 8.8). Hence if $x$ or $y$ is in $W-\{1,s_1\}$, then 
$c_xc_y\in\ci$ and $h_{x,y,s_1}=0$. Using
$$h_{1,1,s_1}=0,h_{1,s_1,s_1}=h_{s_1,1,s_1}=1,h_{s_1,s_1,s_1}=v^{L_1}+v^{-L_1}$$
we see that $\aa(s_1)=L_1$. Thus, (a),(b),(c) are established.

\subhead 13.11\endsubhead
In this subsection we assume that we are in the setup of 7.1 with $4\le m<\iy$ and 
$L_2>L_1$. By 7.8, we have
$$h_{2_{m-1},2_{m-1},2_{m-1}}=(-1)^{(m-2)/2} v^{(mL_2-(m-2)L_1)/2}
+\text{ strictly smaller powers of } v.$$
Hence $\aa(2_{m-1})\ge(mL_2-(m-2)L_1)/2$.

One can show that the function $\aa:W@>>>\bn$ is given as follows:

$\aa(1)=0$,

$\aa(1_1)=L_1, \aa(2_1)=L_2$,
 
$\aa(1_{m-1})=L_2,\aa(2_{m-1})=(mL_2-(m-2)L_1)/2$,

$\aa(2_m)=m(L_1+L_2)/2$,

$\aa(1_k)=\aa(2_k)=L_2$ if $1<k<m-1$.
\nl
This remains true in the case where $L_1=L_2$.

\subhead 13.12\endsubhead
We conjecture that any two-sided cell of $W$ would meet some finite $W_I$; this would imply that
there are only finitely many two-sided cells in $W$. (Compare 18.2.)

On the other hand, the number of left cells in $W$ can be infinite for some non-tame
$W$ with $L=l$ (B\'edard \cite{\BE})

\head 14. Conjectures \endhead
\subhead 14.1\endsubhead
We preserve the setup of 3.1. {\it In this chapter we assume that $W,L$ is bounded, see 13.2.}

For $n\in\bz$ define $\p_n:\ca@>>>\bz$ by $\p_n(\su_{k\in\bz}a_kv^k)=a_n$.

For $z\in W$ we define an integer $\D(z)\ge 0$ by
$$p_{1,z}=n_zv^{-\D(z)}+\text{strictly smaller powers of } v, n_z\in\bz-\{0\}.
\tag a$$
Note that $\D(z)=\D(z\i)$ and $\D(1)=0$, $0<\D(z)\le L(z)$ for $z\ne 1$ (see 5.4). Let 
$$\cd=\{z\in W;\aa(z)=\D(z)\}.$$ 
Clearly, $z\in\cd\imp z\i\in\cd$.

\proclaim{Conjectures 14.2} The following properties hold.

P1. For any $z\in W$ we have $\aa(z)\le\D(z)$.

P2. If $d\in\cd$ and $x,y\in W$ satisfy $\g_{x,y,d}\ne 0$, then $x=y\i$.

P3. If $y\in W$, there exists a unique $d\in\cd$ such that $\g_{y\i,y,d}\ne 0$.

P4. If $z'\le_{\lr}z$ then $\aa(z')\ge\aa(z)$. Hence, if $z'\si_{\lr}z$, then
$\aa(z')=\aa(z)$. 

P5. If $d\in\cd,y\in W,\g_{y\i,y,d}\ne 0$, then $\g_{y\i,y,d}=n_d=\pm 1$.

P6. If $d\in\cd$, then $d^2=1$.

P7. For any $x,y,z\in W$ we have $\g_{x,y,z}=\g_{y,z,x}$.

P8. Let $x,y,z\in W$ be such that $\g_{x,y,z}\ne 0$. Then $x\si_\cl y\i$, 
$y\si_\cl z\i$, $z\si_\cl x\i$.

P9. If $z'\le_\cl z$ and $\aa(z')=\aa(z)$ then $z'\si_\cl z$.

P10. If $z'\le_\car z$ and $\aa(z')=\aa(z)$ then $z'\si_\car z$.

P11. If $z'\le_{\lr}z$ and $\aa(z')=\aa(z)$ then $z'\si_{\lr}z$.

P12. Let $I\sub S$. If $y\in W_I$, then $\aa(y)$ computed in terms of $W_I$ is equal
to $\aa(y)$ computed in terms of $W$.

P13. Any left cell $\G$ of $W$ contains a unique element $d\in\cd$. We have 
$\g_{x\i,x,d}\ne 0$ for all $x\in\G$.

P14. For any $z\in W$ we have $z\si_{\lr}z\i$.

P15. Let $v'$ be a second indeterminate and let $h'_{x,y,z}\in\bz[v',v'{}\i]$ be
obtained from $h_{x,y,z}$ by the substitution $v\m v'$. If $x,x',y,w\in W$ satisfy 
$\aa(w)=\aa(y)$, then 

$\su_{y'}h'_{w,x',y'}h_{x,y',y}=\su_{y'}h_{x,w,y'}h'_{y',x',y}$.
\endproclaim
In \S15-\S17 we will verify the conjectures above in a number of cases.

\subhead 14.3\endsubhead
We consider the following auxiliary statement.

$\tP$. {\it Let $x,y,z,z'\in W$ be such that $\g_{x,y,z\i}\ne 0$, $z'\gt_\cl z$. 
Then there exists $x'\in W$ such that $\p_{\aa(z)}(h_{x',y,z'})\ne 0$. In 
particular, $\aa(z')\ge\aa(z)$.}

In this chapter we will show, that, if P1-P3 and $\tP$ are assumed to be true, then 
P4-P14 are automatically true. The arguments follow \cite{\LC},\cite{\LCC}.

\subhead 14.4\endsubhead 
{\it $\tP$ $\imp$ P4.} Let $z',z$ be as in P4. We can assume that $z\gt_\cl z'$ or 
$z\gt_\car z'$. In the first case, from $\tP$ we get $\aa(z')\ge\aa(z)$. (We can 
find $x,y$ such that $\g_{x,y,z\i}\ne 0$.) In the second case, from $\tP$ we get 
$\aa(z'{}\i)\ge\aa(z\i)$ hence $\aa(z')\ge\aa(z)$.  

\subhead 14.5\endsubhead 
{\it P1,P3 $\imp$ P5.} Let $x,y\in W$. Applying $\t$ to
$c_xc_y=\su_{z\in W}h_{x,y,z}c_z$ gives
$$\su_zh_{x,y,z}p_{1,z}=\su_{x',y'}p_{x',x}p_{y',y}\t(T_{x'}T_{y'})=
\su_{x',y'}p_{x',x}p_{y',y}\d_{x'y',1}=\su_{x'}p_{x',x}p_{x'{}\i,y}$$
hence
$$\su_{z\in W}h_{x,y,z}p_{1,z}=\d_{xy,1}\mod v\i\bz[v\i].\tag a$$
We take $x=y\i$ and note that 
$h_{y\i,y,z}\in v^{\aa(z)}\bz[v\i]$, $p_{1,z}\in v^{-\D(z)}\bz[v\i]$, hence
$$h_{y\i,y,z}p_{1,z}\in v^{\aa(z)-\D(z)}\bz[v\i].$$
The same argument shows that, if $z\in\cd$, then 
$$h_{y\i,y,z}p_{1,z}\in\g_{y\i,y,z\i}n_z+\ca_{<0}.$$
If $z\n\cd$ then, by P1, we have $\aa(z)-\D(z)<0$ so that
$h_{y\i,y,z}p_{1,z}\in\ca_{<0}$. We see that 
$$\su_{z\in W}h_{y\i,y,z}p_{1,z}=\su_{z\in\cd}\g_{y\i,y,z\i}n_z\mod\ca_{<0}.$$
Comparison with (a) gives $\su_{z\in\cd}\g_{y\i,y,z\i}n_z=1$. Equivalently,
$$\su_{z\in\cd}\g_{y\i,y,z}n_z=1.$$
Using this and P3 we see that, in the setup of P5 we have $\g_{y\i,y,d}n_d=1$. Since
$\g_{y\i,y,d},n_d$ are integers, we must have $\g_{y\i,y,d}=n_d=\pm 1$.

\subhead 14.6\endsubhead 
{\it P2,P3 $\imp$ P6.} We can find $x,y$ such that $\g_{x,y,d}\ne 0$. By P2, we have
$x=y\i$ so that $\g_{y\i,y,d}\ne 0$. This implies $\g_{y\i,y,d\i}\ne 0$. (See
13.9(b)). We have $d\i\in\cd$. By the uniqueness in P3 we have $d=d\i$.

\subhead 14.7\endsubhead 
{\it P2,P3,P4,P5 $\imp$ P7.} We first prove the following statement.

(a) {\it Let $x,y,z\in W, d\in\cd$ be such that $\g_{x,y,z}\ne 0$,
$\g_{z\i,z,d}\ne 0$, $\aa(d)=\aa(z)$. Then} $\g_{x,y,z}=\g_{y,z,x}$.
\nl
Let $n=\aa(d)$. From $\g_{x,y,z}\ne 0$ we deduce $h_{x,y,z\i}\ne 0$ hence 
$z\i\le_\car x$, hence $n=\aa(z)=\aa(z\i)\ge \aa(x)$ (see P4). Computing the
coefficient of $c_d$ in two ways, we obtain 
$$\su_{z'}h_{x,y,z'}h_{z',z,d}=\su_{x'}h_{x,x',d}h_{y,z,x'}.$$
Now $h_{z',z,d}\ne 0$ implies $d\le_\car z'$ hence $\aa(z')\le\aa(d)=n$ (see P4); 
similarly, $h_{x,x',d}\ne 0$ implies $d\le_\cl x'$ hence $\aa(x')\le\aa(d)=n$. Thus
we have 
$$\su_{z';\aa(z')\le n}h_{x,y,z'}h_{z',z,d}=
\su_{x';\aa(x')\le n}h_{x,x',d}h_{y,z,x'}.$$
By P2 and our assumptions, the left hand side is 
$$\g_{x,y,z}\g_{z\i,z,d}v^{2n}+\text{strictly smaller powers of } v.$$
Similarly, the right hand side is 
$$\g_{x,x\i,d}\p_n(h_{y,z,x\i})v^{2n}+\text{strictly smaller powers of } v.$$
Hence $\g_{x,x\i,d}\p_n(h_{y,z,x\i})=\g_{x,y,z}\g_{z\i,z,d}\ne 0$. Thus, 
$$\g_{x,x\i,d}\ne 0,\p_n(h_{y,z,x\i})\ne 0.$$
We see that $\aa(x\i)\ge n$. But we have also $\aa(x)\le n$ hence $\aa(x)=n$ and 
$\p_n(h_{y,z,x\i})=\g_{y,z,x}$. Since $\g_{x,x\i,d}\ne 0$, we have (by P5) 
$\g_{x,x\i,d}=\g_{z\i,z,d}$. Using this and 
$\g_{x,x\i,d}\g_{y,z,x}=\g_{x,y,z}\g_{z\i,z,d}$ we deduce $\g_{y,z,x}=\g_{x,y,z}$, 
as required. 

Next we prove the following statement.

(b) {\it Let $z\in W,d\in\cd$ be such that $\g_{z\i,z,d}\ne 0$. Then}
$\aa(z)=\aa(d)$.
\nl
We shall assume that (b) holds whenever $\aa(z)>N_0$ and we shall deduce that it
also holds when $\aa(z)=N_0$. (This will prove (b) by descending induction on 
$\aa(z)$ since $\aa(z)$ is bounded above.) Assume that $\aa(z)=N_0$. From 
$\g_{z\i,z,d}=\pm 1$ we deduce that $h_{z\i,z,d\i}\ne 0$ hence $d\i\le_\cl z\i$ 
hence $\aa(d\i)\ge\aa(z\i)$ (see P4) and $\aa(d)\ge \aa(z)$. Assume that 
$\aa(d)>\aa(z)$, that is, $\aa(d)>N_0$. Let $d'\in\cd$ be such that 
$\g_{d\i,d,d'}\ne 0$ (see P3). By the induction hypothesis applied to $d,d'$ instead
of $z,d$, we have $\aa(d)=\aa(d')$. From $\g_{z\i,z,d}\ne 0$,$\g_{d\i,d,d'}\ne 0$, 
$\aa(d)=\aa(d')$, we deduce (using (a)) that $\g_{z,d,z\i}=\g_{z\i,z,d}$. Hence 
$\g_{z,d,z\i}\ne 0$. It follows that $h_{z,d,z}\ne 0$, hence $z\le_\cl d$, hence 
$\aa(z)\ge\aa(d)$ (see P4). This contradicts the assumption $\aa(d)>\aa(z)$. Hence 
we must have $\aa(z)=\aa(d)$, as required.

We now prove P7. Assume first that $\g_{x,y,z}\ne 0$. Let $d\in\cd$ be such that 
$\g_{z\i,z,d}\ne 0$ (see P3). By (b) we have $\aa(z)=\aa(d)$. Using (a) we then have
$\g_{x,y,z}=\g_{y,z,x}$. Assume next that $\g_{x,y,z}=0$; we must show that 
$\g_{y,z,x}=0$. We assume that $\g_{y,z,x}\ne 0$. By the first part of the proof, we
have
$$\g_{y,z,x}\ne 0\imp\g_{y,z,x}=\g_{z,x,y}\ne 0\imp\g_{z,x,y}=\g_{x,y,z}\ne 0,$$
a contradiction.

\subhead 14.8\endsubhead 
{\it P7 $\imp$ P8.} If $\g_{x,y,z}\ne 0$, then $h_{x,y,z\i}\ne 0$, hence 
$z\i\le_\cl y, z\le_\cl x\i$. By P7 we also have $\g_{y,z,x}\ne 0$ (hence
$x\i\le_\cl z, x\le_\cl y\i$) and $\g_{z,x,y}\ne 0$ (hence
$y\i\le_\cl x, y\le_\cl z\i$). Thus, we have 
$x\si_\cl y\i,y\si_\cl z\i,z\si_\cl x\i$.

\subhead 14.9\endsubhead 
{\it $\tP$,P4,P8 $\imp$ P9.} We can find a sequence $z'=z_0,z_1,\do,z_n=z$ such that
for any $j\in[1,n]$ we have $z_{j-1}\gt_\cl z_j$. By P4 we have 
$\aa(z')=\aa(z_0)\ge\aa(z_1)\ge\do\ge\aa(z_n)=\aa(z)$. Since $\aa(z)=\aa(z')$, we 
have $\aa(z')=\aa(z_0)=\aa(z_1)=\do=\aa(z_n)=\aa(z)$. Thus, it suffices to show 
that, if $z'\gt_\cl z$ and $\aa(z')=\aa(z)$, then $z'\si_\cl z$. Let $x,y\in W$ be 
such that $\g_{x,y,z\i}\ne 0$. By $\tP$, there exists $x'\in W$ such that 
$\p_{\aa(z)}(h_{x',y,z'})\ne 0$. Since $\aa(z')=\aa(z)$, we have 
$\g_{x',y,z'{}\i}\ne 0$. From $\g_{x,y,z\i}\ne 0$,$\g_{x',y,z'{}\i}\ne 0$ we deduce,
using P8, that $y\si_\cl z,y\si_\cl z'$, hence $z\si_\cl z'$.

\subhead 14.10\endsubhead 
{\it P9 $\imp$ P10.} We apply P9 to $z\i, z'{}\i$.

\subhead 14.11\endsubhead 
{\it P4,P9,P10 $\imp$ P11.} We can find a sequence $z'=z_0,z_1,\do,z_n=z$ such that
for any $j\in[1,n]$ we have $z_{j-1}\le_\cl z_j$ or $z_{j-1}\le_\car z_j$. By P4, we
have $\aa(z')=\aa(z_0)\ge\aa(z_1)\ge\do\ge\aa(z_n)=\aa(z)$. Since $\aa(z)=\aa(z')$,
we have $\aa(z')=\aa(z_0)=\aa(z_1)=\do=\aa(z_n)=\aa(z)$. Applying P9 or P10 to 
$z_{j-1},z_j$ we obtain $z_{j-1}\si_\cl z_j$ or $z_{j-1}\si_\car z_j$. Hence 
$z'\si_{\lr}z$.

\subhead 14.12\endsubhead 
{\it P3,P4,P8 for $W$ and $W_I$ $\imp$ P12.} We write $\aa_I:W_I@>>>\bn$ for the
$\aa$-function defined in terms of $W_I$. For $x,y,z\in W_I$, we write 
$h^I_{x,y,z},\g^I_{x,y,z}$ for the analogues of $h_{x,y,z},\g_{x,y,z}$ when $W$ is 
replaced by $W_I$. Let $\ch_I\sub\ch$ be as in 9.9.

Let $d\in\cd$ be such that $\g_{y\i,y,d}\ne 0$. (See P3.) Then
$\p_{\aa(d)}(h_{y\i,y,d\i})\ne 0$. Now $c_{y\i}c_y\in\ch_I$ hence $d\in W_I$ and 
$\p_{\aa(d)}(h^I_{y\i,y,d\i})\ne 0$. Thus, $\aa_I(d\i)\ge\aa(d\i)$. The reverse 
inequality is obvious hence $\aa_I(d)=\aa(d)$. We see that $\g^I_{y\i,y,d}\ne 0$. 
From P8 we see that $y\si_\cl d$ (relative to $W_I$) and $y\si_\cl d$ (relative to
$W$). From P4 we deduce that $\aa_I(y)=\aa_I(d)$ and $\aa(y)=\aa(d)$. It follows 
that $\aa(y)=\aa_I(y)$.

\subhead 14.13\endsubhead 
{\it $\tP$,P2,P3,P4,P6,P8 $\imp$ P13.} If $x\in\G$ then, by P3, there exists 
$d\in\cd$ such that $\g_{x\i,x,d}\ne 0$. By P8 we have $x\si_\cl d\i$ hence 
$d\i\in\G$. By P6, we have $d=d\i$ hence $d\in\G$. It remains to prove the 
uniqueness of $d$. Let $d',d''$ be elements of $\cd\cap\G$. We must prove that 
$d'=d''$. We can find $x',y',x'',y''$ such that $\g_{x',y',d'}\ne 0$,
$\g_{x'',y'',d''}\ne 0$. By P2, we have $x'=y'{}\i,x''=y''{}\i$. By P8, we have 
$y'\si_\cl d'{}\i=d'$ and $y''\si_\cl d''{}\i=d''$, hence $y',y''\in\G$. By the 
definition of left cells, we can find a sequence $y'=x_0,x_1,\do,x_n=y''$ such that
for any $j\in[1,n]$ we have $x_{j-1}\gt_\cl x_j$. Since
$y'\si_\cl y''$, we have $x_j\in\G$ for all $j$. For $j\in[1,n-1]$ let $d_j\in\cd$
be such that $\g_{x_j\i,x_j,d_j}\ne 0$. Let $d_0=d',d_n=d''$. As in the beginning of
the proof, we have $d_j\in\G$ for each $j$. Let $j\in[1,n]$. Since
$x_{j-1}\gt_\cl x_j$, we have (by P8) $\g_{x_j,d_j,x_j\i}\ne 0$. Applying $\tP$ to 
$x_j,d_j,x_j,x_{j-1}$ instead of $x,y,z,z'$, we see that there exists $u$ such that 
$\p_{\aa(x_j)}(h_{u,d_j,x_{j-1}})\ne 0$. Since $x_{j-1}\si x_j$, we have 
$\aa(x_{j-1})=\aa(x_j)$ (see P4), hence 
$\p_{\aa(x_j)}(h_{u,d_j,x_{j-1}})=\g_{u,d_j,x_{j-1}\i}\ne 0$. Using P8, we deduce 
$\g_{x_{j-1}\i,u,d_j}\ne 0$. Using P2 we see that $u=x_{j-1}$ and
$\g_{x_{j-1}\i,x_{j-1},d_j}\ne 0$. We have also
$\g_{x_{j-1}\i,x_{j-1},d_{j-1}}\ne 0$ and by the uniqueness in P3, it follows that 
$d_{j-1}=d_j$. It follows that $d'=d''$, as required.

\subhead 14.14\endsubhead 
{\it P6,P13 $\imp$ P14.} By P13, we can find $d\in\cd$ such that $z\si_\cl d$. Since
$d=d\i$ (see P6), it follows that $z\i\si_\car d$. Thus, $z\si_{\lr}z\i$.

\subhead 14.15\endsubhead 
In this subsection we reformulate conjecture P15, assuming that P4,P9, P10 hold. Let
$\tca=\bz[v,v\i,v',v'{}\i]$ where $v,v'$ are indeterminates. Let $\tch$ be the free
$\tca$-module with basis $e_w (w\in W)$. Let $\ch',c'_w,h'_{x,y,z}$ be obtained from
$\ch,c_w,h_{x,y,z}$ by changing the variable $v$ to $v'$. 

On $\tch$ we have a left $\ch$-module structure given by
$v^nc_ye_w=\su_xv^nh_{y,w,x}e_x$ and a right $\ch'$-module structure given by 
$e_w(v'{}^nc'_y)=\su_xv'{}^nh'_{w,y,x}e_x$. These module structures do not commute
in general. For each $a\ge 0$ let $\tch_{\ge a}$ be the $\tca$-submodule of $\tch$ 
spanned by $\{e_w;\aa(w)\ge a\}$. By P4, this is a left $\ch$-submodule and a right
$\ch'$-submodule of $\tch$. We have
$$\do\tch_{\ge 2}\sub\tch_{\ge 1}\sub\tch_{\ge 0}=\tch$$
and $gr\tch=\op_{a\ge 0}\tch_{\ge a}/\tch_{\ge a+1}$ inherits a left $\ch$-module 
structure and a right $\ch'$-module structure from $\tch$. Clearly, P15 is 
equivalent to the condition that these module structures on $gr\tch$ commute. To 
check this last condition, it is enough to check that

the actions of $c_s, c'_{s'}$ commute on $gr\tch$ for $s,s'\in S$. 
\nl
Let $s,s'\in S$, $w\in W$. A computation using 6.6, 6.7, 8.2 shows 
that $(c_se_w)c'_{s'}-c_s(e_wc'_{s'})$ is $0$ if $sw<w$ or $ws'<w$, while if 
$sw>w,ws'>w$, it is
$$\su_{y;sy<y,ys'<y}(h'_{w,s',y}(v_s+v_s\i)-h_{s,w,y}(v'_{s'}+v'_{s'}{}\i))e_y
+\su_{y;sy<y,ys'<y}\a_ye_y$$
where 
$$\a_y=\su_{y';y's'<y'<sy'}h'_{w,s',y'}h_{s,y',y}
-\su_{y';sy'<y'<y's'}h_{s,w,y'}h'_{y',s',y}.$$
If $y$ satisfies $sy<y,ys'<y$ and either $h'_{w,s',y}$ or $h_{s,w,y}$ is $\ne 0$, 
then $\aa(y)>\aa(w)$. (We certainly have $\aa(y)\ge \aa(w)$ by P4. If we had 
$\aa(y)=\aa(w)$ and $h_{s,w,y}\ne 0$ then by P9 we would have $y\si_\cl w$ hence
$\car(y)=\car(w)$ contradicting $ys'<y,ws'>w$. If we had $\aa(y)=\aa(w)$ and 
$h'_{s,w,y}\ne 0$ then by P10 we would have $y\si_\car w$ hence $\cl(y)=\cl(w)$,
contradicting $sy<y,sw>w$.) Hence, if $sw>w,ws'>w$, we have
$$(c_se_w)c'_{s'}-c_s(e_wc'_{s'})=\su_{y;sy<y,ys'<y,\aa(y)=\aa(w)}\a_ye_y\mod
\tch_{\ge\aa(w)+1}.$$
We see that P15 is equivalent to the following statement.

(a) {\it If $y,w\in W,s,s'\in S$ are such that 
$sw>w,ws'>w,sy<y,ys'<y,\aa(y)=\aa(w)$, then}

$\su_{y';y's'<y'<sy'}h'_{w,s',y'}h_{s,y',y}
=\su_{y';sy'<y'<y's'}h_{s,w,y'}h'_{y',s',y}$.

\head 15. Example: the split case\endhead
\subhead 15.1\endsubhead
We preserve the setup of 3.1. We assume that $L=l$ that is, we are in the split case. From the results on
Soergel bimodules in \cite{\EW}, we see that

(a) $h_{x,y,z}\in\bn[v,v\i]$ for all $x,y,z$ in $W$,

(b) $p_{y,w}\in\bn[v\i]$ for all $y,w$ in $W$.
\nl
In this chapter we assume that $W,l$ is bounded. We will show that
$\tP$ and P1-P3 hold for $W,l$ hence all of P1-P14 
hold for $W,l$; moreover we show that P15 holds. (Note that for P1-P6, the boundedness of $W,l$
will not be needed). The arguments in this chapter are based on \cite{\LCC}.

\subhead 15.2\endsubhead 
From 14.5(a) (which does not depend on the boundedness assumption) we see that for $x,y\in W$ we have 

(a) $\su_{z\in W}h_{x,y,z}p_{1,z}\in\bz[v\i]$.
\nl
From 15.1(a),(b) we see that $h_{x,y,z}p_{1,z}\in\bn[v,v\i]$ for any $z\in W$. Hence in (a) there are no 
cancellations, so that

(b) $h_{x,y,z}p_{1,z}\in\bn[v\i]$ for any $z\in W$.
\nl
Let $z\in W$. By 5.4(a) we have $p_{1,z}=\sum_{j\ge0}e_jv^{-l(z)+j}$ where $e_j\in\ZZ$ are zero for large 
$j$ and $e_0=1$; moreover, by 15.1(b) we have $e_j\in\NN$ for all $j$. Thus,
$v^{-l(z)}h_{x,y,z}+\sum_{j>0}e_jv^{-l(z)+j}h_{x,y,z}\in\bn[v\i]$.
Again in this sum there are no cancellations, hence 
$$v^{-l(z)}h_{x,y,z}\in\bn[v\i].$$ 
(This was proved for finite $W$ in \cite{\LC} and then for general $W$ by Springer (unpublished).)
This shows that the definition of $\aa(z)$ in  13.6 can be given in our case without 
the boundedness assumption. Hence the definition of $\g_{x,y,z}$ in 13.6 can be 
 given in our case without 
the boundedness assumption. Note that the definition of $\D(z)$ in 14.1 can be given without the 
boundedness assumption. Hence the definition of $\cd$ in 14.1 can be given in our case without the 
boundedness assumption. 

We show that {\it P1 holds for $(W,l)$}. 

We fix $z\in W$ and choose $x,y\in W$ so that $\g_{x,y,z\i}\ne 0$. From the  definitions,

(c) $h_{x,y,z}p_{1,z}\in\g_{x,y,z\i}n_zv^{\aa(z)-\D(z)}+\text{strictly smaller powers of } v$
\nl
and the coefficient of $v^{\aa(z)-\D(z)}$ is $\ne 0$. Comparison with (b) gives $\aa(z)-\D(z)\le 0$.

\subhead 15.3\endsubhead 
{\it Proof of P2.} Assume that $x\ne y\i$. From 14.5(a) we see that 

(a) $\su_{z\in W}h_{x,y,z}p_{1,z}\in v\i\bz[v\i]$.
\nl
As in 15.2, this implies (using 15.1(a),(b)) that 

(b) $h_{x,y,z}p_{1,z}\in v\i\bn[v\i]$ for any $z\in W$.
\nl
Assume now that $z=d\i\in\cd$. Then 15.2(c) becomes in our case
$$h_{x,y,z}p_{1,z}\in\g_{x,y,z\i}n_z+v\i\bz[v\i].$$
Comparison with (b) gives $\g_{x,y,z\i}n_z=0$. Since $n_z\ne 0$, we have 
$\g_{x,y,z\i}=0$. This proves P2.

\subhead 15.4\endsubhead 
{\it Proof of P3.} From 14.5(a) we see that 

(a) $\su_{z\in W}h_{y\i,y,z}p_{1,z}\in 1+v\i\bz[v\i]$.
\nl
As in 15.2, this implies (using 15.1(a),(b)) that there is a unique $z$, say $z=d\i$
such that

(b) $h_{y\i,y,d\i}p_{1,d\i}\in 1+v\i\bn[v\i]$ 
\nl
and that

(c) $h_{y\i,y,z}p_{1,z}\in v\i\bn[v\i]$
\nl
for all $z\ne d\i$. For $z=d\i$, 15.2(c) becomes
$$h_{y\i,y,d\i}p_{1,d\i}\in\g_{y\i,y,d}n_{d\i}v^{\aa(d)-\D(d)}
+\text{strictly smaller powers of } v.$$
Here $\aa(d)-\D(d)\le 0$. Comparison with (b) gives $\aa(d)-\D(d)=0$ and \lb
$\g_{y\i,y,d}n_{d\i}=1$. Thus, $d\in\cd$ and $\g_{y\i,y,d}\ne 0$. Thus, the 
existence part of P3 is established.

Assume that there exists $d'\ne d$ such that $d'\in\cd$ and $\g_{y\i,y,d'}\ne 0$. 
For $z=d'{}\i$, 15.2(c) becomes
$$h_{y\i,y,d'{}\i}p_{1,d'{}\i}\in\g_{y\i,y,d'}n_{d'{}\i}+v\i\bz[v\i].$$
Comparison with (c) (with $z=d'{}\i$) gives $\g_{y\i,y,d'}n_{d'{}\i}=0$ hence 
$\g_{y\i,y,d'}=0$, a contradiction. This proves the uniqueness part of P3.

\subhead 15.5\endsubhead
{\it Proof of $\tP$.} We may assume that $z'\ne z$. Then we can find $s\in S$ such 
that $sz'<z',sz>z$ and $h_{s,z,z'}\ne 0$. Since $h_{x,y,z}\ne 0$, we have (by 
13.1(d)) $z\le_\car x$ hence $\cl(x)\sub\cl(z)$ (by 8.6). Since $s\n\cl(z)$, we have
$s\n\cl(x)$, that is, $sx>x$. We have $c_sc_xc_y=\su_up_uc_u$, where
$$p_u=\su_wh_{x,y,w}h_{s,w,u}=\su_{x'}h_{s,x,x'}h_{x',y,u}.$$
In particular, 
$$p_{z'}=\su_wh_{x,y,w}h_{s,w,z'}=h_{x,y,z}h_{s,z,z'}
+\su_{w;w\ne z}h_{x,y,w}h_{s,w,z'}.$$
By 6.5, we have $h_{s,z,z'}\in\bz$ hence
$$\p_n(p_{z'})=\p_n(h_{x,y,z})h_{s,z,z'}+\su_{w;w\ne z}\p_n(h_{x,y,w}h_{s,w,z'})
\tag a$$
for any $n\in\bz$. In particular, this holds for $n=\aa(z)$. By assumption, we have
$\p_n(h_{x,y,z})\ne 0$ and $h_{s,z,z'}\ne 0$; hence, by 15.1(a), we have
$\p_n(h_{x,y,z})>0$ and $h_{s,z,z'}>0$. Again, by 15.1(a) we have 
$\p_n(h_{x,y,w}h_{s,w,z'})\ge 0$ for any $w\ne z$. Hence from (a) we deduce 
$\p_n(p_{z'})>0$. Since $p_{z'}=\su_{x'}h_{s,x,x'}h_{x',y,z'}$, there exists $x'$
such that $\p_n(h_{s,x,x'}h_{x',y,z'})\ne 0$. Since $sx>x$, we see from 6.5 that 
$h_{s,x,x'}\in\bz$ hence 
$$\p_n(h_{s,x,x'}h_{x',y,z'})=h_{s,x,x'}\p_n(h_{x',y,z'}).$$
Thus we have $\p_n(h_{x',y,z'})\ne 0$. This proves $\tP$ in our case.

\subhead 15.6\endsubhead 
Since $\tP$ and P1-P3 are known, we see that P1-P11 and P13,P14 hold in our case 
(see \S14). The same arguments can be applied to $W_I$ where $I\sub S$, hence P1-P11
and P13,P14 hold for $W_I$. By 14.12, P12 holds for $W$. Thus, P1-P14 hold for $W$.

\subhead 15.7\endsubhead 
{\it Proof of P15.} By 14.15, we see that it is enough to prove 14.15(a). Let 
$y,w,s,s'$ be as in 14.15(a). In our case, by 6.5, the equation in 14.15(a) involves
only integers, hence it is enough to prove it after specializing $v=v'$. If in 14.15
we specialize $v=v'$, then the left and right module structures in 14.15 clearly 
commute, since the left and right regular representations of $\ch$ commute. Hence 
the coefficient of $e_y$ in $((c_se_w)c_{s'}-c_s(e_wc_{s'}))_{v=v'}$ is $0$. By the
computation in 14.15, this coefficient is 
$$(h_{w,s',y}-h_{s,w,y})(v+v\i)+\su\Sb y'\\y's'<y'<sy'\eSb h_{w,s',y'}
h_{s,y',y}-\su\Sb y'\\sy'<y'<y's'\eSb h_{s,w,y'}h_{y',s',y}=0.\tag a$$
By 6.5, $h_{s,w,y}$ is the coefficient of $v\i$ in $p_{y,w}$ and
$h_{w,s',y}=h_{s',w\i,y\i}$ is the coefficient of $v\i$ in $p_{y\i,w\i}=p_{y,w}$. 
Thus, $h_{s,w,y}=h_{w,s',y}$ and (a) reduces to the equation in 14.15(a) 
(specialized at $v=v'$). This proves 14.15(a).

\head 16. Example: the quasisplit case\endhead
\subhead 16.1\endsubhead
In this subsection we review some results from \cite{\QG, \S11}.

Let $k$ be an algebraically closed field of characteristic zero. Let $\fC$ be a $k$-linear category, that is
a category in which the space of morphisms between any two objects has a given $k$-vector space structure 
such that composition of morphisms is bilinear and such that finite direct sums exist. A functor from one 
$k$-linear category to another is said to be $k$-linear if it respects the $k$-vector space structures.

Let $\ck(\fC)$ be the Grothendieck group of $\fC$ that is, the free abelian group generated by symbols $[A]$
for each $A\in\fC$ (up to isomorphism) with relations $[A\op B]=[A|+[B]$ for any $A,B\in\fC$.
Let $\nn$ be an integer $\ge1$.
A $k$-linear functor $M\m M^\sh$, $\fC@>>>\fC$ is said to be {\it $\nn$-periodic} if $(\sh)^\nn:\fC@>>>\fC$
is the identity functor. Assuming that such a functor is given we define a new $k$-linear category 
$\fC_\sh$ as follows. The objects of $\fC_\sh$ are pairs $(A,\ph)$ where $A\in\fC$ and $\ph:A^\sh@>>>A$ is an
isomorphism in $\fC$ such that the composition 
$$A^{\sh^\nn}@>\ph^{\sh^{\nn-1}}>>A^{\sh^{\nn-1}}@>\ph^{\sh^{\nn-2}}>>\do@>>>A^\sh@>\ph>>A$$
is the identity map of $A$. Let $(A,\ph)$, $(A',\ph')$ be two objects of $\fC_\sh$. We define a $k$-linear 
map $\Hom_{\fC}(A,A')@>>>\Hom_{\fC}(A,A')$ by $f\m f^!:=\ph'f^\sh\ph\i$. Note that 
the $\nn$-th iteration of ${}^!$ applied to $f$ is $1$. By definition,
$\Hom_{\fC_\sh}((A,\ph),(A',\ph'))=\{f\in\Hom_{\fC}(A,A');f=f^!\}$, a $k$-vector space. 
The direct sum of two objects $(A,\ph)$, $(A',\ph')$ is $(A\op A',\ph\op\ph')$. Clearly, if 
$(A,\ph)\in\fC_\sh$, then $(A,\z\ph)\in\fC_\sh$ for any $\z\in k$ such that $\z^\nn=1$.
An object $(A,\ph)$ of $\fC_\sh$ is said to be {\it traceless} if there exists an 
object $B$ of $\fC$, an integer $t\ge2$ dividing $\nn$ and an isomorphism 
$$A\cong B\op B^\sh\op\do\op B^{\sh^{t-1}}$$ under which $\ph$ corresponds to an isomorphism 
$$B\op B^\sh\op\do\op B^{\sh^{t-1}}@>\si>>B\op B^\sh\op\do\op B^{\sh^{t-1}}$$
which carries the summand of $B^{\sh^j}$ onto the summand $B^{\sh^j}$ for $1\le j\le t-1$) and the
summand $B^{\sh^t}$ onto the summand $B$.

Let $\co$ be the subring of $k$ consisting of all $\ZZ$-linear combinations of $\nn$-th roots of $1$. We
associate to $\fC$ and $\sh$ an $\co$-module $\ck_\sh(\fC)$. By definition $\ck_\sh(\fC)$ is the 
$\co$-module generated by symbols $[B,\ph]$ one for each isomorphism class of objects $(B,\ph)$ of $\fC_\sh$
subject to the following relations:

(a) $[B,\ph]+[B',\ph']=[B\op B',\ph\op\ph']$;

(b) $[B,\ph]=0$ if $(B,\ph)$ is traceless;

(c) $[B,\z\ph]=\z[B,\ph]$ if $\z\in k$ satisfies $\z^\nn=1$.

\subhead 16.2\endsubhead
Let $\tW,\tS$ be a Coxeter group. (The set $\tS$ of simple reflections of $\tW$ is assumed to be finite.) 
We view $\tS$ as a subset of $\tW$. For any $I\sub\tS$ we denote by $\tW_I$ the subgroup of $\tW$ generated 
by $I$. Let $\tl:\tW@>>>\NN$ be the length function of $\tW$. Let $\le$ be the 
standard partial order on $\tW$. Let $\io:\tW@>>>\tW$ be an automorphism such that $\io(\tS)=\tS$. We fix an
integer $\nn\ge1$ such that $\io^\nn=1$. Let $W=\{w\in\tW;\io(w)=w\}$. 
Let $S$ be the set of $\io$-orbits $I$ on
$\tS$ such that $\tW_I$ is finite; for such $I$ let $w_0^I$ be the longest element of $\tW_I$ (note that
$w_0^I\in W$). According to Theorems A.8 and A.9 in the Appendix, $W$ is a Coxeter group on the set of 
generators $\{w_0^I;I\in S\}$ and the restriction of $\tl$ to $W$ is a weight function $L:W@>>>\NN$. Let 
$l:W@>>>\NN$ be the length function of $W,S$. We then say that $W,L$ is in the quasisplit case. We define 
$$\ch,T_x,c_x,p_{y,x},f_{x,y,z},h_{x,y,z}$$
in terms of $W,S,L$ as in 3.2, 5.2, 5.3 (here $x,y,z\in W$). Let 
$$\tch,\ti T_x,\tc_x,\ti p_{y,x},\tf_{x,y,z},\ti h_{x,y,z}$$
be the analogous objects defined in terms of $\tW,\tS,\tl$ (here $x,y,z\in\tW$).

Let $\fh_\RR$ be a reflection representation of $\tW$ over the real numbers $\RR$, as in \cite{\EW, 1.1}; 
for any $s\in\tS$ we fix a linear form $\a_s:\fh_\RR@>>>\RR$ whose kernel is equal to the fixed point set of 
$s:\fh_\RR@>>>\fh_\RR$. Let $\fh=\CC\ot_\RR\fh$; we extend $\a_s$ to a $\CC$-linear function $\fh@>>>\CC$ 
denoted again by $\a_s$. 
Let $R$ be the algebra of polynomial functions $\fh@>>>\CC$ with the $\ZZ$-grading in which
linear functions $\fh@>>>\CC$ have degree $2$. Note that $\tW$ acts naturally on $R$; we write this action 
as $w:r\m{}^wr$ and for $s\in\tS$ we set $R^s=\{r\in R;{}^sr=r\}$, a subalgebra of $R$. Let 
$R^{>0}=\{r\in R;r(0)=0\}$. We can assume that there exists a $\CC$-linear map $\x\m\io(\x)$ of the dual 
space of $\fh$ into itself whose $\nn$-th power is $1$ and is such that $\io(w\x)=\io(w)(\io(\x))$ for 
$w\in W,x\in\fh$ and such that $\io(\a_s)=\a_{\io(s)}$ for $s\in\tS$. It induces an algebra automorphism 
$r\m\io(r)$ of $R$. 

Let $\car$ be the category whose objects are $\ZZ$-graded $(R,R)$-bimodules in which for $M,M'\in\car$,
$\Hom_\car(M,M')$ is the space of  homomorphisms of $(R,R)$-bimodules $M@>>>M'$ compatible with the
$\ZZ$-gradings. For $M\in\car$ and $n\in\ZZ$, the shift $M[n]$ is the object of $\car$ equal in degree $i$ to
$M$ in degree $i+n$. For $M,M'$ in $\car$ we set $MM'=M\ot_RM'$; this is naturally an object of $\car$. For 
$M,M'$ in $\car$ we set 
$$M'{}^M=\op_{n\in\ZZ}\Hom_\car(M,M'[n]),$$
viewed as an object of $\car$ with $(rf)(m)=f(rm)$, $(fr)(m)=f(mr)$ for $m\in M,f\in M'{}^M,r\in R$. For any
$M\in\car$ we set $\un{M}=M/MR^{>0}=M\ot_R\CC$ where $\CC$ is identified with $R/R^{>0}$. We view $\un{M}$ 
as a $\ZZ$-graded $\CC$-vector space. 

For $s\in\tS$ let $B_s=R\ot_{R^s}R[1]\in\car$.
More generally, for any $x\in\tW$, Soergel \cite{\SO, 6.16} shows that there is an object $B_x$ of $\car$
(unique up to isomorphism) such that $B_x$ is an indecomposable direct summand of $B_{s_1}B_{s_2}\do B_{s_q}$
for some/any reduced expression $w=s_1s_2\do s_q$ ($s_i\in\tS$) and such that $B_x$ is not a direct summand 
of $B_{s'_1}B_{s'_2}\do B_{s'_p}$ whenever $s'_1,\do,s'_p\in\tS,p<q$.

Let $\tC$ be the full subcategory of $\car$ whose objects are isomorphic to finite direct sums of shifts of 
objects of the form $B_x$ for various $x\in\tW$. Let $C$ be the full subcategory of $\car$ whose objects are
isomorphic to finite direct sums of objects of the form $B_x$ for various $x\in\tW$.

Let $x\in\tW$. From \cite{\EW} it follows that $\Hom_\car(B_x,B_x)=\CC$ 
and from \cite{\SO, 6.16} it follows that $\dim\un{R_x^{B_x}}_{l(x)}=1$. Thus, as noted in \cite{\LV, 2.2},
$\un{R_x^{B_x}}_{l(x)}\ot_\CC B_x$ is an object of $C$ isomorphic to $B_x$ and
defined up to unique isomorphism
(even though $B_x$ was defined only up to non-unique isomorphism). From now on we will use  the notation 
$B_x$ for this new object; this agrees with the earlier description of $B_s$.

From \cite{\SO} it follows that for $M,M'\in\tC$ we have $MM'\in\tC$.
For any $x\in\tW$ let $R_x$ be the object of $\car$ such that $R_x=R$ as a left $R$-module and such that for
$m\in R_x,r\in R$ we have $mr=({}^xr)m$. The following result appears in \cite{\SO, 6.15}:

(a) {\it For any $M\in\tC$, $R_x^M$ is a finitely generated graded free right $R$-module; hence 
$\dim_\CC\un{R_x^M}<\iy$.}
\nl
We regard $\ck(\tC)$ as an $\ca$-module by $v^n[M]=[M[-n]]$ for $M\in\tC,n\in\ZZ$.
Note that $\ck(\tC)$ is an associative $\ca$-algebra with product defined by $[M][M']=[MM']$ for 
$M\in\tC,M'\in\tC$. From \cite{\SO, 1.10, 5.3} we see that 

(b) {\it the assignment $M\m\sum_{y\in W,i\in\ZZ}\dim\un{R_y^M)}_iv^{-i+\tl(y)}T_y$ defines an $\ca$-algebra 
isomorphism $\c:\ck(\tC)@>\si>>\tch$.}
\nl
From \cite{\EW, Theorem 1.1} it follows that for $x\in\tW$ we have
$$\c(B_x)=\ti c_x.\tag c$$

\subhead 16.3\endsubhead
For $M\in\car$ let $M^\sh$ be the object of $\car$ which is equal to $M$ as a graded $\CC$-vector space, but
left (resp. right) multiplication by $r\in R$ on $M^\sh$ equals left (resp. right) multiplication by $\io(r)$
on $M$. If $f:M_1@>>>M_2$ is a morphism in $\car$ then $f$ can be also viewed as a morphism 
$M^\sh_1@>>>M^\sh_2$ in $\car$. Clearly, $M\m M^\sh$ is a $\CC$-linear $\nn$-periodic functor 
$\car@>>>\car$. Hence $\car_\sh$ is well defined, see 16.1. If $M_1,M_2\in\car$ then the identity maps gives 
an identification $M^\sh_1M^\sh_2=(M_1M_2)^\sh$ as objects in $\car$.

Let $s\in\tS$. We define a $\CC$-linear isomorphism $\o_s:B_s[-1]@>\si>>B_{\io(s)}[-1]$ given by
$x\ot_{R^s}y\m\io(x)\ot_{R^{\io(s)}}\io(y)$ for $x,y\in R$. We have $\o_s(rfr')=\io(r)\o_s(f)\io(r')$ for
$r,r'\in R,f\in B_s[-1]$. Hence $\o_s$ can be viewed as an isomorphism $B_s[-1]@>\si>>B_{\io(s)}^\sh[-1]$ in
$\car$ or as an isomorphism $B_s@>\si>>B_{\io(s)}^\sh$ in $\car$.
Now let $x\in\tW$ and let $s_1s_2\do s_n$ be a reduced expression for $x$. Since $B_x$ is an indecomposable 
direct summand of $B_{s_1}B_{s_2}\do B_{s_n}$ (and $n$ is minimal with this property) we see that $B^\sh_x$ 
is an indecomposable direct summand of 
$$(B_{s_1}B_{s_2}\do B_{s_n})^\sh=B^\sh_{s_1}B^\sh_{s_2}\do B^\sh_{s_n}\cong 
B_{\io\i(s_1)}B_{\io\i(s_2)}\do B_{\io\i(s_n)}$$
(and $n$ is minimal with this property) hence by \cite{\SO, 6.16} we have 
$$B^\sh_x\cong B_{\io\i(x)}.\tag a$$. 
In particular we have $B^\sh_x\in\tC$. It follows that $M\in C\imp M^\sh\in C$ and $M\in\tC\imp M^\sh\in\tC$.
Note that $M\m M^\sh$ are $\CC$-linear, $\nn$-periodic functors $C@>>>C$ and $\tC@>>>\tC$. 
Hence $C_\sh,\tC_\sh$ are defined as in 16.1 and $\ck_\sh(C)$, $\ck_\sh(\tC)$ are well defined $\co$-modules. 

Now if $x\in W$ then from (a) we have that there exists an isomorphism $\ph:B^\sh_x@>\si>> B_x$ in $C$.
Replacing $\ph$ by $c\ph$ for a suitable $c\in\CC^*$ we can assume that $(B_x,\ph)\in C_\sh$. (We use that
$\End(B_x)=\CC$, see \cite{\EW}.)

\subhead 16.4\endsubhead
For $x\in\tW$ we define $\ff_x:R_x^\sh@>>>R_{\io\i(x)}$ by $r\m\io\i(r)$. This is an isomorphism in $\car$.
Now assume that $x\in W$; then $\ff_x:R_x^\sh@>>>R_x$ and $(R_x,\ff_x)\in\car_\sh$; thus,
$(R_x[i],\ff_x[i])\in\car_\sh$ for any $i\in\ZZ$. Hence, if $(M,\ph)\in\tC_\sh$ and $i\in\ZZ$, then 
$f\m f^!$, $\Hom_\car(M,R_x[i])@>>>\Hom_\car(M,R_x[i])$ is defined as in 16.1.
Taking direct sum over $i\in\ZZ$ we obtain a map $f\m f^!$, $R_x^M@>>>R_x^M$ whose $\nn$-th iteration is $1$.
Since $R^{>0}$ is $\io$-stable we see that $f\m f^!$ maps $R_x^MR^{>0}$ into iself hence it induces an
$\CC$-linear map $\un{R_x^M}@>>>\un{R_x^M}$ and (for any $i$) an $\CC$-linear map 
$\un{R_x^M}_i@>>>\un{R_x^M}_i$ (whose $\nn$-th power is $1$), denoted by $\cz_{x,\ph,i}^M$. Let
$$\e^x_i(M,\ph)=\tr_\CC(\cz_{x,\ph,i}^M,\un{R_x^M}_i)\in\co.$$
We now take $M=B_x$ (still assuming $x\in W$) so that $(B_x,\ph)\in C_\sh$ for some $\ph$. 
Then $\un{R_x^{B_x}}_{\tl(x)}=\CC$ hence $\e^x_{\tl(x)}(B_x,\ph)$ is an $\nn$-th root of $1$ in $\CC$. We 
can normalize $\ph:B_x^\sh@>>>B_x$ uniquely so that $\e^x_{\tl(x)}(B_x,\ph)=1$. We shall denote this 
normalized $\ph$ by $\ph_x$.

Next we note that if $x,x'\in\tW$ then by \cite{\EW}, $\Hom_\car(B_{x'},B_x)$ is $\CC$ if $x=x'$ and is $1$
if $x\ne x'$. It follows that $C$ is a semisimple abelian category and the $B_x$ are its simple objects.
Using this and \cite{\QG, 11.1.8} we deduce that

(a) $\ck_\sh(C)$ is the free $\co$-module with basis $\{[B_x,\ph_x];x\in\tW\}$.

\subhead 16.5\endsubhead
Let $\co'=\co[v,v\i]$ where $v$ is an indeterminate. We view $\ck_\s(\tC)$ as an $\co'$-module with 
$v^n[M,\ph]=[M[-n],\ph]$ for $(M,\ph)\in\tC_\sh$, $n\in\ZZ$. We have the following result.

(a) {\it The $\co'$-linear map $q:\co'\ot_\co\ck_\sh(C)@>>>\ck_\sh(\tC)$ given by 
$v^n\ot[M,\ph]\m[M[-n],\ph]$ is an isomorphism.}
\nl
The map $q$ is clearly well defined. To prove 
that it is surjective we shall use the functors $M\m\t_{\le i}M$ from $\tC$ to $\tC$ (resp. $M\m\ch^iM$ from
$\tC$ to $C$) defined in \cite{\EW, 6.2}. (Here $i\in\ZZ$.) These define in an obvious way functors 
$\tC_\sh@>>>\tC_\sh$ (resp. $\tC_\sh@>>>C_\sh$) denoted again by $\t_{\le i}$ (resp. $\ch^i$). Let 
$(M,\ph)\in\tC_\ph$. From the definition we have an exact sequence in $\tC$ (with morphisms in $\tC_\sh$)
$$0@>>>\t_{\le i-1}M@>e>>\t_{\le i}M@>e'>>\ch^iM[-i]@>>>0$$
which is split but the splitting is not necessarily given by morphisms in $\tC_\sh$. Thus there exist 
morphisms 
$$\t_{\le i-1}M@<f<<\t_{\le i}M@<f'<<\ch^iM[-i]$$
in $\tC$ such that $e'f'=1$, $fe=1$, $f'e'+ef=1$. Now $f^!,f'{}^!$ are defined as in 1.1 and, since $e^!=e$,
$e'{}^!=e'$ (notation of 16.1), we have $e'f'{}^!=1$, $f^!e=1$, $f'{}^!e'+ef^!=1$ hence, setting 
$\tf=(f+f^!+(f^1)^!+\do)/\nn$, $\tf'=(f'+f'{}^!+(f'{}^!)^!+\do)/\nn$ (the last two sums have $\nn$ terms) we
have $e'\tf'=1$, $\tf e=1$, $\tf'e'+e\tf=1$ and $\tf^!=\tf$, $\tf'{}^!=\tf'$. Thus we obtain a new splitting
of the exact sequence above which is given by morphisms in $\tC_\sh$. It follows that
$$(\t_{\le i}M,\ph)\cong(\t_{\le i-1}M,\ph)\op(\ch^iM[-i],\ph)$$
in $\tC_\sh$ (the maps $\ph$ are induced by $M^\sh@>>>M$). Hence 
$[\t_{\le i}M,\ph]=[\t_{\le i-1}M,\ph]+[\ch^iM[-i],\ph]$ in $\ck_\sh(\tC)$. Since 
$[M,\ph]=[\t_{\le i}M,\ph]$ for $i\gg0$ and $0=[\t_{\le i}M,\ph]$ for $-i\gg0$ we deduce that
$[M,\ph]=\sum_i[\ch^iM[-i],\ph]$. This proves the surjectivity of $q$.

We define $\ck(\tC_\sh)@>>>\co'\ot\ck(C_\sh)$ by $[M,\ph]\m\sum_{n\in\ZZ}v^{-n}[\ch^nM,\ph_n]$ where 
$\ph_n$ is induced by $\ph$. This clearly induces a homomorphism $q':\ck_\sh(\tC)@>>>\ca\ot\ck_\sh(C)$ which
satisfies $q'q=1$. It follows that $q$ is injective, completing the proof of (a).

\mpb

Using (a) and 16.4(a) we see that 

(b) $\ck_\sh(\tC)$ is a free $\co'$-module with basis $\{[B_x,\ph_x];x\in W\}$, (notation of 16.4).

\subhead 16.6\endsubhead
Let $\cn$ be the free $\co'$-module with basis $\{b_x;x\in W\}$. For any $(M,\ph)\in\tC_\sh$ and any 
$y\in W$ we set
$$\e^y(M,\ph)=\sum_{i\in\ZZ}\e^y_i(M,\ph)v^{-i}\in\co'.$$ 
The homomorphism $\ck(\tC_\sh)@>>>\cn$, 
$$[M,\ph]\m\sum_{y\in W}\e^y(M,\ph)v^{\tl(y)}b_y,$$ 
clearly factors through an $\co'$-module homomorphism 
$$\c':\ck_\sh(\tC)@>>>\cn.\tag a$$
We show:
$$\c'\text{ is an isomorphism}.\tag b$$
For $x\in W$ let $\tb_x=\c'([B_x,\ph_x])$. We can write $\tb_x=\sum_{y\in W}f_{y,x}b_y$ where 
$f_{y,x}\in\co'$ are zero for all but finitely many $y$. In view of 16.5(b), to prove (b) it is enough to 
show:

(c) {\it Let $y\in W$. If $y\not\le x$ then $f_{y,x}=0$. If $y\le x$ then 
$f_{y,x}=\dot p_{y,x}(v)$ where $\dot p_{x,x}=1$ and $\dot p_{y,x}\in v\i\co[v\i]$ if $y<x$.}
\nl
Assume that $f_{y,x}\ne0$. Then for some $i$ we have $\e^y_i(B_x,\ph_x)\ne0$ hence $\un{R_y^{B_x}}\ne0$.
Using 16.4(b),(c), we deduce that the coefficient of $T_y$ in $c_x$ is nonzero; thus we have $y\le x$, as
required. Next we assume that $y\le x$. We have $f_{y,x}=\sum_i\e^y_i(B_x,\ph_x)v^{-i+\tl(y)}$ hence it is 
enough to show that 

$\e^y_i(B_x,\ph_x)\ne0$ implies $-i+\tl(y)\le0$, with strict inequality unless $x=y$. 
\nl
Now $\e^y_i(B_x,\ph_x)\ne0$ implies $\un{R_y^{B_x}}_i\ne0$. Hence it is enough to show that 

$\un{R_y^{B_x}}_i\ne0$ implies $-i+\tl(y)\le0$, with strict inequality unless $x=y$. 
\nl
By 16.2(b),(c), we have 
$$\sum_{j\in\ZZ}\dim\un{R_y^{B_x}}_ju^{-j+\tl(y)}=\ti p_{y,x}(u)$$
and it remains to use that $\ti p_{x,x}=1$ and $\ti p_{y,x}\in v\i\ZZ[v\i]$ if $y<x$. This proves (c) hence 
also (b).

\subhead 16.7\endsubhead
If $(M,\ph)\in\tC_\sh$ and $(M',\ph')\in\tC_\sh$ then $(MM',\ph\ot\ph')$ is again an object of $\tC_\sh$.
Note that if $(M,\ph)$ or $(M',\ph')$ is traceless then $(MM',\ph\ot\ph')$ is again traceless. It follows
easily that $((M,\ph),(M',\ph'))\m(MM',\ph\ot\ph')$ defines an $\co'$-bilinear map
$\ck_\sh(\tC)\T\ck_\sh(\tC)@>>>\ck_\sh(\tC)$ which makes $\ck_\sh(\tC)$ into an associative $\co'$-algebra
with $1$. (The unit element is $[B_1,\ph_1]$.) Via the isomorphism 16.6(a),(b), we obtain an associative 
$\co'$-algebra structure (with $1=b_1$) on $\cn$. One can show that the following identities hold in the 
algebra $\cn$:
$$b_wb_w'=b_{ww'}\text{ if }w,w'\in W, l(ww')=l(w)+l(w');\tag a$$
$$(b_{w_0^I}-v^{\tl(w_0^I)})(b_{w_0^I}+v^{-\tl(w_0^I)})=0\text{ for any }I\in S.\tag b$$
Thus $\cn$ can be identified with the $\co'$-algebra $\co'\ot_\ca\ch$ where $\ch$ is as in 16.2
in such a way that $b_w$ corresponds to $T_w\in\ch$ for any $w\in W$.

\subhead 16.8\endsubhead
For $M\in\tC$ let $D(M)\in\tC$ be the ``dual'' of $M$ defined as in \cite{\SO, 5.9}. Now 
$(M,\ph)\m(D(M),D(\ph)\i)$ induces a ring homomorphism $\bar{}:\ck_\sh(\tC)@>>>\ck_\sh(\tC)$ which maps 
$[B_x,\ph_x]$ to itself for any $x\in W$ and is semilinear with respect to the ring involution 
$\bar{}:\co'@>>>\co'$ given by $v^n\m v^{-n}$ and $\z\m\z\i$ for any $\z\in\co$ such that $\z^\nn=1$. Via 
the isomorphism 16.6(a),(b), this becomes an $\co'$-semilinear ring homomorphism $\bar{}:\cn@>>>\cn$ such 
that $\ov{\tb_x}=\tb_x$ for any $x\in W$. From the results in 16.6 for any $I\in S$ we have
$\tb_{w_0^I}=b_{w_0^I}+z_I$ where $z_I\in v\i\co[v\i]$ satisfies
$\ov{b_{w_0^I}+z_I}=b_{w_0^I}+z_I$ that is $\ov{b_{w_0^I}}=b_{w_0^I}+z_I-\ov{z_I}$. 
Applying $\bar{}$ to the equation 
$b_{w_0^I}^2=(v^n-v^{-n})\b_{w_0^I}+1$
where $n=\tl(w_0^I)$ (see 16.7(b)) we obtain
$$(b_{w_0^I}+z_I-\ov{z_I})^2=(v^{-n}-v^n)(b_{w_0^I}+z_I-\ov{z_I})+1$$
hence
$$(v^n-v^{-n})b_{w_0^I}+1+2(z_I-\ov{z_I})b_{w_0^I}+(z_i-\ov{z_I})^2=(v^{-n}-v^n)(b_{w_0^I}+z_I-\ov{z_I})+1.$$
We deduce that $z_I-\ov{z_I}=v^{-n}-v^n$ that is $\ov{z_I-v^{-n}}=z_I-v^{-n}$. Since $n>0$ we have
$z_I-v^{-n}\in\co[v\i]$ hence $z_I-v^{-n}=0$. Thus we have
$$\tb_{w_0^I}=b_{w_0^I}+v^{-\tl(w_0^I)}.$$
Since this holds for every $I\in S$ we see that under our identification $\cn=\co'\ot_\ca\ch$, 
$\bar{}:\cn@>>>\cn$ corresponds to $\bar{}:\ch@>>>\ch$ (as in 4.1) extended semilinearly to 
$\co'\ot_\ca\ch$. For
$w\in W$, both $\tb_w$ and $c_w$ are fixed by $\bar{}$ and (by results in 16.6) their difference is in
$\sum_{y\in W;y<w}v\i\ZZ[v\i]T_y$; it follows that $\tb_w=c_w$. In particular we have $\tb_w\in\ch$.

\subhead 16.9\endsubhead
Let $x,y\in W$ be such that $y\le x$ and let $n\in\ZZ$. We show:

(a) If the coefficient of $v^n$ in $p_{y,x}$ is $\ne0$ then the coefficient of $v^n$ in $\tp_{y,x}$ is 
$\ne0$;

(b) if the coefficient of $v^n$ in $\tp_{y,x}$ is $1$ then the coefficient of $v^n$ in $p_{y,x}$ is $\pm1$.
\nl
In the setup of (a), the coefficient of $v^n$ in $f_{y,x}$ is $\ne0$ (notation of 16.6). Hence 
$\e^y_{\tl(y)-n}(B_x,\ph_x)\ne0$ (see 16.6) so that
$$\tr_\CC(\cz_{y,\ph_x,\tl(y)-n}^{B_x},\un{R_y^{B_x}}_{\tl(y)-n})\ne0$$ 
(notation of 16.4); in particular we have $\un{R_y^M}_{\tl(y)-n}\ne0$. From 16.2(b) we see that $v^n$ appears 
in the coefficient of $T_y$ in $\c(B_x)$ with $\ne$ coefficient; from 16.2(c) we deduce that $v^n$ appears in
$\tp_{y,x}$ with $\ne$ coefficient. This proves (a).

In the setup of (b), using 16.2(b),(c) we see that $\dim\un{R_y^M}_{\tl(y)-n}=1$. Hence (with notation of 
16.4), $\tr_\CC(\cz_{y,\ph_x,\tl(y)-n}^{B_x},\un{R_y^{B_x}}_{\tl(y)-n})$ is the trace of a linear 
transformation of finite order of a one dimensional vector space so that it is a root of $1$. Thus,  
$\e^y_{\tl(y)-n}(B_x,\ph_x)0$ is a root of $1$ so that by 16.6, the coefficient of $v^n$ in $f_{x,y}$ is a 
root of $1$ and the coefficient of $v^n$ in $p_{y,x}$ is a root of $1$. Since $p_{y,x}$ has integer
coefficients, the coefficient of $v^n$ in $p_{y,x}$ is $\pm1$. This proves (b).

\subhead 16.10\endsubhead
Let $x,y,z$ in $W$ and let $n\in\ZZ$. We form $(B_xB_y,\ph_x\ot\ph_y)\in\tC_\sh$. Now $\ph_x\ot\ph_y$ induces
an isomorphism $\ps_n:(\ch^n(B_xB_y))^\sh@>>>\ch^n(B_xB_y))$ so that $(\ch^n(B_xB_y),\ps_n)\in\tC_\sh$
(notation of 16.5). Let $V_{x,y,z}^n=\Hom_\car(B_z,\ch^n(B_xB_y))$. We can find a linear isomorphism 
$\th:V_n@>>>V_n$ of finite order such that under the obvious imbedding $V_{x,y,z}^n\ot B_z@>>>\ch^n(B_xB_y)$,
$\th\ot\ph_z$ is compatible with $\ps_n$. From the definitions, the coefficient of $v^n$ in $h_{x,y,z}$ is 
equal to $\tr(\th,V_{x,y,z}^n)$. We show:

(a) If the coefficient of $v^n$ in $h_{x,y,z}$ is $\ne0$ then the coefficient of $v^n$ in $\ti h_{x,y,z}$ is 
$\ne0$;

(b) if the coefficient of $v^n$ in $\ti h_{x,y,z}$ is $1$ then the coefficient of $v^n$ in $h_{x,y,z}$ is 
$\pm1$.
\nl
In the setup of (a) we have $\tr(\th,V_{x,y,z}^n)\ne0$ hence $V_{x,y,z}^n\ne0$. Thus $B_z$ appears with 
nonzero multiplicity in $\ch^n(B_xB_y)$. From the definitions we see that the coefficient of $v^n$ in 
$\ti h_{x,y,z}$ is $\ne0$. Thus (a) holds.

In the setup of (b) we have $\dim V_{x,y,z}^n=1$. Since $\th:V_{x,y,z}^n@>>>V_{x,y,z}^n$ has finite order,
it follows that $\tr(\th,V_{x,y,z}^n)$ is a root of $1$. Hence the coefficient of $v^n$ in $h_{x,y,z}$ is a 
root of $1$; but that coefficient is an integer hence it is $\pm1$. Thus (b) holds.

\subhead 16.11\endsubhead
We show that for $x,y,z$ in $W$ we have
$$v^{-L(z)}h_{x,y,z}\in\ZZ[v\i].\tag a$$
We must show that if $n\in\ZZ$ and the coefficient of $v^n$ in $h_{x,y,z}$ is $\ne0$ then $n\le L(z)$. By 
16.10(a), the coefficient of $v^n$ in $\ti h_{x,y,z}$ is $\ne0$; hence by 15.2(b) (applied to $\tW,\tl$) we 
have $n\le\tl(z)$ that is $n\le L(z)$, as required.
Using (a) we see that definition of $\aa(z)\in\NN$, $\D(z)\in\NN$ (relative to $L$) and $\g_{x,y,z}\in\ZZ$ 
(for $x,y,z$ in $W$) and $\cd$ as in \S13 makes sense for $W,L$ even without the assumption 
that $W,L$ is bounded. We shall denote by $\ti\aa(z)$, $\ti\D(z)$, $\ti\g_{x,y,z}$ (for $x,y,z$ in $\tW$) and
$\tcd$ the analogous objects defined in terms of $\tW,\tl$, see 15.2.

We now make the additional assumption that $\tW,\tl$ is bounded. We show:

(b) {\it $W,L$ is bounded.}
\nl
By assumption there exists $N\ge0$ such that for all $x,y,z$ in $\tW$ we have
$v^{-N}\ti h_{x,y,z}\in\ZZ[v\i]$. In particular, if $x,y,z\in W$ and $n$ is such that the coefficient of 
$v^n$ in $\ti h_{x,y,z}$ is nonzero then $n\le N$; using 16.10(a) we deduce that, if $x,y,z\in W$ and $n$ is 
such that the coefficient of $v^n$ in $h_{x,y,z}$ is nonzero then $n\le N$, so that 
$v^{-N}h_{x,y,z}\in\ZZ[v\i]$. This proves (b).

Under the assumption that $\tW,\tl$ is bounded, the results of \S15 are applicable to $\tW,\tl$;
in this chapter we will show that, under the same assumption, P1-P15 hold for $W,L$.

\proclaim{Lemma 16.12} For $z\in W$ we have $\aa(z)=\taa(z)$ and $\ti\D(z)\le\D(z)$.
\endproclaim
We can find $x,y\in W$ such that $\p_{\aa(z)}(h_{x,y,z})\ne 0$. By 16.10(a) we have
\lb$\p_{\aa(z)}(\ti h_{x,y,z})\ne 0$. Hence $\aa(z)\le\taa(z)$. By P3,P5 for $\tW$,
there is a unique $d\in\tcd$ such that $\ti\g_{z\i,z,d}=\pm 1$. The uniqueness of 
$d$ implies that $d$ is fixed by $u$. Thus $d\in W$. By P7 for $\tW$, we have 
$\ti\g_{z,d,z\i}=\pm 1$. Hence $\p_{\taa(z)}(\ti h_{z,d,z})=\pm 1$. By 16.10(b), we 
have $\p_{\taa(z)}(h_{z,d,z})=\pm 1$. Hence $\taa(z)\le\aa(z)$ so that
$\taa(z)=\aa(z)$.

By definition, we have $\p_{-\D(z)}(p_{1,z})\ne 0$. Using 16.9(a), we deduce that
\lb$\p_{-\D(z)}(\tp_{1,z})\ne 0$. Hence $-\D(z)\le-\ti\D(z)$. The lemma is proved.

\proclaim{Lemma 16.13} We have $\cd=\tcd\cap W$.
\endproclaim
Let $d\in\cd$. We have $\aa(d)=\D(d)$. Using 16.12, we deduce $\taa(d)=\D(d)$. By P1
for $\tW$, we have $\taa(d)\le\ti\D(d)$. Hence $\D(d)\le\ti\D(d)$. Using 16.12, we
deduce $\D(d)=\ti\D(d)$ so that $\ti\D(d)=\taa(d)$ and $d\in\tcd$.

Conversely, let $d\in\tcd\cap W$. We have $\taa(d)=\ti\D(d)$. Using 16.12 we deduce
$\aa(d)=\ti\D(d)$. By P5 for $\tW$, we have $\p_{-\ti\D(d)}(\tp_{1,d})=\pm 1$. Using
16.9(b) we deduce $\p_{-\ti\D(d)}(p_{1,d})=\pm 1$. Hence $-\ti\D(d)\le-\D(d)$. Using
16.12 we deduce $\D(d)=\ti\D(d)$ so that $\D(d)=\aa(d)$ and $d\in\cd$. The lemma is 
proved.

\proclaim{Lemma 16.14} (a) Let $x,y,z\in W$ be such that $\g_{x,y,z}\ne 0$. Then 
$\ti\g_{x,y,z}\ne 0$. 

(b) Let $x,y,z\in W$ be such that $\ti\g_{x,y,z}=\pm 1$. Then $\g_{x,y,z}=\pm 1$.
\endproclaim
In the setup of (a) we have $\p_{\aa(z\i)}(h_{x,y,z\i})\ne 0$. Using 16.12 we deduce
that\lb $\p_{\taa(z\i)}(h_{x,y,z\i})\ne 0$. Using 16.10(a), we deduce that
$\p_{\taa(z\i)}(\ti h_{x,y,z\i})\ne 0$. Hence $\ti\g_{x,y,z}\ne 0$. 

In the setup of (b) we have $\p_{\taa(z\i)}(\ti h_{x,y,z\i})=\pm 1$. Using 16.12, we
deduce $\p_{\aa(z\i)}(\ti h_{x,y,z\i})=\pm 1$. Using 16.10(b), we deduce
$\p_{\aa(z\i)}(h_{x,y,z\i})=\pm 1$. Hence $\g_{x,y,z}=\pm 1$. 

\subhead 16.15\endsubhead 
{\it Proof of P1.} By 16.12 and P1 for $\tW$, we have 
$\aa(z)=\taa(z)\le\ti\D(z)\le\D(z)$, hence $\aa(z)\le\D(z)$.

\subhead 16.16\endsubhead 
{\it Proof of P2.} In the setup of P2, we have (by 16.14) $\ti\g_{x,y,d}\ne 0$ and 
$d\in\tcd$ (see 16.13). Using P2 for $\tW$, we deduce $x=y\i$.

\subhead 16.17\endsubhead 
{\it Proof of P3.} Let $y\in W$. By P3 for $\tW$, there is a unique $d\in\tcd$ such
that $\ti\g_{y\i,y,d}\ne 0$. By the uniqueness of $d$, we have $u(d)=d$ hence
$d\in W$. Using P5 for $\tW$, we see that $\ti\g_{y\i,y,d}=\pm 1$. Using 16.14, we
deduce $\g_{y\i,y,d}=\pm 1$. Since $d\in\cd$ by 16.13, the existence part of P3 is 
established. Assume now that $d'\in\cd$ satisfies $\g_{y\i,y,d'}\ne 0$. Using 16.14,
we deduce $\ti\g_{y\i,y,d'}\ne 0$. Since $d'\in\tcd$ by 16.13, we can use the 
uniqueness in P3 for $\tW$ to deduce that $d=d'$. Thus P3 holds for $W$.

\subhead 16.18\endsubhead 
{\it Proof of P4.} We may assume that there exists $s\in S$ such that 
$h_{s,z,z'}\ne 0$ or $h_{z,s,z'}\ne 0$. In the first case, using 16.10(a), we deduce
$\ti h_{s,z,z'}\ne 0$. Hence $z'\le_\cl z$ (in $\tW$) and using P4 for $\tW$, we
deduce that $\taa(z')\ge\ti\aa(z)$. Using now 16.12, we see that $\aa(z')\ge\aa(z)$.
The proof in the second case is entirely similar.

\subhead 16.19\endsubhead 
Now P5 is proved as in 14.5; P6 is proved as in 14.6; P7 is proved as in 14.7; P8 is
proved as in 14.8; P12 is proved as in 14.12.

\subhead 16.20\endsubhead 
{\it Proof of P13.} If $z'\gt_\cl z$ in $W$, then there exists $s\in S$ such that
$h_{s,z,z'}\ne 0$ hence, by 16.10(a), $\ti h_{s,z,z'}\ne 0$, hence $z'\le_\cl z$ in 
$\tW$. It follows that 

(a) $z'\le_\cl z$ (in $W$) implies $z'\le_\cl z$ (in $\tW$). 
\nl
Hence 

(b) $z'\si_\cl z$ (in $W$) implies $z'\si_\cl z$ (in $\tW$). 
\nl
Thus any left cell of $W$ is contained in a left cell of $\tW$. 

In the setup of P13, let $\ti\G$ be the left cell of $\tW$ containing $\G$. Let
$x\in\G$. By P3 for $W$, there exists $d\in\cd$ such that $\g_{x\i,x,d}\ne 0$. By 
P8 for $W$, we have $x\si_\cl d\i$ hence $d\i\in\G$. Using P6 we have $d=d\i$, hence
$d\in\G$. It remains to prove the uniqueness of $d$. Let $d',d''$ be elements of 
$\cd\cap\G$. We must prove that $d'=d''$. Now $d',d''$ belong to $\ti\G$ and, by
16.13, are in $\tcd$. Using P13 for $\tW$, it follows that $d'=d''$. Thus P13 holds 
for $W$.

\proclaim{Lemma 16.21} Let $x,y\in W$. We have $x\si_\cl y$ (in $W$) if and only if
$x\si_\cl y$ (in $\tW$).
\endproclaim
If $x\si_\cl y$ (in $W$) then $x\si_\cl y$ (in $\tW$), by 16.20(b). 

Assume now that $x\si_\cl y$ (in $\tW$). Let $d,d'\in\cd$ be such that $x\si_\cl d$
(in $W$) and $y\si_\cl d'$ (in $W$); see P13. By the first line of the proof we have
$x\si_\cl d$ (in $\tW$) and $y\si_\cl d'$ (in $\tW$). Hence $d\si_\cl d'$ (in 
$\tW$). Since $d,d'\in\tcd$, we deduce (using P13 for $\tW$) that $d=d'$. It follows
that $x\si_\cl y$ (in $W$). The lemma is proved.

\subhead 16.22\endsubhead 
{\it Proof of P9.} We assume that $z'\le_\cl z$ (in $W$) and $\aa(z')=\aa(z)$. By
16.20(a), it follows that $z'\le_\cl z$ (in $\tW$) and, using 16.12, that 
$\taa(z')=\taa(z)$. Using now P9 in $\tW$, it follows that $z'\si_\cl z$ (in $\tW$).
Using 16.21, we deduce that $z'\si_\cl z$ (in $W$).

\subhead 16.23\endsubhead 
Now P10 is proved as in 14.10; P11 is proved as in 14.11; P14 is proved as in 14.14.

\subhead 16.24\endsubhead 
We sketch a proof of P15 in our case.

A refinement of the proof of P15 given in 14.15, 15.7 provides, for any $w,y,x,x'$ in
$\tW$ and any $k$, an isomorphism of vector spaces
$$\op_{j+j'=k}\op_{y'\in\tW}V^{j'}_{w,x',y'}\ot V^j_{x,y',y}@>\si>>
\op_{j+j'=k}\op_{y'\in\tW}V^j_{x,w,y'}\ot V^{j'}_{y',x',y}.$$
which (assuming that $\taa(w)=\taa(y)$) restricts to an isomorphism 
$$\op_{y'\in\tW}V^{j'}_{w,x',y'}\ot V^j_{x,y',y}@>\si>>
\op_{y'\in\tW}V^j_{x,w,y'}\ot V^{j'}_{y',x',y}$$
for any $j,j'$ such that $j+j'=k$.

Assuming now that $w,y,x,x'\in W$, we can take traces of $u$ in both sides; we
deduce
$$\su_{y'\in W}\p_{j'}(h_{w,x',y'})\p_j(h_{x,y',y})=
\su_{y'\in W}\p_j(h_{x,w,y'})\p_{j'}(h_{y',x',y})$$
(the summands corresponding to $y'\in\tW-W$ do not contribute to the trace) or
equivalently
$$\su_{y'\in W}h'_{w,x',y'}h_{x,y',y}=\su_{y'\in W}h_{x,w,y'}h'_{y',x',y},$$
as required. 

\head 17. Example: the infinite dihedral case\endhead
\subhead 17.1\endsubhead
In this chapter we preserve the setup of 7.1. We assume that $m=\iy$ and that 
$L_2>L_1$. We will show that P1-P15 hold in this case.

Let $\z=v^{L_2-L_1}+v^{L_1-L_2}$. For $a\in\{1,2\}$, let $f_a=v^{L_a}+v^{-L_a}$. For
$m,n\in\bz$ we define $\d_{m<n}$ to be $1$ if $m<n$ and to be $0$ otherwise.

\subhead 17.2\endsubhead
From 7.5, 7.6 we have for all $k'\in\bn$:

$c_1c_{2_{k'}}=c_{1_{k'+1}}$,

$c_2c_{1_{k'}}=c_{2_{k'+1}}+\d_{k'>1}\z c_{2_{k'-1}}+\d_{k'>3}c_{2_{k'-3}}$.

\proclaim{Proposition 17.3} For $k\ge 0,k'\ge 1$ we have

(a) $c_{2_{2k+1}}c_{2_{k'}}=f_2\su_{u\in[0,k];2u\le k'-1}c_{2_{2k+k'-4u}}$,

(b) $c_{1_{2k+2}}c_{2_{k'}}=f_2\su_{u\in[0,k];2u\le k'-1}c_{1_{2k+k'+1-4u}}$.
\endproclaim
Assume that $k=0$. Using 17.2 we have $c_2c_{2_{k'}}=f_2c_{2_{k'}}$.

Assume now that $k=1$. Using 17.2, we have $c_{2_3}=c_2c_1c_2-\z c_2$. Using this 
and 17.2, we have
$$\align&c_{2_3}c_{2_{k'}}=c_2c_1c_2c_{2_{k'}}-\z c_2c_{2_{k'}}=
f_2c_2c_{1_{k'+1}}-f_2\z c_{2_{k'}}\\&=
f_2c_{2_{k'+2}}+f_2\z c_{2_{k'}}+\d_{k'>2}f_2c_{2_{k'-2}}
-f_2\z c_{2_{k'}}=f_2c_{2_{k'+2}}+\d_{k'>2}f_2c_{2_{k'-2}},\endalign$$
as required. We prove the equality in (a) for fixed $k'$, by induction on $k$. The 
cases $k=0,1$ are already known. If $k=2$ then using 17.2, we have
$c_{2_5}=c_2c_1c_{2_3}-\z c_{2_3}-c_{2_1}$. Using this, 17.2, and the induction 
hypothesis, we have
$$\align&c_{2_5}c_{2_{k'}}=c_2c_1c_{2_3}c_{2_{k'}}-\z c_{2_3}c_{2_{k'}}
-c_{2_1}c_{2_{k'}}\\&=f_2c_2c_1c_{2_{k'+2}}+\d_{k'>2}f_2c_2c_1c_{2_{k'-2}}
-\z f_2c_{2_{k'+2}}-\d_{k'>2}\z f_2c_{2_{k'-2}}-f_2c_{2_{k'}}\\&
=f_2c_2c_{1_{k'+3}}+\d_{k'>2}f_2c_2c_{1_{k'-1}}-\z f_2c_{2_{k'+2}}-\d_{k'>2}
\z f_2c_{2_{k'-2}}-f_2c_{2_{k'}}\\&=f_2c_{2_{k'+4}}+f_2\z c_{2_{k'+2}}+f_2
c_{2_{k'}}+\d_{k'>2}f_2c_{2_{k'}}+\d_{k'>2}f_2\z c_{2_{k'-2}}+\d_{k'>4}f_2
c_{2_{k'-4}}\\&-\z f_2c_{2_{k'+2}}-\d_{k'>2}\z f_2c_{2_{k'-2}}-f_2c_{2_{k'}}
=f_2c_{2_{k'+4}}+\d_{k'>2}f_2c_{2_{k'}}+\d_{k'>4}f_2c_{2_{k'-4}},\endalign$$
as required. A similar argument applies for $k\ge 3$. This proves (a).

(b) is obtained by multiplying both sides of (a) by $c_1$ on the left. The 
proposition is proved.

\proclaim{Proposition 17.4} For $k\ge 0,k'\ge 1$, we have 

(a) $c_{2_{2k+1}}c_{1_{k'}}=\su_{u\in[0,2k+2]}p_uc_{2_{k'+2k+1-2u}}$,

(b) $c_{1_{2k+2}}c_{1_{k'}}=\su_{u\in[0,2k+2]}p_uc_{1_{k'+2k+2-2u}}$,

(c) $c_{1_{k'}\i}c_{2_{2k+1}}=\su_{u\in[0,2k+2]}p_uc_{2_{k'+2k+1-2u}\i}$,

(d) $c_{1_{k'}\i}c_{1_{2k+2}\i}=\su_{u\in[0,2k+2]}p_uc_{1_{k'+2k+1-2u}\i}$,

(e) $c_{2_{2k+2}}c_{1_{k'}}=\su_{u\in[0,2k+2]}f_1p_uc_{2_{k'+2k+1-2u}}$,    

(f) $c_{1_{2k+3}}c_{1_{k'}}=\su_{u\in[0,2k+2]}f_1p_uc_{1_{k'+2k+2-2u}}$,

(g) $c_{1_1}c_{1_{k'}}=f_1c_{1_{k'}}$,
\nl
where 

$p_0=1$, $p_{2k+2}=\d_{k'>2k+3}$,

$p_u=\d_{k'>u}\z$ for $u=1,3,5,\do,2k+1$,

$p_u=\d_{k'>u-1}+\d_{k'>u+1}$ for $u=2,4,6,\do,2k$.
\endproclaim
We prove (a). For $k=0$ the equality in (a) is
$c_2c_{1_{k'}}=c_{2_{k'+1}}+\d_{k'>1}\z c_{2_{k'-1}}+\d_{k'>3}c_{2_{k'-3}}$
which is contained in 17.2. Assume now that $k=1$. Using 
$c_{2_3}=c_2c_1c_2-\z c_2$ and 17.2, we have
$$\align&c_{2_3}c_{1_{k'}}=c_2c_1c_2c_{1_{k'}}-\z c_2c_{1_{k'}}=c_2c_1
c_{2_{k'+1}}+\d_{k'>1}\z c_2c_1c_{2_{k'-1}}+\d_{k'>3}c_2c_1c_{2_{k'-3}}\\&-
\z c_{2_{k'+1}}-\d_{k'>1}\z^2 c_{2_{k'-1}}-\d_{k'>3}\z c_{2_{k'-3}}=c_2
c_{1_{k'+2}}+\d_{k'>1}\z c_2c_{1_{k'}}\\&+\d_{k'>3}c_2c_{1_{k'-2}}-\z 
c_{2_{k'+1}}-\d_{k'>1}\z^2 c_{2_{k'-1}}-\d_{k'>3}\z c_{2_{k'-3}}\\&=
c_{2_{k'+3}}+\z c_{2_{k'+1}}+\d_{k'>1}c_{2_{k'-1}}+\d_{k'>1}\z c_{2_{k'+1}}
+\d_{k'>1}\z^2c_{2_{k'-1}}+\d_{k'>3}\z c_{2_{k'-3}}\\&+\d_{k'>3}
c_{2_{k'-1}}+\d_{k'>3}\z c_{2_{k'-3}}+\d_{k'>5}c_{2_{k'-5}}-\z c_{2_{k'+1}}
-\d_{k'>1}\z^2c_{2_{k'-1}}-\d_{k'>3}\z c_{2_{k'-3}}\\&=c_{2_{k'+3}}+
\d_{k'>1}\z c_{2_{k'+1}}+(\d_{k'>1}+\d_{k'>3})c_{2_{k'-1}}+\d_{k'>3}\z 
c_{2_{k'-3}}+\d_{k'>5}c_{2_{k'-5}},\endalign$$
as required.

We prove the equality in (a) for fixed $k'$, by induction on $k$. The cases $k=0,1$
are already known. Assume now that $k=2$. Using 
$c_{2_5}=c_2c_1c_{2_3}-\z c_{2_3}-c_{2_1}$, 17.2, and the case $k=1$, we have
$$\align&c_{2_5}c_{1_{k'}}=c_2c_1c_{2_3}c_{1_{k'}}-\z c_{2_3}c_{1_{k'}}
-c_{2_1}c_{1_{k'}}\\&=c_2c_1c_{2_{k'+3}}+\d_{k'>1}\z c_2c_1c_{2_{k'+1}}+
(\d_{k'>1}+\d_{k'>3})c_2c_1c_{2_{k'-1}}+\d_{k'>3}\z c_2c_1c_{2_{k'-3}}\\&
+\d_{k'>5}c_2c_1c_{2_{k'-5}}-\z c_{2_{k'+3}}-\d_{k'>1}\z^2 c_{2_{k'+1}}-
(\d_{k'>1}+\d_{k'>3})\z c_{2_{k'-1}}\\&-\d_{k'>3}\z^2 c_{2_{k'-3}}-
\d_{k'>5}\z c_{2_{k'-5}}-c_{2_{k'+1}}-\d_{k'>1}\z c_{2_{k'-1}}-\d_{k'>3}
c_{2_{k'-3}}\\&=c_2c_{1_{k'+4}}+\d_{k'>1}\z c_2c_{1_{k'+2}}+(\d_{k'>1}+
\d_{k'>3})c_2c_{1_{k'}}+\d_{k'>3}\z c_2c_{1_{k'-2}}\\&+\d_{k'>5}c_2
c_{1_{k'-4}}-\z c_{2_{k'+3}}-\d_{k'>1}\z^2 c_{2_{k'+1}}-(\d_{k'>1}+
\d_{k'>3})\z c_{2_{k'-1}}\\&-\d_{k'>3}\z^2 c_{2_{k'-3}}-\d_{k'>5}\z 
c_{2_{k'-5}}-c_{2_{k'+1}}-\d_{k'>1}\z c_{2_{k'-1}}-\d_{k'>3}c_{2_{k'-3}}\\&=
c_{2_{k'+5}}+\z c_{2_{k'+3}}+c_{2_{k'+1}}+\d_{k'>1}\z c_{2_{k'+3}}+
\d_{k'>1}\z^2 c_{2_{k'+1}}+\d_{k'>1}\z c_{2_{k'-1}}\\&+(\d_{k'>1}+
\d_{k'>3})c_{2_{k'+1}}+(\d_{k'>1}+\d_{k'>3})\z c_{2_{k'-1}}+2\d_{k'>3}
c_{2_{k'-3}}+\d_{k'>3}\z c_{2_{k'-1}}\\&+\d_{k'>3}\z^2c_{2_{k'-3}}+
\d_{k'>5}\z c_{2_{k'-5}}+\d_{k'>5}c_{2_{k'-3}}+\d_{k'>5}\z c_{2_{k'-5}}+
\d_{k'>7}c_{2_{k'-7}}\\&-\z c_{2_{k'+3}}-\d_{k'>1}\z^2c_{2_{k'+1}}-
(\d_{k'>1}+\d_{k'>3})\z c_{2_{k'-1}}-\d_{k'>3}\z^2c_{2_{k'-3}}-\d_{k'>5}
\z c_{2_{k'-5}}\\&-c_{2_{k'+1}}-\d_{k'>1}\z c_{2_{k'-1}}-\d_{k'>3}
c_{2_{k'-3}}\\&=c_{2_{k'+5}}+\d_{k'>1}\z c_{2_{k'+3}}+(\d_{k'>1}+\d_{k'>3})
c_{2_{k'+1}}+\d_{k'>3} c_{2_{k'-3}}+\d_{k'>3}\z c_{2_{k'-1}}\\&+\d_{k'>5}
c_{2_{k'-3}}+\d_{k'>5}\z c_{2_{k'-5}}+\d_{k'>7}c_{2_{k'-7}}\\&=c_{2_{k'+5}}+
\d_{k'>1}\z c_{2_{k'+3}}+(\d_{k'>1}+\d_{k'>3})c_{2_{k'+1}}\\&+\d_{k'>3}\z
c_{2_{k'-1}}+(\d_{k'>3}+\d_{k'>5})c_{2_{k'-3}}+\d_{k'>5}\z c_{2_{k'-5}}+
\d_{k'>7}c_{2_{k'-7}}.\endalign$$
A similar argument applies for $k\ge 4$. This proves (a).

(b) is obtained by multiplying both sides of (a) by $c_1$ on the left; (c),(d) are 
obtained by applying $h\m h^\flat$ to both sides of (a),(b). We prove (e). We have
$$c_{2_{2k+2}}c_{1_{k'}}=c_{2_{2k+1}}c_1c_{1_{k'}}=f_1c_{2_{2k+1}}c_{1_{k'}}$$
and the last expression can be computed from (a). This proves (e). Similarly, (f) 
follows from (b); (g) is a special case of 6.6. The proposition is proved.

\subhead 17.5\endsubhead
From 7.4,7.6 we see that the function $\D:W@>>>\bn$ has the following values:

$\D(2_{2k})=kL_1+kL_2$,

$\D(2_{2k+1})=-kL_1+(k+1)L_2$,

$\D(1_1)=L_1$,

$\D(1_{2k+1})=(k-1)L_1+kL_2$, if $k\ge 1$,

$\D(1_{2k})=kL_1+kL_2$.

It follows that P1 holds and that $\cd$ consists of $2_0=1_0,2_1,1_1,1_3$. Thus, P6
holds.

The formulas in 17.3, 17.4 determine $h_{x,y,z}$ for all $x,y,z$ except when $x=1$ 
or $y=1$, in which case $h_{1,y,z}=\d_{y,z}$, $h_{x,1,z}=\d_{x,z}$. From these 
formulas we see that the triples $(x,y,d)$ with $d\in\cd$, $\g_{x,y,d}\ne 0$ are:

$(2_{2k+1},2_{2k+1},2_1)$, $(1_{2k+2},2_{2k+2},1_3)$, $(1_1,1_1,1_1)$,

$(1,1,1)$, $(2_{2k+2},1_{2k+2},2_1)$,$(1_{2k+3},1_{2k+3},1_3)$, 
\nl
where $k\ge 0$. This implies that P2,P3 hold. From the results in 8.8 we see that
P4,P9,P13 hold. From 14.5 we see that P5 holds. From 14.7 we see that P7 holds. From
14.8 we see that P8 holds. From 14.10 we see that P10 holds. From 14.11 we see that
P11 holds. From 14.12 we see that P12 holds. From 14.14 we see that P14 holds. 

We now verify P15 in our case. With the notation in 14.15, it is enough to show 
that, if $a,b\in\{1,2\}$, $w\in W$, $s_aw>w,ws_b>w$, then  
$$(c_ae_w)c'_b-c_a(e_wc'_b)\in\tch_{\ge\aa(w)+1}.$$
Here $c_a=c_{s_a},c'_b=c'_{s_b}$. If $a$ or $b$ is $1$, then from 17.2 we have
$(c_ae_w)c'_b-c_a(e_wc'_b)=0$. Hence we may assume that $a=b=2$ and $w=1_{2k+1}$. 
Using 17.2 we have
$$\align&c_2(e_{1_{2k+1}}c'_2)=c_2(e_{1_{2k+2}}+\d_{k>0})\z'e_{1_{2k}}
+\d_{k>1}e_{1_{2k-2}}\\&=e_{2_{2k+3}}+\z e_{2_{2k+1}}+\d_{k>0}e_{2_{2k-1}}
+\d_{k>0}\z'e_{2_{2k+1}}+\d_{k>0}\z\z'e_{2_{2k-1}}\\&+\d_{k>1}\z'
e_{2_{2k-3}}+\d_{k>1}e_{2_{2k-1}}+\d_{k>1}\z e_{2_{2k-3}}+\d_{k>2}
e_{2_{2k-5}}\\&=e_{2_{2k+3}}+\z e_{2_{2k+1}}+\d_{k>0}\z'e_{2_{2k+1}}+
\d_{k>0}e_{2_{2k-1}}+\d_{k>1}e_{2_{2k-1}}\\&+\d_{k>0}\z\z'e_{2_{2k-1}}+
\d_{k>1}(\z+\z')e_{2_{2k-3}}+\d_{k>2}e_{2_{2k-5}}.\endalign$$
Similarly,
$$\align&(c_2e_{1_{2k+1}})c'_2=e_{2_{2k+3}}+\z' e_{2_{2k+1}}+\d_{k>0}\z 
e_{2_{2k+1}}+\d_{k>0})e_{2_{2k-1}}+\d_{k>1}e_{2_{2k-1}}\\&+\d_{k>0}\z\z'
e_{2_{2k-1}}+\d_{k>1}(\z+\z')e_{2_{2k-3}}+\d_{k>2}e_{2_{2k-5}}.\endalign$$
Hence 
$$c_2(e_{1_{2k+1}}c'_2)-(c_2e_{1_{2k+1}})c'_2=(\z-\z')(1-\d_{k>0})e_{2_{2k+1}}.$$
If $k>0$, the right hand side is zero. Thus we may assume that $k=0$. In this case,
$$c_2(e_{1_1}c'_2)-(c_2e_{1_1})c'_2=(\z-\z')e_{2_1}.$$
We have $\aa(1_1)=L_1<L_2=\aa(2_1)$. This completes the verification of P15 in our
case.

\head 18. The ring $J$\endhead
\subhead 18.1\endsubhead
{\it In this chapter we assume that $W,L$ is bounded and that P1-P15 in \S14 are valid.} 
In particular the results of this chapter are applicable if we are in the split
case (see \S15) or more generally in the quasisplit case (see \S16) with $W,L$ bounded.

\proclaim{Theorem 18.2} Assume that $W$ is tame.

(a) $W$ has only finitely many left cells. 

(b) $W$ has only finitely many right cells. 

(c) $W$ has only finitely many two-sided cells.

(d) $\cd$ is a finite set.
\endproclaim
We prove (a). Since $\aa(w)$ is bounded above it is enough to show that, for any
$a\in\bn$, $\aa\i(a)$ is a union of finitely many left cells. By P4, $\aa\i(a)$ is a
union of left cells. Let $\ch^1$ be the $\bz$-algebra $\bz\ot_\ca\ch$ where $\bz$ is
regarded as an $\ca$-algebra via $v\m 1$. We write $c_w$ instead of $1\ot c_w$. For
any $a'\ge 0$ let $\ch^1_{\ge a'}$ be the subgroup of $\ch^1$ spanned by 
$\{c_w;\aa(w)\ge a'\}$ (a two-sided ideal of $\ch^1$, by P4). We have a direct sum 
decomposition 
$$\ch^1_{\ge a}/\ch^1_{\ge a+1}=\op_\G E_\G\tag e$$
where $\G$ runs over the left cells contained in $\aa\i(a)$ and $E_\G$ is generated
as a group by the images of $c_w,w\in\G$; these images form a $\bz$-basis of $E_\G$.
Now $\ch^1_{\ge a}/\ch^1_{\ge a+1}$ inherits a left $\ch^1$-module structure from 
$\ch^1$ and (by P9) each $E_\G$ is a $\ch^1$-submodule. Since $W$ is tame, there 
exists a finitely generated abelian subgroup $W_1$ of finite index of $W$. Now 
$\ch^1=\bz[W]$ contains $\bz[W_1]$ as a subring. Since 
$\ch^1_{\ge a}/\ch^1_{\ge a+1}$ is a subquotient of $\ch^1$ (a finitely generated 
$\bz[W_1]$-module) and $\bz[W_1]$ is a noetherian ring, it follows that 
$\ch^1_{\ge a}/\ch^1_{\ge a+1}$ is a finitely generated $\bz[W_1]$-module. Hence in
the direct sum decomposition (e) with only non-zero summands, the number of summands
must be finite. This proves (a).

Since any right cell is of the form $\G\i$ where $\G$ is a left cell, we see that
(b) follows from (a). Since any two-sided cell is a union of left cells, we see that
(c) follows from (a). From P16 we see that (d) follows from (a). The theorem is 
proved.

\subhead 18.3\endsubhead
Let $J$ be the free abelian group with basis $(t_w)_{w\in W}$. We set
$$t_xt_y=\su_{z\in W}\g_{x,y,z\i}t_z.$$
The sum is finite since $\g_{x,y,z\i}\ne 0\imp h_{x,y,z}\ne 0$ and this implies that
$z$ runs through a finite set (for fixed $x,y$). We show that this defines an 
(associative) ring structure on $J$. We must check the identity 
$$\su_z\g_{x,y,z\i}\g_{z,u,u'{}\i}=\su_w\g_{y,u,w\i}\g_{x,w,u'{}\i}\tag a$$
for any $x,y,u,u'\in W$. From P8,P4 we see that both sides of (a) are $0$ unless 
$$\aa(x)=\aa(y)=\aa(u)=\aa(u')=a\tag b$$
for some $a\in\bn$. Hence we may assume that (b) holds. By P8,P4, in the first sum 
in (a) we may assume that $\aa(z)=a$ and in the second sum in (a) we may assume that
$\aa(w)=a$. The equation $(c_xc_y)c_u=c_x(c_yc_u)$ in $\ch$ implies
$$\su_zh_{x,y,z}h_{z,u,u'}=\su_wh_{y,u,w}h_{x,w,u'}.\tag c$$
If $h_{x,y,z}h_{z,u,u'}\ne 0$ then $u'\le_\car z\le_\car x$ hence, by P4,
$\aa(u')\ge\aa(z)\ge\aa(x)$ and $\aa(z)=a$. Hence in the first sum in (c) we may 
assume that $\aa(z)=a$. Similarly in the second sum in (c) we may assume that 
$\aa(w)=a$. Taking the coefficient of $v^{2\aa(z)}$ in both sides of (c) we find 
(a).

For any commutative ring $A$ with $1$ we set $J_A=A\ot J$; this is the free $A$-module with basis
$\{t_x;x\in W\}$. It is naturally an $A$-algebra.

If $\cd$ is finite, the algebra $J_A$ has a unit element $\su_{d\in\cd}n_dt_d$. Here $n_d=\pm 1$ is as in 
14.1(a), see P5. Let us check that 
$t_x\su_dn_dt_d=t_x$ for $x\in W$. This is equivalent to the identity 
$\su_dn_d\g_{x,d,z\i}=\d_{z,x}$. By P7 this is equivalent to 
$\su_dn_d\g_{z\i,x,d}=\d_{z,x}$. This follows from P2,P3,P5. The equality 
$(\su_dn_dt_d)t_x=t_x$ is checked in a similar way.

If $\cd$ is not necessarily finite, then $J_A$ has only a generalized unit element in the sense 
the elements $t_d (d\in\cd)$ of $J_A$ satisfy $t_dt_{d'}=\d_{d,d'}$ for $d,d'\in\cd$ and 
$\sum_{d,d'\in\cd}t_dJ_At_{d'}=J_A$.

For any subset $X$ of $W$, let $J^X_A$ be the $A$-submodule of $J_A$ generated by 
$\{t_x;x\in X\}$. (When $A=\ZZ$ we write $J^X$ instead of $J^X_\ZZ$.)
If $\boc$ is a two-sided cell of $W,L$ then, by P8, $J^\boc_A$ is a
subalgebra of $J_A$ and $J_A=\op_\boc J_A^\boc$ is a direct sum decomposition of $J_A$ as an algebra. 
If $\cd$ is finite then $J^\boc_A$ has a unit element $\su_{d\in\cd\cap\boc}n_dt_d$. 
Similarly, if $\G$
is a left cell of $W,L$ then $J^{\G\cap\G\i}_A$ is a subalgebra of $J_A$ with unit element
$n_dt_d$ where $d\in\cd\cap\G$. 

\proclaim{Proposition 18.4} Assume that we are in the setup of 15.1. Let $x,y\in W$.

(a) The condition $x\si_\cl y$ is equivalent to the condition that $t_xt_{y\i}\ne 0$
and to the condition that, for some $u$, $t_y$ appears with $\ne 0$ coefficient in 
$t_ut_x$.

(b) The condition $x\si_\car y$ is equivalent to the condition that
$t_{x\i}t_y\ne 0$ and to the condition that, for some $u$, $t_y$ appears with 
$\ne 0$ coefficient in $t_xt_u$.

(c) The condition $x\si_{\lr}y$ is equivalent to the condition that $t_xt_ut_y\ne 0$
for some $u$ and to the condition that, for some $u,u'$, $t_y$ appears with $\ne 0$
coefficient in $t_{u'}t_xt_u$.
\endproclaim
Let $J^+=\su_z\bn t_z$. By 15.1(a) we have $J^+J^+\sub J^+$.

We prove (a). The second condition is equivalent to $\g_{x,y\i,u}\ne 0$ for some 
$u$; the third condition is equivalent to $\g_{u,x,y\i}\ne 0$ for some $u$. These 
conditions are equivalent by P7. 

Assume that $\g_{x,y\i,u}\ne 0$ for some $u$. Using P8 we deduce that $x\si_\cl y$.

Assume now that $x\si_\cl y$. Let $d\in\cd$ be such that $x\si_\cl d$. Then we have
also $y\si_\cl d$. By P13 we have $\g_{x\i,x,d}\ne 0$, $\g_{y\i,y,d}\ne 0$. Hence 
$\g_{x\i,x,d}=1,\g_{y\i,y,d}=1$. Hence $t_{x\i}t_x\in t_d+J^+$, 
$t_{y\i}t_y\in t_d+J^+$. Since $t_dt_d=t_d$, it follows that 
$t_{x\i}t_xt_{y\i}t_y\in t_dt_d+J^+=t_d+J^+$. In particular, $t_xt_{y\i}\ne 0$. This
proves (a).

The proof of (b) is entirely similar.

We prove (c). Using the associativity of $J$ we see that the third condition on
$x,y$ is a transitive relation on $W$. Hence to prove that the first condition
implies the third condition we may assume that either $x\si_\cl y$ or $x\si_\car y$,
in which case this follows from (a) or (b). The fact that the third condition 
implies the first condition also follows from (a),(b). Thus the first and third 
condition are equivalent. 

Assume that $t_xt_ut_y\ne 0$ for some $u$. By (a),(b) we then have 
$x\si_\cl u\i,u\i\si_\car y$. Hence $x\si_{\lr}y$. 

Conversely, assume that $x\si_{\lr}y$. By P14 we have $x\si_{\lr}y\i$. By the
earlier part of the proof, $t_{y\i}$ appears with $\ne 0$ coefficient in 
$t_{u'}t_xt_u$ for some $u,u'$. We have $t_{u'}t_xt_u\in at_{y\i}+J^+$ where $a>0$.
Hence $t_{u'}t_xt_ut_y\in at_{y\i}t_y+J^+$. Since $t_{y\i}t_y$ has a coefficient $1$
and the other coefficients are $\ge 0$, it folows that $t_{u'}t_xt_ut_y\ne 0$. Thus,
$t_xt_ut_y\ne 0$. We see that the first and second conditions are equivalent. The 
proposition is proved.

\subhead 18.5\endsubhead
Assume now that we are in the setup of 7.1 with $m=\iy$ and $L_2>L_1$. From the 
formulas in 17.3,17.4 we can determine the multiplication table of $J$. We find

$t_{2_{2k+1}}t_{2_{2k'+1}}=\su_{u\in[0,\tk]}t_{2_{2k+2k'+1-4u}}$,

$t_{1_{2k+3}}t_{1_{2k'+3}}=\su_{u\in[0,\tk]}t_{1_{2k+2k'+3-4u}}$,

$t_{2_{2k+1}}t_{2_{2k'+2}}=\su_{u\in[0,\tk]}t_{2_{2k+2k'+2-4u}}$,

$t_{1_{2k+3}}t_{1_{2k'+2}}=\su_{u\in[0,\tk]}t_{1_{2k+2k'+2-4u}}$,

$t_{2_{2k+2}}t_{1_{2k'+3}}=\su_{u\in[0,\tk]}t_{2_{2k+2k'+2-4u}}$,

$t_{2_{2k+2}}t_{1_{2k'+2}}=\su_{u\in[0,\tk]}t_{2_{2k+2k'+1-4u}}$,

$t_{1_{2k+2}}t_{2_{2k'+1}}=\su_{u\in[0,\tk]}t_{1_{2k+2k'+2-4u}}$,

$t_{1_{2k+2}}t_{2_{2k'+2}}=\su_{u\in[0,\tk]}t_{1_{2k+2k'+3-4u}}$,

$t_{1_1}t_{1_1}=t_{1_1}$,

$t_1t_1=t_1$;
\nl
here $k,k'\ge 0$ and $\tk=\min(k,k')$. All other products are $0$. 

Let $R$ be the free abelian group with basis $(b_k)_{k\in\bn}$. We regard $R$ as a 
commutative ring with multiplication 
$$b_kb_{k'}=\su_{u\in[0,\min(k,k')]}b_{k+k'-2u}.$$
Let $J_0=\su_{w\in W-\{1,1_1\}}\bz t_w$. The formulas above show that 
$J=J_0\op\bz t_1\op\bz t_{1_1}$ (direct sum of rings) and that the ring $J_0$ is 
isomorphic to the ring of $2\T 2$ matrices with entries in $R$, via the isomorphism
defined by:
$$t_{2_{2k+1}}\m\left(\sm b_k&0\\0&0\esm\right),\qua
t_{1_{2k+3}}\m\left(\sm 0&0\\0&b_k\esm\right),\qua    
t_{2_{2k+2}}\m\left(\sm 0&b_k\\0&0\esm\right),\qua  
t_{1_{2k+2}}\m\left(\sm 0&0\\b_k&0\esm\right).$$
Note that $R$ is canonically isomorphic to the representation ring of $SL_2(\bc)$ 
with its canonical basis consisting of irreducible representations.

\subhead 18.6\endsubhead
Assume that we are in the setup of 7.1 with $m=\iy$ and $L_2=L_1$. By methods 
similar (but simpler) to those of \S17 and 18.5, we find
$$t_{2_{2k+1}}t_{2_{2k'+1}}=\su_{u\in[0,2\min(k,k')]}t_{2_{2k+2k'+1-2u}}.$$
Let $J^1$ be the subring of $J$ generated by $t_{2_{2k+1}},k\in\bn$. While, in 18.5,
the analogue of $J^1$ was isomorphic to $R$ as a ring with basis, in the present 
case, $J^1$ is canonically isomorphic to $R'$, the subgroup of $R$ generated by 
$b_k$ with $k$ even. (Note that $R'$ is a subring of $R$, naturally isomorphic to 
the representation ring of $PGL_2(\bc)$.) 

\subhead 18.7\endsubhead
In the setup of 7.1 with $m=4$ and $L_2=2,L_1=1$ (a special case of the situation in
\S16), we have
$$J=\bz t_1\op \bz t_{1_1}\op J_0\op\bz t_{2_3}\op\bz t_{2_4}$$
(direct sum of rings) where $J_0$ is the subgroup of $J$ generated by
$t_{2_1},t_{2_2},t_{1_2},t_{1_3}$. The ring $J_0$ is isomorphic to the ring of
$2\T 2$ matrices with entries in $\bz$, via the isomorphism defined by:
$$t_{2_1}\m\left(\sm 1&0\\0&0\esm\right),\qua
t_{1_3}\m\left(\sm 0&0\\0&1\esm\right),\qua t_{2_2}\m\left(\sm 0&1\\0&0\esm\right),
\qua t_{1_2}\m\left(\sm 0&0\\1&0\esm\right).$$
Moreover, $t_1,t_{1_1},t_{2_4}$ are idempotent. On the other hand,
$$t_{2_3}t_{2_3}=-t_{2_3}.$$
Notice the minus sign! (It is a special case of the computation in 7.8.)

\subhead 18.8 \endsubhead
{\it Until the end of 18.12 we assume that $\cd$ is finite.} 
In the following result, $n_d=\pm1$ (for $d\in\cd$) is as in 14.1(a) (see P.5). 

\proclaim{Theorem 18.9}  The $\ca$-linear maps $\ph:\ch@>>>J_\ca$, $\ph':\ch@>>>J_\ca$ given by
$$\ph(c_x^\da)=\sum_{z\in W,d\in\cd;\aa(d)=\aa(z)}h_{x,d,z}n_dt_z\qua (x\in W),$$
$$\ph'(c_x^\da)=\sum_{z\in W,d\in\cd;\aa(d)=\aa(z)}h_{d,x,z}n_dt_z\qua (x\in W),$$
are homomorphisms of $\ca$-algebras with $1$.
\endproclaim
Note that $\ph'$ is the composition of the algebra isomorphism 
$\ch@>\si>>\ch^{opp}$ given by $c_w\m c_{w\i}$ (see 3.4) with $\ph:\ch^{opp}@>>>(J_\ca)^{opp}$ and with
$(J_\ca)^{opp}@>\si>>J_\ca$ given by $t_w\m t_{w\i}$. Hence it is enough to prove the statement of the
theorem concerning $\ph$. Consider the equality

(a) $\su_wh_{x_1,x_2,w}h'_{w,x_3,y}=\su_wh_{x_1,w,y}h'_{x_2,x_3,w}$
\nl
(see P15) with $\aa(x_2)=\aa(y)=a$. In the left hand side we may assume that 
$y\le_\car w\le_\cl x_2$ hence (by P4) $\aa(y)\ge\aa(w)\ge\aa(x_2)$, hence 
$\aa(w)=a$. Similarly in the right hand side we may assume that $\aa(w)=a$. Picking
the coefficient of $v'{}^a$ in both sides of (a) gives

(b) $\su_wh_{x_1,x_2,w}\g_{w,x_3,y\i}=\su_wh_{x_1,w,y}\g_{x_2,x_3,w\i}$.
\nl
Let $x,x'\in W$. The desired identity 
$\ph(c_x^\da c_{x'}^\da)=\ph(c_x^\da)\ph(c_{x'}^\da)$ is equivalent to 
$$\su\Sb w\in W,d\in\cd\\\aa(d)=a'\eSb h_{x,x',w}h_{w,d,u}n_d=
\su\Sb z,z'\in W,d,d'\in\cd\\\aa(d)=\aa(z)\\\aa(d')=\aa(z')\eSb h_{x,d,z}
h_{x',d',z'}\g_{z,z',u\i}n_dn_{d'}$$
for any $u\in W$ such that $\aa(u)=a'$. In the right hand we may assume that 

$\aa(d)=\aa(z)=\aa(d')=\aa(z')=a'$ 
\nl
(by P8,P4). Hence the right hand side can be rewritten (using (b)):
$$\align&\su\Sb z'\in W,d,d'\in\cd\\\aa(d)=\aa(d')=\aa(z')=a'\eSb h_{x',d',z'}
\su_{z; \aa(z)=a'}h_{x,d,z}\g_{z,z',u\i}n_dn_{d'}\\&
=\su\Sb z'\in W,d,d'\in\cd\\ \aa(d)=\aa(d')=\aa(z')=a'\eSb h_{x',d',z'}
\su_{w;\aa(w)=a'}h_{x,w,u}\g_{d,z',w\i}n_dn_{d'}.\endalign$$
By P2,P3,P5, this equals
$$\su_{z'\in W,d'\in\cd;\aa(d')=\aa(z')=a'}h_{x',d',z'}h_{x,z',u}n_{d'}$$
which by the identity $(c_xc_{x'})c_{d'}=c_x(c_{x'}c_{d'})$ equals
$$\su_{w\in W,d'\in\cd;\aa(d')=a'}h_{x,x',w}h_{w,d',u}n_{d'}.$$
Thus $\ph$ is compatible with multiplication.

Next we show that $\ph$ is compatible with the unit elements of the two algebras. 
An equivalent statement is that for any $z\in W$ such that $\aa(z)=a$, the sum 
$\su_{d\in\cd;\aa(d)=a}h_{1,d,z}n_d$ equals $n_z$ if $z\in\cd$ and is $0$ if 
$z\n\cd$. This is clear since $h_{1,d,z}=\d_{z,d}$.

\subhead 18.10\endsubhead
If we identify the $\ca$-modules $\ch$ and $J_\ca$ via $c_w^\da\m t_w$, the
obvious left $J_\ca$-module structure on $J_\ca$ becomes the left $J_\ca$-module 
structure on $\ch$ given by
$$t_x*c_w^\da=\su_{z\in W}\g_{x,w,z\i}c_z^\da$$
Let $\ch_a=\op_{w;\aa(w)=a}\ca c_w^\da,\ch_{\ge a}=\op_{w;\aa(w)\ge a}\ca c_w^\da$. 
We have $t_x*c_w^\da\in\ch_{\aa(w)}$ for all $x,w$. We show that for any $h\in\ch,w\in W$ we have
$$hc_w^\da=\ph(h)*c_w^\da\mod\ch_{\ge\aa(w)+1}.\tag a$$
Indeed, we may assume that $h=c_x^\da$. Using 18.9(b), we have
$$\align&\ph(c_x^\da)*c_w^\da=\su\Sb d\in\cd,z\\ \aa(d)=\aa(z)\eSb 
h_{x,d,z}n_dt_z*c_w^\da\\&=\su\Sb d\in\cd,z,u\\\aa(d)=\aa(z)\eSb 
h_{x,d,z}\g_{z,w,u\i}\hn_dc_u^\da
\\&=\su\Sb d\in\cd,t,u\\ \aa(d)=\aa(w)=\aa(u)\eSb 
h_{x,t,u}\g_{d,w,t\i}n_dc_u^\da\\&=\su\Sb u\\ \aa(w)=\aa(u)\eSb 
h_{x,w,u}c_u^\da=c_x^\da c_w^\da\mod\ch_{\ge\aa(w)+1},\endalign$$
as required.

\subhead 18.11\endsubhead
Let $\ca@>>>R$ be a ring homomorphism of $\ca$ into a commutative ring $R$ with
$1$. Let $\ch_R=R\ot_\ca\ch$,
$\ch_{R,\ge a}=R\ot_\ca\ch_{\ge a}$. Then $\ph$ extends to a homomorphism of 
$R$-algebras $\ph_R:\ch_R@>>>J_R$. The $J_\ca$-module in 18.10 extends to a 
$J_A$-module structure on $\ch_R$ denoted again by $*$. From 18.10(a) we deduce

(a) $hc_w^\da=\ph_R(h)*c_w^\da\mod\ch_{R,\ge\aa(w)+1}$ for any $h\in\ch_R,w\in W$.

\proclaim{Proposition 18.12} (a) If $N$ is a bound for $W,L$, then 
$(\ker\ph_R)^{N+1}=0$.

(b) If $R=R_0[v,v\i]$ where $R_0$ is a commutative ring with $1$, $v$ is an
indeterminate and $\ca@>>>R$ is the obvious ring homomorphism, then $\ker\ph_R=0$.
\endproclaim
We prove (a). If $h\in\ker\ph_R$ then by 18.11(a), we have
$h\ch_{R,\ge a}\sub\ch_{R,\ge a+1}$ for any $a\ge 0$. Applying this repeatedly, we
see that, if $h_1,h_2,\do,h_{N+1}\in\ch$, we have
$h_1h_2\do h_{N+1}\in\ch_{R,\ge N+1}=0$. This proves (a).

We prove (b). Let $h=\su_xp_xc_x^\da\in\ker\ph_R$ where $p_x\in R$. Assume that 
$h\ne 0$. Then $p_x\ne 0$ for some $x$. We can find $a\ge 0$ such that 
$p_x\ne 0\imp\aa(x)\ge a$ and $X=\{x\in W;p_x\ne 0,\aa(x)=a\}$ is non-empty. We can
find $b\in\bz$ such that $p_x\in v^b\bz[v\i]$ for all $x\in X$ and such that 
$X'=\{x\in X;\p_b(p_x)\ne 0\}$ is non-empty. Let $x_0\in X'$. We can find $d\in\cd$
such that $\g_{x_0,d,x_0\i}=\g_{x_0\i,x_0,d}\ne 0$. We have 
$hc_d^\da=\su_xp_xc_x^\da c_d^\da$. If $\aa(x)>a$, then 
$c_x^\da c_d^\da\in\ch_{R,\ge a+1}$. Hence
$hc_d^\da=\su_{x\in X}p_xc_x^\da c_d^\da\mod\ch_{R,\ge a+1}$. Since $\ph_R(h)=0$,
from 18.11(a) we have $hc_d^\da\in\mod\ch_{R,\ge a+1}$. It follows that 
$\su_{x\in X}p_xc_x^\da c_d^\da\in\ch_{R,\ge a+1}$. In particular the coefficient of
$c_{x_0}^\da$ in $\su_{x\in X}p_xc_x^\da c_d^\da$ is $0$. In other words, 
$\su_{x\in X}p_xh_{x,d,x_0}=0$. The coefficient of $v^{a+b}$ in the last sum is
$$\su_{x\in X}\p_b(p_x)\g_{x,d,x_0\i}=\p_b(p_{x_0})\g_{x_0,d,x_0\i}$$
and this is on the one hand $0$ and on the other hand is non-zero since
$\p_b(p_{x_0})\ne 0$ and $\g_{x_0,d,x_0\i}\ne 0$, by the choice of $x_0,d$. This 
contradiction completes the proof.

\subhead 18.13\endsubhead
We fix a commutative ring $A$ with $1$.
We will show that, without assuming that $\cd$ is finite, $J_A$ can be imbedded naturally in a larger 
$A$-algebra which has a unit element.

Let $\tJ_A$ be the set of formal sums $\sum_{w\in W}f(w)t_w$ where $f:W@>>>A$ is any function. We regard 
$\tJ_A$ as an $A$-module in an obvious way. For a function $f:W@>>>A$, the support of $f$ is 
$\supp(f)=\{w\in W;f(w)\ne0\}$. Note that $J_A$ may be identified with the 
$A$-submodule of $\tJ_A$ consisting of all $\sum_{w\in W}f(w)t_w$ such that $f:W@>>>A$ has finite support.

We say that $f:W@>>>A$ is {\it left (resp. right) admissible} if $\supp(f)$ has finite intersection with 
any left (resp. right) cell in $W$. 
If $\ti f:W@>>>A$ is given by $\ti f(z)=f(z\i)$ then clearly $f$ is left admissible if and only if $\ti f$
is right admissible.

For two functions $f,f':W@>>>A$ we try to define $f'':W@>>>A$ by 
$f''(z)=\sum_{x,y\in W}f(x)f'(y)\g_{x,y,z\i}$; the sum may be infinite in general hence may not make sense. 
We show:

(a) {\it If both $f,f'$ are left admissible then $f''$ is well defined and left admissible. If both $f,f'$ 
are right admissible then $f''$ is well defined and right admissible.}
\nl
Assume first that $f,f'$ are left admissible. To prove the first sentence of (a) it is enough to show that 
for any $d\in\cd$, the set
$$\{(x,y,z)\in W^3;f(x)\ne0,f'(y)\ne0,z\si_\cl d,\g_{x,y,z\i}\ne0\}$$ 
is finite. Using P8 we see that this set is contained in
$$\align&U=\\&\{(x,y,z)\in W^3;f(x)\ne0,f'(y)\ne0,z\si_\cl d,x\si_\cl y\i, y\si_\cl z, l(z)\le l(x)+l(y)\}.\endalign$$
(Note that if $\g_{x,y,z\i}\ne0$ then $h_{x,y,z}\ne0$ hence $l(z)\le l(x)+l(y)$, see 13.1.) Hence it is enough 
to show that $U$ is finite. Let 

$F=\{y\in W;f'(y)\ne0,y\si_\cl d\}$,

$F'=\{x\in W;x\si_\cl d',f(x)\ne0, x\si_\cl h\i \text{ for some }h\in F\}$,

$F''=\{z\in W;l(z)\le l(x)+l(y)\text{ for some }x\in F',y\in F\}$.
\nl
Now $F$ is finite since $f'$ is left admissible. Hence $F'$ is finite. It follows that $F''$ is finite (we
use that $f$ is left admissible). Then $F''$ must be also finite since there are only finitely many elements
of fixed length in $W$. If $(x,y,z)\in U$ then $y\in F$ hence $x\in F'$. We have clearly $z\in F''$. Thus, 
$U\sub F''\T F\T F'''$ so that $U$ is finite. This proves the first sentence in (a).
 
Next we assume that $f,f'$ are right admissible. To prove the second sentence of (a) it is enough to show 
that for any $d\in\cd$, the set
$$\align&\{(x,y,z)\in W^3;f(x)\ne0,f'(y)\ne0,z\si_{\car}d,\g_{x,y,z\i}\ne0\}\\&
=\{(x,y,z)\in W^3;\ti f(x\i)\ne0,\ti f'(y\i)\ne0,z\i\si_\cl d,\g_{x,y,z\i}\ne0\}\endalign$$
is finite. Since $\tf,\tf'$ are left admissible, it is enough to show, by the first part of the proof
applied to $\ti f',\ti f$ instead of $f,f'$, that $\g_{x,y,z\i}\ne0$ implies $\g_{y\i,x\i,z}\ne0$. But the 
identity $h_{x,y,z}=h_{y\i,x\i,z\i}$ implies $\g_{x,y,z\i}=\g_{y\i,x\i,z}$; this completes the proof of (a).

\mpb

Let $\che{}J_A$ (resp. $J\che{}_A$) be the set of formal sums $\sum_{w\in W}f(w)t_w\in\tJ_A$ where 
$f:W@>>>A$ is a left (resp. right) admissible function. Note that $\che{}J_A$ and $J\che{}_A$ are 
$A$-submodules of $\tA_J$. We define $A$-algebra structures on $\che{}J_A$ and on $J\che{}_A$ by
$$(\sum_{x\in W}f(x)t_x)(\sum_{y\in W}f'(y)t_y)=\sum_{z\in W}f''(z)t_z$$
with $f,f',f''$ as in (a). We show:

(b) {\it These algebra structures are associative.}
\nl
We must show that if either each of $f_1,f_2,f_3:W@>>>A$ is left admissible or each of
$f_1,f_2,f_3:W@>>>A$ is right admissible, then for any $u_1\in W$ we have
$$\align&\sum_{x,y,z,u\in W}f_1(x)f_2(y)f_3(z)\g_{x,y,u\i}\g_{u,z,u_1\i}\\&=
\sum_{x,y,z,u'\in W}f_1(x)f_2(y)f_3(z)\g_{y,z,u'{}\i}\g_{x,u',u_1\i}.\endalign$$
(Both sums are finite by repeated application of (a).) It is enough to show that for any $x,y,z,u_1$ in $W$ 
we have
$$\sum_u\g_{x,y,u\i}\g_{u,z,u_1\i}=\sum_{u'}\g_{y,z,u'{}\i}\g_{x,u',u_1\i}$$
This follows from the associativity of $J_A$.

Note that both algebras $\che{}J_A$, $J\che{}_A$ have a unit element, namely 
$1=\sum_{d\in\cd}n_dt_d$; this is checked, using P2,P3,P5,P7, in the same way as in 18.3.

Let $\che J_A=\che{}J_A\cap J\che{}_A$. Thus $\che J_A$ consists of all formal sums 
$\sum_{w\in W}f(w)t_w\in\tJ_A$ such that $f:W@>>>A$ is both left admissible and right admissible. Note 
that $\che J_A$ is an $A$-subalgebra of both $\che{}J_A$ and of $J\che{}_A$ with unit element
$1=\sum_{d\in\cd\cap c}n_dt_d$. Note also that $J_A$ is a subalgebra of $\che J_A$. 
Now the map $\sum_{x\in W}f(x)t_x\m\sum_{x\in W}f(x\i)t_x$ defines algebra isomorphisms
$\che{}J_A@>\si>>(J\che{}_A)^{opp}$,
$\che J_A@>\si>>(\che J_A)^{opp}$, $J_A@>\si>>(J_A)^{opp}$ where the upperscript ${}^{opp}$ 
denotes the opposed algebra.

If $\cd$ is finite (which is the case if $W$ is tame, see 18.2), we have $\che{}J_A=J\che{}_A=\che J_A=J_A$. 

\subhead 18.14\endsubhead
We show:

(a) {\it Let $w\in W$. The function $W@>>>\ca$, $z\m \sum_{d\in\cd;\aa(d)=\aa(z)}h_{w,d,z}n_d$ is left admissible. 
The function $W@>>>\ca$, $z\m \sum_{d\in\cd;\aa(d)=\aa(z)}h_{d,w,z}n_d$ is right admissible.}
\nl
To prove the first assertion of (a) it is enough to show that for any $d\in\cd$, the set
$\{z\in W;z\si_\cl d,h_{w,d,z}\ne0\}$ is finite. This set is
contained in $\{z\in W;l(z)\le l(w)+l(d)\}$ (see 13.1) which is clearly finite. This proves 
the first assertion of (a). The second assertion of (a) is proved in a similar way.
\mpb

From (a) we see that the $\ca$-linear maps
$\ph:\ch@>>>\che{}J_\ca$, $\ph':\ch@>>>J\che{}_\ca$ given by
$$\ph(c_x^\da)=\sum_{d\in\cd,z\in W;\aa(z)=\aa(d)}h_{x,d,z}n_dt_z,\tag b$$
$$\ph'(c_x^\da)=\sum_{z\in c,d\in\cd;\aa(d)=\aa(z)}h_{d,x,z}n_dt_z\tag c$$
are well defined. The proof of the following result is essentially the same as that of Theorem 18.9.

(d) {\it $\ph:\ch@>>>\che{}J_\ca$ and $\ph':\ch@>>>J\che{}_\ca$ are homomorphisms of $\ca$-algebras 
with $1$.}
\nl
Now let $R,\ca@>>>R$ be as in 18.2(b), let $\ch_R$ be as in 18.11; we define an $R$-algebra homomorphism 
$\ph_R:\ch_R@>>>\che{}J_R$ by the same formula as $\ph:\ch@>>>\che{}J_R$ in which $h_{x,d,z}$ is viewed as 
an element of $R$. A proof similar to that of 18.12(b) shows that

(e) $\ker\ph_R=0$.

\subhead 18.15\endsubhead
In the remainder of this chapter we assume that $L=l$ and we fix a two-sided cell $\boc$ of $W$.
We shall study a categorical version $C^\boc$ of the ring $J^\boc$.
To do this, we shall use the theory of Soergel modules as in 16.2. We shall use the notation of 16.2 with $\tW=W$.
Thus, $R,R^{>0},\car,\tC,C, B_x (x\in W)$ are defined as in 16.2. For $M,M'\in\car$, $MM'\in\car$, $M^{M'}\in\car$ 
are defined as in 16.2. If $L\in\tC$ and $j\in\ZZ$ we write $L^j\in C$ for what in \cite{\EW, 6.2} is denoted 
by $\ch^j(L)$. (The fact that $L^j$ is well defined follows from results  of \cite{\SO} and \cite{\EW}.)
For any subset $X$ of $\boc$, let $C^X$ be the full subcategory of $C$ whose objects are isomorphic to 
finite direct sums of objects of the form $B_x (x\in X)$. For any $L\in C$ there is a 
unique direct sum decomposition $L=\un L\op L'$ where $\un L\in C^\boc$ and $L'$ is a direct sum of objects 
of the form $B_x (x\n\boc)$.
(The uniquenes of this direct sum decomposition follows from the results of \cite{\SO} and \cite{\EW}.)
For $M\in C^\boc$ we have $M=\op_{z\in\boc}E^M_z\ot B_z$ where $E^M_z$ are well defined finite
dimensional $\CC$-vector spaces which are $0$ for all but finitely many $z$.

Let $a$ be the value of the $\aa$-function on $\boc$. By arguments similar to those in \cite{\TENS} and 
making use of the results in \cite{\EW} we see that for $L,L'\in C^\boc$ we have $\un{(LL')^j}=0$ if $j>a$
and $L,L'\m L\unb L':=\un{(LL')^a}$ defines a monoidal structure on $C^\boc$. (For three objects
$L,L',L''$ of $C^\boc$ we have $(L\unb L')\unb L''=L\unb(L'\unb L'')=\un{(LL'L'')^{2a}}$.)

For $x,y\in\in\boc$ we have 
$$B_x\unb B_y=\op_{z\in\boc}V_{x,y,z\i}\ot B_z$$
where $V_{x,y,z\i}=E^{B_x\unb B_y}_z$ are
canonically defined $\CC$-vector spaces which are $0$ for all but finitely many $z$. Note that
$$\dim V_{x,y,z\i}=\g_{x,y,z\i}.\tag a$$

For $x\in\boc$ let $d_x$ be the unique element of $\cd\cap\boc$ such that $x\si_\cl d_x$.

If $x\in\boc$ and $d=d_{x\i}$, we have canonically
$$(\che V_{d,d,d}\ot B_d)\unb B_x=B_x.$$
(We shall denote the dual space of a $\CC$-vector space $V$ by $\che V$.)        
Indeed, by (a), it is enough to show that 
$\che V_{d,d,d}\ot V_{d,x,x\i}=\CC$. From $(B_d\unb B_d)\unb B_x=B_d\unb(B_d\unb B_x)$ we deduce using (a)
that $V_{d,d,d}\ot V_{d,x,x\i}=V_{d,x,x\i}\ot V_{d,x,x\i}$ where $V_{d,x,x\i},V_{d,d,d}$ are $1$-dimensional
hence $V_{d,d,d}=V_{d,x,x\i}$ and the desired equality follows. Note also that if 
$d'\in\cd\cap\boc, d'\ne d_{x\i}$, then $(\che V_{d',d',d'}\ot B_{d'})\unb B_x=0$.
Similarly, if $x\in\boc$ and $d=d_x$, we have canonically
$$B_x\unb(\che V_{d,d,d}\ot B_d)=B_x.$$
(We use the equality $\che V_{d,d,d}\ot V_{x,d,x\i}=\CC$ which follows from
$B_x\unb(B_d\unb B_d)=(B_x\unb B_d)\unb B_d$.) Moreover,
if $d'\in\cd\cap\boc$ and $d'\ne d_x$, then $B_x\unb(\che V_{d',d',d'}\ot B_{d'})=0$.
Thus, $\op_{d\in\cd\boc}(\che V_{d,d,d}\ot B_d)$ plays the role of a unit object for the monoidal category
$C^\boc$, although it does not belong to $C^\boc$ (unless $\cd\cap\boc$ is finite).

From (a) we see that if $\G,\G',\G''$ are left cells contained in $\boc$ and
$L\in C^{\G'{}\i\cap\G}$, $L'\in C^{\G\i\cap\G''}$, then $L\unb L'\in C^{\G'{}\i\cap\G''}$.
In particular $L,L'\m L\unb L'$ defines a monoidal structure on $C^{\G\i\cap\G}$. This monoidal structure
admits a unit object, namely $\che V_{d,d,d}\ot B_d$, where $d\in\cd\cap\G$.

\subhead 18.16\endsubhead
We show that for $x,y,z\in\boc$ we have canonically
$$V_{x,y,z}=V_{y\i,x\i,z\i}.\tag a$$
For any $M\in\car$ let $M^\sh$ be the object of $\car$ which is equal to $M$ as a graded $\CC$-vector space, 
but left (resp. right) multiplication by $r\in R$ on $M^\sh$ equals right (resp. left) multiplication by 
$r$ on $M$. In \cite{\LV, 3.1} it is shown that $M\in C$ implies $M^\sh\in C$ and $M\in\tC$ implies 
$M^\sh\in\tC$; it follows that for $M\in\tC$ and $j\in\ZZ$ we have canonically $(M^\sh)^j=(M^j)^\sh$.
More precisely in {\it loc.cit.} it is shown that, if $x\in W$, then we have $B_x^\sh\cong B_{x\i}$.
Let $\ph:B_x^\sh@>>>B_{x\i}$ be an isomorphism. It is well defined up to multiplication by a number in 
$\CC^*$. We show that there is a canonical choice for it.
For any $M\in\car$ we have an obvious isomorphism $(R_x)^M@>>>(R_x^\sh)^{M^\sh}$ (identity map)
of $\CC$-vector spaces.
Despite the fact that this is not necessarily an isomorphism in $\car$, it induces for any $i$
an isomorphism $\un{(R_x)^M}_i@>>>\un{(R_x^\sh)^{M^\sh}}_i$ of $\CC$-vector spaces.
(We use that $R^{>0}(R_x^M)=(R_x^M)R^{>0}$, see \cite{\LV, 3.2}.) Taking $M=B_x$, $i=l(x)$, we may thus 
identify $\un{(R_x)^{B_x}}_{l(x)}@>>>\un{(R_x^\sh)^{B_x^\sh}}_{l(x)}$ as $\CC$-vector spaces.
We now identify $R_x^\sh=R_{x\i}$ as in  \cite{\LV, 3.2} and we identify $B_x^\sh$ with $B_{x\i}$ via $\ph$.
We obtain an identification 
of $\un{(R_x)^{B_x}}_{l(x)}$ with $\un{(R_{x\i})^{B_{x\i}}}_{l(x)}$ that is of $\CC$ with $\CC$. This is 
multiplication by some $\l\in\CC-\{0\}$. By replacing $\ph$ be a nonzero scalar multiple we can
achieve that $\l=1$. This gives a canonical choice for $\ph$; we denote it by $\ph_x$.
We can now identify $B_x^\sh=B_{x\i}$ via $\ph_x$.

For $M,M'\in\tC$ we have canonically $(MM')^\sh=(M')^sh M^\sh$. Hence for $x,y\in\boc$ we have canonically
$$\align&\op_{z\in\boc}V_{y\i,x\i,z\i}B_z=\un{(B_{y\i}B_{x\i})^a}\\&=
\un{(B_y^\sh B_x^\sh)^a}=\un{((B_xB_y)^\sh)^a}=\un{((B_xB_y)^a)^\sh}\\&
=\op_{z\in\boc}(V_{x,y,z}B_{z\i})\sh=\op_{z\in\boc}V_{x,y,z}B_z.\endalign$$
Now (a) follows.

\subhead 18.17\endsubhead
In the remainder of this chapter we assume that $W,S$ is an affine Weyl group, see 1.15 and that
$I$ is a subset of $S$ such that the group $W_I$ generated by $I$ is finite of maximum possible order.
We have $W=W_I\ct$, $W_I\cap\ct=\{1\}$, where $\ct$ is the normal subgroup of $W$ defined in 
1.16. Let $w_0^I$ be the longest element of $W_I$. 
We assume that $\boc$ is the two-sided cell of $W$ containing $w_0^I$.
There is a unique automorphism $w\m w^*$ of $W$ such that $x^*=w_0^Ixw_0^I$ for all $x\in W_I$ and
$y^*=w_0^Iy\i w_0^I$ for all $y\in\ct$. This automorphism maps $S$ onto $S$, $I$ onto $I$ and $\boc$ onto
$\boc$. We shall assume, as we may, that (with notation in \cite{\LV, 2.1}), there is an involutive
automorphism $e\m e^*$ of the dual space $\che\fh$ of the reflection representation $\fh$ of $W$ such that 
$(we)^*=w^*e^*$ for all $w\in W,e\in\che\fh$ and $(\a_s)^*=\a_{s^*}$ for all $s\in S$.

From the definitions we see that for any $x,y,z\in\boc$ we have canonically
$$V_{x^*,y^*,z^*}=V_{x,y,z}.\tag a$$
Let $\G=\{w\in W; w\text{ has maximal length in }wW_I$. According to
\cite{\LC, 8.5}, $\G$ is a left cell of $W$. It is clearly contained in $\boc$. 

The set $\G\i\cap\G$ is the set of all $w\in W$ such that $w$ has maximal length in $W_IwW_I$. Hence
each $W_I,W_I$ double coset in $W$ contains a unique element of $\G\i\cap\G$, see 9.15(e).

By \cite{\BAR, 8.2}, if $w\in W$ has maximal length in its $W_I,W_I$ double coset then $w^*=w\i$.
In particular, we have 
$$w^*=w\i\text{ for all }w\in\G\i\cap\G.\tag b$$

\subhead 18.18\endsubhead
We show that for $x,y,z$ in $\G\i\cap\G$ we have canonically 
$$V_{x,y,z}=V_{y,x,z}.\tag a$$
The method of proof has some common features with one in \cite{\LX, 3.5}

Using 18.17(a), 18.17(b) and 18.16(a) we see that we have canonically
$$V_{x,y,z}=V_{x^*,y^*,z^*}=V_{x\i,y\i,z\i}=V_{y,x,z}.$$
This proves (a). From (a) we deduce that for $x,y$ in $\G\i\cap\G$ we have canonically 
$$B_x\unb B_y=B_y\unb B_x\tag b$$ 
in $C^{\G\i\cap\G}$. 
Now let $M,M'$ be two objects of $C^{\G\i\cap\G}$. We show that we have canonically
$$M\unb M'=M'\unb M\tag c$$
in $C^{\G\i\cap\G}$. We have canonically
$$M=\op_{x\in\G\i\cap\G}E_{M,x}\ot B_x,\qua M'=\op_{y\in\G\i\cap\G}E_{M',y}\ot B_y.$$
Hence, using (b), we have
$$\align&M\unb M'=\op_{x,y\in\G\i\cap\G}(E_{M,x}\ot E_{M',y})\ot B_x\unb B_y\\&=
\op_{x,y\in\G\i\cap\G}(E_{M',y}\ot E_{M,x})\ot B_y\unb B_x=M'\unb M.\endalign$$
We see that the monoidal category $C^{\G\i\cap\G}$ has a natural commutativity constraint.

\subhead 18.19\endsubhead
Let $d$ be the unique element of $\cd\cap\G$.
We define a contravariant functor $D:C^{\G\i\cap\G}@>>>C^{\G\i\cap\G}$ by
$$M\m DM=\op_{z\in\G\i\cap\G}\che E^M_z\ot\che V_{z,z\i,d}\ot\che V_{d,d,d}\ot B_{z\i}\in C^{\G\i\cap\G}.$$
For $M\in C^{\G\i\cap\G}$ we have  
$$DDM=\op_{z\in\G\i\cap\G}E^M_z\ot V_{z,z\i,d}\ot V_{d,d,d}\ot\che V_{z\i,z,d}\ot\che V_{d,d,d}\ot B_z=M.$$
Here we have used that $V_{z,z\i,d}\ot V_{d,d,d}\ot\che V_{z\i,z,d}\ot\che V_{d,d,d}=\CC$
since, by 18.15(a), $V_{d,d,d}$ and $V_{z,z\i,d}=V_{z\i,z,d}$ are $1$-dimensional (note that by 18.17(a) and
18.17(b) we have $V_{z,z\i,d}=V_{z,z^*,d}=V_{z^*,z,d^*}=V_{z\i,z,d}$).

For $x\in\G\i\cap\G$ we have
$$\align&B_x\unb D(B_x)=\che V_{x,x\i,d}\ot\che V_{d,d,d}\ot B_x\unb B_{x\i}\\&
=\op_{z\in\G\i\cap\G}\che V_{x,x\i,d}\ot\che V_{d,d,d}\ot V_{x,x\i,z\i}\ot B_z\\&=
\che V_{x,x\i,d}\ot\che V_{d,d,d}\ot V_{x,x\i,d}\ot B_d\op M_1=M_0\op M_1\endalign$$
where 
$$M_0=\che V_{d,d,d}\ot B_d,$$
$$M_1=\op_{z\in\G\i\cap\G;z\ne d}\che V_{x,x\i,d}\ot\che V_{d,d,d}\ot V_{x,x\i,z\i}B_z.$$
(We have again used that $V_{x,x\i,d}$ is $1$-dimensional.)
Thus we have obvious morphisms $M_0@>j_x>>B_x\unb D(B_x)@>j'_x>>M_0$ where $M_0$ is the unit object
of $C^{\G\i\cap\G}$. Now let $M\in C^{\G\i\cap\G}$. We have 
$$M\unb DM=\op_{x,x'\in\G\i\cap\G}(E^M_x\ot\che E^M_{x'})\ot B_x\unb DB_{x'}$$
Let $j:M_0@>>>M\unb DM$ be the morphism whose $x,x'$ component is $0$ if $x\ne x'$ and is
$a_x\ot j_x$ when $x=x'$ (here $a_x:\CC@>>>E^M_x\ot\che E^M_x$ is the obvious imbedding).
Let $j':M\unb DM@>>>M_0$ be the morphism whose $x,x'$ component is $0$ if $x\ne x'$ and is
$a'_x\ot j'_x$ when $x=x'$ (here $a'_x:E^M_x\ot\che E^M_x@>>>\CC$ is the obvious projection).
The morphisms $M_0@>j>>M\unb DM@>j'>>M_0$ provide a rigid structure for the monoidal category
$C^{\G\i\cap\G}$.

\subhead 18.20\endsubhead
Let $G$ be the simple adjoint group over $\CC$ of type dual to that of $W,S$. Let $Rep G$ be the tensor
category of finite dimensional rational representations of $G$. Note that the simple objects of $Rep G$
are naturally indexed by the elements of $W$ which have maximal length
in their $W_I,W_I$ double coset, hence they are indexed by $\G\i\cap\G$. For $x\in\G\i\cap\G$ let $\cv_x$ be
the corresponding simple object of $Rep G$. Now \cite{\SING, Cor.8.7} can be interpreted as follows:

{\it For any $x,y\in\G\i\cap\G$ we have
$\cv_x\ot\cv_{x'}\cong\op_{z\in\G\i\cap\G}\cv_z^{\op\g_{x,y,z\i}}.$}
\nl
Using 1.1(a) this can be restated as follows:

(a) {\it For any $x,y\in\G\i\cap\G$ we have
$\cv_x\ot\cv_{x'}\cong\op_{z\in\G\i\cap\G} V_{x,y,z\i}\ot\cv_z.$}
\nl
Using (a) we see that the rigid symmetric monoidal category $C^{\G\i\cap\G}$ satisfies the assumptions in 
Deligne's theorem \cite{\DE, 0.6}: the finite $\ot$-generation and the property \cite{\DE, 0.5(i)(b)}
follow from the analogous statements for the category $Rep G$ where they are obvious. We deduce that 
$C^{\G\i\cap\G}$ is equivalent as a tensor category to the category of representations of some supergroup.
One can show that this supergroup is in fact a group isomorphic to $G$ (we omit the proof).

\subhead 18.21\endsubhead
In this subsection we make no assumption on $\boc$.
Let $\G$ be the set of all $x\in\boc$ such that $x$ has minimal length in $xW_I$. According to \cite{\LX}, 
$\G$ is exactly one left cell of $W$. By methods similar to those in 18.18, 18.19 we see that the monoidal 
category $C^{\G\i\cap\G}$ is rigid and has a natural commutativity constraint.

\head 19. Algebras with trace form\endhead  
\subhead 19.1\endsubhead
Let $R$ be a field and let $A$ be an associative $R$-algebra with $1$ of finite
dimension over $R$. We assume that $A$ is semisimple and split over $R$ and that we
are given a {\it trace form} on $A$ that is, an $R$-linear map $\t:A@>>>R$ such that
$(a,a')=\t(aa')=\t(a'a)$ is a non-degenerate (symmetric) $R$-bilinear form
$(,):A\T A@>>>R$. Note that $(aa',a'')=(a,a'a'')$ for all $a,a',a''$ in $A$. Let 
$\Mod A$ be the category whose objects are left $A$-modules of finite dimension over
$R$. We write $E\in\Ir A$ for "$E$ is a simple object of $\Mod A$".

Let $(a_i)_{i\in I}$ be an $R$-basis of $A$. Define an $R$-basis $(a'_i)_{i\in I}$ 
of $A$ by $(a_i,a'_j)=\d_{ij}$. Then 

(a) $\su_ia_i\ot a'_i\in A\ot A$ is independent of the choice of $(a_i)$.

\proclaim{Proposition 19.2} (a) We have $\su_i\t(a_i)a'_i=1$. 

(b) If $E\in\Ir A$, then $\su_i\tr(a_i,E)a'_i$ is in the centre of $A$. It acts on
$E$ as a scalar $f_E\in R$ times the identity and on $E'\in\Ir A$, not isomorphic to
$E$, as zero. Moreover, $f_E$ does not depend on the choice of $(a_i)$.

(c) One can attach uniquely to each $E\in\Ir A$ a scalar $g_E\in R$ (depending only
on the isomorphism class of $E$), so that

$\su_Eg_E\tr(a,E)=\t(a)$ for all $a\in A$,
\nl
where the sum is taken over all $E\in\Ir A$ up to isomorphism. 

(d) For any $E\in\Ir A$ we have $f_Eg_E=1$. In particular, $f_E\ne 0, g_E\ne 0$.

(e) If $E,E'\in\Ir A$, then $\su_i\tr(a_i,E)\tr(a'_i,E')$ is $f_E\dim E$ if $E,E'$
are isomorphic and is $0$, otherwise.
\endproclaim
Let $A=\op_{n=1}^tA_n$ be the decomposition of $A$ as a sum of simple algebras. Let
$\t_n:A_n@>>>R$ be the restriction of $\t$. Then $\t_n$ is a trace form for $A_n$, 
whose associated form is the restriction of $(,)$ and $(A_n,A_{n'})=0$ for 
$n\ne n'$. Hence we can choose $(a_i)$ so that each $a_i$ is contained in some $A_n$
and then $a'_i$ will be contained in the same $A_n$ as $a_i$.

We prove (a). From 19.1(a) we see that $\su_i\t(a_i)a'_i$ is independent of the 
choice of $(a_i)$. Hence we may choose $(a_i)$ as in the first paragraph of the 
proof. We are thus reduced to the case where $A$ is simple. In that case the 
assertion is easily verified.

We prove (b).  From 19.1(a) we see that $\su_i\tr(a_i,E)a'_i$ is independent of the
choice of $(a_i)$. Hence we may choose $(a_i)$ as in the first paragraph of the 
proof. We are thus reduced to the case where $A$ is simple. In that case the 
assertion is easily verified.

We prove (c). It is enough to note that $a\m\tr(a,E)$ form a basis of the space of 
$R$-linear functions $A@>>>R$ which vanish on all $aa'-a'a$ and $\t$ is such a 
function. 

We prove (d). We consider the equation in (c) for $a=a_i$, we multiply both sides by
$a'_i$ and sum over $i$. Using (a), we obtain
$$\su_i\su_Eg_E\tr(a_i,E)a'_i=\su_i\t(a_i)a'_i=1.$$
Hence $\su_Eg_E\su_i\tr(a_i,E)a'_i=1$. By (b), the left hand side acts on a
$E'\in\Ir A$ as a scalar $g_{E'}f_{E'}$ times the identity. This proves (d).

(e) follows immediately from (b). The proposition is proved.

\subhead 19.3\endsubhead
Now let $A'$ be a semisimple subalgebra of $A$ such that $\t'$, the restriction of
$\t$ to $A'$ is a trace form of $A'$. (We do not assume that the unit element 
$1_{A'}$ of $A'$ coincides to the unit element $1$ of $A$.) If $E\in\Mod A$ then 
$1_{A'}E$ is naturally an object of $\Mod A'$. Hence if $E'\in\Ir A'$, then the 
multiplicity $[E':1_{A'}E]$ of $E'$ in $1_{A'}E'$ is well defined.

Note that, if $a'\in A'$, then $\tr(a',1_{A'}E)=\tr(a',E)$.

\proclaim{Lemma 19.4} Let $E'\in\Ir A'$. We have $g_{E'}=\su_E[E':1_{A'}E]g_E$, sum
over all $E\in\Ir A$ (up to isomorphism).
\endproclaim
By the definition of $g_{E'}$, it is enough to show that

(a) $\su_{E'}\su_E[E':1'E]g_E\tr(a',E')=\t(a')$
\nl
for any $a'\in A'$. Here $E'$ (resp. $E$) runs over the isomorphism classes of
simple objects of $\Mod A'$ (resp. $\Mod A$). The left hand of (a) is
$$\su_Eg_E\su_{E'}[E':1'E]\tr(a',E')=\su_Eg_E\tr(a',1'E)=\su_Eg_E\tr(a',E)
=\t(a').$$
This completes the proof.

\head 20. The function $\aa_E$\endhead
\subhead 20.1\endsubhead
In this chapter we assume that $W$ is finite (hence $W,L$ is automatically bounded) and that
P1-P15 are satisfied.

The results of \S19 will be applied in the following cases.

(a) $A=\ch_\bc, R=\bc$. Here $\ca@>>>\bc$ takes $v$ to $1$. We identify $\ch_\bc$ 
with the group algebra $\bc[W]$ by $w\m T_w$ for all $w$. It is well known that 
$\bc[W]$ is a semisimple split algebra. We take $\t$ so that $\t(x)=\d_{x,1}$ for
$x\in W$. Then the bases $(x)$ and $(x\i)$ are dual with respect to $(,)$. 

We will say "$W$-module" instead of "$\bc[W]$-module". We will write $\Mod W,\Ir W$ 
instead of $\Mod\bc[W],\Ir\bc[W]$.

(b) $A=J_\bc$, $R=\bc$. Since $\bc[W]$ is semisimple, we see from 18.12(a) that the
kernel of $\ph_\bc:\bc[W]@>>>J_\bc$ is $0$ so that $\ph_\bc$ is injective. Since 
$\dim\bc[W]=\dim J_\bc=\sh W$ it follows that $\ph_\bc$ is an isomorphism. In 
particular $J_\bc$ is a semisimple split algebra. We take $\t:J_\bc@>>>\bc$ so that
$\t(t_z)$ is $n_z$ if $z\in\cd$ and $0$, otherwise. Then $(t_x,t_y)=\d_{xy,1}$. The
bases $(t_x)$ and $(t_{x\i})$ are dual with respect to $(,)$.

(c) $A=\ch_{\bc(v)}, R=\bc(v)$. Here $\ca@>>>\bc$ takes $v$ to $v$. The homomorphism
$\ph_{\bc(v)}:\ch_{\bc(v)}@>>>J_{\bc(v)}$ is injective. This follows from 18.12(b),
using the fact that injectivity is preserved by tensoring with a field of fractions.
Since $\ch_{\bc(v)},J_{\bc(v)}$ have the same dimension, it follows that 
$\ph_{\bc(v)}$ is an isomorphism. Since $J_{\bc(v)}=\bc(v)\ot J_\bc$, and $J_\bc$ is
semisimple, split, it follows that $J_{\bc(v)}$ is semisimple, split, hence 
$\ch_{\bc(v)}$ is semisimple, split. We take $\t:\ch_{\bc(v)}@>>>\bc(v)$ so that 
$\t(T_w)=\d_{w,1}$. The bases $(T_x)$ and $(T_{x\i})$ are dual with respect to 
$(,)$.

{\it Remark.} The argument above shows also that,

(d) if $R=R_0(v)$, with $R_0$ an arbitrary field and $\ca@>>>R$ carries $v$ to $v$,
then $\ph_R:\ch_R@>>>J_R$ is an isomorphism;

(e) if $R$ in 18.11 is a field of characteristic $0$ then $\ph_R:\ch_R@>>>J_R$ is an
isomorphism if and only if $\ch_R$ is a semisimple $R$-algebra.

\subhead 20.2\endsubhead
For any $E\in\Mod W$ we denote by $E_\sp$ the corresponding $J_\bc$-module. Thus, 
$E_\sp$ coincides with $E$ as a $\bc$-vector space and the action of $j\in J_\bc$ 
on $E_\sp$ is the same as the action of $\ph_\bc\i(j)$ on $E$. The $J_\bc$-module
structure on $E_\sp$ extends in a natural way to a $J_{\bc(v)}$-module structure on
$E_v=\bc(v)\ot_\bc E_\sp$. We will also regard $E_v$ as an $\ch_{\bc(v)}$-module via
the algebra isomorphism $\ph_{\bc(v)}:\ch_{\bc(v)}@>\si>>J_{\bc(v)}$. If $E$ is
simple, then $E_\sp$ and $E_v$ are simple.

Let $E\in\Ir W$. Then $E_\sp$ is a simple $J^\boc_\bc$-module for a unique two-sided
cell $\boc$ of $W$. Then for any $x\in\boc$, we write $E\si_{\lr}x$. If
$E,E'\in\Ir W$, we write $E\si_{\lr}E'$ if for some $x\in W$ we have $E\si_{\lr}x$,
$E'\si_{\lr}x$.

\subhead 20.3\endsubhead
There is the following direct relationship between $E$ and $E_v$ (without going 
through $J$):
$$\tr(x,E)=\tr(T_x,E_v)|_{v=1} \text{ for all } x\in W.$$
Indeed, it is enough to show that $\tr(c_x,E)=\tr(c_x,E_v)|_{v=1}$. Both 
sides are equal to $\su_{z\in W,d\in\cd}\g_{x,d,z\i}\hn_z\tr(t_z,E_\sp)$.

\subhead 20.4\endsubhead
Assume that $E\in\Ir W$. We have 

(a) $(f_{E_v})_{v=1}\dim(E)=\sh W$.
\nl
Indeed, setting $v=1$ in $\su_{x\in W}\tr(T_x,E_v)\tr(T_{x\i},E_v)=f_{E_v}\dim(E)$
gives
$$\su_{x\in W}\tr(x,E)\tr(x\i,E)=(f_{E_v})_{v=1}\dim(E).$$
The left hand side equals $\sh W$; (a) follows.

\subhead 20.5\endsubhead
Let $I\sub S$, let $E'\in\Ir W_I$ and let $E\in\Ir W$. We have 

(a) $[E'_v:E_v]=[E':E]$.
\nl
The right hand side is $\sh(W_I)\i\su_{x\in W_I}\tr(x,E')\tr(x\i,E)$. The left hand 
side is 
$$f_{E'_v}\i\dim(E')\i\su_{x\in W_I}\tr(T_x,E'_v)\tr(T_{x\i},E_v).$$ 
Since this is a constant, it is equal to its value for $v=1$. Hence it is equal to
$$(f_{E'_v}\i)_{v=1}\dim(E')\i\su_{x\in W_I}\tr(x,E'_v)\tr(x\i,E_v).$$
Thus it is enough to show that $(f_{E'_v})_{v=1}\dim(E')=\sh(W_I)$. But this is a
special case of 20.4(a). 

\proclaim{Proposition 20.6} Let $E\in\Ir W$.

(a) There exists a unique integer $\aa_E\ge 0$ such that $\tr(T_x,E_v)\in v^{-\aa_E}\bc[v]$ for all $x\in W$
and $\tr(T_x,E_v)\n v^{-\aa_E+1}\bc[v]$ for some $x\in W$.

(b) For $x\in W$ we have $\tr(T_x,E_v)=\sg(x)v^{-\aa_E}\tr(t_x,E_\sp)\mod v^{-\aa_E+1}\bc[v]$.

(c) Let $\boc$ be the two-sided cell such that $E_\sp\in\Ir J^\boc_\bc$. Then $\aa_E=\aa(z)$ for any 
$z\in\boc$.
\endproclaim
Let $a=\aa(z)$ for any $z\in\boc$. Let $x\in W$. By definition,
$$\tr(c_x^\da,E_v)=\su_{z\in W,d\in\cd;\aa(d)=\aa(z)}h_{x,d,z}n_d\tr(t_z,E_\sp).$$
In the last sum we have $\tr(t_z,E_\sp)=0$ unless $z\in\boc$ in which case 
$\aa(z)=a$ and $h_{x,d,z}=\ov h_{x,d,z}=\g_{x,d,z\i}v^{-a}\mod v^{-a+1}\bz[v]$. Thus we have
$$\tr(c_x^\da,E_v)=\sum_{z\in W,d\in\cd}\g_{x,d,z\i}n_d\tr(t_z,E_\sp)v^{-a}\mod v^{-a+1}\bc[v].$$
For each $z$ in the last sum we have $\sum_{d\in\cd}\g_{x,d,z\i}n_d=\d_{x,z}$. This gives
$$\tr(c_x^\da,E_v)=\tr(t_x,E_\sp)v^{-a}\mod v^{-a+1}\bc[v].\tag d$$
We have $T_x=\sum_{y;y\le x}q'_{y,x}c_y$ with $q'_{y,x}$ as in 10.1. Hence
$\sg(x)\ov T_x=T_x^\da=\su_{y;y\le x}q'_{y,x}c_y^\da$. Applying $\bar{}$ gives
$\sg(x)T_x=\su_{y;y\le x}\ov q'_{y,x}c_y^\da$. Hence
$$\tr(T_x,E_v)=\sg(x)\su_{y;y\le x}\ov q'_{y,x}\tr(c_y^\da,E_v).$$
Using (d) together with $\ov q'_{x,x}=1$, $\ov q'_{y,x}\in v\bz[v]$ for $y<x$ (see 10.1), we deduce
$$\tr(T_x,E_v)=\sg(x)\tr(t_x,E_\sp)v^{-a}\mod v^{-a+1}\bc[v].$$
Since $E_\sp\in\Ir J_\bc$, we have $\tr(t_x,E_\sp)\ne 0$ for some $x\in W$. The proposition follows.

\proclaim{Corollary 20.7} $f_{E_v}=f_{E_\sp}v^{-2\aa_E}+\text{strictly higher powers of } v$. 
\endproclaim  
Using 19.2(e) for $\ch_{\bc(v)}$ and $J_\bc$, we obtain
$$\align &f_{E_v}\dim E=\su_x\tr(T_x,E_v)\tr(T_{x\i},E_v) \\&\in
\su_x\tr(t_x,E_\sp)\tr(t_{x\i},E_\sp)v^{-2\aa_E}+v^{-2\aa_E+1}\bc[v]\\&
=f_{E_\sp}\dim Ev^{-2\aa_E}+v^{-2\aa_E+1}\bc[v].\endalign$$
The corollary follows.

\medpagebreak

Let $\bar{}:\bc[v,v\i]@>>>\bc[v,v\i]$ be the $\bc$-algebra homomorphism given by $v^n\m v^{-n}$ for all $n$.

\proclaim{Corollary 20.8} For any $h\in\ch$ we have $\tr(\bar h,E_v)=\ov{\tr(h,E_v)}$.
\endproclaim
We can assume that $h=c_x^\da$ where $x\in W$. As in the proof of 20.6 we have
$$\tr(c_x^\da,E_v)=\sum_{z\in W,d\in\cd;\aa(d)=\aa(z)}h_{x,d,z}n_d\tr(t_z,E_\sp).$$
Hence it suffices to note that $\ov{h_{x,d,z}}=h_{x,d,z}$ for all $d,z$ in the last sum. 

\medpagebreak

For $E\in\Mod W$ we write $E^\da$,$E_v^\da,E^\da_\sp$ instead of $E\ot\sg,(E\ot\sg)_v,(E\ot\sg)_\sp$.

\proclaim{Lemma 20.9} Let $E\in\Ir W$. For any $x\in W$ we have
$$\tr(T_x,E^\da_v)=(-1)^{l(x)}\ov{\tr(T_x,E_v)}.$$
\endproclaim
There is a unique a $\bc(v)$-algebra involution ${}^\da:\ch_{\bc(v)}@>>>\ch_{\bc(v)}$ extending 
${}^\da:\ch@>>>\ch$ (see 3.5). Let $(E_v)^\da$ be the $\ch_{\bc(v)}$-module with underlying vector space 
$E_v$ such that the action of $h\in\cv_{\bc(v)}$ on $(E_v)^\da$ is the same as the action of $h^\da$ on 
$E_v$. Clearly, $(E_v)^\da\in\Ir\ch_{\bc(v)}$. For $x\in W$ we have 
$$\tr(T_x,(E_v)^\da)=(-1)^{l(x)}\tr(T_{x\i}\i,E_v)=(-1)^{l(x)}\ov{\tr(T_x,E_v)}.$$
(The last equation follows from 20.8.) Setting $v=1$ we obtain
$$\tr(T_x,(E_v)^\da)|_{v=1}=(-1)^{l(x)}\tr(x,E)=\tr(x,E^\da).$$
Using 20.3, we deduce that $(E_v)^\da\cong E^\da_v$ in $\Mod\ch_{\bc(v)}$. The lemma follows.

\proclaim{Proposition 20.10} For any $x\in W$ we have
$$\tr(T_x,E_v)=\tr(t_x,E^\da_\sp)v^{\aa_{E^\da}}+\text{ strictly lower powers of } v.$$
\endproclaim
By 20.9 and 20.6 we have
$$\align&\tr(T_x,E_v)=\sg(x)\ov{\tr(T_x,E^\da_v)}
=\ov{\tr(t_x,E^\da_\sp)v^{-\aa_{E^\da}}+\text{ strictly higher powers of } v}\\&
=\tr(t_x,E^\da_\sp)v^{\aa_{E^\da}}+\text{ strictly lower powers of } v.\endalign$$
The proposition is proved.

\proclaim{Corollary 20.11} $f_{E_v}=f_{E^\da_\sp}v^{2\aa_{E^\da}}+\text{ strictly lower powers of } v$.
\endproclaim
Using 20.10 we have
$$\align &f_{E_v}\dim E=\su_x\tr(T_x,E_v)\tr(T_{x\i},E_v) \\&\in\su_x\tr(t_x,
E^\da_\sp)\tr(t_{x\i},E^\da_\sp)v^{2\aa_{E^\da}}+v^{2\aa_{E^\da}-1}\bc[v\i]\\&
=f_{E^\da_\sp}\dim Ev^{2\aa_{E^\da}}+v^{2\aa_{E^\da}-1}\bc[v\i].\endalign$$

\proclaim{Lemma 20.12} Let $E'\in\Ir W_I$. With notation of 19.2, we have
$$g_{E'_v}=\su_{E;E\in\Ir W}[E':E]g_{E_v}.$$
\endproclaim
We apply 19.4 with $A=\ch_{\bc(v)}$ and $A'$ the analogous algebra for $W_I$ instead
of $W$, identified naturally with a subspace of $A$. (In this case the unit elements
of the two algebras are compatible hence $1_{A'}E_v=E_v$.) It remains to use 20.5(a).

\proclaim{Lemma 20.13} Let $E\in\Ir W$.

(a) For any $x\in W$, $\tr(t_{x\i},E_\sp)$ is the complex conjugate of 
$\tr(t_x,E_\sp)$.

(b) $f_{E_\sp}$ is a strictly positive real number.
\endproclaim
We prove (a). Let $\la,\ra:E\T E@>>>\bc$ be a positive definite hermitian form. We
define $\la,\ra':E_\sp\T E_\sp@>>>\bc$ by 
$\la e,e'\ra'=\su_{z\in W}\la t_ze,t_ze'\ra$.
This is again a positive definite hermitian form on $E_\sp$. We show that
$$\la t_xe,e'\ra'=\la e,t_{x\i}e'\ra'$$
for all $e,e'$. This is equivalent to
$$\su_{y,z}\g_{z,x,y\i}\la t_ye,t_ze'\ra=\su_{y,z}\g_{y,x\i,z\i}\la t_ye,t_ze'\ra$$
which follows from $\g_{z,x,y\i}=\g_{y,x\i,z\i}$. We see that $t_{x\i}$ is the 
adjoint of $t_x$ with respect to a positive definite hermitian form. (a) follows.

We prove (b). By 19.2(e) we have $f_{E_\sp}\dim(E)=\su_x\tr(t_x,E_\sp)\tr(t_{x\i},E_\sp)$. The right hand 
side of this equality is a real number $\ge 0$, by (a). Hence so is the left hand side. Now
$f_{E_\sp}\ne 0$ by 19.2(d) and (b) follows.

\proclaim{Proposition 20.14} Let $E'\in\Ir W_I$.

(a) For any $E\in\Ir W$ such that $[E':E]\ne 0$ we have $\aa_{E'}\le\aa_E$.

(b) We have $g_{E'_\sp}=\su [E':E]g_{E_\sp}$, sum over all $E\in\Ir W$ (up to 
isomorphism) such that $\aa_E=\aa_{E'}$.
\endproclaim
Let $X$ be the set of all $E$ (up to isomorphism) such that $[E':E]\ne 0$ and such 
that $\aa_E$ is minimum, say equal to $a$. Assume first that $a<\aa_{E'}$. Using 
19.2(d) we rewrite 20.12 in the form 

(c) $v^{-2a}f_{E'_v}\i=\su_E[E':E]v^{-2a}f_{E_v}\i$.
\nl
By 20.7, we have 

(d) $(v^{-2\aa_E}f_{E_v}\i)|_{v=0}=f_{E_\sp}\i$,
$(v^{-2\aa_{E'}}f_{E'_v}\i)|_{v=0}=f_{E'_\sp}\i$, 
\nl
hence by setting $v=0$ in (c) we obtain
$$0=\su_{E\in X}[E':E]f_{E_\sp}\i.$$
The right hand side is a real number $>0$ by 20.8(b). This is a contradiction. Thus
we must have $a\ge\aa_{E'}$ and (a) is proved.

We now rewrite (c) in the form

(e) $v^{-2\aa_{E'}}f_{E'_v}\i=\su_E[E':E]v^{-2\aa_{E'}}f_{E_v}\i$.
\nl
Using (d) and (a) we see that, setting $v=0$ in (e) gives
$$f_{E'_\sp}\i=\su_{E;\aa_E=\aa_{E'}}[E':E]f_{E_\sp}\i.$$
This proves (b).

\subhead 20.15\endsubhead
Let $K(W)$ be the $\bc$-vector space with basis indexed by the $E\in\Ir W$ (up to 
isomorphism). If $\tE\in\Mod W$ we identify $\tE$ with the element
$\su_E[E:\tE]E\in K(W)$ ($E$ as above). We define a $\bc$-linear map 
$\jj_{W_I}^W:K(W_I)@>>>K(W)$ by

$\jj_{W_I}^W(E')=\su_E[E':E]E$,
\nl
sum over all $E\in\Ir W$ (up to isomorphism) such that $\aa_E=\aa_{E'}$; here 
$E'\in\Ir W_I$. We call this {\it truncated induction}.

Let $I''\sub I'\sub S$. We show that the following transitivity formula holds:
$$\jj_{W_{I'}}^W\jj_{W_{I''}}^{W_{I'}}=\jj_{W_{I''}}^W:K(W_{I''})@>>>K(W).\tag a$$
Let $E''\in\Ir W_{I''}$. We must show that
$$[E'':E]=\su_{E';\aa_{E'}=\aa_{E''}}[E'':E'][E':E]$$
for any $E''\in\Ir W_{I''},E\in\Ir W$ such that $\aa_{E''}=\aa_E$; in the sum we 
have $E'\in\Ir W_{I'}$. Clearly,
$$[E'':E]=\su_{E'}[E'':E'][E':E].$$
Hence it is enough to show that, if $[E'':E'][E':E]\ne 0$, then we automatically 
have
$\aa_{E'}=\aa_{E''}$. By 2.10(a) we have $\aa_{E''}\le\aa_{E'}\le\aa_E$. Since
$\aa_{E''}=\aa_E$, the desired conclusion follows.

\subhead 20.16\endsubhead
For any $x\in W$ we set 
$$\g_x=\su_{E;E\in\Ir W}\tr(t_x,E_\sp)E\in K(W).$$
We sometimes write $\g_x^W$ instead of $\g_x$, to emphasize dependence on $W$. Note
that $\g_x$ is a $\bc$-linear combination of $E$ such that $E\si_{\lr}x$. Hence, if
$E,E'$ appear with $\ne 0$ coefficient in $\g_x$ then $E\si_{\lr}E'$.

\proclaim{Proposition 20.17} If $x\in W_I$, then $\g_x^W=\jj_{W_I}^W(\g_x^{W_I})$.
\endproclaim
An equivalent statement is
$$\tr(t_x,E_\sp)=\su_{E';\aa_E=\aa_{E'}}\tr(t_x,E'_\sp)[E':E]\tag a$$
for any $E\in\Ir W$; in the sum we have $E'\in\Ir W_I$. Clearly, we have
$$v^{\aa_E}\tr(T_x,E_v)=\su_{E';E'\in\Ir W_I}v^{\aa_E}\tr(T_x,E'_v)[E':E].\tag b$$
In the right hand side we may assume that $\aa_{E'}\le\aa_E$. Using this and 20.6, 
we see that setting $v=0$ in (b) gives (a). The proposition is proved. 

\proclaim{Lemma 20.18} (a) We have $\aa_{\sg}=L(w_0)$.

(b) We have $f_{\sg_\sp}=1$.

(c) We have $\g_{w_0}=\sg$. 
\endproclaim
$\sg_v$ is the one dimensional $\ch_{\bc(v)}$-module on which $T_x$ acts as 
$\sg(x)v^{-L(x)}$. (This follows from 20.3.) From 20.6(b) we see that 
$\aa_{\sg}=L(w_0)$ and that $\tr(t_{w_0},\sg_\sp)=1$. This proves (a). To prove (c)
it remains to show that, if $\tr(t_{w_0},E_\sp)\ne 0$ ($E$ simple) then $E\cong\sg$.
This assumption shows, by 20.6(c), that $E_\sp\in\Ir J^\boc_\bc$ where $\boc$ is 
the two-sided cell such that $\sg_\sp\in\Ir J^\boc_\bc$. Since 
$\tr(t_{w_0},\sg_\sp)=1$, we have $w_0\in\boc$. From 13.8 it follows that $\{w_0\}$
is a two-sided cell. Thus $\boc=\{w_0\}$ and $J^\boc_\bc$ is one dimensional. Hence
it cannot have more than one simple module. Thus, $E\cong\sg$. This yields (c) and 
also (b). The lemma is proved.

\subhead 20.19\endsubhead
Assume that $I,I'$ form a partition of $S$ such that $W=W_I\T W_{I'}$. If 
$E\in\Ir W_I$ and $E'\in\Ir W_{I'}$, then $E\bxt E'\in\Ir W$. From the definitions,
$$\aa_{E\bxt E'}=\aa_E+\aa_{E'}, f_{(E\bxt E')_\sp}=f_{E_\sp}f_{E'_\sp}.$$
Moreover, if $x\in W_I,x'\in W_{I'}$, then
$$\g_{xx'}^W=\g_x^{W_I}\bxt\g_{x'}^{W_{I'}}.$$

\subhead 20.20\endsubhead
Until the end of 20.23 we assume that $w_0$ is in the centre of $W$. 
Then, for any $E\in\Ir W$, $w_0$ acts on $E$ as $\e_E$ times identity where $\e_E=\pm 1$.
Now $E\m\e_EE$ extends to a $\bc$-linear involution $\z:K(W)@>>>K(W)$.

\proclaim{Lemma 20.21} Let $E\in\Ir W$. For any $x\in W$ we have

$\tr(T_{w_0x},E_v)=\e_E v^{-\aa_E+\aa_{E^\da}}\ov{\tr(T_x,E_v)}$.
\endproclaim
Since $w_0$ is in the centre of $W$, $T_{w_0}$ is in the centre of $\ch_{\bc(v)}$ 
hence it acts on $E_v$ as a scalar $\l\in\bc(v)$ times the identity. Now 
$\tr(T_x,E_v)\in\bc[v,v\i]$ and $\tr(T_x\i,E_v)\in\bc[v,v\i]$. In particular, 
$\l\in\bc[v,v\i]$ and $\l\i\in\bc[v,v\i]$. This implies $\l=cv^n$ where $c\in\bc$.
For $v=1$, $\l$ becomes $\e_E$. Hence $\l=\e_Ev^n$ for some $n$. We have 
$$\tr(T_{w_0x},E_v)=\tr(T_{w_0}T_{x\i}\i,E_v)=\l\tr(T_{x\i}\i,E_v)=
\l\ov{\tr(T_x,E_v)}.$$
We have
$$\su_x\tr(T_{w_0x},E_v)\tr(T_{x\i w_0},E_v)=
\l^2\su_x\ov{\tr(T_x,E_v)}\ov{\tr(T_{x\i},E_v)}$$
hence $f_{E_v}\dim(E)=\l^2\ov{f_{E_v}}\dim(E)$ so that $f_{E_v}=v^{2n}\ov{f_{E_v}}$.
By 20.9, we have
$$\su_x\ov{\tr(T_x,E_v)}\ov{\tr(T_{x\i},E_v)}=
\su_x\tr(T_x,E^\da_v)\tr(T_{x\i},E^\da_v)$$
hence $\ov{f_{E_v}}=f_{E^\da_v}$. We see that $f_{E_v}=v^{2n}f_{(E^\da)_v}$. 
Comparing the lowest terms we see that
$-2\aa_E=2n-2\aa_{E^\da}$ hence $n=-\aa_E+\aa_{E^\da}$ and that 

(a) $f_{E_\sp}=f_{E^\da_\sp}$.

\proclaim{Lemma 20.22} 
$v^{\aa_E}\tr(T_{w_0x},E_v)=\e_E(-1)^{l(x)}v^{\aa_{E^\da}}\tr(T_x,E^\da_v)$.
\endproclaim
We combine 20.8, 20.21.

\proclaim{Lemma 20.23} For any $x\in W$ we have $\g_{xw_0}=\sg(x)\z(\g_x)\ot\sg$.
\endproclaim
An equivalent statement is
$$\tr(t_{xw_0},E_\sp)=\sg(x)\tr(t_x,E^\da_\sp)\e_{E^\da}$$
for any $E\in\Ir W$. Setting $v=0$ in the identity in 20.22 gives 
$$\sg(xw_0)\tr(t_{w_0x},E_\sp)=\e_E\tr(t_x,E^\da_\sp).$$
It remains to show that $\e_{E^\da}=\e_E\sg(w_0)$. This is clear.

\subhead 20.24\endsubhead
By the Cayley-Hamilton theorem, any element $r\in J$ satisfies an equation of the 
form $r^n+a_1r^{n-1}+\do+a_n=0$ where $a_i\in\bz$. (We use that the structure 
constants of $J$ are integers.) This holds in particular for $r=t_x$ where $x\in W$.
Hence for any $\ce\in\Ir J_\bc$, $\tr(t_x,\ce)$ is an algebraic integer. If $R$ is a
subfield of $\bc$ such that the group algebra $R[W]$ is split over $R$, then $J_R$
is split over $R$ and it follows that for $x,\ce$ as above,  $\tr(t_x,\ce)$ is an 
algebraic integer in $R$. In particular, if we can take $R=\bq$, then 
$\tr(t_x,\ce)\in\bz$.

\head 21. Study of a left cell\endhead
\subhead 21.1\endsubhead
In this chapter we preserve the setup of 20.1. Let $\G$ be a left cell of $W,L$. Let
$d$ be the unique element in $\G\cap\cd$. The $\ca$-submodule 
$\su_{y\in\G}\ca c_y^\da$ of $\ch$ can be regarded as an $\ch$-module by the rule
$c_x^\da\cdot c_w^\da=\su_{z\in\G}h_{x,y,z}c_z^\da$ with $x\in W,y\in W$. By change
of scalars ($v\m 1$) this gives rise to an $\ch_\bc=\bc[W]$-module $[\G]$. On the 
other hand, $J^\G_\bc=\op_{y\in\G}\bc t_y$ is a left ideal in $J_\bc$ by 14.2(P8).

\proclaim{Lemma 21.2} The $\bc$-linear isomorphism $t_y\m c_y^\da$ for $y\in\G$
is an isomorphism of $J_\bc$-modules $J^\G_\bc@>\si>>[\G]_\sp$.
\endproclaim
We have $\G\sub X_a=\{w\in W;\aa(x)=a\}$ for some $a\in\bn$. The $\ca$-submodule 
$\su_{y\in X_a}\ca c_y^\da$ of $\ch$ can be regarded as an $\ch$-module by the rule
$$c_x^\da\cdot c_w^\da=\su_{z\in X_a}h_{x,y,z}c_z^\da$$
with $x\in W,y\in W$. By change of scalars ($v\m 1$) this gives rise to an 
$\ch_\bc=\bc[W]$-module $[X_a]$. On the other hand, 
$J^{X_a}_\bc=\op_{y\in X_a}\bc t_y$ is a left (even two-sided) ideal in $J_\bc$. The
$\bc$-linear map in the lemma extends by the same formula to a $\bc$-linear 
isomorphism $J^{X_a}_\bc@>\si>>[X_a]_\sp$. It is enough to show that this is 
$J_\bc$-linear. This follows from the computation in 18.10. The lemma is proved.

\proclaim{Lemma 21.3} Let $\ce\in\Ir J_\bc$. The $\bc$-linear map 
$u:\Hom_{J_\bc}(J^\G_\bc,\ce)@>>>t_d\ce$ given by $\xi\m\xi(n_dt_d)$ is an 
isomorphism.
\endproclaim
$u$ is well defined since $\xi(n_dt_d)=t_d\xi(t_d)\in t_d\ce$. We define a linear 
map in the opposite direction by $e\m[j\m je]$. It is clear that this is the inverse
of $u$. (We use that $jn_dt_d=j$ for $j\in J^\G_\bc$.) The lemma is proved.

\proclaim{Proposition 21.4} We have $\g_d=n_d\su_E[E:[\G]]E$ (sum over all 
$E\in\Ir W$ up to isomorphism).
\endproclaim
An equivalent statement is that $\tr(n_dt_d,E_\sp)=[E:[\G]]$, for $E$ as above. By
21.2, we have $[E:[\G]]=[E_\sp:J^\G_\bc]$. Hence it remains to show that 
$\tr(n_dt_d,\ce)=[\ce:J^\G_\bc]$ for any $\ce\in\Ir J_\bc$. Since 
$\ce=\op_{d'\in\cd}n_{d'}t_{d'}\ce$ and $n_dt_d:\ce@>>>\ce$ is the projection to the
summand $n_dt_d\ce$, we see that $\tr(n_dt_d,\ce)=\dim(t_d\ce)$. It remains to show
that $\dim(t_d\ce)=[\ce:J^\G_\bc]$. This follows from 21.3.

\proclaim{Proposition 21.5} $[\G]^\da,[\G w_0]$ are isomorphic in $\Mod W$.
\endproclaim
We may identify $[\G]^\da$ with the $W$-module with $\bc$-basis $(e_y)_{y\in\G}$ 
where $s\in S$ acts by $e_y\m-e_y+\su_{z\in\G}h_{s,y,z}e_z$. 

On the other hand we may identify $[\G w_0]$ with the $W$-module with $\bc$-basis 
$(e'_{yw_0})_{y\in\G}$ where $s\in S$ acts by 
$e'_{yw_0}\m e'_{yw_0}-\su_{z\in\G}h_{s,yw_0,zw_0}e'_{zw_0}$.

The $W$-module dual to $[\G]^\da$ has a $\bc$-basis $(e''_y)_{y\in\G}$ (dual to 
$(e_y)$) in which the action of $s\in S$ is given by
$e''_y\m-e''_y+\su_{z\in\G}h_{s,z,y}e''_z$. We define a $\bc$-isomorphism between 
this last space and $[\G w_0]$ by $e''_y\m\sg(y)e'_{yw_0}$ for all $y$. We show that
this comutes with the action of $W$. It suffices to show that for any $s\in S$, we 
have

(a) $-h_{s,z,y}=\sg(y)\sg(z)h_{s,yw_0,zw_0}$ for all $z\ne y$ and

(b) $1-h_{s,y,y}=-1+h_{s,yw_0,yw_0}$ for all $y$.
\nl
We use 6.6. Assume first that $sz>z$. If $sy>y$ and $y\ne z$, both sides of (a) are
$0$. If $sy<y<z$ then (a) follows from 11.6. If $y=sz$ then  both sides of (a) are 
$-1$. If $sy<y$ but $y\not<z$ or $y\ne sz$ then both sides of (a) are $0$.

Assume next that $sz<z$. If $z\ne y$ then both sides of (a) are $0$.

If $sy>y$, both sides of (b) are $1$. If $sy<y$, both sides of (b) are $-1$. Thus 
(a),(b) are verified. Since $[\G]^\da$ and its dual are isomorphic in $\Mod W$ (they
are defined over $\bq$), the lemma follows.

\proclaim{Corollary 21.6} Let $E\in\Ir W$ and let $\boc$ be the two-sided cell of 
$W$ such that $E_\sp\in\Ir J^{\boc}_\bc$. Then $E^\da_\sp\in\Ir J^{\boc w_0}_\bc$.
\endproclaim
Replacing $\G$ by $\boc$ in the definition of $[\G]$ we obtain a $W$-module 
$[\boc]$. Then 21.2, 21.5 hold with $\G$ replaced by $\boc$ with the same proof. Our
assumption implies (by 21.2 for $\boc$) that $E$ appears in the $W$-module $[\boc]$.
Using 21.5 for $\boc$ we deduce that $E^\da$ appears in the $W$-module $[\boc w_0]$.
Using 21.2 for $\boc w_0$, we deduce that $E^\da_\sp$ appears in the $J_\bc$-module
$J^{\boc w_0}_\bc$. The corollary follows. 

\proclaim{Corollary 21.7} Let $E,E'\in\Ir W$ be such that $E\si_{\lr}E'$. Then
$E^\da\si_{\lr}E'{}^\da$.
\endproclaim
By assumption there exists a two-sided cell $\boc$ such that 
$E_\sp,E'_\sp\in\Ir J^{\boc}_\bc$. By 21.6, 
$E^\da_\sp,E'{}^\da_\sp\in\Ir J^{\boc w_0}_\bc$. The corollary follows.

\subhead 21.8\endsubhead
The results of \S19 are applicable to $A$, the $\bc$-subspace $J^{\G\cap\G\i}_\bc$
of $J_\bc$ spanned by $\G\cap\G\i$ and $R=\bc$. This is a $\bc$-subalgebra of 
$J_\bc$ with unit element $n_dt_d$. 
In 21.9 we will show that $J^{\G\cap\G\i}_\bc$ is
semisimple. It is then clearly split. We take $\t:J^{\G\cap\G\i}_\bc@>>>\bc$ so that
$\t(t_x)=n_d\d_{x,d}$. (This is the restriction of $\t:J_\bc@>>>\bc$.) We have
$(t_x,t_y)=\d_{xy,1}$. The bases $(t_x)$ and $(t_{x\i})$ (where $x$ runs through 
$\G\cap\G\i$) are dual with respect to $(,)$. 

\subhead 21.9\endsubhead
We show that the $\bc$-algebra $J^{\G\cap\G\i}_\bc$ is semisimple. It is enough to
prove the analogous statement for the $\bq$-algebra $A'$, the $\bq$-span of 
$\G\cap\G\i$ in $J_\bq$. We define a $\bq$-bilinear pairing $(|):A'\T A'\to\bq$ by
$(t_x|t_y)=\d_{x,y}$  for $x,y\in\G\cap\G\i$. Let $j\m\ti j$ be the $\bq$-linear map
$A'\to A'$ given by $\ti t_x=t_{x\i}$ for all $x$. We show that 
$$(j_1j_2|j_3)=(j_2|\ti j_1j_3)\tag a$$
for all $j_1,j_2,j_3$ in our ring. We may assume that $j_1=t_x,j_2=t_y,j_3=t_z$. 
Then (a) follows from 
$$\g_{x,y,z\i}=\g_{x\i,z,y\i}.$$
Now let $I$ be a left ideal of $A'$. Let $I^\pe=\{a\in A';(a|I)=0\}$. Since $(|)$ is
positive definite, we have $A'=I\op I^\pe$. From (c) we see that $I^\pe$ is a left 
ideal. This proves that $A'$ is semisimple.

The same proof could be used to show directly that $J_\bc$ is semisimple.

\proclaim{Proposition 21.10} Let $E,E'\in\Ir W$,
$N=\su_{x\in\G\cap\G\i}\tr(t_x,E_\sp)\tr(t_{x\i},E'_\sp)$. Then 
$N=f_{E_\sp}[E:[\G]]$ if $E,E'$ are isomorphic and $N=0$, otherwise.
\endproclaim
If $\ce\in\Ir J_\bc$, then $t_d\ce$ is either $0$ or in $\Ir J^{\G\cap\G\i}_\bc$.
Moreover, $\ce\m t_d\ce$ defines a bijection between the set of simple 
$J_\bc$-modules (up to isomorphism) which appear in the $J_\bc$-module $J^\G_\bc$ 
and the set of simple $J^{\G\cap\G\i}_\bc$-modules (up to isomorphism). We then have
$\dim(t_d\ce)=[\ce:J_\bc^\G]$. For $j\in J^{\G\cap\G\i}_\bc$ we have 
$\tr(j,\ce)=\tr(j,t_d\ce)$. If $t_dE_\sp=0$ or $t_dE'_\sp=0$, then $N=0$ and the 
result is clear. If $t_dE_\sp\ne 0$ and $t'_dE_\sp\ne 0$ then, by 19.2(e), we see 
that $N=f_{t_dE_\sp}[E_\sp:J^\G_\bc]$ if $E,E'$ are isomorphic and to $0$, 
otherwise. It remains to show that $f_{t_dE_\sp}=f_{E_\sp}$, 
$[E:[\G]]=[E_\sp:J^\G_\bc]$ and the analogous equalities for $E'$. Now 
$f_{t_dE_\sp}=f_{E_\sp}$ follows from 19.4 applied to 
$(A',A)=(J^{\G\cap\G\i}_\bc,J_\bc)$; the equality $[E:[\G]]=[E_\sp:J^\G_\bc]$ 
follows from 21.2. The proposition is proved.

\head 22. Constructible representations\endhead
\subhead 22.1\endsubhead
In this chapter we preserve the setup of 20.1. 

We define a class $\C(W)$ of $W$-modules (relative to our fixed $L:W@>>>\bn$) by 
induction on $\sh S$. If $\sh S=0$ so that $W=\{1\}$, $\C(W)$ consists of the unit
representation. Assume now that $\sh S>0$. Then $\C(W)$ consists of the $W$-modules
of the form $\jj_{W_I}^W(E')$ or $\jj_{W_I}^W(E')\ot\sg$ for various subsets 
$I\sub S,I\ne S$ and various $E'\in\C(W_I)$. (If we restrict ourselves to $I$ such
that $\sh(S-I)=1$ we get the same class of $W$-modules, by the transitivity of 
truncated induction.) The $W$-modules in $\C(W)$ are said to be the {\it 
constructible representations} of $W$.

Now the unit representation of $W$ is constructible (it is obtained by truncated 
induction from the unit representation of the subgroup with one element). Hence 
$\sg\in\C(W)$.

\proclaim{Lemma 22.2} If $E\in\C(W)$, then there exists a left cell $\G$ of $W$ such
that $E=[\G]$.
\endproclaim
We argue by induction on $\sh S$. If $\sh S=0$ the result is obvious. Assume now 
that $\sh S>0$. Let $E\in\C(W)$.

{\it Case 1}. $E=\jj_{W_I}^W(E')$ where $I\sub S,I\ne S$ and $E'\in\C(W_I)$. By the
induction hypothesis there exists a left cell $\G'$ of $W_I$ such that $E'=[\G']$. 
Let $d\in\G'\cap\cd$. By 21.4 we have $\g_d^{W_I}=[\G']=E'$. By 20.17 we have 
$E=\jj_{W_I}^W(E')=\jj_{W_I}^W(\g_d^{W_I})=\g_d^W$. Let $\G$ be the left cell of $W$
that contains $d$. By 21.4 we have $\g_d^W=[\G]$. Hence $E=[\G]$.

{\it Case 2}. $E=\jj_{W_I}^W(E')\ot\sg$ where $I\sub S,I\ne S$ and $E'\in\C(W_I)$.
Then by Case 1, $E\ot\sg=[\G]$ for some left cell $\G$ of $W$. By 21.5 we have 
$E=[\G w_0]$. The lemma is proved.

\proclaim{Proposition 22.3} For any $E\in\Ir W$ there exists a constructible 
representation of $W$ which contains a simple component isomorphic to $E$.
\endproclaim
The general case can be easily reduced to the case where $W$ is irreducible. Assume
now that $W$ is irreducible. If $L=al$ for some $a>0$, the constructible
representations of $W$ are listed in \cite{\CG} and the proposition is easily 
checked. (See also the discussion of types $A,D$ in 22.5, 22.26.) In the cases where
$W$ is irreducible but $L$ is not of the form $al$, the constructible 
representations are described later in this chapter and this yields the proposition
in all cases.

\subhead 22.4\endsubhead
Let $W=\fS_n$ be the group of permutations of $1,2,\do,n$. We regard $W$ as a 
Coxeter group with generators
$$s_1=(1,2),s_2=(2,3),\do,s_{n-1}=(n-1,n),$$
(transpositions). We take $L=al$ where $a>0$.

The simple $W$-modules (up to isomorphism) are in 1-1 correspondence with the 
partitions $\a=(\a_1\ge\a_2\ge\do)$ such that $\a_N=0$ for large $N$ and 
$\su_i\a_i=n$. The correspondence (denoted by $\a\m\p_\a$) is defined as follows.
Let $\a$ be as above, let $(\a'_1\ge\a'_2\ge\do)$ be the partition dual to $\a$. Let
$\p_\a$ be the simple $W$-module whose restriction to $\fS_{\a_1}\T\fS_{\a_2}\do$ 
contains $1$ and whose restriction to $\fS_{\a'_1}\T\fS_{\a'_2}\do$ contains the 
sign representation. We have (a consequence of results of Steinberg):
$$f_{(\p_\a)_v}=v^{-\su_i2\bin{\a'_i}{2}}+\text{strictly higher powers of } v.$$
It follows that 

(a) $\aa_{\p_\a}=\su_ia\bin{\a'_i}{2}$ and $f_{(\p_\a)_\sp}=1$.

\proclaim{Lemma 22.5} In the setup of 22.4, a $W$-module is constructible if and 
only if it is simple.
\endproclaim
For any sequence $\b=(\b_1,\b_2,\do)$ in $\bn$ such that $\b_N=0$ for large $N$ and
$\su_i\b_i=n$, we set
$$I_\b=\{s_i;i\in[1,n-1],i\ne\b_1,i\ne\b_1+\b_2,\do\}.$$
From 22.4(a) we see easily that, if $\b$ is the same as $\a'$ up to order, then

(a) $\jj_{W_{I_\b}}^W(\sg)=\p_\a$. 
\nl
Since the $\sg\in\C(W_{I_\b})$, it follows that $\p_\a\in\C(W)$. Thus any simple
$W$-module is constructible.

We now show that any constructible representation $E$ of $W=\fS_n$ is simple. We may
assume that $n\ge 1$ and that the analogous result is true for any $W_{I'}\ne W$. We
may assume that $E=\jj_{W_{I_\b}}^W(C)$ where $\b$ is as above, $W_{I_\b}\ne W$ and
$C\in\C(W_{I_\b})$. By the induction hypothesis, $C$ is simple. Since the analogue
of (a) holds for $W_{I_\b}$ (instead of $W$) we have 
$C=\jj_{W_{I_{\b'}}}^{W_{I_\b}}(\sg)$ for some $\b'$ such that 
$W_{I_{\b'}}\sub W_{I_\b}$. By the transitivity of truncated induction we have 
$E=\jj_{W_{I_{\b'}}}^W(\sg)$. Hence, by (a), for $\b'$ instead of $\b$, $E$ is 
simple. The lemma is proved.

\subhead 22.6\endsubhead
We now develop some combinatorics which is useful for the verification of 22.3 for 
$W$ of classical type.

Let $a>0,b\ge 0$ be integers. We can write uniquely $b=ar+b'$ where $r,b'\in\bn$ and
$b'<a$. Let $N\in\bn$. Let $\cm_{a,b}^N$ be the set of multisets 
$\tZ=\{\tz_1\le\tz_2\le\do\le\tz_{2N+r}\}$ of integers $\ge 0$ such that 

(a) if $b'=0$, there are at least $N+r$ distinct entries in $\tZ$, no entry is 
repeated more than twice and all entries of $\tZ$ are divisible by $a$;

(b) if $b'>0$, all inequalities in $\tZ$ are strict and $N$ entries of $\tZ$ are
divisible by $a$ and $N+r$ entries of $\tZ$ are equal to $b'$ modulo $a$.
\nl
The entries which appear in $\tZ$ exactly once are called the {\it singles} of 
$\tZ$; they form a set $Z$. The other entries of $\tZ$ are called the {\it doubles}
of $\tZ$.

For example, the multiset $\tZ^0$ whose entries are (up to order)
$$0,a,2a,\do,(N-1)a,b',a+b',2a+b',\do,(N+r-1)a+b'$$
belongs to $\cm_{a,b}^N$. Clearly, the sum of entries of $\tZ$ minus the sum of 
entries of $\tZ^0$ is $\ge 0$ and divisible by $a$, hence it is equal to $an$ for a
well defined $n\in\bn$ said to be the {\it rank} of $\tZ$. We have
$$\su_{k=1}^{2N+r}\tz_k=an+aN^2+N(b-a)+a\bin{r}{2}+b'r.$$
Note that $\tZ^0$ has rank $0$. Let $\cm_{a,b;n}^N$ be the set of multisets of rank
$n$ in $\cm_{a,b}^N$. We define an (injective) map $\cm_{a,b}^N@>>>\cm_{a,b}^{N+1}$
by 
$$\{\tz_1\le\tz_2\le\do\le\tz_{2N+r}\}\m
\{0,b',\tz_1+a\le\tz_2+a\le\do\le\tz_{2N+r}+a\}.$$
This restricts for any $n\in\bn$ to an (injective) map 

(c) $\cm_{a,b;n}^N@>>>\cm_{a,b;n}^{N+1}$.
\nl
It is easy to see that, for fixed $n$, $\sh(\cm_{a,b;n}^N)$ is bounded as $N\to\iy$,
hence the maps (c) are bijections for large $N$. Let $\cm_{a,b;n}$ be the inductive
limit of $\cm_{a,b;n}^N$ as $N\to\iy$ (with respect to the maps (c)).

\subhead 22.7\endsubhead
Let $\Sy_{a,b;n}^N$ be the set consisting of all tableaux (or {\it symbols})
$$\align\l_1&,\l_2,\do,\l_{N+r}\\&\mu_1,\mu_2,\do,\mu_N\tag a\endalign$$
where $\l_1<\l_2<\do<\l_{N+r}$ are integers $\ge 0$, congruent to $b'$ modulo $a$, 
$\mu_1,\mu_2,\do,\mu_N$ are integers $\ge 0$, divisible by $a$ and
$$\su_i\l_i+\su_j\mu_j=an+aN^2+N(b-a)+a\bin{r}{2}+b'r.$$
If we arrange the entries of $\L$ in a single row, we obtain a multiset
$\tZ\in\cm_{a,b;n}^N$. This defines a (surjective) map 
$\p_N:\Sy_{a,b;n}^N@>>>\cm_{a,b;n}^N$. We define an (injective) map 

(b) $\Sy_{a,b;n}^N@>>>\Sy_{a,b;n}^{N+1}$
\nl
by associating to (a) the symbol
$$\align b',\l_1+a&,\l_2+a,\do,\l_{N+r}+a\\&0,\mu_1+a,\mu_2+a,\do,\mu_N+a.\endalign
$$
This is compatible with the map $\cm_{a,b}^N@>>>\cm_{a,b}^{N+1}$ in 22.6 (via 
$\p_N,\p_{N+1}$).

Since for fixed $n$, $\sh(\Sy_{a,b;n}^N)$ is bounded as $N\to\iy$, the maps (b) are 
bijections for large $N$. Let $\Sy_{a,b;n}$ be the inductive limit of 
$\Sy_{a,b;n}^N$ as $N\to\iy$ (with respect to the maps (b)).

\subhead 22.8\endsubhead
Let $\tZ=\{\tz_1\le\tz_2\le\do\le\tz_{2N+r}\}\in\cm_{a,b;n}^N$.
Let $t$ be an integer which is large enough so that the multiset

(a) $\{at+b'-\tz_1,at+b'-\tz_2,\do,at+b'-\tz_{2N+r}\}$
\nl
is contained in the multiset

(b) $\{0,a,2a,\do,ta, b',a+b',2a+b',\do,ta+b'\}$
\nl
and let $\btZ$ be the complement of (a) in (b). Then $\btZ\in\cm_{a,b}^{t+1-N-r}$.
The sum of entries of $\btZ$ is
$$\align&\su_{k\in[0,t]}(2ka+b')-(at+b')(2N+r)+\su_h\tz_h\\&=
at(t+1)+(t+1)b'-(at+b')(2N+r)+an+aN^2+N(b-a)+a\bin{r}{2}+b'r\\&
=an+a(t+1-N-r)^2+(t+1-N-r)(b-a)+a\bin{r}{2}+b'r.\endalign$$
Thus, $\btZ$ has rank $n$.

We define a bijection $\p_N\i(\tZ)@>\si>>\p_{t+1-N-r}\i(\btZ)$ by $\L\m\ov\L$ where
$\L$ is as in 22.7(a) and $\ov\L$ is
$$\align&\{b',a+b',2a+b',3a+b',\do,ta+b'\}\\&
-\{at+b'-\mu_1,at+b'-\mu_2,\do,at+b'-\mu_N\}\\&\{0,a,2a,3a,\do,ta\}
-\{at+b'-\l_1,at+b'-\l_2,\do,at+b'-\l_{N+r}\}.\endalign$$ 

\subhead 22.9\endsubhead
Let $W=W_n$ be the group of permutations of $1,2,\do,n,n',\do,2',1'$ which commute 
with the involution $i\m i', i'\m i$. We regard $W_n$ as a Coxeter group with 
generators $s_1,s_2,\do,s_n$ given as products of transpositions by 

$s_1=(1,2)(1',2'),s_2=(2,3)(2',3'),\do,s_{n-1}=(n-1,n)((n-1)',n')$,

$s_n=(n,n')$.

\subhead 22.10\endsubhead
A permutation in $W$ defines a permutation of the $n$ element set consisting of the
pairs $(1,1'),(2,2'),\do,(n,n')$. Thus we have a natural homomorphism of $W_n$ onto
$\fS_n$, the symmetric group in $n$ letters. Define a homomorphism 
$\c_n:W_n@>>>\pm 1$ by

$\c_n(\s)=1$ if $\{\s(1),\s(2),\do,\s(n)\}\cap\{1',2',\do,n'\}$ has even cardinal,

$\c_n(\s)=-1$, otherwise.
\nl
The simple $W$-modules (up to isomorphism) are in 1-1 correspondence with the 
ordered pairs $\a,\b$ where $\a=(\a_1\ge\a_2\ge\do)$ and $\b=(\b_1\ge\b_2\ge\do)$ 
are partitions such that $\a_N=\b_N=0$ for large $N$ and $\su_i\a_i+\su_j\b_j=n$.
The correspondence (denoted by $\a,\b\m E^{\a,\b}$) is defined as follows. Let 
$\a,\b$ be as above, let $(\a'_1\ge\a'_2\ge\do)$ be the partition dual to $\a$ and 
let $(\b'_1\ge\b'_2\ge\do)$ be the partition dual to $\b$. Let 
$k=\su_i\a_i,l=\su_j\b_j$. Let $\p_\a$ be the simple $\fS_k$-module defined as in
22.4 and let $\p_\b$ be the analogously defined simple $\fS_l$-module. We regard 
$\p_\a,\p_\b$ as simple modules of $W_k,W_l$ via the natural homomorphisms 
$W_k@>>>\fS_k,W_l@>>>\fS_l$ as above. We identify $W_k\T W_l$ with the subgroup of 
$W$ consisting of all permutations in $W$ which map $1,2,\do,k,k',\do,2',1'$ into 
itself hence also map $k+1,k+2,\do,n,n',\do,(k+2)',(k+1)'$ into itself. Consider the
representation $\p_\a\ot(\p_\b\ot\c_l)$ of $W_k\T W_l$. We induce it to $W$; the
resulting representation of $W$ is irreducible; we denote it by $E^{\a,\b}$. 

We fix $a>0,b\ge 0$ and we write $b=ar+b'$ as in 22.6.

Let $\a,\b$ be as above. Let $N$ be an integer such that $\a_{N+r+1}=0$,
$\b_{N+1}=0$. (Any large enough integer satisfies these conditions.) We set
$$\l_i=a(\a_{N+r-i+1}+i-1)+b',(i\in[1,N+r]),\qua\mu_j=a(\b_{N-j+1}+j-1),(j\in[1,N]).
$$    
We have $0\le\l_1<\l_2<\do<\l_{N+r}$, $0\le\mu_1<\mu_2<\do<\mu_N$. Let $\L$ denote 
the tableau 22.7(a). It is easy to see that $\L\in\Sy_{a,b;n}^N$. Moreover, if $N$ 
is replaced by $N+1$, then $\L$ is replaced by its image under
$\Sy_{a,b;n}^N@>>>\Sy_{a,b;n}^{N+1}$ (see 22.7). Let $[\L]=E^{\a,\b}$. Note that 
$[\L]$ depends only on the image of $\L$ under the canonical map 
$\Sy_{a,b;n}^N@>>>\Sy_{a,b;n}$. In this way, we see that 

{\it the simple $W$-modules are naturally in bijection with the set $\Sy_{a,b;n}$.}
\nl
For $i\in[1,N]$ we have $a(\a_{N-i+1}+i-1)+b=a(\a_{N+r-i-r+1}+i+r-1)+b'=\l_{i+r}$. 

If $N$ is large we have $\l_i=a(i-1)+b'$ for $i\in[1,r]$ and $\mu_j=a(j-1)$ for 
$j\in[1,r]$.

\subhead 22.11\endsubhead
Let $q,y$ be indeterminates. With the notation in 22.10, let
$$H_\a(q)=q^{-\su_i\bin{\a'_i}{2}}\prod_{i,j}\fra{q^{\a_i+\a'_j-i-j+1}-1}{q-1},$$
$$G_{\a,\b}(q,y)=q^{-\su_i\a'_i\b'_i/2}\prod_{i,j}(q^{\a_i+\b'_j-i-j+1}y+1);$$
both products are taken over all $i\ge 1,j\ge 1$ such that $\a_i\ge j,\a'_j\ge i$. 

Define a weight function $L:W@>>>\bn$ by $L(s_1)=L(s_2)=\do=L(s_{n-1})=a$, 
$L(s_n)=b$. We now assume that both $a,b$ are $>0$. We also assume that $a,b$ are 
such that $W,L$ satisfies the assumptions of 18.1. Then $f_{E^{\a,\b}_v}$ is defined
in terms of this $L$.

\proclaim{Lemma 22.12 (Hoefsmit \cite{\HO})} We have
$$f_{E^{\a,\b}_v}=H_\a(v^{2a})H_\b(v^{2a})G_{\a,\b}(v^{2a},v^{2b})G_{\b,\a}
(v^{2a},v^{-2b}).$$
\endproclaim

We will rewrite the expression above using the following result.

\proclaim{Lemma 22.13} Let $N$ be a large integer. We have
$$H_\a(q)=q^{\su_{i\in[1,N-1]}\bin{i}{2}}
\fra{\prod_{i=1}^N \prod_{h\in[1,\a_{N-i+1}+i-1]}\fra{q^h-1}{q-1}}
{\prod_{1\le i<j\le N}\fra{q^{\a_{N-j+1}+j-1}-q^{\a_{N-i+1}+i-1}}{q-1}},$$
$$\align&G_{\a,\b}(q,y)G_{\b,\a}(q,y\i)=q^{\su_{i\in[1,N-1]}i^2}
(\sqrt y+\sqrt y\i)^N\\&\T\fra{\prod_{i=1}^N\prod_{h\in[1,\a_{N-i+1}+i-1]}(q^hy+1)
\prod_{j=1}^N\prod_{h\in[1,\b_{N-j+1}+j-1]}(q^hy\i+1)}{\prod_{i,j\in[1,N]}
(q^{\a_{N-i+1}+i-1}\sqrt y+q^{\b_{N-j+1}+j-1}\sqrt y\i)}.\endalign$$
\endproclaim
The proof is by induction on $n$. We omit it.

\proclaim{Proposition 22.14} (a) If $b'=0$ then $f_{[\L]_\sp}$ is equal to $2^d$ 
where $2d+r$ is the number of singles in $\L$. If $b'>0$ then $f_{[\L]_\sp}=1$.

(b) We have $\aa_{[\L]}=A_N-B_N$ where
$$A_N=\su_{i\in[1,N+r],j\in[1,N]}\min(\l_i,\mu_j)
+\su_{1\le i<j\le N+r}\min(\l_i,\l_j)+\su_{1\le i<j\le N}\min(\mu_i,\mu_j),$$
$$\align&B_N=\su_{i\in[1,N+r],j\in[1,N]}\min(a(i-1)+b',a(j-1))\\&
+\su_{1\le i<j\le N+r}\min(a(i-1)+b',a(j-1)+b')
+\su_{1\le i<j\le N}\min(a(i-1),a(j-1)).\endalign$$
\endproclaim
It is enough to prove (a) assuming that $N$ is large. Since
$$A_{N+1}-A_N=a(N+r)N+a\bin{N}{2}+a\bin{N+r}{2}+b'N+b'(N+r)=B_{N+1}-B_N,$$
we have $A_{N+1}-B_{N+1}=A_N-B_N$ hence it is enough to prove (b) assuming that $N$
is large. In the remainder of the proof we assume that $N$ is large.

For $f,f'\in\bc(v)$ we write $f\cong f'$ if $f'=fg$ with $g\in\bc(v)$, $g|_{v=0}=1$.
Using 22.12, 22.13, we see that
$$\align&f_{[\L]_v}\cong\prod_{i\in[1,N]}(v^{2a-2b}+1)(v^{4a-2b}+1)\do
(v^{2\mu_i-2b}+1)\\&(v^b+v^{-b})^Nv^{2a\su_{i=1}^{N-1}
(2i^2-i)}\prod_{i,j\in[1,N]}(v^{2\l_{i+r}-b}+v^{2\mu_j-b})\i\\&\prod_{1\le i<j\le N}
(v^{2\l_{j+r}-2b}-v^{2\l_{i+r}-2b})\i\prod_{1\le i<j\le N}
(v^{2\mu_j}-v^{2\mu_i})\i\endalign$$
hence 
$$f_{[\L]_v}=2^dv^{-K}+\text{strictly higher powers of $v$}$$
where $d=0$ if $b'>0$,
$$\align&d=\sh(j\in[1,N]:b\le\mu_j)-\sh(i,j\in[1,N]:\l_{i+r}=\mu_j)\\&
=N-\sh(i\in[1,r],j\in[1,N]:(i-1)a=\mu_j)-\sh(i,j\in[1,N]:\l_{i+r}=\mu_j)
\\&=N-\sh(i\in[1,r],j\in[1,N]:\l_i=\mu_j)-\sh(i,j\in[1,N]:\l_{i+r}=\mu_j)
\\&=N-\sh(i\in[1,N+r],j\in[1,N]:\l_i=\mu_j)=(\sh \text{ singles}-r)/2,\endalign$$
if $b'=0$,
$$\align&-K=-bN+2a\su_{i\in[1,N-1]}(2i^2-i)
+\su_{j\in[1,N]}\su\Sb k\in [1,r]\\ak\le\mu_j\eSb(2ak-2b)\\&
-\su_{i,j\in[1,N]}(-b+2\min(\l_{i+r},\mu_j))
-\su_{1\le i<j\le N}(-2b+2\min(\l_{i+r},\l_{j+r}))\\&
-\su_{1\le i<j\le N}2\min(\mu_i,\mu_j)=-bN+2a\su_{i\in[1,N-1]}(2i^2-i)+2bN^2-bN\\&
+\su_{j\in[1,N]}\su_{k\in[1,r],ak\le\mu_j}(2ak-2b)
-\su_{i,j\in[1,N]}2\min(\l_{i+r},\mu_j)\\&
-\su_{1\le i<j\le N}2\min(\l_{i+r},\l_{j+r})
-\su_{1\le i<j\le N}2\min(\mu_i,\mu_j)\\&
=\su_{j\in[1,N]}\su_{k\in[1,r],ak\le\mu_j}(2ak-2b)
-\su_{i,j\in[1,N]}2\min(\l_{i+r},\mu_j)\\&
-\su_{1\le i<j\le N}2\min(\l_{i+r},\l_{j+r})
-\su_{1\le i<j\le N}2\min(\mu_i,\mu_j)+\bst.\endalign$$
(We will generally write $\bst$ for an expression which depends only on $a,b,N$.) We
have
$$\align&\su\Sb j\in[1,N]\\k\in[1,r]\\ak\le\mu_j\eSb(2ak-2b)=
\su\Sb j\in[1,r]\\k\in[1,r]\\ak\le\mu_j\eSb(2ak-2b)
+\su\Sb j\in[r+1,N]\\k\in[1,r]\\ak\le\mu_j\eSb(2ak-2b)\\&
=\su\Sb j\in[1,r]\\k\in[1,r]\\ak\le a(j-1)\eSb(2ak-2b)
+\su\Sb j\in[r+1,N]\\k\in[1,r]\eSb(2ak-2b)=\bst,\endalign$$
hence
$$\align&-K\\&=-2(\su_{i,j\in[1,N]}\min(\l_{i+r},\mu_j)-\su_{1\le i<j\le N}
\min(\l_{i+r},\l_{j+r})-\su_{1\le i<j\le N}\min(\mu_i,\mu_j))\\&+\bst.\endalign$$
We have
$$\align&\su_{i\in[1,r],j\in[1,N]}\min(\l_i,\mu_j)
=\su_{i\in[1,r],j\in[1,N]}\min(a(i-1)+b',\mu_j)\\&
=\su_{i\in[1,r],j\in[1,r]}\min(a(i-1)+b',a(j-1))
+\su_{i\in[1,r],j\in[r+1,N]}\min(a(i-1)+b',\mu_j)\\&
=\su_{i\in[1,r],j\in[1,r]}\min(a(i-1)+b',a(j-1))
+\su_{i\in[1,r],j\in[r+1,N]}(a(i-1)+b')=\bst,\endalign$$
hence
$$\su_{i,j\in[1,N]}\min(\l_{i+r},\mu_j)=\su_{i\in[1,N+r],j\in[1,N]}\min(\l_i,\mu_j)
+\bst.$$
We have
$$\align&\su_{1\le i<j\le N+r}\min(\l_i,\l_j)=
\su_{1\le i<j\le N}\min(\l_{i+r},\l_{j+r})+\su_{i\in[1,r]}\l_i(N+r-i)\\&
=\su_{1\le i<j\le N}\min(\l_{i+r},\l_{j+r})+\su_{i\in[1,r]}(a(i-1)+b')(N+r-i)\\&
=\su_{1\le i<j\le N}\min(\l_{i+r},\l_{j+r})+\bst.\endalign$$
We see that
$$-K=-2A_N+\bst.\tag c$$
In the special case where $\a=\b=(0\ge 0\ge\do)$ we have $K=0$. On the other hand, 
by (c), we have $0=-2B_N+\bst$ where $\bst$ is as in (c). Hence in general we have 
$-K=-2A_N+2B_N$. This proves the proposition, in view of 20.11 and 20.21(a).

\subhead 22.15\endsubhead
We identify $\fS_k\T W_l$ ($k+l=n)$ with the subgroup of $W$ consisting of all 
permutations in $W$ which map $\{1,2,\do,k\}$ into itself (hence also map
$\{1',2',\do,k'\}$ and $\{k+1,\do,n,n',\do,(k+1)'\}$ into themselves. This is a 
standard parabolic subgroup of $W$. We consider an irreducible representation of 
$\fS_k\T W_l$ of the form $\sg_k\bxt[\L']$ where $\sg_k$ is the sign representation
of $\fS_k$ and $\L'\in\Sy_{a,b;l}^N$. We may assume that $\L'$ has at least $k$ 
entries. We want to associate to $\L'$ a symbol in $\Sy_{a,b;n}^N$ by increasing 
each of the $k$ largest entries in $\L'$ by $a$. It may happen that the set of $r$
largest entries of $\L'$ is not uniquely defined but there are two choices for it.
(This can only happen if $b'=0$.) Then the same procedure gives rise to two distinct
symbols $\L^I,\L^{II}$ in $\Sy_{a,b;n}^N$.

\proclaim{Lemma 22.16} (a) $g_{(\sg_k\ot[\L'])_\sp}=g_{[\L']_\sp}$ is equal to
$g_{[\L]_\sp}$ or to $g_{[\L^I]_\sp}+g_{[\L^{II}]_\sp}$.

(b) $\aa_{\sg_k\ot[\L']}=a\bin{k}{2}+\aa_{[\L']}$ is equal to $\aa_{[\L]}$ or to
$\aa_{[\L^I]}=\aa_{[\L^{II}]}$.
\endproclaim
$\L$, if defined, has the same number of singles as $\L'$. Moreover, $\L^I$ (and
$\L^{II}$), if defined, has one more single than $\L'$. Hence (a) follows from 
22.14(a) using 20.18, 20.19.

By 22.14(b), the difference $\aa_{[\ti\L]}-\aa_{[\L]}$ (where $\ti\L$ is either $\L$
or $\L^I$ or $\L^{II}$) is $a$ times the number of $i<j$ in $[1,k]$. Thus, it is 
$a\bin{k}{2}$. Hence (a) follows from 20.18, 20.19. The lemma is proved.

\proclaim{Lemma 22.17} $\jj_{\fS_k\T W_l}^W(\sg_k\ot[\L'])$ equals $[\L]$ or 
$[\L^I]+[\L^{II}]$.
\endproclaim
By a direct computation (involving representations of symmetric groups) we see that:

(a) if $\L$ is defined then $[[\L']:[\L]]\ge 1$;

(b) if $\L^I,\L^{II}$ are defined then $[[\L']:[\L^I]]\ge 1$ and
$[[\L']:[\L^{II}]]\ge 1$.
\nl
In the setup of (a) we have (by 20.14(b)):

$g_{[\L']_\sp}=\su_{E;\aa_E=\aa_{E'}}[[\L']:E]g_{E_\sp}$ hence using 22.16(a) we
have

$$g_{[\L]_\sp}=\su_{E;\aa_E=\aa_{E'}}[[\L']:E]g_{E_\sp}.\tag c$$
By 22.16(b), $E=[\L]$ enters in the last sum and its contribution is
$\ge g_{[\L]_\sp}$; the contribution of the other $E$ is $\ge 0$ (see 20.13(b)). 
Hence (c) forces $[[\L']:[\L]]=1$ and $[[\L']:E]=0$ for all other $E$ in the sum. In
this case the lemma follows.

In the setup of (b) we have (by 20.14(b)) 
$g_{[\L']_\sp}=\su_{E;\aa_E=\aa_{E'}}[[\L']:E]g_{E_\sp}$ hence, by 22.16(a),
$$g_{[\L^I]_\sp}+g_{[\L^{II}]_\sp}=\su_{E;\aa_E=\aa_{E'}}[[\L']:E]g_{E_\sp}.\tag d$$
By 22.16(b), $E=[\L^I]$ and $E=[\L^{II}]$ enter in the last sum and their
contribution is $\ge g_{[\L^I]_\sp}+g_{[\L^{II}]_\sp}$; the contribution of the
other $E$ is $\ge 0$ (see 20.13(b)). Hence (d) forces 
$[[\L']:[\L^I]]=[[\L']:[\L^{II}]]=1$ and $[[\L']:E]=0$ for all other $E$ in the sum.
The lemma follows.

\proclaim{Lemma 22.18}$[\L]\ot\sg=[\ov\L]$. (Notation of 22.14.)
\endproclaim
This can be reduced to a known statement about the symmetric group. We omit the
details.

\subhead 22.19\endsubhead
Let $Z$ be a totally ordered finite set $z_1<z_2<\do<z_M$. For any $r\in[0,M]$ such
that $r=M\mod 2$ let $\uZ_r$ be the set of subsets of $Z$ of cardinal $(M-r)/2$. An
involution $\io:Z@>>>Z$ is said to be $r$-admissible if the following hold:

(a) $\io$ has exactly $r$ fixed points;

(b) if $M=r$, there is no further condition; if $M>r$, there exist two consecutive 
elements $z,z'$ of $Z$ such that $\io(z)=z',\io(z')=z$ and the induced involution of
$Z-\{z,z'\}$ is $r$-admissible.
\nl
Let $\In_r(Z)$ be the set of $r$-admissible involutions of $Z$. To $\io\in\In_r(Z)$
we associate a subset $\cs_\io$ of $\uZ_r$ as follows: a subset $Y\sub Z$ is in 
$\cs_\io$ if it contains exactly one element in each non-trivial $\io$-orbit. 
Clearly, $\sh(\cs_\io)=2^{p_0}$ where $p_0=(M-r)/2$. 
(In fact, $\cs_\io$ is naturally an affine space over the field $\bF_2$.)

\proclaim{Lemma 22.20} Assume that $p_0>0$. Let $Y\in\uZ_r$.

(a) We can find two consecutive elements $z,z'$ of $Z$ such that exactly one of
$z,z'$ is in $Y$.

(b) There exists $\io\in\In_r(Z)$ such that $Y\in\cs_\io$.

(c) Assume that for some $k\in[0,p_0-1]$, $z_1,z_2,\do,z_k$ belong to $Y$ but 
$z_{k+1}\n Y$. Let $l$ be the smallest number such that $l>k$ and $z_l\in Y$. There
exists $\io\in\In_r(Z)$ such that $Y\in\cs_\io$ and $\io(z_l)=z_{l-1}$.
\endproclaim
We prove (a). Let $z_k$ be the smallest element of $Y$. If $k>1$ then we can take 
$(z,z')=(z_{k-1},z_k)$. Hence we may assume that $z_1\in Y$. Let $z_{k'}$ be the 
next smallest element of $Y$. If $k'>2$ then we can take $(z,z')=(z_{k'-1},z_{k'})$.
Continuing like this we see that we may assume that $Y=\{z_1,z_2,\do,z_{p_0}\}$. 
Since $p_0<M$, we may take $(z,z')=(z_{p_0},z_{p_0+1})$. 

We prove (b). Let $z,z'$ be as in (a). Let $Z'=Z-\{z,z'\}$ with the induced order.
Let $Y'=Y\cap Z'$. If $p_0\ge 2$ then by induction on $p_0$ we may assume that there
exists $\io'\in\In_r(Z')$ such that $Y'\in\cs_{\io'}$. Extend $\io'$ to an 
involution $\io$ of $Z$ by $z\m z',z'\m z$. Then $\io\in\In_r(Z)$ and $Y\in\cs_\io$.
If $p_0=1$, 
define $\io:Z@>>>Z$ so that $z\m z',z'\m z$ and $\io=1$ on $Z-\{z,z'\}$. 
Then $\io\in\In_r(Z)$ and $Y\in\cs_\io$.

We prove (c). We have $l\ge k+2$. Hence $z_{l-1}\n Y$. Let $(z,z')=(z_{l-1},z_l)$.
We continue as in the proof of (b), except that instead of invoking an induction 
hypothesis, we invoke (b) itself.

\subhead 22.21\endsubhead
Assume that $M>r$. We consider the graph whose set of vertices is $\uZ_r$ and in
which two vertices $Y\ne Y'$ are joined if there exists $\io\in\In_r(Z)$ such that
$Y\in\cs_\io,Y'\in\cs_\io$.

\proclaim{Lemma 22.22} This graph is connected.
\endproclaim
We show that any vertex $Y=\{z_{i_1},z_{i_2},\do,z_{i_{p_0}}\}$ 
is in the same connected
component as $Y_0=\{z_1,z_2,\do,z_{p_0}\}$. We argue by induction on 
$m_Y=i_1+i_2+\do+i_{p_0}$. If $m_Y=1+2+\do+p_0$ then 
$Y=Y_0$ and there is nothing to 
prove. Assume now that $m>1+2+\do+p_0$ so that $Y\ne Y_0$. Then the assumption of 
Lemma 22.20(c) is satisfied. Hence we can find $l$ such that $z_l\in Y,z_{l-1}\n Y$
and $\io\in\In_r(Z)$ such that $Y\in\cs_\io$ and $\io(z_l)=z_{l-1}$. Let 
$Y'=(Y-\{z_l\})\cup\{z_{l-1}\}$. Then $Y'\in\cs_\io$ hence $Y,Y'$ are joined in our
graph. We have $m_{Y'}=m_Y-1$ hence by the induction hypothesis $Y',Y_0$ are in the
same connected component. It follows that $Y,Y_0$ are in the same connected 
component. The lemma is proved.

\subhead 22.23\endsubhead
Assume that $b'=0$. Let $\tZ\in\cm_{a,b;n}^N$. Let $Z$ be the set of singles of 
$\tZ$. Each set $Y\in\uZ_r$ gives rise to a symbol $\L_Y$ in $\p_N\i(\tZ)$: the
first row of $\L_Y$ consists of $Z-Y$ and one element in each double of $\tZ$; the
second row consists of $Y$ and one element in each double of $\tZ$. For any 
$\io\in\In_r(Z)$ we set
$$c(\tZ,\io)=\op_{Y\in\cs_\io}[\L_Y]\in\Mod W.$$ 

\proclaim{Proposition 22.24} (a) In the setup of 22.23, let $\io\in\In_r(Z)$. Then
$c(\tZ,\io)\in\C(W)$.

(b) All constructible representations of $W$ are obtained as in (a).
\endproclaim
We prove (a) by induction on $n$. If $n=0$ the result is clear. Assume that
$n\ge 1$. We may assume that $0$ is not a double of $\tZ$. Let $at$ be the largest
entry of $\tZ$.

(A) Assume that there exists $i,0\le i<t$, such that $ai$ does not appear in $\tZ$.
Then $\tZ$ is obtained from $\tZ'\in\cm_{a,b;n-k}^N$ with $n-k<n$ by increasing each
of the $k$ largest entries by $a$ and this set of largest entries is unambiguously 
defined. The set $Z'$ of singles of $\tZ'$ is naturally in order preserving 
bijection with $Z$. Let $\io'$ correspond to $\io$ under this bijection. By the 
induction hypothesis, $c(\tZ',\io')\in\C(W_{n-k})$. Since, by 22.5, the sign 
representation $\sg_k$ of $\fS_k$ is constructible, it follows that
$\sg_k\bxt c(\tZ',\io')\in\C(\fS_k\T W_{n-k})$. Using 22.17, we have 
$$\jj_{\fS_k\T W_{n-k}}^W(\sg_k\bxt c(\tZ',\io'))=c(\tZ,\io)$$
hence $c(\tZ,\io)\in\C(W)$.

(B) Assume that there exists $i, 0<i\le t$ such that $ai$ is a double of $\tZ$. Let
$\btZ$ be as in 22.8 (with respect to our $t$). Then $0$ is not a double of $\btZ$ 
and the largest entry of $\btZ$ is $at$. Let $\bZ$ be the set of singles of $\btZ$.
We have $\bZ=\{at-z;z\in Z\}$. Thus $\bZ,Z$ are naturally in (order reversing) 
bijection under $j\m at-j$. Let $\io'\in\In_r(\bZ)$ correspond to $\io$ under this
bijection. Since $at-ai$ does not appear in $\btZ$, (A) is applicable to $\btZ$. 
Hence $c(\btZ,\io')\in\C(W)$. By 22.18 we have $c(\btZ,\io')\ot\sg=c(\tZ,\io)$ hence
$c(\tZ,\io)\in\C(W)$. 

(C) Assume that we are not in case (A) and not in case (B). Then 
$\tZ=\{0,a,2a,\do,ta\}=Z$. We can find $ia,(i+1)a$ in $Z$ such that $\io$ 
interchanges $ia,(i+1)a$ and induces on $Z-\{ia,(i+1)a\}$ an $r$-admissible
involution $\io_1$. We have
$$\tZ'=\{0,a,2a,\do,ia,ia,(i+1)a,(i+2)a,\do,(t-1)a\}\in\cm_{a,b;n-k}^N$$
with $n-k<n$. The set of singles of $\tZ'$ is 
$$Z'=\{0,a,2a,\do,(i-1)a,(i+1)a,\do,(t-1)a\}.$$
It is in natural (order preserving) bijection with $Z-\{ia,(i+1)a\}$. Hence $\io_1$
induces $\io'\in\In_r(Z')$. By the induction hypothesis we have
$c(\tZ',\io')\in\C(W_{n-k})$. Hence $\sg_k\bxt c(\tZ',\io')\in\C(\fS_k\T W_{n-k})$ 
where $\sg_k$ is as in (A). By 22.17 we have 
$$\jj_{\fS_k\T W_{n-k}}^W(\sg_k\bxt c(\tZ',\io'))=c(\tZ,\io)$$
hence $c(\tZ,\io)\in\C(W)$. This proves (a).

We prove (b) by induction on $n$. If $n=0$ the result is clear. Assume now that
$n\ge 1$. By an argument like the ones used in (B) we see that the class of 
representations of $W$ obtained in (a) is closed under $\ot\sg$. Therefore, to show
that $C\in\C(W)$ is obtained in (a), we may assume that
$C=\jj_{\fS_k\T W_{n-k}}^W(C')$ for some $k>0$ and some $C'\in\C(\fS_k\T W_{n-k})$.
By 22.5 we have $C'=E\bxt C''$ where $E$ is a simple $\fS_k$-module and 
$C''\in\C(W_{n-k})$. Using 22.5(a) we have
$E=\jj_{\fS_{k'}\T\fS_{k''}}^{\fS_k}(\sg\bxt E')$ where $k'+k''=k,k'>0$ 
and $E'$ is a simple $\fS_{k''}$-module. Let 
$\tC=\jj_{\fS_{k''}\T W_{n-k}}^{W_{n-k'}}(E'\ot C')\in\C(W_{n-k'})$. Then 
$C=\jj_{\fS_{k'}\T W_{n-k'}}^W(\sg_{k'}\ot\tC)$. By the induction hypothesis, $\tC$
is of the form described in (a). Using an argument as in (A) or (C) we deduce that 
$C$ is of the form described in (a). The proposition is proved.

\proclaim{Proposition 22.25}Assume that $b'>0$. 

(a) Let $E\in\Ir W$. Then $E\in\C(W)$.

(b) All constructible representations of $W$ are obtained as in (a).
\endproclaim
We prove (a). We may assume that $E=[\L]$ where $\L\in\Sy_{a,b;n}^N$ does not
contain both $0$ and $b'$. We argue by induction on $n$. If $n=0$ the result is
clear. Assume now that $n\ge 1$.

(A) Assume that either (1) there exist two entries $z,z'$ of $\L$ such that $z'-z>a$
and there is no entry $z''$ 
of $\L$ such that $z<z''<z'$, or (2) there exists an entry
$z'$ of $\L$ such that $z'\ge a$ and there is no entry $z''$ of $\L$ such that
$z''<z'$. Let $\L'$ be the symbol obtained from $\L$ by substracting $a$ from each 
entry $\tz$ of $\L$ such that $\tz\ge z'$ and leaving the other entries of $\L$ 
unchanged. Then $\L'\in\Sy_{a,b;n-k}^N$ with $n-k<n$. By the induction hypothesis, 
$[\L']\in\C(W_{n-k})$. Since, by 22.5, the sign representation $\sg_k$ of $\fS_k$ is
constructible, it follows that $\sg_k\bxt[\L']\in\C(\fS_k\T W_{n-k})$. Using 22.17,
we have $\jj_{\fS_k\T W_{n-k}}^W(\sg_k\bxt[\L'])=[\L]$ hence $[\L]\in\C(W)$.

(B) Assume that there exist two entries $z,z'$ of $\L$ such that $0<z'-z<a$. Let $t$
be the smallest integer such that $at+b'\ge\l_i$ for all $i\in[1,N+r]$ and 
$at\ge\mu_j$ for all $j\in[1,N]$. Let $\ov\L\in\Sy_{a,b;n}^{t+1-N-r}$ be as in 22.8
with respect to this $t$. Then $\ov\L$ does not contain both $0$ and $b'$. Now (A)
is applicable to $\ov\L$. Hence $[\ov\L]\in\C(W)$. By 22.18 we have 
$[\ov\L]\ot\sg=[\L]$ hence $[\L]\in\C(W)$. 

(C) Assume that we are not in case (A) and not in case (B). Then the entries of $\L$
are either $0,a,2a,\do,ta$ or $b',a+b',2a+b',\do,ta+b'$. This cannot happen for
$n\ge 1$. This proves (a).

The proof of (b) is entirely similar to that of 22.24(b). The proposition is proved.

\subhead 22.26\endsubhead 
We now assume that $n\ge 2$ and that $W'=W'_n$ is the kernel of $\c_n:W_n@>>>\pm 1$
in 22.10. We regard $W'_n$ as a Coxeter group with generators $s_1,s_2,\do,s_{n-1}$
as in 22.9 and $s'_n=(n-1,n')((n-1)',n)$ (product of transpositions). Let 
$L:W'@>>>\bn$ be the weight function given by $L(w)=al(w)$ for all $w$. Here $a>0$.

For $\L\in\Sy_{a,0}^N$ we denote by $\L^{tr}$ the symbol whose first (resp. second)
row is the second (resp. first) row of $\L$. We then have $\L^{tr}\in\Sy_{a,0}^N$. 
From the definitions we see that the simple $W_n$-modules $[\L],[\L^{tr}]$ have the
same restriction to $W'$; this restriction is a simple $W'$-module $[\un\L]$ if 
$\L\ne\L^{\tr}$ and is a direct sum of two non-isomorphic simple $W'$-modules 
$[{}^I\un\L]$,$[{}^{II}\un\L]$ if $\L=\L^{\tr}$. In this way we see that

{\it the simple $W'$-modules are naturally in bijection with the set of orbits of 
the involution of $\Sy_{a,0;n}$ induced by $\L\m\L^{tr}$ except that each fixed 
point of this involution corresponds to two simple $W'$-modules.}

Let $\tZ\in\cm_{a,0;n}^N$. Let $Z$ be the set of singles of $\tZ$. Assume first that
$Z\ne\em$. Each set $Y\in\uZ_0$ gives rise to a symbol $\L_Y$ in $\Sy_{a,0;n}^N$: 
the first row of $\L_Y$ consists of $Z-Y$ and one element in each double of $\tZ$; 
the second row consists of $Y$ and one element in each double of $\tZ$. For any 
$\io\in\In_0(Z)$ we define $c(\tZ,\io)\in\Mod W$ by 
$$c(\tZ,\io)\op c(\tZ,\io)=\op_{Y\in\cs_\io}[\L_Y]\in\Mod W.$$
Note that $Y$ and $Z-Y$ have the same contribution to the sum. A proof entirely
similar to that of 22.24 shows that $c(\tZ,\io)\in\C(W)$. Moreover, if $Z=\em$ and
$\L=\L^{tr}\in\Sy_{a,1;n}^N$ is defined by $\p_N(\L)=\tZ$, then 
$[{}^I\un\L]\in\C(W)$ and $[{}^{II}\un\L]\in\C(W)$. All constructible 
representations of $W$ are obtained in this way.
 
\subhead 22.27\endsubhead 
Assume that $W$ is of type $F_4$ and that the values of $L:W@>>>\bn$ on $S$ are
$a,a,b,b$ where $a>b>0$.

{\it Case 1}. Assume that $a=2b$. There are four simple $W$-modules
$\r_1,\r_2,\r_8,\r_9$ (subscript equals dimension) with $\aa=3b$. Then
$$\r_1\op\r_8,\r_2\op\r_9,\r_8\op\r_9\in\C(W).$$
(They are obtained by $\jj$ from the $W_I$ of type $B_3$ with parameters $a,b,b$.) 

The simple $W$-modules $\r_1^\da,\r_2^\da,\r_8^\da,\r_9^\da$ have $\aa=15b$ and
$$\r_1^\da\op\r_8^\da,\r_2^\da\op\r_9^\da,\r_8^\da\op\r_9^\da\in
\C(W).$$
There are five simple $W$-modules $\r_{12},\r_{16},\r_6,\r'_6,\r_4$ (subscript
equals dimension) with $\aa=7b$. Then
$$\r_4\op\r_{16},\r_{12}\op\r_{16}\op\r_6,\r_{12}\op\r_{16}\op\r'_6\in\C(W).$$
All $12$ simple $W$-modules other than the $13$ listed above, are constructible. All
constructible representations of $W$ are thus obtained.

{\it Case 2}. Assume that $a\n\{b,2b\}$. The simple $W$-modules 
$\r_{12},\r_{16},\r_6,\r'_6,\r_4$ in Case 1 now have $\aa=3a+b$ and
$$\r_4\op\r_{16},\r_{12}\op\r_{16}\op\r_6,\r_{12}\op\r_{16}\op\r'_6\in\C(W).$$
All $20$ simple $W$-modules other than the $5$ listed above, are constructible. All
constructible representations of $W$ are thus obtained.

\subhead 22.28\endsubhead 
Assume that $W$ is of type $G_2$ and that the values of $L:W@>>>\bn$ on $S$ are
$a,b$ where $a>b>0$. Let $\r_2,\r'_2$ be the two $2$-dimensional simple $W$-modules.
They have $\aa=a$ and $\r_2\op\r'_2$ is constructible. All $4$ simple $W$-modules 
other than the $2$ listed above, are constructible. All constructible 
representations of $W$ are thus obtained.

\subhead 22.29\endsubhead 
Let $\cl$ be the set of all weight functions $L:W@>>>\bn$ such that $L(s)>0$ for all
$s\in S$. We assume that $P1-P15$ in \S14 hold for any $L\in\cl$. For $L,L'\in\cl$ 
we write $L\si L'$ if the constructible representations of $W$ with respect to $L$ 
are the same as those with respect to $L'$. This is an equivalence relation on 
$\cl$. From the results in this chapter we see that any equivalence class for $\si$
contains some $L$ which is attached to some $(G,F,\cp,\be)$ as in 0.3.

We expect that the constructible representations of $W$ are exactly the 
representations of $W$ carried by the left cells of $W$ (for fixed $L\in\cl$). (For
$L=l$ this holds by \cite{\CG}. For $W$ of type $F_4$ and general $L$ this holds by
\cite{\GE}.) This would imply that for $L,L'\in\cl$ we have $L\si L'$ if and only if
the representations of $W$ carried by the left cells of $W$ with respect to $L$ are
the same as those with respect to $L'$.

\head 23. Two-sided cells\endhead
\subhead 23.1\endsubhead
We preserve the setup of 20.1.
We define a graph $\cg_W$ as follows. The vertices of $\cg_W$ are the simple
$W$-modules up to isomorphism. Two non-isomorphic simple $W$-modules are joined in 
$\cg_W$ if they both appear as components of some constructible representation of 
$W$. Let $\fra{\cg_W}{\si}$ 
be the set of connected components of $\cg_W$. The connected
components of $\cg_W$ are determined explicitly by the results in \S22 for $W$ 
irreducible.

For example, in the setup of 22.4,22.5 we have $\cg_W=\fra{\cg_W}{\si}$. In 
the setup of 22.24, $\fra{\cg_W}{\si}$ 
is naturally in bijection with $\cm_{a,b;n}$. (Here, 22.22 is 
used). In the setup of 22.25, we have $\cg_W=\fra{\cg_W}{\si}$.

We show that:

(a) {\it if $E,E'$ are in the same connected component of $\cg_W$ then}
$E\si_{\lr}E'$.
\nl
We may assume that both $E,E'$ appear in some constructible representation of $W$. 
By 22.2, there exists a left cell $\G$ such that $[E:[\G]]\ne 0$, $[E':[\G]]\ne 0$.
By 21.2, we have $[E_\sp:J^\G_\bc]\ne 0$, $[E'_\sp:J^\G_\bc]\ne 0$. Hence 
$E\si_{\lr}E'$, as desired.

\subhead 23.2\endsubhead
Let $c_W$ be the set of two-sided cells of $W,L$. Consider the (surjective) map
$\Ir W@>>>c_W$ which to $E$ associates the two-sided cell $\boc$ such that 
$E\si_{\lr}x$ for $x\in\boc$. By 23.1 this induces a (surjective) map

(a) $\o_W:\fra{\cg_W}{\si}@>>>c_W$.
\nl
We conjecture that $\o_W$ is a bijection. This is made plausible by:

\proclaim{Proposition 23.3} Assume that $W,L$ is split. Then $\o_W$ is a bijection.
\endproclaim
Let $E,E'\in\Ir W$ be such that $E\si_{\lr}E'$. By 22.3, we can find constructible
representations $C,C'$ such that $[E:C]\ne 0,[E':C']\ne 0$. By 22.2, we can find 
left cells $\G,\G'$ such that $C=[\G],C'=[\G']$. Then $[E:[\G]]\ne 0$,
$[E':[\G']]\ne 0$. Let $d\in\cd\cap\G,d'\in\cd\cap\G'$. Since $\g_d=[\G]$ and 
$[E:[\G]\ne 0$, we have $E\si_{\lr}d$. Similarly, $E'\si_{\lr}d'$. Hence 
$d'\si_{\lr}d'$. By 18.4(c), there exists $u\in W$ such that $t_dt_ut_{d'}\ne 0$. 
(Here we use the splitness assumption.) Note that $j\m jt_ut_{d'}$ is a 
$J_\bc$-linear map $J^\G_\bc@>>>J^{\G'}_\bc$. This map is non-zero since it takes
$t_d$ to $t_dt_ut_{d'}\ne 0$. Thus, $\Hom_{J_\bc}(J^\G_\bc,J^{\G'}_\bc)\ne 0$. Using
21.2, we deduce that $\Hom_W([\G],[\G'])\ne 0$. Hence there exists $\tE\in\Ir W$
such that $\tE$ is a component of both $[\G]=C$ and $[\G']=C'$. Thus, both $E,\tE$ 
appear in $C$ and both $\tE,E'$ appear in $C'$. Hence $E,E'$ are in the same 
connected component of $\cg_W$. The proposition is proved.

\subhead 23.4\endsubhead
Assume now that $W,S,L,\tW,\io$ are as in 16.2 and $\tW$ is an irreducible Weyl group. 

Let $c_{\tW}^!$ be the set of all $\io$-stable two-sided cells of $\tW$. Let 
$c_{\tW}^\star$ be the set of all two-sided cells of $\tW$ which meet $W$. We have
$c_{\tW}^\star\sub c_{\tW}^!\sub c_{\tW}$. Let $f:c_W@>>>c_{\tW}^\star$ be the map 
which attaches to a two-sided cell of $W$ the unique two-sided cell of $\tW$ 
containing it; this map is well defined by 16.20(b) and is obviously surjective. 

\proclaim{Proposition 23.5} In the setup of 23.4, $\o_W$ is a bijection and
$f:c_W@>>>c_{\tW}^\star$ is a bijection.
\endproclaim
Since $\o_W,f$ are surjective, the composition $f\o_W:\fra{\cg_W}{\si}@>>>c_{\tW}^\star$
is surjective. Hence it is enough to show that this composition is injective. For 
this it suffices to check one of the two statements below:

(a) $\sh(\fra{\cg_W}{\si})=\sh(c_{\tW}^\star)$;

(b) the composition 
$\fra{\cg_W}{\si}@>f\o_W>>c_{\tW}^\star\sub c_{\tW}@>f'>>\bn\op\bn$ 
(where $f'(\boc)=(\aa(x),\aa(xw_0))$ for $x\in\boc$) is injective.
\nl
Note that the value of the composition (b) at $E$ is $(\aa_E,\aa_{E^\da})$.

{\it Case 1.} $W$ is of type $G_2$ and $\tW$ is of type $D_4$. Then (b) holds: the
composition (b) takes distinct values $(0,12),(1,7),(3,3),(7,1),(12,0)$ on the $5$ 
elements of $\fra{\cg_W}{\si}$. 

{\it Case 2.} $W$ is of type $F_4$ and $\tW$ is of type $E_6$. Then again (b) holds.

{\it Case 3.} $W$ is of type $B_n$ with $n\ge 2$ and $\tW$ is of type $A_{2n}$ or
$A_{2n+1}$. Then $\io$ is conjugation by the longest element $\tw_0$ of $\tW$. We show that (a) holds.

Let $Y$ be the set of all $E\in\Ir\tW$ (up to isomorphism) such that
$\tr(\tw_0,E)\ne 0$. Let $Y'$ be the set of all $E'\in\Ir W$ (up to isomorphism). By
23.4 and 23.1 we have a natural bijection between $c_{\tW}$ and the set of 
isomorphism classes of $E\in\Ir\tW$. If $\boc\in c_{\tW}$ corresponds to $E$, then
the number of fixed points of $\io$ on $\boc$ is clearly $\pm\dim(E)\tr(\tw_0,E)$.
Hence $\sh(c_{\tW}^\star)=\sh Y$. From 23.1 we have $\sh(\fra{\cg_W}{\si})=\sh Y$.
Hence to show (a) it suffices to show that $\sh Y=\sh Y'$. But this is shown in 
\cite{\LU}.

{\it Case 4.} Assume that $\tW$ is of type $D_n$ and $W$ is of type $B_{n-1}$ with
$n\ge 3$. We will show that (a) holds. We change notation and write $W'$ instead of
$\tW$, $W'{}^\io$ instead of $W$. Then $W'$ is as in 22.26 and we may assume that
$\io:W'@>>>W'$ is conjugation by $s_n$ (as in 22.26). Let $\cm_{1,0;n}^{N,!}$ be the
set of all elements in $\cm_{1,0;n}^N$ whose set of singles is non-empty. Let 
$$\cm_{1,0;n}^!=\lim_{N\to\iy}\cm_{1,0;n}^{N,!}.$$
By 22.26 and 23.3, $c_{W'}^!$ is naturally in bijection with $\cm_{1,0;n}^!$. By 
23.1, $\fra{\cg_{W'{}^\io}}{\si}$ is naturally in bijection with $\cm_{1,2;n-1}$. The 
identity map is clearly a bijection $\cm_{1,2;n-1}^N@>\si>>\cm_{1,0;n}^{N+1,!}$. It
induces a bijection $\cm_{1,2;n-1}@>\si>>\cm_{1,0;n}^!$. Hence to prove that
$\sh(\fra{\cg_{W'{}^\io}}{\si})\le\sh(c_{W'}^\star)$ it suffices to prove that
$\sh(\cm_{1,0;n}^!)=\sh(\cm_{1,0;n}^\star)$. In other words, we must show that 

(c) {\it any $\io$-stable two-sided cell of $W'$ meets $W'{}^\io$.}
\nl
Now 22.26 and 23.3 provide an inductive procedure to obtain any $\io$-stable two-sided
cell of $W'$. Namely such a cell is obtained by one of two procedures:

(i) we consider a $\io$-stable two-sided cell in a parabolic subgroup of type $\fS_k\T D_{n-k}$ (where 
$n-k\in[2,n-1]$) and we attach to it the unique two-sided cell of $W'$ that contains it;

(ii) we take a two-sided cell obtained in (i) and multiply it on the right by the longest element of $W'$.
\nl
Since we may assume that (c) holds when $n$ is replaced by $n-k\in[2,n-1]$, we see 
that the procedures (i) and (ii) yield only two-sided cells that contain $\io$-fixed
elements. This proves (c). The proposition is proved.

\head 24. Virtual cells\endhead
\subhead 24.1\endsubhead
In this chapter we preserve the setup of 20.1. 

A {\it virtual cell} of $W$ (with respect to $L:W@>>>\bn$) is an element of $K(W)$ 
of the form $\g_x$ (see 20.16) for some $x\in W$.

\proclaim{Lemma 24.2}Let $x\in W$ and let $\G$ be the left cell containing $x$. 

(a) If $\g_x\ne 0$ then $x\in\G\cap\G\i$. 

(b) $\g_x$ is a $\bc$-linear combination of $E\in\Ir W$ such that $[E:[\G]]\ne 0$.
\endproclaim
Assume that $\g_x\ne 0$. Then there exists $\ce\in\Ir J_\bc$ such that
$\tr(t_x,\ce)\ne 0$. We have $\ce=\op_{d\in\cd}t_d\ce$ and $t_x:\ce@>>>\ce$ maps the
summand $t_d\ce$ (where $x\si_\cl d$) into $t_{d'}\ce$, where $d'\si_\cl x\i$ and all
other summands to $0$. Since $\tr(t_x,\ce)\ne 0$, we must have 
$t_d\ce=t_{d'}\ce\ne 0$ hence $d=d'$ and $x\si_\cl x\i$. This proves (a).

We prove (b). Let $d\in\cd\cap\G$. Assume that $E\in\Ir W$ appears with $\ne 0$ 
coefficient in $\g_x$. Then $\tr(t_x,E_\sp)\ne 0$. As we have seen in the proof of
(a), we have $t_dE_\sp\ne 0$. Using 21.3,21.2, we deduce $[E_\sp:J_\bc^\G]\ne 0$ and
$[E_\sp:[\G]_\sp]\ne 0$. Hence $[E:[\G]]\ne 0$. The lemma is proved.

\subhead 24.3\endsubhead
In the remainder of this chapter we will give a number of explicit computations of 
virtual cells.

\proclaim{Lemma 24.4}In the setup of 22.10, $w_0$ acts on $[\L]$ as multiplication
by

$\e_{[\L]}=(-1)^{\su_j(a\i\mu_j-j+1)}$.
\endproclaim
Using the definitions we are reduced to the case where $k=n$ or $l=n$. If $k=n$ we 
have $\e_{[\L]}=1$ since $[\L]$ factors through $\fS_n$ and the longest element of
$W_n$ is in the kernel of $W_n@>>>\fS_n$. Similarly, if $l=n$ we have 
$\e_{[\L]}=\e_{\c_n}=(-1)^n$. The lemma is proved.

\proclaim{Proposition 24.5} Assume that we are in the setup of 22.23. Let 
$\io\in\In_r(Z)$ and let $\k:\cs_\io@>>>\bF_2$ be an affine-linear function. Let

$c(\tZ,\io,\k)=\su_{Y\in\cs_\io}(-1)^{\k(Y)}[\L_Y]\in K(W)$.
\nl
There exists $x\in W$ such that $\g_x=\pm c(\tZ,\io,\k)$.
\endproclaim
To some extent the proof is a repetition of the proof of 22.24(a), but we have to 
keep track of $\k$, a complicating factor.

We argue by induction on the rank $n$ of $\tZ$. If $n=0$ the result is clear. Assume
now that $n\ge 1$. We may assume that $0$ is not a double of $\tZ$. Let $at$ be the
largest entry of $\tZ$.

(A) Assume that there exists $i, 0\le i<t$, such that $ai$ does not appear in $\tZ$.
Then $\tZ$ is obtained from a multiset $\tZ'$ of rank $n-k<n$ by increasing each of
the $k$ largest entries by $a$ and this set of largest entries is unambiguously 
defined. The set $Z'$ of singles of $\tZ'$ is naturally in bijection with $Z$. 

Let $\io',\k'$ correspond to $\io,\k$ under this bijection. By the induction
hypothesis, there exists $x'\in W_{n-k}$ such that 
$\g_{x'}^{W_{n-k}}=\pm c(\tZ',\io',\k')\in K(W_{n-k})$. Let $w_{0,k}$ be the longest
element of $\fS_k$. Then
$$\g_{w_{0,k}x'}^{\fS_k\T W_{n-k}}=\g_{w_{0,k}}^{\fS_k}\bxt\g_{x'}^{W_{n-k}}
=\sg_k\bxt\g_{x'}^{W_{n-k}}\tag a$$
and
$$\align&\g_{w_{0,k}x'}^W=\jj_{\fS_k\T W_{n-k}}^W(\g_{w_{0,k}x'}^{\fS_k\T W_{n-k}})
\\&=\jj_{\fS_k\T W_{n-k}}^W(\sg_k\bxt\g_{x'}^{W_{n-k}})=\pm\jj_{\fS_k\T W_{n-k}}^W
(\sg_k\bxt c(\tZ',\io',\k'))=\pm c(\tZ,\io,\k),\tag b\endalign$$
as required.

(B) Assume that there exists $i, 0<i\le t$ such that $ai$ is a double of $\tZ$. Let
$\btZ$ be as in 22.8 (with respect to our $t$). Then $0$ is not a double of $\btZ$ 
and the largest entry of $\btZ$ is $at$. Let $\bZ$ be the set of singles of $\btZ$.
We have $\bZ=\{at-z;z\in Z\}$. Thus $\bZ,Z$ are naturally in (order reversing) 
bijection under $j\m at-j$. Let $\io'\in\In_r(\bZ)$ correspond to $\io$ under this
bijection and let $\k':\cs_{\io'}@>>>\bF_2$ correspond to $\k$ under this bijection.
Define $\k'':\cs_{\io'}@>>>\bF_2$ by $\k''(Y)=\k'(Y)+\su_{y\in Y}a\i y$ (an 
affine-linear function). Since $at-ai$ does not appear in $\btZ$, (A) is applicable
to $\btZ$. Hence there exists $x'\in W$ such that $\g_{x'}=\pm c(\btZ,\io',\k'')$. 
By 20.23, 22.18, 24.4, we have
$$\align&\g_{x'w_0}\ot\sg=(-1)^{l(x')}\z(\g_{x'})=\pm\z(c(\btZ,\io',\k''))\ot\sg
\\&=\pm c(\btZ,\io',\k')\ot\sg=\pm c(\tZ,\io,\k),\endalign$$
as desired.

(C) Assume that we are not in case (A) and not in case (B). Then 
$\tZ=\{0,a,2a,\do,ta\}=Z$. We can find $ia,(i+1)a$ in $Z$ such that $\io$ 
interchanges $ia,(i+1)a$ and induces on $Z-\{ia,(i+1)a\}$ an $r$-admissible
involution $\io_1$. 

(C1) Assume first that $\k(Y)=\k(Y*\{ia,(i+1)a\})$ for any $Y\in\cs_\io$. ($*$ is
symmetric difference.) Let 
$$\tZ'=\{0,a,2a,\do,ia,ia,(i+1)a,(i+2)a,\do,(t-1)a\}.$$
This has rank $n-k<n$. The set of singles of $\tZ'$ is 
$$Z'=\{0,a,2a,\do,(i-1)a,(i+1)a,\do,(t-1)a\}.$$
It is in natural (order preserving) bijection with $Z-\{ia,(i+1)a\}$. Hence $\io_1$
induces $\io'\in\In_r(Z')$. We have an obvious surjective map of affine spaces 
$\p:\cs_\io@>>>\cs_{\io'}$ and $\k$ is constant on the fibres of this map. Hence 
there is an affine-linear map $\k':\cs_{\io'}@>>>\bF_2$ such that $\k=\k'\p$. By the
induction hypothesis, there exists $x'\in W_{n-k}$ such that 
$\g_{x'}^{W_{n-k}}=\pm c(\tZ',\io',\k')\in K(W_{n-k})$. Let $w_{0,k}$ be the longest
element of $\fS_k$. Then (a), (b) hold and we are done.

(C2) Assume next that $\k(Y)\ne\k(Y*\{ia,(i+1)a\})$ for some (or equivalently any)
$Y\in\cs_\io$. We have
$$\btZ=\{0,0,a,a,2a,2a,\do,ta,ta\}-\{at-0,at-a,\do,at-at\}=\tZ=Z.$$
Let $\io'\in\In_r(Z)$ correspond to $\io$ under the bijection $z\m ta-z$ of $Z$ with
itself; let $\k':\cs_{\io'}@>>>\bF_2$ correspond to $\k$ under this bijection. Let 
$\k'':\cs_{\io'}@>>>\bF_2$ be given by $\k''(Y)=\k'(Y)+\su_{y\in Y}a\i y$ (an
affine-linear function). Note that $\io'$ interchanges $(t-i-1)a,(t-i)a$ and induces
on $Z-\{(t-i-1)a,(t-i)a\}$ an $r$-admissible involution. We show that for any 
$Y\in\cs_{\io'}$ we have 

$\k''(Y)=\k''(Y*\{(t-i-1)a,(t-i)a\})$,
\nl
or equivalently 

$\k'(Y)=\k'(Y*\{(t-i-1)a,(t-i)a\})+1$.
\nl
This follows from our assumption $\k(Y)=\k(Y*\{ia,(i+1)a\})+1$ for any 
$Y\in\cs_\io$. We see that case (C1) applies to $\io',\k''$ so that there exists 
$x'\in W$ with $\g_{x'}=\pm c(\tZ,\io',\k'')$. By 20.23, 22.18, 24.4, we have 
$$\align&\g_{x'w_0}=(-1)^{l(x')}\z(\g_{x'})\ot\sg=\pm\z(c(\btZ,\io',\k''))\ot\sg\\&
=\pm c(\btZ,\io',\k')\ot\sg=\pm c(\tZ,\io,\k),\endalign$$
as desired. The proposition is proved.

\subhead 24.6\endsubhead
Assume that we are in the setup of 22.27. By 22.27,
$$\r_4+\r_{16},\r_{12}+\r_{16}+\r_6,\r_{12}+\r_{16}+\r'_6$$
are constructible, hence (by 22.2, 21.4) are of the form $n_d\g_d$ for suitable
$d\in\cd$, hence are $\pm$ virtual cells. 

Let $d\in\cd$ be such that $n_d\g_d=\r_{12}+\r_{16}+\r_6$. Let $\G$ be the left cell
that contains $d$. Recall (21.4) that $[\G]=A\op B\op C$ where 
$A=\r_{12},B=\r_{16},C=\r_6$. By the discussion in 21.10 we see that 
$J_\bc^{\G\cap\G\i}$ has exactly three simple modules (up to isomorphism), namely
$t_dA_\sp,t_dB_\sp,t_dC_\sp$, and these are $1$-dimensional. Since $J^{\G\cap\G\i}$
is a semisimple algebra (21.9), it follows that it is commutative of dimension $3$.
Hence $\G\cap\G\i$ consists of three elements $d,x,y$. Let $p_A,p_B,p_C$ denote the
traces of $t_x$ on $A_\sp,B_\sp,C_\sp$ respectively. Let $q_A,q_B,q_C$ denote the 
traces of $t_y$ on $A_\sp,B_\sp,C_\sp$ respectively. By 20.26, 
$p_A,p_B,p_C,q_A,q_B,q_C$ are integers. Recall that the traces of $n_dt_d$ on
$A_\sp,B_\sp,C_\sp$ are $1,1,1$ respectively. Since $f_{A_\sp},f_{B_\sp},f_{C_\sp}$
are $6,2,3$ we see that the orthogonality formula 21.10 gives 
$$1+p_A^2+q_A^2=6,1+p_B^2+q_B^2=2, 1+p_C^2+q_C^2=3,$$
$$1+p_Ap_B+q_Aq_B=0,1+p_Ap_C+q_Aq_C=0,1+p_Bp_C+q_Bq_C=0.$$
Solving these equations with integer unknowns we see that there exist 
$\e,\e'\in\{1,-1\}$ so that (up to interchanging $x,y$) we have
$$(p_A,q_A)=(2\e,\e'),(p_B,q_B)=(0,-\e'),(p_C,q_C)=(-\e,\e').$$
Then $\e\g_x=2\r_{12}-\r_6$, $\e'\g_y=\r_{12}-\r_{16}+\r_6$. Hence 

$2\r_{12}-\r_6,\r_{12}-\r_{16}+\r_6$ are $\pm$ virtual cells.
\nl
The same argument shows that 
$2\r_{12}-\r_{6'},\r_{12}-\r_{16}+\r_{6'}$ are $\pm$ virtual cells. 
A similar (but simpler) argument shows that
$\r_4-\r_{16}$ is $\pm$ a virtual cell. 

Assume now that we are in the setup of 22.27 (Case 1). By 22.27,
$$\r_1+\r_2,\r_1+\r_8,\r_2+\r_9,\r_8+\r_9,\r_1^\da+\r_2^\da,\r_1^\da+\r_8^\da,
\r_2^\da+\r_9^\da,\r_8^\da+\r_9^\da,$$
are constructible, hence by 22.2, 21.4 are of the form $n_d\g_d$ for suitable
$d\in\cd$, hence are $\pm$ virtual cells. By an argument similar to that above (but
simpler) we see that
$$\r_1-\r_2,\r_1-\r_8,\r_2-\r_9,\r_8-\r_9,\r_1^\da-\r_2^\da,\r_1^\da-\r_8^\da,
\r_2^\da-\r_9^\da,\r_8^\da-\r_9^\da,$$
are $\pm$ virtual cells. 

\subhead 24.7\endsubhead
Assume that we are in the setup of 22.29. By 22.29, $\r_2+\r'_2$ is constructible,
hence by 22.2, 21.4, is of the form $n_d\g_d$ for some $d\in\cd$, hence is $\pm$ a 
virtual cell. As in 24.6, we see that $\r_2-\r'_2$ is $\pm$ a virtual cell. 

\head 25. Relative Coxeter groups\endhead
\subhead 25.1\endsubhead
Let $W,S$ be a Coxeter group and let $u\in A_W$ (see 1.17). We assume that $W$ is a
Weyl group or an affine Weyl group. Let $J$ be a $u$-stable subset of $S$ such that
$W_J$ is finite (that is, $J\ne S$ when $W$ is infinite). Let 
$U:W@>>>\{\text{permutations of } R\}$ be as in 1.5. Let $\cw$ be the set of all
$w\in W$ such that $U(w)$ carries $\{(1,s);s\in J\}$ onto itself. (A subgroup of 
$W$.) Alternatively, 
$$\cw=\{w\in W; wW_J=W_Jw, w \text{ has minimal length in } wW_J=W_Jw\}.$$
Let $K$ be the set of all $u$-orbits $k$ on $S-J$ such that $W_{J\cup k}$ is finite.
(In the case where $W$ is infinite, $K$ consists of all $u$-orbits on $S-J$ if
$\sh(u\bsl(S-J))\ge 2$ and $K=\em$ if $\sh(u\bsl(S-J))=1$.) We assume that $J$ is 
{\it $u$-excellent} in the following sense: for any $k\in K$ we have 
$w_0^{J\cup k}Jw_0^{J\cup k}=J$. 

For $k\in K$ we have $w_0^{J\cup k}w_0^Jw_0^{J\cup k}=w_0^J$ hence
$$\t_k:=w_0^{J\cup k}w_0^J=w_0^Jw_0^{J\cup k}$$
satisfies $\t_k^2=1$.

If $k\in K$ then $U(w_0^{J\cup k})$ maps $\{(1,s);s\in J\cup k\}$ onto
$\{(-1,s);s\in J\cup k\}$. It also maps $\{(\pm 1,s);s\in J\}$ onto
$\{(\pm 1,s);s\in J\}$. Hence it maps $\{(1,s);s\in J\}$ onto $\{(-1,s);s\in J\}$. 
Similarly, $U(w_0^J)$ maps $\{(-1,s);s\in J\}$ onto $\{(1,s);s\in J\}$. Hence 
$U(\t_k)=U(w_0^J)U(w_0^{J\cup k})$ maps $\{(1,s);s\in J\}$ onto $\{(1,s);s\in J\}$.
Thus, $\t_k\in\cw$. More precisely, $\t_k\in\cw^u$, the fixed point set of 
$u:\cw@>>>\cw$.

The following result is proved in \cite{\COX} assuming that $W$ is a Weyl group (see
\cite{\LPI} for the case where $W$ is an affine Weyl group).

(a) {\it $\cw^u$ is a Coxeter group on the generators $\{\t_k;k\in K\}$. Moreover, 
if $W$ is a Weyl group then $\cw^u$ is a Weyl group; if $W$ is an affine Weyl group
and $\sh(u\bsl(S-J))\ge 2$ then $\cw^u$ is an affine Weyl group; if $W$ is an affine
Weyl group and $\sh(u\bsl(S-J))=1$ then $\cw^u=\{1\}$.}

\subhead 25.2\endsubhead
We now strengthen our assumption on $J$ by assuming that there exists an adjoint 
reductive group $G_J$ defined over $\bF_q$ whose Coxeter graph is $J$ (a full 
subgraph of the Coxeter graph of $W$), such that $u:J@>>>J$ is induced by the 
Frobenius map of $G_J$ and that $G_J(\bF_q)$ 
admits a unipotent cuspidal representation $E$; let $\boc_0$ be the two-sided cell 
of $W_J$ (with the weight function given by length) corresponding to this unipotent 
representation in the classification \cite{\RED}. 
The function $\aa:W@>>>\bn$ (see 13.6) (defined in terms of the 
weight function given by the length) takes a constant value $a$ on $\boc_0$ and a 
constant value $a_k$ on $\boc_0\t_k$ for $k\in K$ (see 9.13, P.11, 15.6). The 
function
$\{\t_k;k\in K\}@>>>\bz$ given by $\t_k\m a_k-a$ takes equal values at two elements
$\t_k,\t_{k'}$ that are conjugate in $\cw^u$ (case by case check) hence it is the 
restriction of a weight function $L:\cw^u@>>>\bz$. This weight function takes $>0$ 
values on $\{\t_k;k\in K\}$. Let $\aa_L:\cw^u@>>>\bn$ be the function defined like 
$\aa:W@>>>\bn$ (see 13.6) in terms of $\cw^u$ (instead of $W$) and the weight
function just defined. Define $\aa':\cw^u@>>>\bn$ by $\aa'(x)=\aa(yx)$ where $y$ is
any element of $\boc_0$. 
This is independent of the choice of $y$, by 9.13, P.11, 15.6.

\proclaim{Conjecture 25.3} (a) $\aa_L=\aa'$.

(b) Let $\boc$ be a two-sided cell of $\cw^u$ (relative to the weight function $L$)
as in 25.2. There exists a (necessarily unique) two-sided cell $\ti\boc$ of $W$ 
(relative to the weight function given by length) such that $yx\in\ti\boc$ for any
$y\in\boc_0,x\in\boc$. Moreover the map $\boc\m\ti\boc$ is injective.
\endproclaim
This would reduce the problem of computing the two-sided cells of $\cw^u$ (relative
to the weight function $L$) to the analogous problem for $W$ (relative to the weight
function given by length).

\head 26. Representations\endhead
\subhead 26.1\endsubhead
Let $W,S$ be an affine Weyl group and let $u\in A_W$ (see 1.17). Let $J$ be a 
$u$-stable subset of $S$ with $J\ne S$. Let $\cu(J)$ be the set of isomorphism 
classes of unipotent cuspidal representations of $G_J(\bF_q)$ (as in 25.2). Note 
that $\cu(J)$ is independent of the choice of $G_J$. Let $E\in\cu(J)$. Let 
$\ch(W,J,E)$ be the Iwahori-Hecke algebra attached to $\cw^u$ (defined as in 25.1 in
terms of $W,S,J$) and to the weight function $L:\cw^u@>>>\bn$ (defined as in 25.2).
Let $\Om$ be as in 1.18. Let 
$$\Om^u=\{a\in\Om; ua=au\},\Om^{u,J}=\{a\in\Om^u;a(J)=J\}.$$
If $a\in\Om^{u,J}$ then $a:W@>>>W$ restricts to an
automorphism of $\cw^u$ as a Coxeter group; this automorphism is compatible with the
weight function $L:\cw^u@>>>\bn$ hence it induces an automorphism of the algebra 
$\ch(W,J,E)$. Hence we may form a semidirect product algebra 
$\ch(W,J,E)\ot_\ca\ca[\Om^{u,J}]$ where $\ca[\Om^{u,J}]$ is the group algebra of 
$\Om^{u,J}$ over $\ca$. 

Let $v_0\in\bc^*$ be such that $v_0=1$ or $v_0$ is not a root of $1$. Let 
$$(\ch(W,J,E)\ot_\ca\ca[\Om^{u,J}])_{v_0}$$
be the $\bc$-algebra obtained from 
$\ch(W,J,E)\ot_\ca\ca[\Om^{u,J}]$ by the change of scalars $\ca@>>>\bc,v\m v_0$. Let
$$\ci=\sqc\Ir(\ch(W,J,E)\ot_\ca\ca[\Om^{u,J}])_{v_0}$$
where $\Ir$ stands for the set of isomorphism classes of simple modules of an 
algebra and the disjoint union is taken over all $(J,E)$ as above modulo the action
of $\Om^u$.

On the other hand, let $\cg$ be a connected, simply connected almost simple 
reductive group over $\bc$, of type "dual" to that of $W$.
 Let $A(\cg)$ be the group of automorphisms of $\cg$ modulo the 
group of inner automorphisms of $\cg$. There is a natural action of $A(\cg)$ on 
$\cg$ (well defined up to conjugacy) and we form the semidirect product $\tcg$ of
$\cg$ and $A(\cg)$ via this action. Note that $\cg$ may be identified with the 
identity component of $\tcg$. Let $\cj$ be the set of all pairs $(C,\ce)$ where $C$
is a $\cg$-conjugacy class in $\tcg$ and $\ce$ is an irreducible $\cg$-equivariant 
local system on $C$.

\proclaim{Theorem 26.2} There is a natural bijection $\ci\lra\cj$.
\endproclaim
This is shown in \cite{\LPI},\cite{\LP}.
Using this bijection we may transfer the partition of 
$\ci$ into pieces indexed by the various $(J,E)$ into a partition of $\cj$ into
pieces again indexed by the various $(J,E)$. This partition can be described purely
in terms of the geometry of $\tcg$ (see \cite{\LP}).

\head 27. A new realization of Hecke algebras\endhead
\subhead 27.1\endsubhead
Let $G,F,\cp,\be,W,S,J,\cw,u,\do$ be as in 0.3. Let $H=H(G^F,\cp^F,\be)$. In this 
chapter we give a new realization of the Hecke algebra $H$ as a function space. We 
will identify $\bbq=\bc$ (where $l$ is a prime number invertible in $\bF_q$) via
some field isomorphism.

Let $P_0\in\cp^F$. Let $L$ be an $F$-stable Levi subgroup of $P_0$, $NL$  
the normalizer of $L$ in $G$, $Z_L$ the centre of $L$. Let $M=NL/Z_L$. We 
have canonically $NL/L=\cw$ hence $M=\sqc_{w\in\cw}M_w$ where $M_w$ is the inverse 
image of $w$ under the obvious map $NL/Z_L@>>>NL/L$. We have $M_1=L/Z_L=L_{ad}$. The
conjugation action defines an (injective) homomorphism $M@>>>Aut(L)$ which 
restricts, for any $w\in\cw$, to an 

(a) isomorphism of $M_w$ onto an $L_{ad}$-coset $Aut(L)_w$ in $Aut(L)$. 
\nl
By known properties of unipotent representations, there is a unique 
$L_{ad}^F$-module structure on $\be_{P_0}$ that extends the given $P_0^F$-module 
structure on $\be_{P_0}$ via the obvious homomorphism $P_0^F@>>>L_{ad}^F$. We choose
an $M^F$-module structure $\io:M^F@>>>GL(\be_{P_0})$ on $\be_{P_0}$ that extends 
this $L_{ad}^F$-module structure. (This exists by known properties of unipotent 
cuspidal representations.)

Let $w\in\cw$ and let $\co_w$ be the corresponding good $G$-orbit on $\cp\T\cp$. For
$(P_1,P_2)\in\co_w$ let $\bP_2@>\psi^{P_2}_{P_1}>>\bP_1$ be the unique isomorphism 
which takes the image of any $x\in P_1\cap P_2$ under $P_1\cap P_2@>>>\bP_2$ to the
image of $x\in P_1\cap P_2$ under $P_1\cap P_2@>>>\bP_1$. Then the composition
$$\bP_0@>\Ad(g_2)>>\bP_2@>\psi^{P_2}_{P_1}>>\bP_1@>\Ad(g_1\i)>>\bP_0$$
where $g_1,g_2\in G$, $g_1P_0g_1\i=P_1, g_2P_0g_2\i=P_2$, may be regarded as an 
element of $Aut(L)_w$ (we identify $\bP_0=L$). This corresponds under (a) to an 
element $\a_{g_1,g_2}\in M_w$.

\subhead 27.2\endsubhead
Assume now that $F(w)=w$. Define ${}^w\ph\in H$ as follows: if $(P_1,P_2)\in\co_w$ 
then $({}^w\ph)^{P_2}_{P_1}:\be_{P_2}@>>>\be_{P_1}$ is the composition
$$\be_{P_2}@>g_2\i>>\be_{P_0}@>\io(\a_{g_1,g_2})>>\be_{P_0}@>g_1>>\be_{P_1}$$
where $g_1,g_2\in G^F$, $g_2P_0g_2\i=P_2,g_1P_0g_1\i=P_1$; if $(P_1,P_2)\n\co_w$ 
then $({}^w\ph)^{P_2}_{P_1}:\be_{P_2}@>>>\be_{P_1}$ is $0$. 
($({}^w\ph)^{P_2}_{P_1}$ is independent of the choices of $g_1,g_2$.) 

\subhead 27.3\endsubhead
For $w$ as in 27.2 we have $\co_{w\i}=\{(P_2,P_1)\in\cp\T\cp;(P_1,P_2)\in\co_w\}$.
Let 

$\cu=\{P_1\in\cp^F;(P_0,P_1)\in\co_w\}$.
\nl
Then $\sh\cu=q^{l(w)}$ where $l$ is
length in $W$. The composition $({}^w\ph)({}^{w\i}\ph)$ has as $(P_0,P_0)$-component
the sum over all $P_1\in\cu$ of the compositions
$$\be_{P_0}@>\io(\a_{g_1,1})>>\be_{P_0}@>g_1>>\be_{P_1}@>g_1\i>>\be_{P_0}
@>\io(\a_{1,g_1})>>\be_{P_0}$$
where $g_1\in G^F$, $g_1P_0g_1\i=P_1$, that is $q^{l(w)}$ times the identity map of 
$\be_{P_0}$. Thus,

(a) $({}^w\ph)({}^{w\i}\ph)=q^{l(w)}({}^1\ph)+\text{ linear combination of
${}^{w'}\ph$ with } w'\ne 1$.

\subhead 27.4\endsubhead
Let $w,w'\in\cw^u$ be such that $l(ww')=l(w)+l(w')$. (Here $l$ is length in
$W$.) Then 

(a) $(P_1,P_2)\in\co_w,(P_2,P_3)\in\co_{w'}\imp(P_1,P_3)\in\co_{ww'}$,

(b) if $(P_1,P_3)\in\co_{ww'}$ then there is a unique $P_2\in\cp$ such that
$(P_1,P_2)\in\co_w,(P_2,P_3)\in\co_{w'}$.
\nl
If $P_1,P_2,P_3$ are as in (a) we have
$\psi^{P_3}_{P_1}=\psi^{P_2}_{P_1}\psi^{P_3}_{P_2}:\be_{P_3}@>>>\be_{P_1}$. From the
definitions we see that

(c) $({}^w\ph)({}^{w'}\ph)={}^{ww'}\ph$.

\subhead 27.5\endsubhead
For $w\in\cw^u$, ${}^w\ph$ is a basis element of $H_{\co_w}$. If
$w=t_{k_1}t_{k_2}\do t_{k_r}$ is a reduced expression in $\cw^u$ (see 0.3) then
$l(w)=l(t_{k_1})+l(t_{k_2})+\do+l(t_{k_r})$ (where $l$ is as in 27.4) and 
$T_w=T_{\t_{k_1}}T_{\t_{k_2}}\do T_{\t_{k_r}}$ (notation of 0.3) is a well defined 
basis element of $H_{\co_w}$ independent of the reduced expression. Hence 
${}^w\ph=x_wT_w$ where 

(a) $x_w=x_{\t_{k_1}}x_{\t_{k_2}}\do x_{\t_{k_r}}$
\nl
with $x_w\in\bc^*$ for all $w\in\cw^u$. From 27.4(c) we see that ${}^1\ph$ is the 
unit element of $H$. Hence ${}^1\ph=T_1$. By 27.3(a), we have 
$({}^{\t_k}\ph)({}^{\t_k}\ph)=q^{l(\t_k)}({}^1\ph)+\do$ hence

(b) $x_{\t_k}^2T_{\t_k}^2=q^{l(\t_k)}T_1+\text{ linear combination of
$T_{w'}$ with } w'\ne 1$.
\nl
On the other hand, by 0.3(d) we have
$T_{\t_k}^2=(q^{N_k/2}-q^{-N_k/2})T_{\t_k}+T_1$. Comparing with (b) we see that 
$x_{\t_k}^2=q^{l(\t_k)}$ hence $x_{\t_k}=\e_kq^{l(\t_k)/2}$ where 
$\e_k\in\{1,-1\}$. From (a) we see that for $w\in\cw^u$ we have 
$x_w=\e_wq^{l(w)/2}$ where $w\m\e_w$ is a function $\cw^u@>>>\{1,-1\}$ satisfying
$\e_k\e_{k'}\e_k\do=\e_{k'}\e_k\e_{k'}\do$ for $k\ne k'$ (both products have a 
number of terms equal to the order of $\t_k\t_{k'}$ in $\cw^u$). It follows that 
$w\m\e_w$ is a group homomorphism $\cw^u@>>>\{1,-1\}$. Since $M^F/L_{ad}^F=\cw^u$, 
we may regard $\e$ as a homomorphism $M^F@>>>\bc^*$ which is trivial on $L_{ad}^F$.

\subhead 27.6\endsubhead
Replacing $\io:M^F@>>>GL(\be_{P_0})$ by its tensor product with $\e:M^F@>>>\bc^*$ we
obtain a new homomorphism $\io_0:M^F@>>>GL(\be_{P_0})$. If we now redefine ${}^w\ph$
in terms of $\io_0$ rather than $\io$, then the $\e$-factors disappear and we have 

(a) ${}^w\ph=q^{l(w)/2}T_w, w\in\cw^u$.
\nl
Note that $\io_0$ is uniquely determined by property (a) and by its restriction to 
$L_{ad}^F$.

\subhead 27.7\endsubhead
Let $D=\dim\be_{P_0}$. Let $Y$ be the set of all triples $(P,P',gU_P)$ where 
$P,P'\in\cp$ and $gU_P\in G/U_P$ is such that $gPg\i=P'$ (hence $gU_P=U_{P'}g$). Now
$Y$ is naturally defined over $\bF_q$, with Frobenius map 

$F:(P,P',gU_P)@>>>(F(P),F(P'),F(g)U_{F(P)})$. 
\nl
Let $Y_0$ be the set of all triples $(P,P',gU_P^F)$ where $P,P'\in\cp^F$ and 
$gU_P^F\in G^F/U_P^F$ is such that $gPg\i=P'$ (hence $gU_P^F=U_{P'}^Fg$). We have a
bijection $Y_0@>\si>>Y^F$ given by $(P,P',gU_P^F)\m(P,P',gU_P)$.

Let $\fB$ be the vector space of all functions $f:Y_0@>>>\bc$. We define a 
multiplication $\fB\T\fB@>>>\fB$, $f',f''\m f'*f''$ by
$$(f'*f'')(P,P',gU_P^F)=\fra{D}{\sh\bP^F}\su_{\tP,g'U_P^F,g''U_{\tP}^F}
f'(P,\tP,g'U_P^F)f''(\tP,P',g''U_{\tP}^F)\tag a$$
where the sum is taken over all 
$$\tP\in\cp^F,g'U_P^F\in G^F/U_P^F,g''U_{\tP}^F\in G^F/U_{\tP}^F$$
such that 
$$g'Pg'{}\i=\tP,g''\tP g''{}\i=P',g''g'\in U_{P'}^Fg=gU_P^F.$$
Equivalently,
$$(f'*f'')(P,P',gU_P^F)=\fra{D}{\sh\bP^F}\sh(U_P^F)\i\su_{\tP,g'}
f'(P,\tP,g'U_P^F)f''(\tP,P',gg'{}\i U_{\tP}^F)\tag b$$
where the sum is taken over all $\tP\in\cp^F,g'\in G^F$ such that $g'Pg'{}\i=\tP$.
With this multiplication, $\fB$ becomes an associative algebra. 

Define $\k:H@>>>\fB$ by $\ph\m\k(\ph)$ where $\k(\ph)(P,P',gU_P^F)$ is the trace of
the composition
$$\be_P@>g>>\be_{P'}@>\ph^{P'}_P>>\be_P.$$ 
(This is independent of the choice of $g$ in its $U_P^F$-coset; $\ph^{P'}_P$ is as 
in 0.1.) 

\proclaim{Lemma 27.8} $\k:H@>>>\fB$ is an algebra homomorphism. 
\endproclaim
Let $\ph,\ph'\in H$ and let $(P,P',gU_P^F)\in Y_0$.  Then 
$$\align&\k(\ph\ph')(P,P',gU_P^F)=\tr(\be_P@>g>>\be_{P'}@>(\ph\ph')^{P'}_P>>\be_P)
\\&=\su_{\tP}\tr(\be_P@>g>>\be_{P'}@>\ph'{}^{P'}_{\tP}>>\be_{\tP}
@>\ph^{\tP}_P>>\be_P).\endalign$$
On the other hand, 
$$\align&\k(\ph)*\k(\ph')(P,P',gU_P^F)\\&
=\sh(U_P^F)\i\fra{D}{\sh\bP^F}\su_{\tP,g';g'Pg'{}\i=\tP}
\k(\ph)(P,\tP,g'U_P^F)\k(\ph')(\tP,P',gg'{}\i U_{\tP}^F)\\&
=\fra{D}{\sh P^F}\su_{\tP,g';g'Pg'{}\i=\tP}\tr(\be_P@>g'>>\be_{\tP}
@>\ph^{\tP}_P>>\be_P)\tr(\be_{\tP}@>gg'{}\i>>\be_{P'}@>\ph^{P'}_{\tP}>>\be_{\tP})\\&
=\fra{D}{\sh P^F}\su_{\tP,g';g'Pg'{}\i=\tP}
\tr(\be_P@>g'>>\be_{\tP}@>\ph^{\tP}_P>>\be_P)
\tr(\be_P@>g>>\be_{P'}@>\ph^{P'}_{\tP}>>\be_{\tP}@>g'{}\i>>\be_P).\endalign$$
It is then enough to show that for any $\tP\in\cp^F$, we have
$$\align&\tr(\be_P@>g>>\be_{P'}@>\ph'{}^{P'}_{\tP}>>\be_{\tP}@>\ph^{\tP}_P>>\be_P)
\\&=\fra{D}{\sh P^F}\su\Sb g'\in G^F\\g'Pg'{}\i=\tP\endSb\tr(\be_P@>g'>>\be_{\tP}
@>\ph^{\tP}_P>>\be_P)\tr(\be_P@>g>>\be_{P'}@>\ph^{P'}_{\tP}>>\be_{\tP}
@>g'{}\i>>\be_P).\endalign$$
Let $\g\in G^F$ be such that $\g P\g\i=\tP$. We rewrite the equality to be proved 
using the substitution $g'=\g h$:
$$\tr(\be_P@>AB>>\be_P)=\fra{D}{\sh(P^F)}
\su_{h\in P^F}\tr(\be_P@>Ah>>\be_P)\tr(\be_P@>h\i B>>\be_P)\tag a$$ 
where $A$ is the composition $\be_P@>\g>>\be_{\tP}@>\ph^{\tP}_P>>\be_P$ and $B$ is
the composition $\be_P@>g>>\be_{P'}@>\ph^{P'}_{\tP}>>\be_{\tP}@>\g\i>>\be_P$. (Then
$AB$ is the composition $\be_P@>g>>\be_{P'}@>\ph^{P'}_{\tP}>>\be_{\tP}@>\ph^{\tP}_P
>>\be_P$.) Now (a) follows immediately from the Schur orthogonality relations for 
the matrix coefficients of the irreducible representation of $P^F$ on $\be_P$. The 
lemma is proved.

\subhead 27.9\endsubhead
Let $w\in\cw^u$. Let $f_w:Y_0@>>>\bc$ be the image of $q^{-l(w)/2}({}^w\ph)$ 
(defined as in 27.2 in terms of $\io_0$) under $\k:H@>>>\fB$.

If $(P_1,P_2,gU_{P_1}^F)\in Y_0,(P_1,P_2)\n\co_w$ then $f_w(P_1,P_2,gU_{P_1}^F)=0$.

If $(P_1,P_2,gU_{P_1}^F)\in Y_0,(P_1,P_2)\in\co_w$ then $f_w(P_1,P_2,gU_{P_1}^F)$ is
$q^{-l(w)/2}$ times the trace of the composition
$$\be_{P_1}@>g_2\i g>>\be_{P_0}@>\io_0(\a_{g_1,g_2})>>\be_{P_0}@>g_1>>\be_{P_1}$$
where $g_1,g_2\in G^F$, $g_2P_0g_2\i=P_2,g_1P_0g_1\i=P_1$; here we may assume that
$g_1=g\i g_2$ hence
$$f_w(P_1,P_2,gU_{P_1}^F)=q^{-l(w)/2}\tr(\be_{P_0}@>\io_0(\a_{g\i g_2,g_2})>>
\be_{P_0})$$
where $g_2\in G^F$, $g_2P_0g_2\i=P_2$ (with $\a_{g\i g_2,g_2}$ as in 27.1).

In particular, if $(P_1,P_2,gU_{P_1}^F)\in Y_0,P_1\ne P_2$ then 
$f_1(P_1,P_2,gU_{P_1}^F)=0$; if $P_1\in\cp^F,g\in P_1^F$ then 
$$f_1(P_1,P_1,gU_{P_1}^F)=\tr(\be_{P_1}@>g>>\be_{P_1}).$$ 
Thus, $f_1$ is not identically zero.

Here are some properties of the functions $f_w$ which follow immediately from the
corresponding properties of the functions ${}^w\ph$ using 27.8.

(a) $(f_{\t_k}-q^{-N_k/2}f_1)(f_{\t_k}+q^{N_k/2}f_1)=0$ for all $k$,

(b) $f_wf_{w'}=f_{ww'}, w,w'\in\cw^u, l(ww')=l(w)+l(w')$.
\nl
Let $\bH=\k(H)$. This is a subalgebra of $\fB$ generated as a vector space by 
$\{f_w;w\in\cw^u\}$. From (b) we see that $f_1f_w=f_wf_1=f_w$ for all $w\in\cw^u$, 
hence $f_1$ is the unit element of the algebra $\bH$. From (a) we see that 
$f_{\t_k}$ is invertible in this algebra for any $k$ and then from (b) we see that 
$f_w$ is invertible in this algebra for any $w\in W$. Since $f_1\ne 0$ we have 
$f_w\ne 0$ for any $w\in\cw^u$.  Now the $f_w$ have disjoint supports. (The support
of $f_w$ is contained in $Y_w^F$ where
$Y_w=\{(P_1,P_2,gU_{P_1}^F)\in Y;(P_1,P_2)\in\co_w\}$.) It follows that the elements
$f_w (w\in\cw^u)$ are linearly independent in the vector space $\bH$. Hence 
$\k:H@>>>\bH$ is an isomorphism of algebras.

We have thus obtained a new model $\bH$ for the Hecke algebra $H$ as the vector
space of functions $f:Y_0@>>>\bc$ spanned by the functions $f_w (w\in\cw^u)$ with
multiplication $*$ as in 27.7.

\subhead 27.10\endsubhead
Let 

$Z=\{(P',gU_{P_0});P'\in\cp,g\in G/U_{P_0};gP_0g\i=P'\}$.
\nl
Let $w\in\cw^u$. Let

$Z_w=\{(P',gU_{P_0})\in Z,(P_0,P')\in\co_w\}$.
\nl
We have a morphism $\l:Z_w@>>>Aut(L)_w=M_w$ where $\l(P',gU_{P_0})$ is the 
composition
$$L=\bP_0@>\Ad(g)>>\bP'@>\psi^{P'}_{P_0}>>\bP_0=L.$$
Note that $\l$ is a smooth morphism with connected fibres.

Now $Z_w,M_w$ are naturally defined over $\bF_q$ with Frobenius maps $F$ and $\l$ 
commutes with $F$. Hence $\l$ restricts to a map 
$$Z_w^F@>\l>>M_w^F.\tag a$$
Define $f^0_w:Z_w^F@>>>\bc$ by 
$$f^0_w(P',gU_{P_0})=f_w(P_0,P',gU_{P_0}).$$
(Here $g\in G^F/U_{P_0}^F$.) We have 
$$f^0_w=q^{-l(w)/2}\l^*(\c(\io_0)_w)$$
where $\c(\io_0)_w:M_w^F@>>>\bc$ is the character of $\io_0:M^F@>>>GL(\be_{P_0})$ 
restricted to $M_w^F$.

\subhead 27.11\endsubhead
The obvious homomorphism $Aut(L)@>>>Aut(L_{ad})$ defines for any $w\in\cw^u$ an
isomorphism of $M_w=Aut(L)_w$ with a connected component of the reductive algebraic
group $Aut(L_{ad})$ with identity component $L_{ad}$. Hence we have the notion of 
character sheaf on $M_w$ (see \cite{\ICS}). Let $\hM_w$ be the set of isomorphism 
classes of character sheaves on $M_w$. Let $A\in\hM_w$. Then $A$ is 
$L_{ad}$-equivariant for the conjugation action of $M_w$. Since $\l$ is smooth with
connected fibres of fixed dimension, a suitable shift of $\l^*(A)$ is a simple
perverse sheaf $\tA$ on $Z_w$. Let $\tA^\sh$ be the unique simple perverse sheaf on
$Z$, whose support is the closure in $Z$ of the support of $\tA$ and which satisfies
$\tA^\sh|_{Z_w}=\tA$.

Let $\hM_w^F$ be the set of all $A\in\hM_w$ such that $F^*A\cong A$. For any 
$A\in\hM_w^F$ we choose an isomorphism $\ph:F^*A@>\si>>A$. There are induced 
isomorphisms $\ph:F^*\tA@>\si>>\tA,\ph:F^*\tA^\sh@>\si>>\tA^\sh$. Let
$$
\c_{A,\ph}:M_w^F@>>>\bbq,\c_{\tA,\ph}:Z_w{}^F@>>>\bbq,\c_{\tA^\sh,\ph}:Z^F@>>>\bbq$$
be the corresponding characteristic functions (alternating sums of traces of 
Frobenius at stalks of cohomology sheaves at various $F$-fixed points). We have 
$c_{\tA,\ph}=(-1)^N\l^*(c_{A,\ph})$ where $N=\dim Z_w-\dim M_w$ and
$c_{\tA,\ph}=c_{\tA^\sh,\ph}|_{Z_w^F}$.

It is known \cite{\ICS} that the functions $\c_{A,\ph}$ (where $A$ runs through 
$\hM_w^F$) form a basis for the vector space of functions $M_w^F@>>>\bc$ that are
constant on the orbits of $M^0{}^F$ (acting on $M_w^F$ by conjugation). Hence  
$$\c(\io_0)_w=\su_{A\in\hM_w^F}\xi_A\c_{A,\ph}$$
where $\xi_A\in\bc$ are uniquely determined. Applying $\l^*$ to both sides we deduce
$$q^{l(w)/2}f_w^0=\su_A\xi_A(-1)^N\c_{\tA,\ph}.$$
Hence 
$$f^0_w=q^{-l(w)/2}(-1)^N\su_A\xi_A\c_{\tA^\sh,\ph}|_{Z_w^F}.$$
The following conjecture provides a geometric interpretation of the polynomials
$p_{y,w}$ (see 5.3) attached to the Coxeter group $\cw^u$ with its weight function 
$L:\cw^u@>>>\bn$.

\proclaim{Conjecture 27.12} Assume that $y\in\cw^u$. We have
$$q^{-l(w)/2}(-1)^N\su_A\xi_A\c_{\tA^\sh,\ph}|_{Z_y^F}=p_{y,w}|_{v=\sqrt q}f_y^0,$$
$$\sum_A\xi_A\c_{\tA^\sh,\ph}|_{Z^F-\cup_{y\in\cw^u}Z_y^F}=0.$$
\endproclaim

\subhead 27.13\endsubhead
We now consider the special case where $\cp$ is the set of Borel subgroups of $G$
and $\be$ is the trivial vector bundle $\bc$. Then $\cw=W$. In this case the 
homomorphism $\io_0$ is trivial. For $w\in\cw^u=W^u$ and 
$(P_1,P_2,gU_{P_1}^F)\in Y_0$ we have

$f_w(P_1,P_2,gU_{P_1}^F)=0$ if $(P_1,P_2)\n\co_w$,
 
$f_w(P_1,P_2,gU_{P_1}^F)=q^{-l(w)/2}$, if $(P_1,P_2)\in\co_w$.
\nl
In particular, the functions in $\k(H)$ do not depend on the third coordinate
$gU_{P_1}^F$ which can therefore be omitted. For $f',f''$ in $\k(H)$ we have

$(f'*f'')(P,P')=\su_{\tP\in\cp^F}f'(P,\tP)f''(\tP,P')$.
\nl
In the present case, Conjecture 27.12 states that $p_{y,w}|_{v=\sqrt q}$ is (up to 
normalization) the restriction to $Z_y^F$ of the characteristic function of the 
intersection cohomology sheaf of the closure of $Z_w$ in $Z$. Equivalently, if for 
$w\in W^u$ we set

$\cp_w=\{P'\in\cp;(P_0,P')\in\co_w)$,
\nl
then $p_{y,w}|_{v=\sqrt q}$ is (up to normalization) the restriction to $\cp_y^F$ of
the characteristic function of the intersection cohomology sheaf of the closure of 
$\cp_w$ in $\cp$. 

This property is known to be true; it is proved in \cite{\KLL} in the case where 
$u=1$ on $W$ and is stated in the general case in \cite{\LU}. 

\subhead 27.14\endsubhead
Let $\bg,F,\cp,\be,W,S,J,\cw,u,\do$ be as in 0.6. Let $H=H(\bg^F,\cp^F,\be)$. 
Everything in 27.1-27.13 extends to this case (we replace $G$ by $\bg$ throughout) 
with the following modifications. In the definition of $\fB$ (see 27.7) we must now
restrict ourselves to functions $f:Y_0@>>>\bc$ such that 
$$\{(P,P')\in\cp^F\T\cp^F; f(P,P',gU_P^F)\ne 0 \text{ for some } g\in\bg^F\}$$ 
is contained in the union of finitely many $\bg$-orbits on $\cp\T\cp$. Also, when 
defining the multiplication $*$ in 27.7 only the definition 27.7(a) makes now sense
(in 27.7(b) the quantity $\sh(U_P^F)$ is infinite hence does not make sense). In 
27.13 one should use Iwahori subgroups instead of Borel subgroups.

\head Appendix\endhead
\subhead A.1\endsubhead
Let $\tW$ be the Coxeter group associated to a finite set $\tS$ and to the Coxeter matrix
$(m_{s,s'})_{s,s'\in\tS}$. We view $\tS$ as a subset of $\tW$; for any $I\sub\tS$ we denote by $\tW_I$ the 
subgroup of $\tW$ generated by $I$. For any $I\sub\tS$ let $[I]$ the full subgraph of the Coxeter graph of 
$\tW$ with set of vertices $I$. Let $\tl:\tW@>>>\NN$ be the length function of $\tW$. Let $\le$ be the 
standard partial order on $\tW$.

For any $m\in\ZZ_{\ge1}\cup\{\iy\}$ we set $\k_m=2\cos(\pi/m)\in\RR$.

Let $\EE$ be the $\RR$-vector space with basis $\{e_s;s\in\tS\}$. Following \cite{\BO,Ch.V, 4.1}, for any 
$s\in\tS$ we define a linear map $\s_s:\EE@>>>\EE$ by $\s_s(e_{s'})=e_{s'}+\k_{m_{s,s'}}e_s$ for all 
$s'\in\tS$. According to \cite{\BO, Ch.V, 4.3} there is a unique group homomorphism $\s:\tW@>>>GL(\EE)$ 
(``reflection representation'' of $\tW$) such that $\s(s)=\s_s$ for all $s\in\tS$.  
Let $(,):\EE\T\EE@>>>\RR$ be the symmetric bilinear form given by $(e_s,e_{s'})=-\k_{m_{s,s'}}/2$ for 
$s,s'\in\tS$. Note that for any $w\in\tW$, $\s(w)$ is an isometry of this form.

Assume that we are given a group automorphism $\t:\tW@>>>\tW$ such that $\t(\tS)=\tS$; we have necessarily 
$m_{\t(s),\t(s')}=m_{s,s'}$ for any $s,s'$ in $\tS$. We define a vector space isomorphism $\t:\EE@>\si>>\EE$
by $\t(e_s)=e_{\t(s)}$ for any $s\in\tS$. For any $s,s'\in\tS$ we have 
$\t(\s_s(e_{s'})=\s_{\t(s)}(\t(e_{s'}))$. Hence for any $w\in\tW$ we have
$\t\s(w))=\s(\t(w))\t:\EE@>>>\EE$.

Let $S$ be the set of $\t$-orbits $I$ on $\tS$ such that $\tW_I$ is finite; for any $I\in S$ let $w_0^I$ be 
the longest element of the finite Coxeter group $\tW_I$. Let ${}'\tW$ be the subgroup of $\tW$ generated by 
$\{w_0^I;I\in S\}$. Let $W=\{w\in\tW;\t(w)=w\}$. We show:

(a) $W={}'\tW$.
\nl
In the proof we shall make use of the following fact.

(b) If $I\sub\tS$ and $w\in\tW$ satisfies $sw<w$ for any $s\in I$ then $\tl(yw)=\tl(w)-\tl(y)$ for any 
$y\in\tW_I$. In particular, $y\m\tl(y)$ is bounded on $\tW_I$ so that $\tW_I$ is finite.
\nl
The proof of the first sentence in (b) is identical to the proof of 9.8(d) (with the assumption
that $w$ has maximal length in its $\tW_I$ coset replaced by $sw<w$ for any $s\in I$). The second sentence in
(b) follows from the first since $\tl(y)\le\tl(w)$ for any $y\in\tW_I$.

We prove (a). The inclusion ${}'\tW\sub W$ is obvious. We now prove the reverse inclusion. Let 
$w\in W$. We show that $w\in{}'\tW$ by induction on $\tl(w)$. If $\tl(w)=0$ then $w=1$ and the result is
obvious. Assume now that $\tl(w)\ge1$. We can find $s\in\tS$ such that $sw<w$. Then for any $i$ we have 
$\t^i(sw)<\t^i(w)$ that is $\t^i(s)w<w$. Thus $s'w<w$ for any $s'\in I$ where $I$ is the $\t$-orbit of $s$. 
Using (b) we see that $\tW_I$ is finite and $\tl(w_0^Iw)=\tl(w)-\tl(w_0^I)$. The induction hypothesis is 
applicable to $w_0^Iw$ instead of $w$ and yields $w_0^Iw\in{}'\tW$. It follows that $w\in{}'\tW$. This 
proves (a).

Next we show:

(c) Let $I\ne I'$ in $S$ be such that $\tW_{I\cup I'}$ is infinite. Then the subgroup $\cw$ of $\tW$ 
generated by $w_0^I$ and $w_0^{I'}$ is infinite.
\nl
Note that each element of $\cw$ is fixed by $\t$. Assume that $\cw$ is finite. Then we can find $w\in\cw$ of
maximal length among the elements of $\cw$. If $sw>w$ for some $s\in I$ then for any $i$ we have 
$\t^i(sw)>\t^i(w)$ hence $\t^i(s)w>w$; thus $s'w>w$ for any $s'\in I$. Using 9.7 we deduce that
$\tl(w_0^Iw)=\tl(w_0^I)+\tl(w)>\tl(w)$ contradicting the maximality of $\tl(w)$. Thus we have $sw<w$ for any
$s\in I$. Similarly we have $sw<w$ for any $s\in I'$. Using (b) we see that $\tW_{I\cup I'}$ is finite. This 
contradiction proves (c).

\mpb

Let $I\in S$. We set $m=\max\{m_{s,s'};s\in I,s'\in I\}$. We show that if $m\ge3$ then the following holds.

(d) There exists $i\in\ZZ$ such that $\t^i:I@>>>I$ is a fixed point free involution and $m_{s,s'}=m$ if 
$s,s'\in I,s\ne s'$ are in the same $\t^i$-orbit and $m_{s,s'}=2$ if $s,s'\in I$ are not in the same
$\t^i$-orbit.
\nl
We can find $s_0\in I,s'_0\in I$ such that $m_{s_0,s'_0}=m\ge3$. Let $K\sub I$ be set of vertices of the 
connected component $[K]$ of the Coxeter graph of $\tW_I$ that contains $s_0$ and $s'_0$. If $s,s'\in K$ 
then $s'=\t^i(s)$ for some $i$ and we have necessarily $\t^i(K)=K$ (indeed, $\t^i(K),K$ are sets of vertices
of connected components of the Coxeter graph of $\tW_I$ containing $s'$). Thus the group of automorphisms of 
$[K]$ acts transitively on $K$. Using the known classification of finite Coxeter groups we see that 
$\sh(K)=2$ that is, $K=\{s_0,s'_0\}$. We also see that for some $i\in\ZZ$ we have $\t^i(s_0)=s'_0$, 
$\t^i(s'_0)=s_0$. Since $I$ is a $\t$-orbit, we deduce that $I$ is a disjoint union of $\t^i$-orbits of size
$2$ and $\t^i:I@>>>I$ is an involution; moreover for each $\t^i$-orbit $\{s,s'\}$ on $I$ we have 
$m_{s,s'}=m$. We also see that if $s,s'$ are not in the same $\t^i$-orbit, then $m_{s,s'}=2$. This proves 
(d).

\mpb

Now let $I\ne I'$ in $S$ be such that $\tW_{I\cup I'}$ is finite. Let 

$m=\max\{m_{s,s'};s\in I,s'\in I\}$, $m'=\max\{m_{s,s'};s\in I',s'\in I'\}$,
$\mu=\max\{m_{s,s'};s\in I,s'\in I'\}$. 
\nl
Note that $m<\iy$, $m'<\iy$, $\mu<\iy$. We show that if $\mu\ge3$ then (after possibly interchanging $I,I'$),
(e),(f) below hold.

(e) There exists $p\in\{1,2,3\}$ (with $p=1$ if $\mu>3$) and a $p$-fold covering $u:I@>>>I'$ which commutes 
with the action of $\t$, such that for any $s\in I,s'\in I'$ we have $m_{s,s'}=\mu$ if $s'=u(s)$ and 
$m_{s,s'}=2$ if $s'\ne u(s)$. 

(f) If $p\in\{2,3\}$ then $m\le2,m'\le2$. If $p=1$ and $\mu\ge4$ then $m\le2,m'\le2$. If $p=1$ and $\mu=3$ 
then $(m,m')$ is $(4,2)$ or $(3,2)$ or $(2,2)$ or $(1,1)$.
\nl
We can find $s_0\in I,s'_0\in I'$ such that $m_{s_0,s'_0}=\mu\ge3$. Let $K\sub I$ be the set of vertices of 
the connected component $[K]$ of the Coxeter graph of $\tW_{I\cup I'}$ that contains $s_0$ and $s'_0$. Let 
$Aut[K]$ be the group of automorphisms of the Coxeter graph $[K]$.

If $s_1,s_2\in I\cap J$ then we can find $i\in\ZZ$ such that $\t^i(s_1)=s_2$. Then $K$, $t^i(K)$ are sets of
vertices of connected components of the Coxeter graph of $\tW_{I\cup I'}$ and both contain $s_2$; hence 
$K=t^i(K)$. We see that

(g) for any $s_1,s_2$ in $I\cap K$ there exists an element of $Aut[K]$ which carries $s_1$ to $s_2$. 
Similarly, for any $s'_1,s'_2$ in $I'\cap K$ there exists an element of $Aut[K]$ which carries $s'_1$ to 
$s'_2$. In particular, $Aut[K]$ acts on $K$ with at most two orbits.
\nl
For any $s\in I$ we set $\ck'_s=\{s'\in I';m_{s,s'}\ge3\}$; for any $s'\in I'$ we set 
$\ck_{s'}=\{s\in I;m_{s,s'}\ge3\}$. Let $\ck'=\ck'_{s_0}$, $\ck=\ck_{s'_0}$. We have $\ck'\sub I'\cap K$,
$\ck\sub I\cap K$.

Assume first that $\mu\ge4$. If there exists $s'\in\ck'_{s_0}-\{s'_0\}$, then using the known classification 
of finite Coxeter groups we see that $\mu\in\{4,5\}$. Since $s',s'_0\in I'\cap K$ we can find $i$ such that 
$\t^i\in Aut[K]$ and $\t^i$ carries $s'_0$ to $s'$ hence is nontrivial in $\Aut[K]$. By the classification 
of finite Coxeter groups we deduce that $\mu=4$ and $\t^i$ interchanges $s_0,s'_0$. This contradicts the 
fact that $s_0,s'_0$ are in different $\t$-orbits. Thus $\sh(\ck'_{s_0})=1$. Using the fact that $I$ is a 
$\t$-orbit, we deduce that $\sh(\ck'_s)=1$ for any $s\in I$. Similarly, we have $\sh(\ck_{s'})=1$ for any 
$s'\in I'$. Hence we have a bijection $u:I@>>>I'$ given by $s\m s'$ where $\ck'_{s}=\{s'\},\ck_{s'}=\{s\}$ 
and (e) is verified in this case. Assume now that $s\in I$ satisfies $m_{s,s_0}\ge3$. Then 
$s,s_0\in I\cap K$ hence we can find $j$ such that $\t^j\in Aut[K]$ and $\t^j$ carries $s_0$ to $s$. Thus 
$[K]$ has at least three distinct vertices $s,s_0,s'_0$ with $m_{s_0,s'_0}=4$ and has an automorphism which 
carries $s_0$ to $s$. This is impossible, by the classification of finite Coxeter groups. We see that 
$m\le2$. Similarly we have $m'\le2$. Thus (f) holds in this case.

Next we assume that $\mu=3$. If $m\ge3$ then by (d) we have $m_{s_0,s_1}=m$ for some $s_1\in I$ and by (g) 
some automorphism $r$ of $[K]$ carries $s_0$ to $s_1$. Also, using (g) and the classification of finite 
Coxeter groups we see that $m\in\{3,4\}$ and there exists $s'_1\in I'-\{s'_0\}$ such that 
$K=\{s'_0,s_0,s_1,s'_1\}$ and the edges of $[K]$ are $(s'_0,s_0),(s_0,s_1),(s_1,s'_1)$. Note that $r$ 
interchanges $s_0$ with $s_1$ and $s'_0$ with $s'_1$. In particular $s'_0$ is not connected with any 
$s'\in I'-\{s'_0\}$ in the Coxeter graph of $\tW_{I\cup I'}$ so that $m'\le2$. We also see that in this case 
$\ck=\{s_0\}$, $\ck'=\{s'_0\}$. Similarly, if $m'\ge3$, then $m'\in\{3,4\}$, $m\le2$ and $\ck=\{s_0\}$, 
$\ck'=\{s'_0\}$. Assume now that $\sh(\ck')\ge3$. Since $\ck'\sub I'\cap K$, $Aut[K]$ has at least three 
distinct elements (see (g)). Using the classification of finite Coxeter groups, we deduce that 
$\sh(\ck')=3$, $K=\ck'$, hence $\ck=\{s_0\}$; using the fact that $I$ and $I'$ are $\t$-orbits, we see that 
we have $\sh(\ck'_s)=3$ for any $s\in I$ and $\sh(\ck_{s'})=1$ for any $s'\in I'$, so that if we define 
$u:I'@>>>I$ by $u(s')=s$ where $s$ is such that $s\in\ck_{s'}$, then $u$ satisfies (e) with $I,I'$ 
interchanged and (f) holds as well. 

Thus we can assume that $\sh(\ck')\le2$; similarly we can assume that $\sh(\ck)\le2$. Assume now that 
$\sh(\ck')=2$, $\sh(\ck)\le2$; we write $\ck'=\{s'_0,s'_1\}$. This is not compatible with the inequality 
$m\ge3$, by a previous argument. Thus $m\le2$ and similarly, $m'\le2$. Hence $[K]$ is a graph of type $A_n$,
$n\ge2$. From (g) we see that $n\le3$. Since $\{s_0,s'_0,s'_1\}\sub K$ it follows that $\{s_0,s'_0,s'_1\}=K$
that is $\ck'=K$; hence $\ck=\{s_0\}$. Using the fact that $I$ and $I'$ are $\t$-orbits, we see that we have
$\sh(\ck'_s)=2$ for any $s\in I$ and $\sh(\ck_{s'})=1$ for any $s'\in I'$ so that if we define $u:I'@>>>I$ 
by $u(s')=s$ where $s$ is such that $s\in\ck_{s'}$, then $u$ satisfies (e) with $I,I'$ interchanged and (f) 
holds as well. 

Thus (e),(f) hold if $\sh(\ck')\ge2$; similarly they hold if $\sh(\ck)\ge2$. We may assume therefore that
$\sh(\ck')=\sh(\ck)=1$. Using the fact that $I$ and $I'$ are $\t$-orbits, we see that we have 
$\sh(\ck'_s)=1$ for any $s\in I$ and $\sh(\ck_{s'})=1$ for any $s'\in I'$ so that if we define $u:I'@>>>I$ by
$u(s')=s$ where $s$ is such that $s\in\ck_{s'}$ then $u$ satisfies (e) with $I,I'$ interchanged. Thus (e) is 
proved. It remains to prove (f) in the case where $p=1$ and $\mu=3$. If $m\ge3$ then by an earlier argument 
we have $m\in\{3,4\}$ and $m'\le2$. If $m'\ge3$ then again by an earlier argument we have $m'\in\{3,4\}$, 
$m\le2$. Interchanging $I,I'$ we have again $m\in\{3,4\}$ and $m'\le2$. Thus we can assume that $m\le2$, 
$m'\le2$. Then (f) is clear. This completes the proof of (f).

\subhead A.2\endsubhead
Let $E$ be an $\RR$-vector space with basis $e_1,e_2,\do,e_k,f_1,f_2,\do,f_k$. Let $\mu\in\ZZ_{\ge3}$.
Let $e=(e_1+\do+e_k)/\sqrt{k}$, $f=(f_1+\do+f_k)/\sqrt{k}$. 

For $i\in[1,k]$ we define a linear map $s_i:E@>>>E$ by 

$s_i(e_i)=-e_i$, $s_i(e_j)=e_j$ for $j\ne i$, $s_i(f_i)=f_i+\k_\mu e_i$, $s_i(f_j)=f_j$ for $j\ne i$.
\nl
For $i\in[1,k]$ we define a linear map $t_i:E@>>>E$ by 

$t_i(f_i)=-f_i$, $t_i(f_j)=f_j$ for $j\ne i$, $t_i(e_i)=e_i+\k_\mu f_i$, $t_i(e_j)=e_j$ for $j\ne i$.
\nl
Note that $s_1,\do,s_k$ commute and $t_1,\do,t_k$ commute. We set $\s=s_1s_2\do s_k$, $\ti\s=t_1t_2\do t_k$.
For all $j$ we have

$\s(e_j)=-e_j$, $\s(f_j)=f_j+\k_\mu e_j$, $\ti\s(e_j)=e_j+\k_\mu f_j$, $\ti\s(f_j)=-f_j$.
\nl
Hence 

(a) $\s\ti\s(e_j)=(\k_\mu^2-1)e_j+\k_\mu f_j$, $\s\ti\s(f_j)=-\k_\mu e_j-f_j$ for all $j$.
\nl
We see that 

(b) $\s(e)=-e$, $\s(f)=f+\k_\mu e$, $\ti\s(e)=e+\k_\mu f$, $\ti\s(f)=-f$.
\nl
Note that $E_j=\RR e_j+\RR f_j$ is $\s\ti\s$-stable. From (a) we see that the characteristic polynomial of 
$\s\ti\s$ on $E_j$ is $X^2-(\k_\mu^2-2)X+1$. Hence 

(c) $(\s\ti\s)^m=1$ on $E$.

\subhead A.3\endsubhead
Let $E$ be an $\RR$-vector space with basis $e_1,e_2,\do,e_k,e'_1,e'_2,\do,e'_k,f_1,\do,f_k$. 
Let $e=(e_1+\do+e_k+e'_1+\do+e'_k)/\sqrt{2k}$, $f=(f_1+\do+f_k)/\sqrt{k}$. 
For $i\in[1,k]$ we define a linear map $s_i:E@>>>E$ by 

$s_i(e_i)=-e_i$, $s_i(e_j)=e_j$ for $j\ne i$, $s_i(e'_j)=e'_j$ for all $j$, 

$s_i(f_i)=f_i+e_i$, $s_i(f_j)=f_j$ for $j\ne i$.
\nl
For $i\in[1,k]$ we define a linear map $s'_i:E@>>>E$ by 

$s'_i(e_j)=e_j$ for all $j$, $s'_i(e'_i)=-e'_i$, $s'_i(e'_j)=e'_j$ for $j\ne i$, 

$s'_i(f_i)=f_i+e'_i$, $s'_i(f_j)=f_j$ for $j\ne i$.
\nl
For $i\in[1,k]$ we define a linear map $t_i:E@>>>E$ by 

$t_i(e_i)=e_i+f_i$, $t_i(e_j)=e_j$ for $j\ne i$, 

$t_i(e'_i)=e'_i+f_i$, $t_i(e'_j)=e'_j$ for $j\ne i$, 

$t_i(f_i)=-f_i$, $t_i(f_j)=f_j$ for $j\ne i$.
\nl
Note that $s_1,\do,s_k,s'_1,\do,s'_k$ commute and $t_1,\do,t_k$ commute. We set 
$\s=s_1s_2\do s_ks'_1s'_2\do s'_k$, $\ti\s=t_1t_2\do t_k$. For all $j$ we have

$\s(e_j)=-e_j$, $\s(e'_j)=-e'_j$, $\s(f_j)=f_j+e_j+e'_j$, $\ti\s(e_j)=e_j+f_j$, 
$\ti\s(e'_j)=e'_j+f_j$, $\ti\s(f_j)=-f_j$.
\nl
Hence 

(a) $\s\ti\s(e_j)=f_j+e'_j$, $\s\ti\s(e'_j)=f_j+e_j$, $\s\ti\s(f_j)=-f_j-e_j-e'_j$.
\nl
We see that

(b) $\s(e)=-e$, $\s(f)=f+\sqrt{2}e$, $\ti\s(e)=e+\sqrt{2}f$, $\ti\s(f)=-f$.
\nl
Note that $E_j=\RR e_j+\RR e'_j+\RR f_j$ is $\s\ti\s$-stable.
The characteristic polynomial of $\s\ti\s$ on $E_j$ is $(X^2+1)(X+1)$. Hence 

(c) $(\s\ti\s)^4=1$ on $E$.

\subhead A.4\endsubhead
Let $E$ be an $\RR$-vector space with basis 

$e_1,e_2,\do,e_k,e'_1,e'_2,\do,e'_k,e''_1,\do,e''_k,f_1,\do,f_k$. 
\nl
Let $e=(e_1+\do+e_k+e'_1+\do+e'_k+e''_1+\do+e''_k)/\sqrt{3k}$, $f=(f_1+\do+f_k)/\sqrt{k}$. 

For $i\in[1,k]$ we define a linear map $s_i:E@>>>E$ by 

$s_i(e_i)=-e_i$, $s_i(e_j)=e_j$ for $j\ne i$, $s_i(e'_j)=e'_j$ for all $j$, 

$s_i(e''_j)=e''_j$ for all $j$, $s_i(f_i)=f_i+e_i$, $s_i(f_j)=f_j$ for $j\ne i$.
\nl
For $i\in[1,k]$ we define a linear map $s'_i:E@>>>E$ by 

$s'_i(e_j)=e_j$ for all $j$, $s'_i(e'_i)=-e'_i$, $s'_i(e'_j)=e'_j$ for $j\ne i$, 

$s'_i(e''_j)=e''_j$, for all $j$, $s'_i(f_i)=f_i+e'_i$, $s'_i(f_j)=f_j$ for $j\ne i$.
\nl
For $i\in[1,k]$ we define a linear map $s''_i:E@>>>E$ by 

$s''_i(e_j)=e_j$ for all $j$, $s''_i(e'_j)=e'_j$ for all $j$, $s''_i(e''_i)=-e''_i$, 

$s''_i(e''_j)=e''_j$ for $j\ne i$, $s''_i(f_i)=f_i+e''_i$, $s''_i(f_j)=f_j$ for $j\ne i$.
\nl
For $i\in[1,k]$ we define a linear map $t_i:E@>>>E$ by 

$t_i(e_i)=e_i+f_i$, $t_i(e_j)=e_j$ for $j\ne i$, 

$t_i(e'_i)=e'_i+f_i$, $t_i(e'_j)=e'_j$ for $j\ne i$, 

$t_i(e''_i)=e''_i+f_i$, $t_i(e''_j)=e''_j$ for $j\ne i$, 

$t_i(f_i)=-f_i$, $t_i(f_j)=f_j$ for $j\ne i$.
\nl
Note that $s_1,\do,s_k,s'_1,\do,s'_k,s''_1,\do,s''_k$ commute and $t_1,\do,t_k$ commute. We set 

$\s=s_1s_2\do s_ks'_1s'_2\do s'_ks''_1\do s''_k$, $\ti\s=t_1t_2\do t_k$.
\nl
For all $j$ we have

 $\s(e_j)=-e_j$, $\s(e'_j)=-e'_j$, $\s(e''_j)=-e''_j$, $\s(f_j)=f_j+e_j+e'_j+e''_j$, 

$\ti\s(e_j)=e_j+f_j$, $\ti\s(e'_j)=e'_j+f_j$, $\ti\s(e''_j)=e''_j+f_j$, $\ti\s(f_j)=-f_j$.
\nl
Hence 

(a) $\s\ti\s(e_j)=e'_j+e''_j+f_j$, $\s\ti\s(e'_j)=e_j+e''_j+f_j$, $\s\ti\s(e''_j)=e_j+e'_j+f_j$, 
$\s\ti\s(f_j)=-f_j-e_j-e'_j-e''_j$.
\nl
We see that 

(b) $\s(e)=-e$, $\s(f)=f+\sqrt{3}e$, $\ti\s(e)=e+\sqrt{3}f$, $\ti\s(f)=-f$.

Note that $E_j=\RR e_j+\RR e'_j+\RR e''_j+\RR f_j$ is $\s\ti\s$-stable.
From (a) we see that $(\s\ti\s)^6=1$ on $E_j$. (Its characteristic polynomial is $(X^2-X+1)(X+1)^2$.) Hence 

(c) $(\s\ti\s)^6=1$ on $E$.

\subhead A.5\endsubhead
Let $E$ be an $\RR$-vector space with basis 

$e_1,e_2,\do,e_k,e'_1,e'_2,\do,e'_k,f_1,\do,f_k,f'_1,\do,f'_k$. 
\nl
Let $e=(e_1+\do+e_k+e'_1+\do+e'_k)\sqrt{2}/\sqrt{2k}$, $f=(f_1+\do+f_k+f'_1+\do+f'_k)/\sqrt{2k}$. 

For $i\in[1,k]$ we define a linear map $s_i:E@>>>E$ by 

$s_i(e_i)=-e_i$, $s_i(e_j)=e_j$ for $j\ne i$, $s_i(e'_i)=e'_i+e_i$, 

$s_i(e'_j)=e'_j$ for $j\ne i$, $s_i(f_i)=f_i+e_i$, $s_i(f_j)=f_j$ for $j\ne i$, 

$s_i(f'_j)=f'_j$ for all $j$.
\nl
For $i\in[1,k]$ we define a linear map $s'_i:E@>>>E$ by 

$s'_i(e_i)=e_i+e'_i$, $s'_i(e_j)=e_j$ for $j\ne i$, 

$s'_i(e'_i)=-e'_i$, $s'_i(e'_j)=e'_j$ for $j\ne i$, 

$s'_i(f_j)=f_j$ for all $j$, $s'_i(f'_i)=f'_i+e'_i$, $s'_i(f'_j)=f'_j$ for $j\ne i$.
\nl

For $i\in[1,k]$ we define a linear map $t_i:E@>>>E$ by 

$t_i(e_i)=e_i+f_i$, $t_i(e_j)=e_j$ for $j\ne i$, $t_i(e'_j)=e'_j$ for all $j$, 

$t_i(f_i)=-f_i$, $t_i(f_j)=f_j$ for $j\ne i$, $t_i(f'_j)=f'_j$ for all $j$.
\nl
For $i\in[1,k]$ we define a linear map $t'_i:E@>>>E$ by 

$t'_i(e_i)=e_i+f_i$, $t'_i(e_j)=e_j$ for $j\ne i$, $t'_i(e'_j)=e'_j$ for all $j$, 

$t'_i(f_i)=-f_i$, $t'_i(f_j)=f_j$ for $j\ne i$, $t'_i(f'_j)=f'_j$ for all $j$.
\nl
We set 

$\s=s_1s'_1s_1s_2s'_2s_2\do s_ks'_ks_k$, $\ti\s=t_1t_2\do t_kt'_1t'_2\do t'_k$.
\nl
For all $j$ we have 

$\s(e_j)=-e'_j$, $\s(e'_j)=-e_j$, $\s(f_j)=f_j+e_j+e'_j$, $\s(f'_j)=f'_j+e_j+e'_j$, 

$\ti\s(e_j)=e_j+f_j$, $\ti\s(e'_j)=e'_j+f'_j$, $\ti\s(f_j)=-f_j$, $\ti\s(f'_j)=-f'_j$.
\nl
Hence 

(a) $\s\ti\s(e_j)=e_j+f_j$, $\s\ti\s(e'_j)=e'_j+f'_j$, $\s\ti\s(f_j)=-e_j-e'_j-f_j$, 
$\s\ti\s(f'_j)=-e_j-e'_j-f'_j$.
\nl
We see that 

(b) $\s(e)=-e$, $\s(f)=f+\sqrt{2}e$, $\ti\s(e)=e+\sqrt{2}f$, $\ti\s(f)=-f$.
\nl
Note that $E_j=\RR e_j+\RR e'_j+\RR f_j+\RR f'_j$ is $\s\ti\s$-stable.
From (a) we see that the characteristic polynomial of $\s\ti\s$ on $E_j$ is $(X^2+1)(X^2-1)$. Hence 

(c) $(\s\ti\s)^4=1$ on $E$.

\subhead A.6\endsubhead
Let $E$ be an $\RR$-vector space with basis 

$e_1,e_2,\do,e_k,e'_1,e'_2,\do,e'_k,f_1,\do,f_k,f'_1,\do,f'_k$. 
\nl
Let $e=(e_1+\do+e_k+e'_1+\do+e'_k)\sqrt{2+\sqrt{2}}/\sqrt{2k}$, $f=(f_1+\do+f_k+f'_1+\do+f'_k)/\sqrt{2k}$. 

For $i\in[1,k]$ we define a linear map $s_i:E@>>>E$ by 

$s_i(e_i)=-e_i$, $s_i(e_j)=e_j$ for $j\ne i$, $s_i(e'_i)=e'_i+\sqrt{2}e_i$, $s_i(e'_j)=e'_j$ for $j\ne i$,  

$s_i(f_i)=f_i+e_i$, $s_i(f_j)=f_j$ for $j\ne i$, $s_i(f'_j)=f'_j$ for all $j$.
\nl
For $i\in[1,k]$ we define a linear map $s'_i:E@>>>E$ by 

$s'_i(e_i)=e_i+\sqrt{2}e'_i$, $s'_i(e_j)=e_j$ for $j\ne i$, 

$s'_i(e'_i)=-e'_i$, $s'_i(e'_j)=e'_j$ for $j\ne i$, 

$s'_i(f_j)=f_j$ for all $j$, $s'_i(f'_i)=f'_i+e'_i$, $s'_i(f'_j)=f'_j$ for $j\ne i$.
\nl
For $i\in[1,k]$ we define a linear map $t_i:E@>>>E$ by 

$t_i(e_i)=e_i+f_i$, $t_i(e_j)=e_j$ for $j\ne i$, $t_i(e'_j)=e'_j$ for all $j$, 

$t_i(f_i)=-f_i$, $t_i(f_j)=f_j$ for $j\ne i$, $t_i(f'_j)=f'_j$ for all $j$.
\nl
For $i\in[1,k]$ we define a linear map $t'_i:E@>>>E$ by 

$t'_i(e_i)=e_i+f_i$, $t'_i(e_j)=e_j$ for $j\ne i$, $t'_i(e'_j)=e'_j$ for all $j$, 

$t'_i(f_i)=-f_i$, $t'_i(f_j)=f_j$ for $j\ne i$, $t'_i(f'_j)=f'_j$ for all $j$.
\nl
We set $\s=s_1s'_1s_1s'_1s_2s'_2s_2s'_2\do s_ks'_ks_ks'_k$, $\ti\s=t_1t_2\do t_kt'_1t'_2\do t'_k$.
For all $j$ we have 

$\s(e_j)=-e_j$, $\s(e'_j)=-e'_j$, $\s(f_j)=f_j+2e_j+\sqrt{2}e'_j$, 

$\s(f'_j)=f'_j+\sqrt{2}e_j+2e'_j$, $\ti\s(e_j)=e_j+f_j$, $\ti\s(e'_j)=e'_j+f'_j$, 

$\ti\s(f_j)=-f_j$, $\ti\s(f'_j)=-f'_j$.
\nl
Hence 

(a) $\s\ti\s(e_j)=e_j+\sqrt{2}e'_j+f_j$, $\s\ti\s(e'_j)=\sqrt{2}e_j+e'_j+f'_j$, 

$\s\ti\s(f_j)=-2e_j-\sqrt{2}e'_j-f_j$, $\s\ti\s(f'_j)=-\sqrt{2}e_j-2e'_j-f'_j$.
\nl
We see that 

(b) $\s(e)=-e$, $\s(f)=f+\sqrt{2+\sqrt{2}}e$, $\ti\s(e)=e+\sqrt{2+\sqrt{2}}f$, $\ti\s(f)=-f$.
\nl
Note that $E_j=\RR e_j+\RR e'_j+\RR f_j+\RR f'_j$ is $\s\ti\s$-stable. From (a) we see that the 
characteristic polynomial of $\s\ti\s$ on $E_j$ is $(X^2-\sqrt{2}X+1)(X^2+\sqrt{2}X+1)$. Hence 

(c) $(\s\ti\s)^8=1$ on $E$.

\subhead A.7\endsubhead
For any $I,I'$ in $S$ we define $M_{I,I'}\in\ZZ_{>0}\cup\{\iy\}$ as follows. If $I=I'$ we set $M_{I,I'}=1$.
If $\tW_{I\cup I'}$ is infinite we set $M_{I,I'}=\iy$. Now assume that $I\ne I'$ and $\tW_{I\cup I'}$ is 
finite. Let 

(a) $m=\max\{m_{s,s'};s\in I,s'\in I\}$, $m'=\max\{m_{s,s'};s\in I',s'\in I'\}$,
$\mu=\max\{m_{s,s'};s\in I,s'\in I'\}$.
\nl
If $\mu=2$ we set $M_{I,I'}=2$. 

If $\mu=3$ and $p$ in A.1(e) is $2$ or $3$ we set $M_{I,I'}=2p$.

If $\mu=3$, $p$ in A.1(e) is $1$ and $(m,m')$ is $(4,2)$ or $(2,4)$, we set $M_{I,I'}=8$.

If $\mu=3$, $p$ in A.1(e) is $1$ and $(m,m')$ is $(3,2)$ or $(2,3)$, we set $M_{I,I'}=4$.

If $\mu\ge3$ and $m\le2,m'\le2$ we set $M_{I,I'}=\mu$.
\nl
Let $W'$ be the Coxeter group with generators $g_I$ ($I\in S$) with relations $(g_Ig_{I'})^{M_{I,I'}}=1$ for
any $I,I'\in S$ such that $M_{I\cup I'}<\iy$. We have the following result.

\proclaim{Theorem A.8}The map $g_I\m w_0^I$ ($I\in S$) extends uniquely to a group homomorphism 
$\l:W'@>>>\tW$ which is an isomorphism of $W'$ onto the subgroup $W=\{w\in\tW;\t(w)=w\}$ of $\tW$.
\endproclaim
We show that for any $I,I'\in S$ such that $\tW_{I\cup I'}$ is finite we have

(a) $(w_0^Iw_0^{I'})^{M_{I,I'}}=1$ in $\tW$. 
\nl
If $I=I'$ this is a well known property of the longest element in a finite Coxeter group. Now assume that 
$I\ne I'$. By the injectivity of $\s$, see \cite{\BO, Ch.V,4.4}, it is enough to show that 
$(\s(w_0^Iw_0^{I'}))^{M_{I,I'}}=1$ in $GL(\EE)$. Let $m,m',\mu$ be as in A.7(a). Let $\EE_{I\cup I'}$ be the
subspace of $\EE$ spanned by $\{e_s;s\in I\cup I'\}$. Note that $\s(w_0^Iw_0^{I'})$ leaves stable this 
subspace and induces the identity map on $\EE/\EE_{I\cup I'}$. Moreover, the bilinear form $(,)$ is positive
definite on $\EE_{I\cup I'}$ (see \cite{\BO, Ch.5,4.8}) hence $\EE$ is the direct sum of $\EE_{I\cup I'}$ and
its perpendicular in $\EE$ on which the isometry  $\s(w_0^Iw_0^{I'})$ must act as the identity. Thus it is 
enough to show that $(\s(w_0^Iw_0^{I'}))^{M_{I,I'}}=1$ on the subspace $\EE_{I\cup I'}$.

If $\mu=2$ then $w_0^I, w_0^{I'}$ commute hence 

$(w_0^Iw_0^{I'})^{M_{I,I'}}=(w_0^Iw_0^{I'})^2=(w_0^I)^2(w_0^{I'})^2=1$,
\nl
as required.

If $\mu=3$ and $p$ in A.1(e) is $2$ (resp. $3$) then from A.3(c) (resp. A.4(c)) we see that the $4$-th power
(resp. $6$-th power) of $\s(w_0^Iw_0^{I'})|\EE_{I\cup I'}$ is $1$, as required.

If $\mu=3$ and $p$ in A.1(e) is $1$ and $(m,m')$ is $(4,2)$ or $(2,4)$ (resp. $(m,m')$ is $(3,2)$ or $(2,3)$)
then from A.6(c) (resp. A.5(c)) we see that the $8$-th power (resp. $4$-th power) of 
$\s(w_0^Iw_0^{I'})|\EE_{I\cup I'}$ is $1$, as required.

If $\mu=3$ and $p$ in A.1(e) is $1$ or if $\mu>3$ then from A.2(c) we see that the $\mu$-th power of 
$\s(w_0^Iw_0^{I'})|\EE_{I\cup I'}$ is $1$, as required. This proves (a).

From (a) we see that the map $g_I\m w_0^I$ ($I\in S$) extends (uniquely) to a group homomorphism 
$\l:W'@>>>\tW$. From A.1(a) we see that the image of $\l$ is exactly $W$. It remains to show that $\l$ is 
injective.

For any $I\in S$ we set $\ti\e_I=\sum_{s\in I}e_s$, $\ps_I=sqrt{(\ti e_I,\ti e_I)}\in\RR_{>0}$ (we have 
$(\ti e_I,\ti e_I)\in\RR_{>0}$ by \cite{\BO, Ch.V,4.8,Thm.2}); we also set
$$e_I=\ti e_I/\ps_I.$$
Note that $(e_I,e_I)=1$. Setting $m=\max\{m_{s,s'};s,s'\in I\}$, we have $\ps_I=1$ if $m=1$ and
$$\ps_I=\sqrt{\sh(I)(1-\k_m/2)}$$ 
if $m\ge2$. (If $m\le2$ this is obvious; if $m\ge3$, this follows from A.1(d).) For example, if $m\le2$ we 
have $\ps_I=\sqrt{\sh(I)}$; if $m=3$ we have $\ps_I=\sqrt{\sh(I)}\k_4\i$; if $m=4$ we have 
$\ps_I=\sqrt{\sh(I)}\k_8\i$.

\mpb

We show that for $I,I'\in S$ such that $\tW_{I\cup I'}$ is finite we have

(b) $\s(w_0^I)(e_{I'})=e_{I'}+\k_{M_{I,I'}}e_I$.
\nl
When $I=I'$ this reduces to $\s(w_0^I)(e_I)=-e_I$ which follows easily from the definitions.

We now assume that $I\ne I'$. Let $m,m',\mu$ be as in A.7(a).

If $\mu=2$ we have from the definitions $\s(w_0^I)(e_{I'})=e_{I'}$, as required (since $\k_{M_{I,I'}}=0$).

If $\mu=3$ and $p$ in A.1(e) is $2$ (resp. $3$) then from A.3(b) (resp. A.4(b)) we see that (b) holds. (We 
use that $\k_4=\sqrt{2}$, $\k_6=\sqrt{3}$.)

If $\mu=3$ and $p$ in A.1(e) is $1$ and $(m,m')$ is $(4,2)$ or $(2,4)$ (resp. $(m,m')$ is $(3,2)$ or $(2,3)$)
then from A.6(b) (resp. A.5(b)) we see that (b) holds. (We use that $\k_4=\sqrt{2}$, 
$\k_8=\sqrt{2+\sqrt{2}}$.)

If $\mu=3$ and $p$ in A.1(e) is $1$ or if $\mu>3$ then from A.2(b) we see that (b) holds. This proves (b). 

\mpb

We show that for $I,I'\in S$ such that $\tW_{I\cup I'}$ is infinite we have

(c) $\s(w_0^I)(e_{I'})=e_{I'}+xe_I$, $\s(w_0^{I'})(e_I)=e_I+xe_{I'}$ where $x\in\RR_{\ge2}$; moreover, we 
have $x=-2(e_I,e_{I'})$.
\nl
By A.1(c), the product $w_0^Iw_0^{I'}$ has infinite order in $W$. 

Next we remark that, by a standard argument, if $s_{i_1}s_{i_2}\do s_{i_q}$ is a reduced expression in $\tW$ 
($q\ge1$) then $\s(s_{i_1}s_{i_2}\do s_{i_{q-1}})e_{s_{i_q}}$ is an $\RR_{\ge0}$-linear combination of 
elements $e_{s_{i_r}}$, $r\in[1,q]$ and at least one coefficient is $>0$; moreover if $s_{i_q}$ is 
different from $s_{i_1},s_{i_2},\do,s_{i_{q-1}}$ then the coefficient of $e_{s_{i_q}}$ is $1$.

Now let $s'\in I'$. We have $m_{s,s'}\ge3$ for some $s\in I$. (Otherwise, $w_0^I,w_0^{I'}$ would commute and
$w_0^Iw_0^{I'}$ would have order $2$, which is not the case, since it has infinite order). Hence we can find
a reduced expression $s_{i_1}s_{i_2}\do s_{i_q}$ for $w_0^I$ with $s_{i_r}\in I$ and $s_{i_q}$ such that 
$m_{s_{i_q},s'}\ge3$. We have $\s(s_{i_q})e_{s'}=e_{s'}+ce_{s_{i_q}}$ where $c\in\RR_{>0}$ and
$$\align&\s(w_0^I)e_{s'}=\s(s_{i_1}s_{i_2}\do s_{i_{q}})e_{s'}\\&=
\s(s_{i_1}s_{i_2}\do s_{i_{q-1}})(e_{s'}+ce_{s_{i_q}})\\&=
\s(s_{i_1}s_{i_2}\do s_{i_{q-1}})e_{s'}+c\s(s_{i_1}s_{i_2}\do s_{i_{q-1}})e_{s_{i_q}}.\endalign$$
By the remark above this is equals $e_{s'}$ plus an $\RR_{\ge0}$-linear combination of elements 
$e_{s_{i_r}}$, $r\in[1,q]$ and at least one of the $e_{s_{i_r}}$ appears with coefficient $>0$.
Taking sum over all $s'\in I'$ we see that $\s(w_0^I)e_{I'}$ is equal to $e_{I'}$ plus an 
$\RR_{\ge0}$-linear combination of elements $e_s$, $s\in I$, and at least one of the $e_s$, $s\in I$,
appears with coefficient $>0$. Since $\t:\EE@>>>\EE$ commutes with the action of $\s(w_0^I)$ and it keeps 
fixed $e_{I'}$ it follows that $\s(w_0^I)e_{I'}-e_{I'}$ is fixed by $\t$ hence is of the form $xe_I$ for a 
well defined $x\in\RR_{>0}$. Thus we have $\s(w_0^I)e_{I'}=e_{I'}+xe_I$ with $x\in\RR_{>0}$. Similarly we 
have $\s(w_0^{I'})(e_I)=e_I+ye_{I'}$ with $y\in\RR_{>0}$. Since $\s(w_0^I)e_I=-e_I$, 
$\s(w_0^{I'})e_{I'}=-e_{I'}$, we see that the $2$-dimensional subspace $\ce=\RR e_I+\RR e_{I'}$ is stable 
under $\s(w_0^I)$ and $\s(w_0^{I'})$. 

We have $\s(w_0^I)\s(w_0^{I'})(e_I)=(xy-1)+ye_{I'}$, $\s(w_0^I)\s(w_0^{I'})(e_{I'})=-xe_I-e_{I'}$. We see 
that $\s(w_0^I)\s(w_0^{I'}):\EE_{I,I'}@>>>\EE_{I,I'}$ has determinant $1$ and trace $xy-2$.

Assume that the bilinear form $(,)$ is positive definite on $\ce$. Then $\EE$ is the direct sum of $\ce$ and
$\ce^\perp$, the perpendicular to $\c$ with respect to $(,)$, and both $\s(w_0^I)$ and $\s(w_0^{I'})$ are 
contained in the group $\cg$ of isometries of $(,)$ which preserve $\ce$ and induces the identity map on
$\ce^\perp$. By \cite{\BO, V,4.4,Cor.3} the powers of $w_0^Iw_0^{I'}$ form a discrete subgroup of the group 
of isometries of $(,)$ hence a discrete subgroup of $\cg$. This discrete subgroup is also infinite, as we 
have seen, and $\cg$ is compact. This is a contradiction. 

We see that $(,)$ is not positive definite on $\ce$; however, as we have seen above, we have $(e_I,e_I)=1$. 
It follows that the set of isotropic vectors in $\ce$ is either a union of two lines $L,L'$ (if $(,)$ is 
nondegenerate on $\ce$) or a line $L$ (if $(,)$ is degenerate on $\ce$). In the first case both $L,L'$ must 
be stable under the isometry $\s(w_0^I)\s(w_0^{I'})$ (which has determinant $1$); it follows that 
$\s(w_0^I)\s(w_0^{I'})$ is diagonalizable over $\RR$; hence it has real eigenvalues (with product $1$). In 
the second case $L$ must be stable under the isometry $\s(w_0^I)\s(w_0^{I'})$ hence $\s(w_0^I)\s(w_0^{I'})$ 
has again real eigenvalues (with product $1$). Thus in both cases the trace of $\s(w_0^Iw_0^{I'})$ on $\ce$
must be $\ge2$. In other words we have $xy-2\ge2$ that is $xy\ge4$. Since $\s(w_0^I)$ is an isometry for 
$(,)$ we have $(e_{I'}+xe_I,e_{I'}+xe_I)=(e_{I'},e_{I'})$ hence $2x(e_I,e_{I'})+x^2(e_I,e_I)=0$. Since $x>0$
and $(e_I,e_I)=1$ it follows that $x=-2(e_I,e_{I'})$. Similarly, we have $y=-2(e_I,e_{I'}$. Thus we have 
$x=y$ and $x^2\ge4$ hence $x\ge2$. This proves (c).

Let $\bEE$ be the subspace of $\EE$ spanned by the vectors $e_I;I\in S'$. From (b) we see that the action of
$W'$ on $\EE$ (by $w\m\s(\l(w))$) leaves $\bEE$ stable and the action of the generators $g_I$ is given on 
$\bEE$ by the formulas for the reflection representation of $W'$ with the following modification: when 
$I,I'\in S'$ are such that $M_{I,I'}=\iy$ then 
$$g_I(e_{I'})=e_{I'}+x_{I,I'}e_I, g_{I'}(e_I)=e_I+x_{I,I'}e_{I'}$$ 
where $x_{I,I'}\in\RR_{\ge2}$, while in the actual reflection representation we would have $x_{I,I'}=2$. 
Hence if $w\in W'$ satisfies $\l(w)=1$ then $w$ acts as identity on the (modified) reflection representation
of $W'$. But the proof of faithfulness of the reflection representation in \cite{\BO,Ch.V,4.4} extends to a 
proof of faithfulness of the modified reflection representation. (The only place where the proof must be 
changed is in \cite{\BO,ChV, 4.5, case (a)} which is an easily verified statement about an infinite dihedral
group.) This proves that $\l$ is injective. The theorem is proved.

\proclaim{Theorem A.9}Let $L:W@>>>\NN$ be the restriction to $W$ of the length function $\tl:\tW@>>>\NN$.
Then $L$ is a weight function for $W,S$.
\endproclaim
Let $l:W@>>>\NN$ be the length function of the Coxeter group $W,S$. It is enough to prove the following
statement.

(a) If $w\in W$ and $I\in S$ are such that $l(w_0^Iw)=1+l(w)$ then $\tl(w_0^Iw)=\tl(w_0^I)+\tl(w)$.
\nl
Our assumption implies that $\s(w\i)e_I$ is an $\RR_{\ge0}$-linear combination of elements $e_{I'}\in\bEE$
with $I'\in S$; hence for some $s\in I$, $\s(w\i)e_s$ is an $\RR_{\ge0}$-linear combination of elements 
$e_{s'}\in\EE$ with $s'\in\tS$ (notation as in the proof of A.8). It follows that $\tl(sw)=\tl(w)+1$. Hence 
for any $i$ we have $\tl(\t^i(sw))=\tl(\t^i(w))+1$ that is $\tl(\t^i(s)w)=\tl(w)+1$. Hence we have $s'w>w$ 
for any $s'\in I$. Using 9.7 we deduce that $\tl(w_0^Iw)=\tl(w_0^I)+\tl(w)$, as required.

\enddocument